\documentclass[12pt]{amsart}

\usepackage{amsmath}
\usepackage{amsfonts}
\usepackage{amssymb}
\usepackage{amscd}
\usepackage{graphicx}
\usepackage[abbrev,alphabetic]{amsrefs}
\RequirePackage[dvipsnames,usenames]{color}
\usepackage{soul,xcolor}
\setstcolor{red}
\usepackage{stmaryrd}
\usepackage{mathtools}
\usepackage{booktabs}
\usepackage{multirow}
\usepackage{stmaryrd} 
\newtagform{tiny}{\tiny(}{)}

\colorlet{green}{PineGreen}

\usepackage{mathtools}
\usepackage{hyperref}
\usepackage[margin=1.25in]{geometry}

\usepackage{amsthm}
\usepackage{comment}
\usepackage[all,cmtip]{xy}
\usepackage{tikz-cd}
\usetikzlibrary{cd}

\usepackage[all]{xy}

\usepackage{eqparbox} 
\newcommand{\eqmathbox}[2][M]{\eqmakebox[#1]{$\displaystyle#2$}}

\newcommand\dhxrightarrow[2][]{%
  \mathrel{\ooalign{$\xrightarrow[#1\mkern4mu]{#2\mkern4mu}$\cr%
  \hidewidth$\rightarrow\mkern4mu$}}
}

\makeatletter
\def\@tocline#1#2#3#4#5#6#7{\relax
  \ifnum #1>\c@tocdepth 
  \else
    \par \addpenalty\@secpenalty\addvspace{#2}%
    \begingroup \hyphenpenalty\@M
    \@ifempty{#4}{%
      \@tempdima\csname r@tocindent\number#1\endcsname\relax
    }{%
      \@tempdima#4\relax
    }%
    \parindent\z@ \leftskip#3\relax \advance\leftskip\@tempdima\relax
    \rightskip\@pnumwidth plus4em \parfillskip-\@pnumwidth
    #5\leavevmode\hskip-\@tempdima
      \ifcase #1
       \or\or \hskip 1em \or \hskip 2em \else \hskip 3em \fi%
      #6\nobreak\relax
    \hfill\hbox to\@pnumwidth{\@tocpagenum{#7}}\par
    \nobreak
    \endgroup
  \fi}
\makeatother

\newsavebox{\pullback}
\sbox\pullback{%
\begin{tikzpicture}%
\draw (0,0) -- (1ex,0ex);%
\draw (1ex,0ex) -- (1ex,1ex);%
\end{tikzpicture}}

\newsavebox{\pullbackdl}
\sbox\pullbackdl{%
\begin{tikzpicture}%
\draw (-1ex,0ex) -- (0ex,0ex);%
\draw (0ex,-1ex) -- (0ex,0ex);%
\end{tikzpicture}}

\newsavebox{\pushoutdr}
\sbox\pushoutdr{%
\begin{tikzpicture}%
\draw (-1ex,-1ex) -- (-1ex,0ex);%
\draw (-1ex,0ex) -- (0ex,0ex);%
\end{tikzpicture}}

\newcommand{\cyan}{\color{cyan}}
\newcommand{\blue}{\color{blue}}

\newcommand{\stacksproj}[1]{{\cite[Tag~{#1}]{stacks-project}}}

\newcommand{\rup}[1]{\lceil #1 \rceil}
\newcommand{\rdown}[1]{\lfloor #1 \rfloor}

\renewcommand{\mod}{\ \textrm{mod}\ }
 
\renewcommand{\P}{\mathbb{P}}
\newcommand{\Z}{\mathbb{Z}}
\newcommand{\Q}{\mathbb{Q}}

\newcommand{\fram}{\mathfrak{m}}
\newcommand{\mcO}{\mathcal{O}}

\newcommand{\F}{\mathbb{F}}

\newcommand{\cHom}{\mathcal{H}om}
\newcommand{\cExt}{\mathcal{E}xt}

\newcommand{\bP}{\mathbb{P}}
\newcommand{\bQ}{\mathbb{Q}}

\newcommand{\bZ}{\mathbb{Z}}

\newcommand{\cF}{\mathcal{F}}
\newcommand{\cG}{\mathcal{G}}

\newcommand{\cI}{\mathcal{I}}

\newcommand{\cO}{\mathcal{O}}

\newcommand{\MO}{\mathcal{O}}

\newcommand{\mcF}{\mathcal{F}}

\newcommand{\m}{\mathfrak{m}}

\renewcommand{\q}{q^{\infty}}

\newcommand{\qS}{q^{\infty}_{{\rm uni}}S}

\newcommand{\red}{\mathrm{red}}
\newcommand{\perf}{\mathrm{perf}}

\renewcommand{\Im}{\mathrm{Im}}

\DeclareMathOperator{\image}{Im}

\DeclareMathOperator{\res}{res}
\DeclareMathOperator{\Tr}{Tr}
\DeclareMathOperator{\tr}{tr}

\DeclareMathOperator{\reg}{reg}
\DeclareMathOperator{\adj}{adj}

\DeclareMathOperator{\Supp}{Supp}
\DeclareMathOperator{\Spec}{Spec}

\DeclareMathOperator{\codim}{codim}
\DeclareMathOperator{\Hom}{Hom}
\DeclareMathOperator{\Ext}{Ext}

\DeclareMathOperator{\Ex}{Exc}

\DeclareMathOperator{\Ker}{Ker}

\DeclareMathOperator{\depth}{depth}
\DeclareMathOperator{\Coker}{Coker}

\newcommand{\mydot}{{{\,\begin{picture}(1,1)(-1,-2)\circle*{2}\end{picture}\ }}}

\theoremstyle{plain}
\newtheorem{theorem}{Theorem}[section]
\newtheorem{thm}[theorem]{Theorem}

\newtheorem{proposition}[theorem]{Proposition}
\newtheorem{prop}[theorem]{Proposition}
\newtheorem{lemma}[theorem]{Lemma}
\newtheorem{lem}[theorem]{Lemma}
\newtheorem{corollary}[theorem]{Corollary}
\newtheorem{cor}[theorem]{Corollary}

\newtheorem{claim}[theorem]{Claim}
\newtheorem*{claim*}{Claim}

\theoremstyle{definition}

\newtheorem{dfn}[theorem]{Definition}
\newtheorem{definition}[theorem]{Definition}
\newtheorem{setting}[theorem]{Setting}

\newtheorem{notation}[theorem]{Notation}

\newtheorem*{setup*}{Setup}

\theoremstyle{remark}
\newtheorem{rem}[theorem]{Remark}
\newtheorem{remark}[theorem]{Remark}

\def\commentbox#1{\textcolor{Mahogany}%
{\footnotesize\newline{\color{Mahogany}\fbox{\parbox{\textwidth-15pt}{\textbf{comment: } #1}}}\newline}}



\numberwithin{equation}{theorem}

\title{Quasi-$F^e$-splittings and quasi-$F$-regularity}

\author{Hiromu Tanaka} 
\address{Graduate School of Mathematical Sciences, 
The University of Tokyo, 
3-8-1 Komaba, Meguro-ku, Tokyo 153-8914, JAPAN} 
\email{tanaka@ms.u-tokyo.ac.jp}

\author{Jakub Witaszek} 
\address{Department of Mathematics \\
Fine Hall, Washington Road\\
Princeton University\\  
Princeton NJ 08544-1000, USA}
\email{jwitaszek@princeton.edu}

\author{Fuetaro Yobuko}
\address{Department of Mathematics, Faculty of Science and Technology, Tokyo University of Science, 2641 Yamazaki, Noda, Chiba, 278-8510, Japan}
\email{soratobumusasabidesu@gmail.com}

\begin{document}

\begin{abstract}
We develop the theory of quasi-$F^e$-splittings, quasi-$F$-regularity, and quasi-$+$-regularity. 
\end{abstract}

\maketitle

\setcounter{tocdepth}{2}

\tableofcontents

\section{Introduction}

A key tool in positive characteristic birational geometry and commutative algebra is the $F$-splitting of {an $\F_p$-}scheme $X$, which is a splitting of the Frobenius {homomorphism} 
$\cO_X \to F_*\cO_X$, or equivalently that of $\cO_X \to F^e_*\cO_X$ for every integer $e>0$. One considers also a more restrictive notion of global $F$-regularity requiring that for every effective divisor $D$, there exists $e>0$ such that the composition $\cO_X \to F^e_*\cO_X \hookrightarrow F^e_*(\MO_X(D))$ 
splits.

Recently, the third author introduced in \cite{yobuko19} a new notion, called quasi-$F$-splitting, which shares many properties with $F$-splittings but is satisfied by a broader class of schemes. We say that {an $\F_p$-}scheme $X$ is \emph{$n$-quasi-$F$-split} for an integer $n > 0$ if there exists a dashed arrow making the following diagram
\begin{equation} \label{eq:intro-quasi-F-def}
\begin{tikzcd}
W_n\mathcal{O}_X \arrow[r, "F"] \arrow[d, "R^{n-1}"']
& F_*W_n\mathcal{O}_X \arrow[ld, dashed, "\exists"'] \\
\mathcal{O}_X 
\end{tikzcd}
\end{equation}
commutative.  We say that $X$ is \emph{quasi-$F$-split} if it is $n$-quasi-$F$-split for some integer $n>0$. For a more comprehensive introduction to the theory of quasi-$F$-splittings, we refer to the introductions of \cite{yobuko19} and \cite{KTTWYY1}. 

One of the key limitations of the above definition is that it involves only one power of Frobenius. If one naively replaces $F$ with $F^e$ in (\ref{eq:intro-quasi-F-def}), the notion becomes redundant, being equivalent to the standard $F$-splitting {(Proposition \ref{p-A'-F-split})}. 
The main goal of our article is to introduce proper definitions of quasi-$F^e$-splittings and quasi-$F$-regularity and to verify that they satisfy the comprehensive set of properties one would expect from such notions, especially in the context of birational geometry. This effort builds upon \cite{KTTWYY1}, where such properties were established for standard quasi-$F$-splittings.

Let $Q_{X,n}$ be the pushout of (\ref{eq:intro-quasi-F-def}), which one can check to be an $\cO_X$-module. Then $X$ is $n$-quasi-$F$-split, if the induced map $\cO_X \to Q_{X,n}$ splits, or, equivalently, that the dual map $\Hom_{\cO_X}(Q_{X,n},\cO_X) \to H^0(X, \cO_X)$ is surjective. Similarly, for every integer $e>0$, we can consider the pushout diagram
\[
\begin{tikzcd}
W_n\mathcal{O}_X \arrow[r, "F^e"] \arrow[d, "R^{n-1}"']
& F^e_*W_n\mathcal{O}_X \arrow{d}\\
\mathcal{O}_X  \arrow{r} & Q^e_{X,n}. \arrow[lu, phantom, "\usebox\pushoutdr" , very near start, yshift=0em, xshift=0.6em, color=black]
\end{tikzcd}
\]
\begin{definition} \label{def:intro-quasi-F^e-split}
We say that $X$ is \emph{$n$-quasi-$F^e$-split} if and only if the following map \[
\Hom_{W_n\cO_X}(Q^e_{X,n},W_n\omega_X(-K_X)) \to H^0(X, \cO_X), 
\]
obtained by applying $\Hom_{W_n\cO_X}(-, W_n\omega_X(-K_X))$ to $\cO_X \to Q^e_{X,n}$, is surjective. 

Here, $W_n\omega_X$ is the dualising sheaf of the scheme $W_nX$, and $W_n\omega_X(-K_X)$ is the $S_2$-hull of $W_n\omega_X \otimes_{W_n\cO_X} W_n\cO_X(-K_X)$, where $W_n\cO_X(-K_X)$ is the Teichmuller lift of $\cO_X(-K_X)$.
\end{definition}
\noindent We emphasise that $W_n\omega_X(-K_X)$ is not equal to $W_n\cO_X$ unless $n=1$. The reader should feel free to skip the next two remarks.

\begin{remark}
What is special about the case of $e>1$ is that $Q^e_{X,n}$ is not an $\cO_X$-module any more: it is only a $W_{\min(e,n)}\cO_X$-module {(Proposition \ref{p-W_eO-str})}. In particular, the above definition is not equivalent to the splitting of $\cO_X \to Q^e_{X,n}$; as mentioned before asking for such a splitting would give a redundant notion again, equivalent to the standard $F$-splitting (Proposition \ref{p-A'-F-split}). One can however reformulate Definition \ref{def:intro-quasi-F^e-split} as the existence of a dashed arrow rendering the following diagram
\[
\begin{tikzcd}
W_n\mathcal{O}_X \arrow[r, "F^e"] \arrow[d, "R^{n-1}"']
& F^e_*W_n\mathcal{O}_X \arrow[ldd, dashed, "\exists"'] \\
\mathcal{O}_X \arrow[d] &   \\
W_n\omega_X(-K_X). &
\end{tikzcd}
\]
commutative {(Proposition \ref{p-many-defs})}, where $\cO_X \to W_n\omega_X(-K_X)$ is obtained from the restriction $R^{n-1} \colon W_n\cO_X \to \cO_X$ by applying ${\cHom_{W_n\cO_X}(-,W_n\omega_X(-K_X))}$. In the above, one can equivalently replace $W_n\omega_X(-K_X)$ by $W_e\omega_X(-K_X)$.
\end{remark}

\begin{remark}
Definition \ref{def:intro-quasi-F^e-split} can be also reformulated in terms of local cohomology 
{(Lemma \ref{lem:cohomological-criterion-for-log-quasi-F-splitting})}. For example, if $X = \Spec R$ for a Gorenstein local ring $(R, \m)$ of dimension $d$, then $X$ is $n$-quasi-$F^e$-split if and only if the following map is injective
\[
H^d_{\m}(\cO_X) \to H^d_{W_n\m}(Q^e_{X,n}).
\]
\end{remark}

In contrast to the case of $F$-splittings, 
the higher the $e$ is, the more restrictive the condition of being quasi-$F^e$-split becomes.
For example, a supersingular elliptic curve $X$ is $n$-quasi-$F^e$-split if and only if $e \leq n-1$. 


\begin{theorem}[{Corollary \ref{c-str-CY}}]
    Let $k$ be an algebraically closed field of characteristic $p>0$.
    Let $X$ be a $d$-dimensional smooth proper variety  over $k$ such that $d\geq 2$,  $\omega_X \simeq \cO_X$, and $H^{d-1}(X, \cO_X)=0$.
Let $h$ be the quasi-$F$-split height of $X$\footnote{this means that $X$ is $n$-quasi-$F$-split if and only if $n \geq h$}. 

Then $X$ is $n$-quasi-$F^e$-split if and only if $n \geq eh-e+1$.
\end{theorem}
\noindent See also Corollary \ref{c-abel-var} for the case of abelian varieties. It is an intriguing problem to find all pairs of integers $(n,e)$ for which an explicit variety $X$ is $n$-quasi-$F^e$-split.

In view of the above, we say that:
\begin{itemize}
    \item $X$ is \emph{quasi-$F^\infty$-split} if $\forall_{e>0}\, \exists_{n>0}\,$ $X$ is $n$-quasi-$F^e$-split.
    \item $X$ is \emph{{uniformly} quasi-$F^\infty$-split} if $\exists_{n>0}\,\forall_{e>0}\,$ $X$ is $n$-quasi-$F^e$-split.
\end{itemize}
We emphasise that by the above result, Calabi-Yau varieties are never {uniformly} quasi-$F^\infty$-split {except when they are $F$-split}.

We now move to the definition of quasi-$F$-regularity. Since working with small coefficients often causes problems in the setting of quasi-$F$-splittings\footnote{in many of our theorems for pairs $(X,\Delta)$, we need to assume that $\Delta \geq \{p^i\Delta\}$ for every $i>0$; this assumption is satisfied for example when $\Delta$ has standard coefficients, but not in general}, we also define the notion of quasi-$+$-regularity which is often easier to work with.
\begin{definition}
We say that $X$ is \emph{globally $n$-quasi-$F$-regular} if for every {effective} divisor $D$, there exists {a rational number} ${0 < \epsilon < 1}$ 
such that $(X,\epsilon D)$ is $n$-quasi-$F^{\infty}$-split. 
{We say that $X$}  
is \emph{globally quasi-$F$-regular} if it is globally $n$-quasi-$F$-regular for some $n>0$.
\end{definition}
\begin{definition}
We say that $X$ is \emph{globally $n$-quasi-$+$-regular} if for every finite surjective 
{morphism} $f \colon Y \to X$ from a normal integral  scheme $Y$, the following map 
is surjective: \[
\Hom_{W_n\cO_X}(Q^f_{X,n},W_n\omega_X(-K_X)) \to H^0(X, \cO_X).
\]
Here 
{this} map is obtained by applying $\Hom_{W_n\cO_X}(-, W_n\omega_X(-K_X))$ to the following pushout diagram:
\[
\begin{tikzcd}
W_n\MO_X \arrow[r, "f^*"] \arrow[d, "R^{n-1}"'] & f_*W_n\MO_Y \arrow[d] \\
\MO_X \arrow[r, "\Phi^f_{X, n}"] & Q^f_{X, n}.
\end{tikzcd}
\]
We say {that} $X$ is \emph{globally quasi-$+$-regular} if it is globally $n$-quasi-$+$-regular for some $n>0$.
\end{definition}

\begin{proposition}[{Proposition \ref{p-QFR-to-Q+R}}]
If $X$ is globally quasi-$F$-regular, then it is globally quasi-$+$-regular. 
\end{proposition}

\noindent In the forthcoming paper \cite{KTTWYY3}, we will show that the converse is true for $\bQ$-Gorenstein affine schemes $X$ by building up the theory of quasi-test ideals\footnote{
the equivalence of global $F$-regularity and global $+$-regularity, 
for example for non-$\bQ$-Gorenstein affine schemes or even smooth projective varieties, is a long-standing open problem.}.

As a source of examples, we verify that all one-dimensional log Fano pairs $(\bP^1, \Delta)$ are globally quasi-$F$-regular (Theorem \ref{t-fano-curve-strong}). Moreover, we show the following result generalising \cite[Theorem C]{KTTWYY1}.
\begin{theorem}[{Theorem \ref{t-2dim-klt-QFR2}}]
Let $X$ be an affine klt surface over a perfect field of characteristic $p>0$. Then $X$ is globally quasi-$+$-regular. 
\end{theorem}
However, generalising the main result of \cite{KTTWYY2} to the quasi-$+$-regular case (specifically that three-dimensional $\bQ$-factorial klt singularities in characteristic $p \geq 42$ are quasi-$+$-regular) is a bit more subtle. This will be done in the forthcoming paper \cite{KTTWYY3} with help of quasi-test ideals.

The following diagram summarises the connection between some of the notions defined above:
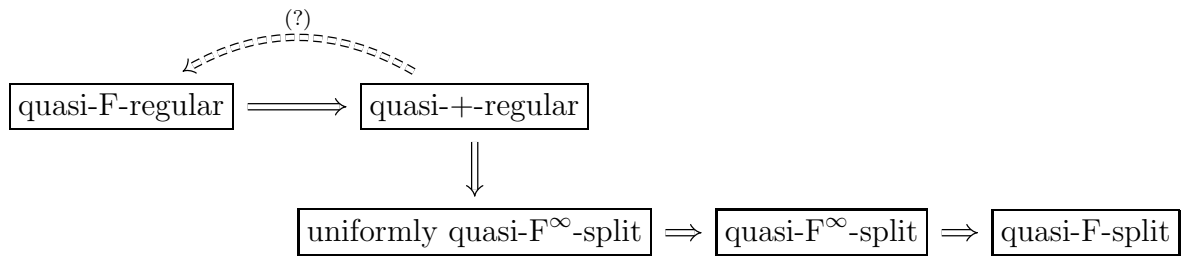
\[
\begin{tikzcd}[column sep = small]
\fbox{\text{quasi-F-regular}} \ar[r, Rightarrow] & \fbox{\text{quasi-+-regular}} \ar[l, bend right = 30, Rightarrow, dashed, "(?)"'] \ar[d, Rightarrow]  & & \\
& \fbox{$\text{{uniformly} quasi-F}^\infty\text{-split}$} \ar[r, Rightarrow] & \fbox{$\text{quasi-F}^\infty\text{-split}$} \ar[r, Rightarrow] & \fbox{\text{quasi-F-split}}
\end{tikzcd}
\]
By changing the order of quantifiers in the above definitions, we also define weaker notions of \emph{feeble quasi-$F$-regularity} and \emph{feeble quasi-$+$-regularity} {(cf.\ 
Definition \ref{d-weak-QFR}, Definition \ref{d-Q+R})}.\\

The key advantage of using quasi-$+$-regularity or quasi-$F^e$-splittings is that one can significantly weaken the assumptions of theorems proven in \cite{KTTWYY1} for quasi-$F$-splittings. This is most apparent in the following version of inversion of adjunction\footnote{Let us remind the reader that, in general, the inversion of adjunction is false for quasi-$F$-splittings (see {\cite[Example 7.7]{KTY22}}), but we could establish it in the relatively projective case when restricting to an anti-ample divisor $D$ (\cite[Corollary 4.12]{KTTWYY1}). Such a case is essential in the study of singularities via plt blow-ups (extractions of Koll\'ar's components).}. 

\begin{theorem}[{Corollary \ref{c-IOA-anti-ample+}}]
\label{thm:intro-quasi-adjunction}
Let $X$ be a normal variety admitting a projective morphism  $\pi: X \rightarrow Z$ to an affine normal variety $Z$ over a perfect field of characteristic $p > 0$. Let $S$ be a normal prime divisor on $X$. 
Suppose that 
$(X, S)$ is a $\Q$-factorial plt pair, 
$\MO_X(D)$ is Cohen-Macaulay for every Weil divisor $D$ on $X$, 
and 
$\MO_S(\Gamma)$ is Cohen-Macaulay for every Weil divisor $\Gamma$ on $S$.
Further assume that 
\begin{enumerate}
    \item $X$ is locally quasi-$+$-regular, 
    \item $A := -(K_X+S)$ and $-S$ are ample, {and}
    \item $(S, B_S)$ is globally quasi-$+$-regular, where $K_S+B_S = (K_X+S)|_S$. 
\end{enumerate}
Then $X$ is globally quasi-$+$-regular over an open neighbourhood of $\pi(S)$.
\end{theorem}
\noindent We also prove an analogous result for 
{uniform and non-uniform} 
quasi-$F^\infty$-splittings  
{(Corollary \ref{c-IOA-anti-ample}, Corollary \ref{c-IOA-anti-ample-uni})}.  

The above result provides two significant advancements in comparison to \cite[{Corollary 4.12}]{KTTWYY1}. First, we do not need to assume that $X$ satisfies relative Kawamata-Viehweg vanishing over $Z$ (briefly speaking, Kawamata-Viehweg vanishing is replaced by Serre's vanishing by taking $e\gg 0$ in the definition of quasi-$F^e$-splitting). Second, the statement is now intrinsic to the category of quasi-$F$-singularities; specifically, in Assumption (1) we only need to assume that $X$ is locally quasi-$+$-regular (or locally quasi-$F^\infty$-split), as opposed to locally $F$-split as in \cite[Corollary 4.12]{KTTWYY1}. Assumption (1) can often be verified by performing another plt blow-up.


We conclude the introduction by listing some other results of this paper:
{
\begin{enumerate}
    \item invariance of quasi-$F^e$-splittings under finite {covers} 
    (Subsection \ref{ss:finite-covers});
    \item criterion for quasi-$F^e$-splittings via Cartier operator (Theorem \ref{thm:higher-Cartier-criterion-for-quasi-F-split});
    \item introduction of quasi-$F^e$-stable sections $q^eS^0(X,L)$ and quasi-$+$-stable sections $qB^0(X,L)$ for a line bundle $L$ (Subsection \ref{ss:quasi-F^e-stable-sections} and Definition \ref{d-def-+-section});
    \item introduction of adjoint variants of the above results and definitions for pairs $(X,S)$ (Subsection \ref{ss:pure-quasi-F-split}, \ref{ss:quasi-F^e-stable-sections}, etc.).
\end{enumerate}}

\noindent A technical difficulty in our paper is establishing that quasi-$F^\infty$-stable sections and quasi-$+$-regular stable sections can be calculated for a single $e>0$ or by  a single finite cover, respectively {(Corollary \ref{cor-stab-q^e} and 
Proposition \ref{prop:qeB0-stabilisation}). 
This is well known for $F$-splittings and $+$-regularity in positive characteristic {in the local case} (see, e.g., \cite{BST15}), but as we could not find a reference which covers all the cases needed for our paper, we derived such results from scratch in the appendix.

\subsection{Future directions}
In the forthcoming 
{paper} \cite{KTTWYY3} joint with Kawakami, Takamatsu, and Yoshikawa, we establish the following results, {amongst} other things.
\begin{enumerate}
    \item quasi-$+$-regularity and quasi-$F$-regularity agree for $\bQ$-Gorenstein singularities;
    \item quasi-$+$-regular $\bQ$-Gorenstein singularities are klt;
    \item quasi-$+$-regular singularities of dimension $d \leq 3$ are Cohen-Macaulay.  
    \item three-dimensional $\bQ$-factorial klt singularities are quasi-$+$-regular in characteristic $p \geq 42$.
\end{enumerate}
We do not know 
{whether} (3) holds in 
dimension $\geq 4$. 

\medskip
\noindent {\bf Acknowledgements.}
The authors are grateful to the anonymous referee for carefully reading the manuscript and pointing out numerous errors. The referee's comments and corrections have significantly improved the paper.

\begin{itemize}
    \item The first author was supported by JSPS KAKENHI Grant numbers JP22H01112 and JP23K03028. 
    \item The second author was supported by NSF research grant DMS-2101897.
    \item The third author was supported by JSPS KAKENHI Grant number JP19K14501 and 24K16893.
\end{itemize}

\section{Preliminaries}

\subsection{Notation}\label{ss:notation}



In this subsection, we summarise notation and basic definitions used in this {article}. 
\begin{enumerate}
\item Throughout the paper, $p$ denotes a prime number and $\F_p := \Z/p\Z$. 
Given an $\F_p$-scheme, we denote by $F \colon X \to X$ the absolute Frobenius morphism. 
\item For an integral scheme $X$, 
we define the {\em function field} $K(X)$ of $X$ 
as the stalk $\MO_{X, \xi}$ at the generic point $\xi$ of $X$. 
\item {An effective Cartier divisor $D \subseteq X$ on a Noetherian scheme $X$ is called \emph{simple normal crossing} if for every $x \in D$, the local ring $\cO_{X,x}$ is regular and there exists a regular system of parameters $x_1,\ldots, x_d$ in the maximal ideal $\m$ of $\cO_{X,x}$ and $1 \leq r \leq d$ such that $D$ is defined by $x_1\cdots x_r$ in $\cO_{X,x}$ (cf.\ \stacksproj{0BI9 and Tag 0BIA}).}
\item {Given an integral normal Noetherian scheme $X$ and a $\Q$-divisor $\Delta$, 
a projective birational morphism $\pi \colon Y \to X$ is called \emph{a log resolution of $(X, \Delta)$} 
if $Y$ is regular and 
the closed subscheme $\mathrm{Exc}(\pi) \cup \Supp \pi^{-1}\Delta$, equipped with the reduced scheme structure, is a simple normal crossing divisor.}
    \item We say that $X$ is a {\em variety} (over a field $k$) if $X$ is an integral scheme 
    which is separated and of finite type over $k$. We say that $X$ is a {\em curve} if $X$ is a variety of dimension one, and that $X$ is a {\em surface} if $X$ is a variety of dimension two.
\item We say that a scheme $X$ is {\em excellent} if it is quasi-compact and 
there exists an affine open cover $X =\bigcup_{i \in I} \Spec A_i$ such that 
 $A_i$ is an excellent ring for every $i \in I$. 
We note that the regular locus of an integral normal excellent scheme $X$ 
is an open dense subset of $X$ which contains all the points of codimension one
(see \stacksproj{07P7 and 0BX2}).
    \item \label{sss-Ffinite} {We say that an $\F_p$-scheme $X$ is {\em $F$-finite} if 
$F \colon X \to X$ is a finite morphism,  
 and  we say that an $\F_p$-algebra $R$ is {\em $F$-finite} if $\Spec R$ is $F$-finite. Such schemes admit many good properties.
\begin{enumerate}
\renewcommand{\labelenumi}{(\roman{enumi})}
    \item If $R$ is an $F$-finite Noetherian $\F_p$-algebra, 
    then it is a homomorphic image of a regular ring of finite Krull dimension \cite[Remark 13.6]{Gab04}; 
    in particular, $R$ is excellent, it admits a dualising complex, and $\dim R < \infty$.  
    \item Being $F$-finite is stable under localisation and ideal-adic completions \cite[Example 9]{Has15}.
    \item If a scheme $X$ is of finite type over an $F$-finite Noetherian $\F_p$-scheme $Y$, then it is also $F$-finite. 
\end{enumerate}}

\item \label{it:essential}
We say that a ring $A$ is {\em essentially of finite type over } a ring $R$ if 
$A$ is isomorphic to $S^{-1}(R[x_1, ..., x_n]/I)$ as an $R$-algebra, 
where $R[x_1, ..., x_n]$ is the polynomial ring over $R$ with variables $x_1, ..., x_n$, 
$I$ is an ideal of $R[x_1, ..., x_n]$, and $S$ is a multiplicatively closed subset of 
$R[x_1, ..., x_n]/I$. 
\item 
Given an integral normal Noetherian scheme $X$ and a $\bQ$-divisor $D$, 
we define the subsheaf $\MO_X(D)$ of the constant sheaf $K(X)$ on $X$ 
{by the following formula}
\[
\Gamma(U, \MO_X(D)) = 
\{ \varphi \in K(X) \mid 
\left({\rm div}(\varphi)+D\right)|_U \geq 0\}
\]
{for every} open subset $U$ of $X$. 
In particular, 
$\cO_X(\rdown{{D}}) = \cO_X({D})$.
\item  
Take $n \in \Z_{>0}$. 
Given an $F$-finite Noetherian $\F_p$-algebra $R$ and an ideal $I$ of $R$, 
we set 
\[
W_nI := \{ (a_0, a_1, ..., a_{n-1}) \in W_nR \mid a_0, a_1, \ldots, a_{n-1} \in I\}. 
\]
It is well known that $W_nI$ is an ideal of $W_nR$. 
Moreover, if $(R, \m)$ is a local ring, then $(W_nR, \sqrt{W_n\m})$ is a local ring. In particular, it holds that 
\[
H^i_{W_n\m}(M) = H^i_{\sqrt{W_n\m}}(M)
\]
for a $W_nR$-module $M$.
\item\label{ss-nota-mod} Given $b \in \Z$ and $n \in \Z_{>0}$, 
the equality $a := b \mod n$ means that $a$ is the integer satisfying $a -b \in n\Z$ and 
$0 \leq a < n$. 
For example,  $a := 13 \mod 5$ means that $a \in \Z$ and $a=3$; 
in particular, $a$ is not defined as an element of $\Z/5\Z$.
\end{enumerate}

Last, we emphasise the following 
non-standard definition used in  the article.

\begin{definition}
We say that an integral normal Noetherian scheme $X$ is \emph{divisorially Cohen-Macaulay} 
if $\MO_X(D)$ is Cohen-Macaulay for every Weil divisor $D$ on $X$. 
\end{definition}
\noindent In particular, divisorially Cohen-Macaulay schemes are Cohen-Macaulay. 


\subsubsection{Dualising complexes and canonical divisors}\label{sss-dualising}

For basic properties {of} dualising complexes, we refer to \stacksproj{0A85} and \cite{Har66}. 

Throughout the article, whenever we consider a dualising complex, we implicitly work over a fixed excellent base scheme, which is denoted by $B$ in this subsection.
We assume that $B$ admits a dualising complex and we implicitly make a choice of {one}, say $\omega^{\mydot}_B$. It {automatically} holds that $\dim B < \infty$ \cite[Corollary V.7.2]{Har66}.

Given an irreducible $S_2$ excellent scheme $X$ and a separated morphism $\pi \colon X \to B$ of finite type, we set $\omega^{\mydot}_X := \pi^! \omega^{\mydot}_B$. Then $\omega^{\mydot}_X$ is a dualising complex of $X$ (cf.\ \stacksproj{0AU5}). In what follows, if a base scheme is not specified, then we implicitly take $X=B$ and $\pi = \mathrm{id}$. 
There exists a unique $e \in \Z$ such that 
$\mathcal H^{-e}(\omega^{\mydot}_X) \neq 0$ and 
$\mathcal H^{i}(\omega^{\mydot}_X) = 0$ for $i< -e$. We say that $\omega^{\mydot}_X$ is \emph{normalised} if $e = \dim X$. We set $\omega_{X} := \mathcal H^{-e}(\omega^{\mydot}_X)$. The sheaf $\omega_X$ is called a \emph{dualising sheaf} of $X$, 
and it is known that  
it is $S_2$ \stacksproj{0AWE}.

Throughout the article, we always assume that the dualising complex $\omega^{\mydot}_B$ on the base scheme $B$ is normalised. 
In this case, given a proper morphism $\pi \colon X \to B$ of irreducible excellent schemes 
such that $\pi$ is surjective or $\dim B = \dim \MO_{B, b}$ for every closed point $b \in B$, 
 the dualizing complex $\omega^{\mydot}_X = \pi^!\omega_B^{\mydot}$ is also  normalised \cite[Subsubsection 2.1.1]{KTTWYY1}. 

Assume that $X$ is integral and normal.  
Then the dualising sheaf $\omega_X$ is invertible on the regular locus of $X$    \stacksproj{0AWX}. 
We fix a Weil divisor $K_X$ such that $\cO_X(K_X) \simeq \omega_X$. 
Any such Weil divisor is called a {\em canonical divisor} on $X$.  If $K_X$ and $K'_X$ are canonical divisors on $X$, 
then there exists an invertible sheaf $L$ on $B$ such that $\MO_X(K_X) \simeq \MO_X(K'_X) \otimes \pi^*L$ \cite[Theorem V.3.1]{Har66}.

\subsubsection{Singularities of minimal model program}\label{ss-mmp-sing}

We will freely use the standard notation in birational geometry, for which we refer to \cite{kollar13} and \cite{KM98}. 
In particular, we will use the abbreviated names for singularities such as klt or plt. 
We say that $(X,\Delta)$ is a {\em log pair} if $X$ is an integral  normal excellent scheme admitting a dualising complex and $\Delta$ is an effective $\bQ$-divisor such that $K_X+\Delta$ is $\bQ$-Cartier. 
In contrast to \cite[Definition 2.8]{kollar13}, we always require that $\Delta$ is effective, and so a klt  pair $(X,\Delta)$ is always a log pair.


\begin{remark} 
 In \cite{kollar13}, it is assumed that all schemes $X$ are of finite type over a regular base scheme $B$. This assumption was only introduced so that $B$, and so $X$,  admit a dualising complex, but the results of \cite{kollar13} go through for general excellent schemes admitting dualising complexes. 
We {point out} that 
the discrepancy $a(E, X, \Delta)$ for a log pair $(X,\Delta)$ and a prime divisor $E$ over $X$ does not depend on the choice of the dualising complex $\omega^{\mydot}_X$.
\end{remark}

\subsection{$S_2$ sheaves}\label{ss-S2}

Let $Y$ be an irreducible excellent $S_2$ scheme admitting a dualising complex $\omega_Y^{\mydot}$. 
Let $\cF$ be a coherent sheaf with $\Supp \cF=Y$. 
We say that  $\cF$ is $S_r$ if $\depth \cF_y \geq \min \{ r, \dim \MO_{Y, y}\}$ 
(note that we have $\dim \cF_y = \dim \MO_{Y, y}$ by $\Supp \cF =Y$). 
It is well known that (1-a) and (1-b) below are equivalent. 
\begin{enumerate}
\item[(1-a)] $\cF$ is $S_1$. 
\item[(1-b)] 
If $U$ and $V$ are non-empty open subsets satisfying $U \supseteq V$, 
then the restriction map $\cF(U) \to \cF(V)$ is injective. 
\end{enumerate}
Moreover,  (2-a) and (2-b) below are equivalent. 
\begin{enumerate}
\item[(2-a)] $\cF$ is $S_2$. 
\item[(2-b)] 
$\cF$ is $S_1$, and 
if $U$ and $V$ are non-empty open subsets satisfying $U \supseteq V$ and 
$\codim_U (U \setminus V) \geq 2$, 
then the restriction map $\cF(U) \to \cF(V)$ is bijective. 
\end{enumerate}
For the dualising sheaf $\omega_Y$ on $Y$ and 
an arbitrary coherent sheaf $\cF$ with $\Supp \cF = Y$, 
we set 
\[
\cF^{**} := \cHom_{\MO_Y}(\cHom_{\MO_Y}(\cF, \omega_Y), \omega_Y)  
\]
which  is called the {\em $S_2$-hull} (or $S_2$-fication) of $\cF$ \cite[Theorem 2]{Kol22}.

\begin{remark} {{In the case of Noetherian integral normal schemes}, we could have written
\[
\cF^{**} := \cHom_{\MO_Y}(\cHom_{\MO_Y}(\cF, \MO_Y), \MO_Y).  
\]
The reason we stated the above definition in this way, is that often for us we will set $Y := W_n X$ for {an integral normal excellent $\F_p$-}scheme $X$. 
In this case, it is more natural to define duality by mapping to $\omega_Y = W_n\omega_X$, or, as will be the case in this paper, to $W_n\omega_X(-K_X)$.}
\end{remark}

For a coherent sheaf $\cF$ with $\Supp \cF =Y$, 
the following properties hold for its $S_2$-hull $\cF^{**}$ {(cf.\ \cite[Definition 13]{Kol22})}. 
\begin{enumerate}
\item $\cF^{**}$ is $S_2$ (Lemma \ref{l-hom-S2}). 
\item We have the natural $\MO_Y$-module homomorphism $\theta: \cF \to \cF^{**}$. 
\item $\cF$ is $S_1$ if and only if $\theta : \cF \to \cF^{**}$ is injective. 
\item $\codim_Y(\cF^{**}/\theta(\cF)) \geq 2$. 
\end{enumerate}



{
\begin{lem}\label{l-hom-S2}
Let $Y$ be an irreducible excellent $S_2$ scheme. 
Take a coherent sheaf $\cF$ and an $S_2$ coherent sheaf $\cG$ satisfying $\Supp \cF =\Supp \cG = Y$. 
Then $\cHom_{\MO_Y}(\cF, \cG)$ is $S_2$. 
\end{lem}
\begin{proof}
Since $\cG$ is $S_1$, it is easy to see that $\cHom_{\MO_Y}(\cF, \cG)$ is $S_1$ as well. Thus it is enough to show that $\cHom_{\cO_Y}(\cF, \cG)(U) \to \cHom_{\cO_Y}(\cF, \cG)(V)$ is surjective for every non-empty open subsets $V \subseteq U \subseteq Y$ such that the complement of the inclusion $i \colon V \hookrightarrow U$ is of codimension {$\geq 2$}.

To this end, take $\varphi_V \in \cHom_{\cO_Y}(\cF, \cG)(V)$. Then the composition
\[
\varphi_U : \cF|_U  \to i_*(\cF|_V) \xrightarrow{i_*{\varphi_V}} i_*(\cG|_V)  \xleftarrow{\simeq} \cG|_U
\]
gives an element lifting $\varphi_V$ under $\cHom_{\cO_Y}(\cF, \cG)(U) \to \cHom_{\cO_Y}(\cF, \cG)(V)$, proving its surjectivity. Here, the isomorphism $i_*(\cG|_V)  \xleftarrow{\simeq} \cG|_U$ follows from $\cG$ being $S_2$. 
\end{proof}}

\subsection{Generalities on restrictions}

\begin{definition} \label{definition:restricting-divisor}
Let $X$ be an integral normal Noetherian scheme, let $S$ be a prime divisor, and let $D$ be a $\bQ$-Cartier $\bQ$-divisor such that $S \not\subseteq \Supp D$. 
For a positive integer $m$ such that $mD$ is Cartier, we define 
\[
D|_S := \frac{1}{m} ( (mD)|_S), 
\]
which is a $\Q$-Cartier $\Q$-divisor on $S$. Note that $D|_S$ is independent of the choice of $m$. 
\end{definition}

\begin{proposition} \label{prop:basic_restriction}
Let $X$ be an integral normal divisorially Cohen-Macaulay excellent scheme admitting a dualising complex and let $S$ be a normal prime divisor. Let $D$ be a $\bQ$-Cartier  $\bQ$-divisor 
such that $S \not\subseteq \Supp D$ and 
$(X,S+\{D\})$ is plt. Then 
there exists a unique $\MO_X$-module homomorphism 
\[
\res \colon \cO_X(\lfloor D \rfloor) 
\to \cO_S(\lfloor D|_S \rfloor) 
\]
such that $\res|_{X\, \setminus\, \Supp D}$ coincides with 
the restriction homomorphism $\rho \colon \MO_{X \,\setminus\, \Supp D} \to \MO_{S\, \setminus\, \Supp D}$. 
Moreover, $\res$ is surjective. 
\end{proposition}

\begin{proof}
See \cite[Proposition 2.29]{KTTWYY1}.    
\end{proof}

\subsection{Matlis duality} \label{ss:Matlis-duality}


In this subsection, we recall the foundations of Matlis duality and refer to \stacksproj{08XG} for details. 

Let $(R, \m)$ be a 
local Noetherian ring admitting a normalised dualising complex $\omega^{\mydot}_R$. 
Let $E$ be the injective hull of the residue field $R/\m$ over $R$ 
(the reader should be aware that $E$ is an injective Artinian $R$-module). 
The operation $(-)^{\vee} := \Hom_R(-, E)$ is called \emph{Matlis duality}. 
Note that $(-)^\vee$ is exact, and so this operation extends canonically to the corresponding derived categories; 
we will denote such an extension by the same symbol $(-)^{\vee}$.
If $(R, \m)$ is complete, then it is an anti-equivalence between the  category of Noetherian $R$-modules and the category of Artinian $R$-modules. 

{A} key property of Matlis duality is that {it turns the local cohomology of a complex into its Grothendieck dual} 
as indicated by the following result. This is a generalisation of Serre duality to the relative setting. {In what follows, $(-)^{\wedge}$ denotes the derived $\m$-adic completion ({see \stacksproj{0922} for the definition in the Noetherian case, cf.\ \stacksproj{091V}}). 

\begin{proposition}[{cf.\ \cite[Tag 0A84 and 0AAK]{stacks-project}}] \label{prop:Matlis-duality}
Let $(R, \m)$ be a local Noetherian ring and 
let $K \in D^b_{\mathrm{coh}}(R)$. Then
\[
    R\Hom_R(K, \omega^{\mydot}_R)^{\wedge} \simeq R\Gamma_{\m}(K)^{\vee}.
\]
In particular, if $K$ is {a finitely generated} $R$-module, then
\[
 \Ext^{-i}_R(K, \omega^{\mydot}_R)^{\wedge} \simeq H^{i}_{\m}(K)^{\vee}.
\]
\end{proposition}
\noindent 
{In the second isomorphism, $(-)^{\wedge}$ denotes the usual $\m$-adic completion, 
but there is no ambiguity here as the derived {$\m$-adic} completion
agrees with the usual {$\m$-adic} completion for {finitely generated} modules over Noetherian rings by \cite[Tag 0EEU(2)]{stacks-project}.} 

In order to make proofs easier to read, we introduce the following notation used throughout the article.
\begin{notation} \label{notation:global-local-cohomology}
Given $i \in \Z$, 
a {Noetherian} local ring $(R,\m)$, a proper morphism of schemes $g \colon Y \to \Spec R$, and a coherent sheaf $\cF$ on $Y$, we define 
\[
H^i_{\m}(Y,\cF) := H^iR\Gamma_{\m}R\Gamma({Y},\cF), 
\]
which is an Artinian $R$-module. 
Here, $R\Gamma({Y},\cF)$ may be identified with  $Rg_*\cF$.
\end{notation}

Given a short exact sequence of coherent sheaves on $Y$, 
the functors $H^i_{\m}(Y, -)$  induce 
a long exact sequence of cohomologies. {Note that} if $Y=\Spec R$ {and $g = {\rm id}$}, then $H^i_{\m}(Y,\cF) $ agrees with the usual local cohomology, and if $R=k$ is a field and $\m = (0)$, then $H^i_{\m}(Y,\cF)  = H^i(Y,\cF)$ is the coherent cohomology.

{The following lemma will be useful in translating problems concerning global sections into problems concerning local cohomology.}


{\begin{lemma}  \label{lem:Matlis-duality-for-highest-cohomology}
Let $(R, \m)$ be {a Noetherian local ring admitting a dualising complex such that $\Spec R$ is irreducible}. 
Let $X$ be {a $d$-dimensional} irreducible scheme  which is proper over $\Spec R$ and let $\cF$ be a coherent sheaf on $X$. Then
\[
{H^d_{\m}(X,\cF)^{\vee}} \simeq \Hom_{\cO_X}(\cF, \omega_X)^{\wedge}.
\]
Moreover, if $X$ is Cohen-Macaulay, then
\[
{H^{i}_{\m}(X,\cF)^{\vee}} \simeq \Ext^{d-i}_{\cO_X}(\cF, \omega_X)^{\wedge}.
\]
\end{lemma}

\begin{proof}
Consider the following isomorphisms
\begin{alignat*}{1}
\big(H^iR\Gamma_{\m}R\Gamma(X,\cF)\big)^{\vee} &\simeq  \big( H^{-i}R\Hom_R(R\Gamma(X,\cF),\omega^{\mydot}_R)\big)^{\wedge} \\
&\simeq \big(H^{-i}R\Gamma \circ  R\cHom_{\cO_X}(\cF,\omega^{\mydot}_X)\big)^{\wedge} \simeq \big(H^{-i}R\Hom_{\cO_X}(\cF,     \omega_X^{\mydot})\big)^{\wedge},
\end{alignat*}
where the first one {holds by} Matlis duality (Proposition \ref{prop:Matlis-duality}) and \cite[Tag 0A06]{stacks-project}, whilst the second one is Grothendieck duality. 

If $X$ is Cohen-Macaulay, then $\omega_X^\mydot=\omega_X[d]$, and so
\[
{H^{i}_{\m}(X,\cF)^{\vee}} \simeq \Ext^{d-i}_{\cO_X}(\cF, \omega_X)^{\wedge}.
\]

For $i=d$, but without assuming Cohen-Macauliness, the higher cohomologies of $\omega_X^{\mydot}$ do not interfere, and so the same conclusion holds (see \cite[Lemma 2.25]{KTTWYY1} or the proof of \cite[Lemma 2.2]{BMPSTWW20}).
\qedhere

\end{proof}

\begin{lem}\label{l-easy-Serre}
Let $(R, \m)$ be {a Noetherian local ring admitting a dualising complex such that $\Spec R$ is irreducible}. 
Let $X$ be {a $d$-dimensional} irreducible scheme  which is proper over $\Spec R$, let $\cF$ be a coherent sheaf on $X$, and let $A$ be an ample Cartier divisor.

Suppose that $X$ and 
{$\cF$} are Cohen-Macaulay. 
Then there exists $s_0 \in \bZ_{>0}$ such that
\[
H^i_{\m}(X, \cF(-sA))=0
\]
for all integers $i$ and $s$ satisfying $i < d$ and $s \geq s_0$, 
{where $\cF(-sA) := \cF \otimes_{\MO_X} \MO_X(-sA)$.} 
\end{lem}

\begin{proof}
Let $s \in \bZ_{>0}$ and $i<d$. 
By Lemma \ref{lem:Matlis-duality-for-highest-cohomology},
\[
H^i_{\m}(X, \cF(-sA))^\vee \simeq \Ext^{d-i}_{\cO_X}(\cF(-sA), \omega_X)^{\wedge}.
\]
Moreover,
\[
\cExt^{d-i}_{\cO_X}(\cF(-sA), \omega_X)^{\wedge} 
{=} 0
\]
for $i<d$, as $\cF(-sA)$ is Cohen-Macaulay (see, e.g.\ Proposition \ref{prop:Matlis-duality}). Therefore, by the Leray spectral sequence for $R\Hom = R\Gamma \circ R\cHom$, we get that
\[
\Ext^{d-i}_{\cO_X}(\cF(-sA), \omega_X)^{\wedge} \simeq H^{d-i}(X, \cHom_{{\MO_X}}(\cF(-sA), \omega_X))^{\wedge}, 
\]
which is zero for $s \gg 0$ by {the} Serre vanishing {theorem}.
\end{proof}
}

\subsection{Witt dualising sheaves and $W_n\omega_X(-K_X)$}\label{ss-Witt-dualising}

Let $X$ be an integral normal {Noetherian} 
$F$-finite 
$\F_p$-scheme. Fix $n \in \Z_{>0}$. 
Recall that the ringed  space $W_nX := (X, W_n\MO_X)$ is a Noetherian scheme over $\Z/p^n\Z$.  
Then $W_n\omega^{\mydot}_X$  denotes the dualising complex (if it exists) and 
$W_n\omega_X$ denotes the dualising sheaf, which is an $S_2$ coherent $W_n\MO_X$-module (cf.\ Subsection \ref{sss-dualising}). 
Throughout this paper, we shall basically  work under the following setting. 


\begin{setting}[{general case}] \label{setting:most-general-foundations-of-log-quasi-F-splitting}
Let $X$ be a $d$-dimensional 
integral normal {Noetherian} 
$F$-finite 
$\F_p$-scheme such that $W_nX$ is an excellent scheme admitting a dualising complex $W_n\omega^{\mydot}_X$ for every $n \in \Z_{>0}$. 
\end{setting}

{In some cases, we might need to assume that our scheme is projective over a fixed affine base as described by Setting \ref{setting:foundations-of-log-quasi-F-splitting} to which we will refer whenever necessary.} 

\begin{setting}[{projective case}] \label{setting:foundations-of-log-quasi-F-splitting}
Let $R$ be {an $F$-finite Noetherian domain of characteristic $p>0$} such that 
$W_nR$ is an excellent ring admitting a dualising complex for every $n \in \Z_{>0}$. 
Let $X$ be a $d$-dimensional integral normal scheme which is projective over $\Spec R$. 
\end{setting}

\begin{remark}
{In the situation of} Setting \ref{setting:most-general-foundations-of-log-quasi-F-splitting}, 
it {automatically} holds that 
 $d = \dim W_nX=\dim X <\infty$ (\cite[Ch.\ V, Corollary 7.2]{Har66}). 
{In the situation of} Setting \ref{setting:foundations-of-log-quasi-F-splitting}, 
we have that $\dim W_nR = \dim R<\infty$. 
\end{remark}

\begin{rem}\label{r-dcpx-in-application}
If $S$ is a complete $F$-finite Noetherian local $\F_p$-algebra and $R$ is a ring essentially of finite type over $S$, 
then $W_nR$ is an excellent ring admitting a dualising complex for every $n \in \Z_{>0}$. 
Indeed, $W_nS$ is a complete Noetherian local ring and $W_nR$ is a ring 
essentially of finite type over $W_nS$ 
{(cf.\ \cite[Scholie 7.8.3(iii)]{EGAIV2}, \cite[Chapter V, \S 10]{Har66})}. 
\end{rem}


\begin{rem}\label{r-n-vs-n-1-omega}
\begin{enumerate}
\item 
In the situation of Setting \ref{setting:most-general-foundations-of-log-quasi-F-splitting}, 
whenever we fix $n$, 
we always assume that the dualising complexes on $W_mX$ and $W_nX$ are compatible for every $m <n$, i.e., 
we define $W_m\omega_X^{\mydot} := i^!W_n\omega_X^{\mydot}$ 
for the induced closed immersion $i : W_mX \hookrightarrow W_nX$. 
The authors do not know whether we can assume this compatibility condition 
without fixing $n$. 
\item 
In the situation of Setting \ref{setting:foundations-of-log-quasi-F-splitting}, 
we consider $\Spec W_nR\ (=W_n(\Spec R))$ as the base scheme of $W_nX$. 
As in (1), when we fix $n \in \Z_{>0}$, 
we always assume that the dualising complexes on $\Spec W_mR$ and $\Spec W_nR$ are compatible for every $m <n$. 
In this case, $W_m\omega_X^{\mydot}$ and $W_n\omega_X^{\mydot}$ 
are also compatible by 
$W_m\omega_X^{\mydot} = \pi_m^!(\omega_{\Spec W_mR}^{\mydot})$ and 
$W_n\omega_X^{\mydot} = \pi_n^!(\omega_{\Spec W_nR}^{\mydot})$ 
(cf.\ Subsubsection \ref{sss-dualising}), 
where $\pi_m : W_mX \to \Spec W_mR$ and $\pi_n : W_nX \to \Spec W_nR$ denote the induced morphisms. 
\end{enumerate}
\end{rem}

In the situation of Setting \ref{setting:most-general-foundations-of-log-quasi-F-splitting}, 
the $S_2$-hull $\mcF^{**}$ of a coherent $W_n\MO_X$-module $\mcF$ is given as follows (Subsection \ref{ss-S2}): 
\begin{eqnarray*}
\mcF^{**}  := \cHom_{W_n\MO_X}(\cHom_{W_n\MO_X}(\mcF, W_n\omega_X), W_n\omega_X).  
\end{eqnarray*}
For a Weil divisor $D$ on $X$, 
we define $W_n\omega_X(D)$ as the $S_2$-hull 
of $W_n\omega_X \otimes_{W_n\MO_X} W_n\MO_X(D)$, 
where $W_n\MO_X(D)$ denotes the Witt divisorial sheaf introduced in \cite[Subsection 3.1]{tanaka22}. 
Hence 
\begin{eqnarray*}
W_n\omega_X(D) &=&  
\cHom_{W_n\MO_X}(\cHom_{W_n\MO_X}(W_n\omega_X(D), W_n\omega_X), W_n\omega_X)\\
&\simeq& j_*(W_n\omega_{X_{\reg}} \otimes_{W_n\MO_{X_{\reg}}} W_n{\MO_{X_{\reg}}}(D|_{X_{\reg}})), 
\end{eqnarray*}
where $j : X_{\reg} \hookrightarrow X$ denotes the open immersion  from the regular locus $X_{\reg}$ of $X$. 
By construction, $W_n\omega_X(D)$ is $S_2$. 
We shall frequently use the case when $D = -K_X$, i.e., $W_n\omega_X(-K_X)$. 
For a coherent $W_n\MO_X$-module $\cF$, we set $\cF^* := \cHom_{W_n\MO_X}(\cF, W_n\omega_X(-K_X))$, which is an $S_2$ coherent $W_n\MO_X$-module 
(Lemma \ref{l-hom-S2}). 
This notation $\cF^*$ is compatible with the notation 
$\cF^{**}$ introduced earlier up to isomorphism, because $\cF^{**} \simeq (\cF^*)^*$ holds as follows: 
\begin{eqnarray*}
\cF^{**} 
&=& \cHom_{W_n\MO_X}(\cHom_{W_n\MO_X}(\cF, W_n\omega_X), W_n\omega_X)\\
&\simeq & \cHom_{W_n\MO_X}(\cHom_{W_n\MO_X}(\cF, W_n\omega_X(-K_X)), W_n\omega_X(-K_X))\\
&=& (\cF^*)^*, 
\end{eqnarray*}
where the isomorphism can be checked by restricting to the regular locus $X_{\reg}$ of $X$. 

\begin{rem}
Take a coherent $W_n\MO_X$-module $\cF$ with $\Supp F=X$. 
For $n \leq m$ and the induced closed immersion $i : W_nX \hookrightarrow W_mX$, we have that 
$\cF$ is an $S_2$ coherent $W_n\MO_X$-module if and only if 
$i_*\cF$ is an $S_2$ coherent $W_m\MO_X$-module. 
\end{rem}



In the situation of Setting \ref{setting:most-general-foundations-of-log-quasi-F-splitting}, we now recall the definition of 
$\underline{p}: W_{n-1}\omega_X \to W_{n}\omega_X$. 
Since $\Ker( p :W_{n}\MO_X \to W_{n}\MO_X) = \Ker(R: W_{n}\MO_X \to W_{n-1}\MO_X)$, 
we obtain the following decomposition 
\begin{equation}\label{e-def-pbar-WO}
p:W_{n}\MO_X \xrightarrow{R} W_{n-1}\MO_X \xrightarrow{\underline{p}} W_{n}\MO_X, 
\end{equation}
where $R: W_{n}\MO_X \to W_{n-1}\MO_X$ is surjective and 
$\underline p:  W_{n-1}\MO_X \to W_{n}\MO_X$ is injective. 

We now apply $\cHom_{W_{n}\MO_X}(-, W_n\omega_X)$ to (\ref{e-def-pbar-WO}). 
We have $\cHom_{W_{n}\MO_X}(W_n\MO_X, W_n\omega_X) \simeq W_n\omega_X$ and 
\begin{eqnarray*}
\cHom_{W_n\MO_X}(W_{n-1}\MO_X, W_n\omega_X) 
&\simeq& 
j_* \cHom_{W_n\MO_{X_{\reg}}}(W_{n-1}\MO_{X_{\reg}}, W_n\omega_{X_{\reg}})\\
&\simeq& j_* \cHom_{W_{n-1}\MO_{X_{\reg}}}(W_{n-1}\MO_{X_{\reg}}, W_{n-1}\omega_{X_{\reg}})\\
&\simeq& \cHom_{W_{n-1}\MO_{X}}(W_{n-1}\MO_X, W_{n-1}\omega_{X})\\
&\simeq& W_{n-1}\omega_X,
\end{eqnarray*}
where $j:X_{\reg} \hookrightarrow X$ denotes the open immersion from the regular locus $X_{\reg}$ of $X$ 
and the second isomorphism follows from the Grothendieck duality 
(cf.\ Remark \ref{r-n-vs-n-1-omega}). 
Thus we obtain a composition of $W_n\MO_X$-module homomorphisms: 
\begin{equation}\label{e-def-pbar-Womega}
p : W_{n}\omega_X \xrightarrow{R \,:=\, \underline{p}^*} 
W_{n-1}\omega_X \xrightarrow{\underline{p} \,:=\, R^*} W_{n}\omega_X, 
\end{equation}
where $\underline{p}^* := \cHom_{W_n\MO_X}(\underline p, W_n\omega_X)$ and $R^*:= \cHom_{W_n\MO_X}(R, W_n\omega_X)$.

Fix a Weil divisor $D$ on $X$. 
We naturally get  the following logarithmic version of (\ref{e-def-pbar-Womega}): 
\begin{equation}\label{e-def-pbar-Womega-log}
p : W_{n}\omega_X(D) \xrightarrow{R} 
W_{n-1}\omega_X(D) \xrightarrow{\underline{p}} W_{n}\omega_X(D). 
\end{equation}
More specifically, if $X$ is regular, then we obtain (\ref{e-def-pbar-Womega-log})  from (\ref{e-def-pbar-Womega}) by taking the tensor product with the invertible $W_n\MO_X$-module $W_n\MO_X(D)$ \cite[Proposition 3.12]{tanaka22}. 
The general case is reduced to this case, because (\ref{e-def-pbar-Womega-log}) is obtained by applying $j_*$ to 
$p : W_{n}\omega_{X_{\reg}}(D|_{X_{\reg}}) \xrightarrow{R} 
W_{n-1}\omega_{X_{\reg}}(D|_{X_{\reg}}) \xrightarrow{\underline{p}} W_{n}\omega_{X_{\reg}}(D|_{X_{\reg}})$, 
where $j : X_{\reg} \hookrightarrow X$ is the open immersion from the regular locus $X_{\reg}$. 
Note that $\underline{p} : W_{n-1}\omega_X(D) \to W_n\omega_X(D)$ is injective by definition. 


\begin{rem}\label{r-Illusie-def}
Assume that $X$ is a smooth variety over a perfect field $k$ of characteristic $p>0$. 
In this case, we have $W_n\omega_X \simeq W_n\Omega_X^N$ for $N := \dim X$ \cite[Theorem 4.1 in page 197]{Eke84}. 
Via this isomorphism, (\ref{e-def-pbar-Womega}) coincides with the sequence given in \cite[Ch. I,  Proposition 3.4]{illusie_de_rham_witt}
\[
p : W_n\Omega_X^N \xrightarrow{R} W_{n-1}\Omega_X^N \xrightarrow{\underline{p}} W_n\Omega_X^N,  
\]
because each of them coincides with  the unique decomposition of the $p$-multiplication map 
into a surjection followed by an injection (cf.\ Proposition \ref{p-non-surje-R}(1)). 
\end{rem}

\begin{rem}\label{r-non-surje-R}
As pointed out in Remark \ref{r-Illusie-def}, 
$R: W_n\omega_X \to W_{n-1}\omega_X$ is surjective when 
$X$ is a smooth variety over a perfect field  of characteristic $p>0$. 
On the other hand, 
$R: W_n\omega_X \to W_{n-1}\omega_X$ is not necessarily surjective in general, whilst $\underline{p}: W_{n-1}\omega_X \to W_n\omega_X$ is always injective. 
Indeed, if $X$ is a {Gorenstein} affine normal surface which is not $F$-split but $n$-quasi-$F$-split for some 
integer $n\geq 2$, then $R: W_{n}\omega_X \to W_{n-1}\omega_X$ is not surjective (Proposition \ref{p-non-surje-R}). 
\end{rem}

\begin{prop}\label{p-non-surje-R}
In the situation of Setting \ref{setting:most-general-foundations-of-log-quasi-F-splitting}, 
assume that $X$ is affine, {Gorenstein}, and $n$-quasi-$F$-split for an integer $n \geq 2$. 
Consider the following conditions 
\begin{enumerate}
\item[(a)] $X$ is $F$-split. 
\item[(b)] $R: W_{n}\omega_X \to W_{n-1}\omega_X$ is surjective. 
\end{enumerate}
Then the following hold. 
\begin{enumerate}
    \item If {\rm (a)} holds, then {\rm (b)} holds. 
    \item If {\rm (b)} holds and $\dim X =2$, then {\rm (a)} holds. 
\end{enumerate}
\end{prop}

\begin{proof}
Since the problem is local, we may assume that $X=\Spec R$ for a local ring $(R, \m)$. 
In what follows, we set $H^i_{\m}(\cF) := H^i_{\m}( \Gamma(X, \cF))$ for a coherent $W_n\MO_X$-module $\cF$. 
Set $d := \dim X = \dim R$. 
The problem is reduced to the case when  $d \geq 2$, 
as otherwise the assertions are clear. 

For $Q_{X, n} :=Q_{X, 0, n} =F_*(W_n\MO_X/pW_n\MO_X)$ and $B_{X, n} := B_{X, 0, n} 
={\rm Coker}(F : W_n\MO_X \to F_*W_n\MO_X)$, 
we have the following exact sequences 
\cite[Remark 3.4]{KTTWYY1}: 
\begin{equation}\label{e1-non-surje-R}
0 \to \MO_X \xrightarrow{\Phi_{X, n}} Q_{X, n} \to B_{X, n} \to 0, 
\end{equation}
\begin{equation}\label{e2-non-surje-R}
0 \to F_*W_{n-1}\MO_X \xrightarrow{\underline{p}} F_*W_n\MO_X \to Q_{X, n}  \to 0. 
\end{equation}
Since $X$ is Cohen-Macaulay, we have $H_{\m}^{d-1}(\MO_X)=0$. 
As $X$ is $n$-quasi-$F$-split, 
$H^d_{\m}(\Phi_{X, n}) : H^d_{\m}(\MO_X) \to H^d_{\m}(Q_{X, n})$ is injective \cite[Lemma 3.13]{KTTWYY1}. 
Then these, together with (\ref{e1-non-surje-R}), imply 
\begin{equation}\label{e3-non-surje-R}
H^{d-1}_{\m}(Q_{X, n}) \simeq H^{d-1}_{\m}(B_{X, n}).
\end{equation}
By applying $\mathcal Hom_{W_n\MO_X}(-, W_n\omega_X)$ to (\ref{e2-non-surje-R}),  we obtain another exact sequence 
{\small 
\[
F_*W_n\omega_X \xrightarrow{\underline{p}^* = R} F_*W_{n-1}\omega_X \to 
\cExt^1_{W_n\MO_X}(Q_{X, n}, W_n\omega_X) \to \cExt^1_{W_n\MO_X}(F_*W_n\MO_X, W_n\omega_X)=0. 
\]
}
Since $Q_{X, n}$ is naturally a coherent $\MO_X$-module \cite[Proposition 3.6]{KTTWYY1}, Grothendieck duality yields  
$\cExt^1_{W_n\MO_X}(Q_{X, n}, W_n\omega_X) \simeq \cExt^1_{\MO_X}(Q_{X, n}, \omega_X)$.  
Therefore, we obtain the equivalence (b) $\Leftrightarrow$ (b') for the condition (b') below: 
\vspace{0.3em}
\begin{enumerate}
\item[(b')] 
$\cExt^1_{\MO_X}(Q_{X, n}, \omega_X)=0$. 
\end{enumerate}
\vspace{0.3em}
Moreover, we get (b) $\Leftrightarrow$ (b') $\Leftrightarrow$ (b'') by applying Matlis duality (Proposition \ref{prop:Matlis-duality}): 
\vspace{0.3em}
\begin{enumerate}
\item[(b'')] 
$H^{d-1}_{\m}( Q_{X, n}) = 0$. 
\end{enumerate}
\vspace{0.3em}
{Since $X$ is Gorenstein,} (a) is equivalent to the following condition (a')
by the exact sequence $0 \to \MO_X \to F_*\MO_X \to B_{X, 1} \to 0$. 
\vspace{0.3em}
\begin{enumerate}
\item[(a')] 
$H^{d-1}_{\m}(B_{X, 1})=0$. 
\end{enumerate}
\vspace{1em}

Let us show (1). 
Assume (a). 
Since  (a') holds, we obtain $H^{d-1}_{\m}(B_{X, n})=0$ by the exact sequence 
\cite[(3.8.2)]{KTTWYY1}: 
\begin{equation}\label{e4-non-surje-R}
0 \to F_*B_{X, n-1} \to B_{X, n} \to B_{X, 1} \to 0. 
\end{equation}
Then (\ref{e3-non-surje-R}) implies (b''). 
This completes the proof of (1). 

Let us show (2). 
Assume $d=\dim X=2$ and (b). 
It follows from (b'') and (\ref{e3-non-surje-R}) that $H^1_{\m}(B_{X, n})=0$. {Observe that $H^0_{\m}(B_{X, 1})=0$ as $H^0_{\m}(F_*\cO_X)=0$ and $H^1_{\m}(\cO_X)=0$.} Thus, 
by (\ref{e4-non-surje-R}), 
we obtain $H^1_{\m}(B_{X, n-1})=0$. 
Applying this argument repeatedly, 
we obtain $H^1_{\m}(B_{X, i})=0$ for every $1 \leq i \leq n$ 
by descending induction on $i$. 
Thus (a') holds. 
This completes the proof of (2). 
\end{proof}

\section{Iterative quasi-F-splitting}\label{s-iterative-QFS}


\subsection{Definition of quasi-$F^e$-splitting}\label{ss-log-QFS}

{Let $X$ be as in Setting \ref{setting:most-general-foundations-of-log-quasi-F-splitting}.} 
Take $e \in \Z_{>0}$ and let $\Delta$ be a (non-necessarily effective) $\bQ$-divisor on $X$. We define a $W_n\MO_X$-module $Q^e_{X,\Delta,n}$ and 
a $W_n\MO_X$-module homomorphism $\Phi^e_{X, \Delta, n}$ by the following pushout diagram of $W_n\MO_X$-module homomorphisms: 
\begin{equation} \label{diagram:quasi-F-split-definition}
\begin{tikzcd}
W_n\cO_X(\Delta) \arrow{r}{F^e} \arrow{d}{R^{n-1}} & F^e_* W_n \cO_X(p^e\Delta)  \arrow{d}\\
\cO_X(\Delta) \arrow{r}{\Phi^e_{X, \Delta, n}}&
\arrow[lu, phantom, "\usebox\pushoutdr" , very near start, yshift=0em, xshift=0.6em, color=black] Q^e_{X, \Delta, n}.
\end{tikzcd}
\end{equation}
For the definition of the Witt divisorial sheaf $W_n\MO_X(\Delta)$, 
we refer to \cite[Subsection 3.1]{tanaka22}. 
We remind the reader that 
$\MO_X(\Delta) = \MO_X(\rdown{\Delta})$, 
but 
it is \emph{not} true in general that $W_n\cO_X(\Delta) = W_n\cO_X(\rdown{\Delta})$.  

We define a $W_n\cO_X$-module $B^e_{X, \Delta, n}$ by 
\begin{align}
\label{eq:definition-of-log-B} B^e_{X, \Delta, n} &:= {\rm Coker}( W_n\MO_X(\Delta) \xrightarrow{F^e} F^e_*W_n\MO_X(p^e\Delta))\\
&\hphantom{:}= F^e_*W_n\MO_X(p^e\Delta)/F^e(W_n\MO_X(\Delta)). \nonumber
\end{align}

\begin{remark}\label{r-Q-B-big-diag}
{The key properties of the construction of $Q^e_{X,\Delta,n}$,  $B^e_{X,\Delta,n}$, and $\Phi^e_{X,\Delta,n}$ may be encapsulated by the following diagram}
\begin{equation}\label{e-big-BC-diagram}
\begin{tikzcd}
& 0 \arrow{d} & 0 \arrow{d} & & \\
&F_*W_{n-1}\MO_X(p\Delta) \arrow{d}{V} \arrow[r,dash,shift left=.1em] \arrow[r,dash,shift right=.1em] & F_*W_{n-1}\MO_X(p\Delta) \arrow{d}{F^eV =VF^e} & & \\
0 \arrow{r} & W_n\MO_X(\Delta) \arrow{r}{F^e} \arrow{d}{R^{n-1}} & F^e_*W_n\MO_X(p^e\Delta) \arrow{r} \arrow{d} & B^e_{X, \Delta, n} \arrow[d,dash,shift left=.1em] \arrow[d,dash,shift right=.1em] \arrow{r}  & 0 \\
0 \arrow{r} & \MO_X(\Delta) \arrow{r}{\Phi^e_{X, \Delta, n}} \arrow{d} &  \arrow[lu, phantom, "\usebox\pushoutdr" , very near start, yshift=0em, xshift=0.6em, color=black] Q^e_{X, \Delta, n} \arrow{r} \arrow{d} & B^e_{X, \Delta, n} \arrow{r} & 0.\\
& 0 & 0 & &
\end{tikzcd}
\end{equation}
{All  the horizontal and vertical sequences are exact, as} $F^e \colon W_n\MO_X(\Delta) \to F^e_*W_n\MO_X(p^e\Delta)$ is injective, $R^{n-1} \colon W_n\MO_X(\Delta) \to \MO_X(\Delta)$ is surjective, and 
(\ref{diagram:quasi-F-split-definition}) is a pushout diagram. 
\end{remark}

One can now check that 
\begin{equation}\label{e-C-tensor}
Q^e_{X,\Delta+D,n} \simeq Q^e_{X,\Delta,n} \otimes_{W_n\MO_X} W_n\cO_X(D)
\end{equation} 
for any Cartier divisor $D$. 
Indeed,
{$W_n\cO_X(D)$ is an invertible $W_n\MO_X$-module, and hence}  
both sides of (\ref{e-C-tensor}) are isomorphic to the pushout of 
\[
F^e_*W_n\MO_X(p^e\Delta+p^eD) \xleftarrow{F^e} W_n\MO_X(\Delta+D) \xrightarrow{R^{n-1}} \MO_X(\Delta+D).
\]

In what follows, we shall often assume that $\rdown{\Delta}=0$, which is equivalent to 
the condition that the coefficients of $\Delta$ are contained in $[0, 1)$. 
In this case, $\cO_X(\Delta) =\MO_X$.

\begin{dfn}\label{d-Psi}
In the situation of Setting  \ref{setting:most-general-foundations-of-log-quasi-F-splitting}, let $\Delta$ be a $\bQ$-divisor on $X$.  
We define a $W_n\MO_X$-module homomorphism 
\[
(\Phi^e_{X, \Delta, n})^* : (Q^e_{X, \Delta, n})^* \to \MO_X(\Delta)^*
\]
by applying $(-)^* := \mathcal Hom_{W_n\cO_X}(-, W_n\omega_X(-K_X))$ to the $W_n\cO_X$-module homomorphism $\Phi^e_{X, \Delta, n} \colon \cO_X(\Delta)  \to Q^e_{X,\Delta,n}$. 

More explicitly, $(\Phi^e_{X, \Delta, n})^*$ can be written as follows: 
\[
(\Phi^e_{X, \Delta, n})^* : 
\mathcal Hom_{W_n\MO_X}(Q^e_{X, \Delta, n}, W_n\omega_X(-K_X)) \to 
\mathcal Hom_{W_n\MO_X}(\MO_X(\Delta), W_n\omega_X(-K_X)).  
\]
Note that 
$\cO_X(\Delta)^* = \mathcal Hom_{W_n\MO_X}(\MO_X(\Delta), W_n\omega_X(-K_X)) 
\simeq \MO_X(-\rdown{\Delta})$. 
\end{dfn}
\vspace{1em}

\noindent 
For the case when $\rdown{\Delta}=0$, 
the following are equivalent. 
\begin{enumerate}
\item[(A)] $\Phi_{X, \Delta, n} : \MO_X \to Q_{X, \Delta, n}$ splits as a ${W_n}\MO_X$-module homomorphism. 
\item[(B)] $H^0(X, (\Phi_{X, \Delta, n})^*) : 
\Hom_{{W_n}\MO_X}(Q_{X, \Delta, n}, \MO_X) \to \Hom_{{W_n}\MO_X}(\MO_X, \MO_X)$ is surjective. 
\end{enumerate}
\vspace{0.5em}
\noindent On the other hand, their iterative versions are not equivalent in general. 
\begin{enumerate}
\item[(A')] $\Phi^e_{X, \Delta, n} : \MO_X \to Q^e_{X, \Delta, n}$ splits as a $W_{n}\MO_X$-module homomorphism.
\item[(B')] $H^0(X, (\Phi^e_{X, \Delta, n})^*) : 
\Hom_{{W_n}\MO_X}(Q^e_{X, \Delta, n}, W_n\omega_X(-K_X)) \to 
\Hom_{{W_n}\MO_X}(\MO_X, W_n\omega_X(-K_X))$ is surjective. 
\end{enumerate}
As Yobuko's original definition is given by (A), it is tempting to adopt (A') as the definition of quasi-$F^e$-splitting. 
However, {condition (A') is too restrictive (e.g., if (A') holds for $\Delta =0$ and $e \geq 2$, then $X$ is $F$-split (Proposition \ref{p-A'-F-split}))}, 
and so  our  definition is (B'), which is weaker  than (A'). 


\begin{dfn}\label{d-IQFS}
In the situation of Setting  \ref{setting:most-general-foundations-of-log-quasi-F-splitting}, let $\Delta$ be a $\bQ$-divisor on $X$. 
Take $n, e \in \Z_{>0}$. 
We say that $(X, \Delta)$ is {\em $n$-quasi-$F^e$-split} 
if $\rdown{\Delta}=0$ and the induced map 
\[
(\Phi_{X, \Delta, n}^e)^* \colon H^0(X,(Q^e_{X, \Delta, n})^*) \to H^0(X,\cO^*_X)
\]
is surjective, {where $\cO^*_X \simeq \cO_X$.}

We call $(X, \Delta)$ {\em quasi-$F^e$-split} if it is $n$-quasi-$F^e$-split for some $n \in \Z_{>0}$. 
We call  $(X, \Delta)$ 
 {\em $n$-quasi-$F^e$-pure} (resp.\ 
 {\em quasi-$F^e$-pure}) 
if 
 there exists an open cover $X = \bigcup_{i \in I} X_i$ such that 
 $(X_i, \Delta|_{X_i})$ is  {$n$-quasi-$F^e$-split} (resp.\ 
 {quasi-$F^e$-split}) for all $i \in I$. 
\end{dfn}

\begin{rem}\label{r-IQFS} 
By definition, we have $Q_{X, \Delta, n} = Q_{X, \Delta, n}^{1}$ (compare Definition \ref{d-IQFS} and \cite[(3.3.1)]{KTTWYY1}). 
In particular, $(X, \Delta)$ is quasi-$F$-split if and only if $(X, \Delta)$ is quasi-$F^1$-split. 
\end{rem}

\begin{dfn}\label{d-IQFS-ht}
In the situation of Setting  \ref{setting:most-general-foundations-of-log-quasi-F-splitting}, let $\Delta$ be a $\bQ$-divisor on $X$. 
For $e \in \Z_{>0}$, we define 
\[
h^e(X, \Delta) \in \Z_{>0} \cup \{ \infty \} 
\]
as follows. 
\begin{itemize}
    \item If $(X, \Delta)$ is not quasi-$F^e$-split, then we set $h^e(X, \Delta) := \infty$. 
    \item 
    If $(X, \Delta)$ is quasi-$F^e$-split, then 
    $h^e(X, \Delta)$ is defined as the smallest positive integer $n$ 
    such that $(X, \Delta)$ is $n$-quasi-$F^e$-split. 
    In this case, we say that $X$ is {\em quasi-$F^e$-split of height} $n$. 
\end{itemize}
\end{dfn}

\noindent When $1<h^e(X, \Delta)<\infty$, 
$X$ is quasi-$F^e$-split of height $n$ if and only if 
$X$ is not $(n-1)$-quasi-$F^e$-split but $n$-quasi-$F^e$-split. 




\begin{rem}\label{r-e-to-e+1}
In the situation of Setting  \ref{setting:most-general-foundations-of-log-quasi-F-splitting}, let $\Delta$ be a $\bQ$-divisor on $X$ with $\rdown{\Delta}=0$. 
Fix  $n \in \Z_{>0}$ and $e \in \Z_{>0}$. 
We have the following factorisations of $W_n\MO_X$-module homomorphisms: 
\[
\Phi^{e+1}_{X, \Delta, n} :
\MO_X \xrightarrow{\Phi^e_{X, \Delta, n}} Q^e_{X, \Delta, n} \to Q_{X, \Delta, n}^{e+1}
\]
\[
(\Phi^{e+1}_{X, \Delta, n})^* : 
(Q_{X, \Delta, n}^{e+1})^* \to 
(Q^e_{X, \Delta, n})^* \xrightarrow{(\Phi^e_{X, \Delta, n})^*}
(\MO_X)^*,
\]
where $\mathcal F^* :=  \mathcal Hom_{W_n\MO_X}( \mathcal F, W_n\omega_X(-K_X))$. 
In particular, if $(X, \Delta)$ is $n$-quasi-$F^{e+1}$-split, 
then $(X, \Delta)$ is $n$-quasi-$F^e$-split. 
In other words, $h^e(X, \Delta) \leq h^{e+1}(X, \Delta)$. 
\end{rem}

\begin{remark}\label{r-big-open-QFS}
In the situation of Setting  \ref{setting:most-general-foundations-of-log-quasi-F-splitting}, 
let $\Delta$ be an effective $\Q$-divisor on $X$ with $\rdown{\Delta}=0$, and let $X'$ be an open subset of $X$ such that ${\rm codim}_X(X \setminus X') \geq 2$. 
By Definition \ref{d-IQFS}, 
$(X, \Delta)$ is $n$-quasi-$F^e$-split {if and only if} 
$(X', \Delta|_{X'})$ is $n$-quasi-$F^e$-split. 
\end{remark}



{In what follows, {we set $Q^e_{X, n} := Q^e_{X, 0, n}$ and $\Phi^e_{X, n} := \Phi^e_{X, 0, n}$.} 
The following result shows that asking for the splitting of 
$\Phi^e_{X, n} \colon \cO_X \to Q^e_{X,n}$ for ${e \geq 2}$ gives a redundant definition.}
 
\begin{prop}\label{p-A'-F-split}
In the situation of Setting  \ref{setting:most-general-foundations-of-log-quasi-F-splitting}, 
take integers $e \geq 2$ and $n \geq 1$. 
Then the following are equivalent. 
\begin{enumerate}
\item $X$ is $F$-split, i.e., $F: \MO_X \to F_*\MO_X$ splits as an $\MO_X$-module homomorphism. 
\item $\Phi^e_{X, n} : \MO_X \to Q^e_{X, n}$ splits as a $W_{n}\MO_X$-module homomorphism.
\end{enumerate}
\end{prop}

\begin{proof}
By Remark \ref{r-Q-B-big-diag}, (1) implies (2). 
Let us show (2) $\Rightarrow$ (1). 
Assume (2). 
We then get $W_n\MO_X$-module homomorphisms: 
\[
{\rm id}: \MO_X \xrightarrow{\Phi^e_{X, n}} Q^e_{X, n} \xrightarrow{\exists \theta} \MO_X. 
\]
By $p\MO_X =0$,  the $W_n\MO_X$-module homomorphism $\theta :  Q^e_{X, n} \to \MO_X$ factors through 
\[
\frac{Q^e_{X, n}}{pQ^e_{X, n}} = 
\frac{F_*^eW_n\MO_X}{pF_*^eW_n\MO_X} = 
F_*^e\frac{W_n\MO_X}{pW_n\MO_X}. 
\]
Therefore, we obtain the following $W_n\MO_X$-module homomorphisms: 
\[
{\rm id}: \MO_X \xrightarrow{\Phi'} F_*^e(W_n\MO_X /pW_n\MO_X) \xrightarrow{\theta'} \MO_X. 
\]
It is enough to prove that $\Phi' : \MO_X \to F_*^e(W_n\MO_X /pW_n\MO_X)$ factors through $F: \MO_X \to F_*\MO_X$, which holds by the following commutative diagram
\[
\begin{tikzcd}
W_n\MO_X \arrow[r, "F"] \arrow[d, ->>, "R^{n-1}"] 
& F_*W_n\MO_X  \arrow[r, "F^{e-1}"] 
\arrow[d, ->>, "F_*R^{n-1} =:\rho"] & F_*^eW_n\MO_X\arrow[d, ->>, "\pi"] \\
\MO_X \arrow[r, "F"] \arrow[rr, bend right, "\Phi'"] & F_*\MO_X \arrow[r, dashed]  
& F_*^e(W_n\MO_X/pW_n\MO_X), 
\end{tikzcd}
\]
{where the {above} dashed arrow exists, {because 
the assumption $e \geq 2$ implies that}}
\[
F^{e-1}(\Ker(\rho)) 
=F^{e-1}(V(  F_*^2W_{n-1}\MO_X )) 
\subseteq pF_*^eW_n\MO_X \subseteq \Ker (\pi). 
\]
\qedhere 


\end{proof}

As the definition of quasi-$F^e$-splitting (Definition \ref{d-IQFS}) shows, 
we are mainly interested in the case when $\rdown{\Delta} =0$. 
However, we shall work under a more general setting, 
as it  will be necessary in our applications. 
{We have the following  cohomological criterion for whether a pair is quasi-$F^e$-split.}


\begin{remark}\label{r-loc-coh-W_nm}
In the situation of Setting  \ref{setting:foundations-of-log-quasi-F-splitting}, assume that $(R, \m)$ is a local ring. 
Recall that $(W_nR, {\sqrt{W_n\m}})$ is a local ring for every $n >0$. 
For $i \in \Z$, integers $n' \geq n >0$, and a coherent $W_n\MO_X$-module $\cF$, 
we have that $\cF$ is naturally a $W_{n'}\MO_X$-module and 
\[
H^i_{W_{n'}\m}(X, \cF) = H^i_{W_n\m}(X, \cF). 
\]
\end{remark}

\begin{lemma} \label{lem:cohomological-criterion-for-log-quasi-F-splitting}
In the situation of Setting  \ref{setting:foundations-of-log-quasi-F-splitting}, {assume that $(R, \m)$ is a local ring and} 
let $\Delta$ be a $\Q$-divisor on $X$. 
Fix $e, n \in \Z_{>0}$. 
Then the following are equivalent, 
where  $(-)^* := \cHom_{W_n\MO_X}(-, W_n\omega_X(-K_X))$. 
\begin{enumerate}
\item $H^0(X, (\Phi^e_{X, \Delta, n})^*) : 
H^0(X, (Q^e_{X, \Delta, n})^*) \to H^0(X, \MO_X(\Delta)^*)$ is surjective. 
\item 
$H^d_{\m}(X, \Phi^e_{X, K_X+\Delta, n}) \colon H^d_{\m}(X, \cO_X(K_X+\Delta)) \to H^d_{\m}(X, Q^e_{X,K_X+\Delta,n})$ is injective. 
\end{enumerate}
\end{lemma}

\noindent 
As a consequence, when $\rdown{\Delta}=0$, 
$(X,\Delta)$ is $n$-quasi-$F^e$-split if and only if 
\[
H^d_{\m}(X, \Phi^e_{X, K_X+\Delta, n}) \colon H^d_{\m}(X, \cO_X(K_X+\Delta)) \to H^d_{\m}(X, Q^e_{X,K_X+\Delta,n})
\]
is injective (cf.\ Definition \ref{d-IQFS}). 
{For the definition of $H^d_{\m}(X, -)$, see Notation \ref{notation:global-local-cohomology}.}  



\begin{proof}
Consider the following exact sequence of $W_nR$-modules: 
\begin{equation}\label{e1:cohomological-criterion-for-log-quasi-F-splitting}
0 \to K \to H^d_{\m}(X, \cO_X(K_X+\Delta)) \xrightarrow{
H^d_{\m}(X, \Phi^e_{X, K_X+\Delta, n})} H^d_{\m}(X, Q^e_{X,K_X+\Delta,n}), 
\end{equation}
where $K$ is the kernel of $H^d_{\m}(X, \Phi^e_{X, K_X+\Delta, n})$. 
By applying Matlis duality $\Hom_{W_nR}(-, E)$ (which is exact), 
Lemma \ref{lem:Matlis-duality-for-highest-cohomology} 
yields an exact sequence 
\begin{equation}\label{e2:cohomological-criterion-for-log-quasi-F-splitting}
0 \leftarrow K^{\vee} \leftarrow 
H^0(X, \MO_X(\Delta)^*)^{\wedge} 
\xleftarrow{H^0(X, (\Phi^e_{X, \Delta, n})^*)^{\wedge}}
H^0(X, (Q^e_{X, \Delta, n})^*)^{\wedge} 
\end{equation}
where $K^{\vee} = \Hom_{W_nR}(K, E)$ and 
we used the following (Lemma \ref{l-hom-S2}):  
\[
\cHom_{W_n\MO_X}(Q^e_{X, K_X+\Delta, n}, W_n\omega_X) 
\simeq 
\cHom_{W_n\MO_X}(Q^e_{X, \Delta, n}, W_n\omega_X(-K_X)) 
= (Q^e_{X, \Delta, n})^*. 
\]
{Then the following is true:}
\begin{align*}
H^0(X, (\Phi^e_{X, \Delta, n})^*)\text{ is surjective}
& \Longleftrightarrow  H^0(X, (\Phi^e_{X, \Delta, n})^*)^{\wedge}\text{ is surjective}\\
&\Longleftrightarrow  H^d_{\m}(X, \Phi^e_{X, K_X+\Delta, n})\text{ is injective}.
\end{align*}
Here the first equivalence 
holds by the fact that $(-)^{\wedge} \simeq (-) \otimes_{W(R)} \widehat{W_nR}$ 
and the completion $W_nR \to \widehat{W_nR}$ is faithfully flat 
(cf.\ \stacksproj{00MA(3)}, \stacksproj{00MC}), 
and the second  one follows from 
(\ref{e1:cohomological-criterion-for-log-quasi-F-splitting}) and 
(\ref{e2:cohomological-criterion-for-log-quasi-F-splitting}). {In the last assertion, we used that $K=0$ if and only if $K^{\vee}=0$, which is true because $K^{\vee\vee} \simeq K$} (cf. \stacksproj{08Z9}). 
\end{proof}


\begin{lemma} \label{l-BQO}
In the situation of Setting \ref{setting:most-general-foundations-of-log-quasi-F-splitting}, let $\Delta$ be a $\bQ$-divisor on $X$. 
Take $e, n \in \Z_{>0}$. 
Then we have the following exact sequence of coherent $W_{n+1}\MO_X$-modules: 
\begin{align} \label{eq:C-restriction-sequence}
0 \to 
F_*B^e_{X, p\Delta, n} \to Q^e_{X,\Delta,n+1} \to F^e_* \cO_X(p^e\Delta) \to 0. 
\end{align}
\end{lemma}

\begin{proof}
The exact sequence (\ref{eq:C-restriction-sequence}) 
is obtained by applying the snake lemma 
to the following commutative diagram in which the horizontal sequence is exact
\vspace{0.3em}
\[
\begin{tikzcd}
& F_*W_{n}\MO_X(p\Delta) \arrow{d}{F^e} \arrow[r,dash,shift left=.1em] \arrow[r,dash,shift right=.1em] & F_*W_{n}\MO_X(p\Delta) \arrow{d}{VF^e} & &\\
0 \arrow{r}&  F_*^{e+1}W_{n}\MO_X(p^{e+1}\Delta) \arrow{r}{V} & F^e_*W_{n+1}\MO_X(p^e\Delta) \arrow{r}{R^n} & F^e_*\MO_X(p^e\Delta) \arrow{r} & 0.
\end{tikzcd}
\]
\end{proof}

\begin{prop}\label{p-W_eO-str}
In the situation of Setting  \ref{setting:most-general-foundations-of-log-quasi-F-splitting}, let $\Delta$ be a $\bQ$-divisor on $X$.  
Then  $B^e_{X, \Delta, n}$ and $Q^e_{X, \Delta, n}$
are naturally coherent $W_e\MO_X$-modules. 
\end{prop}

\begin{proof}
Since $B^e_{X, \Delta, n}$ and $Q^e_{X, \Delta, n}$ are 
coherent $W_n\MO_X$-modules, 
there is nothing to show when $n \leq e$. 
In what follows, we assume 
$n >e$. 
By the following exact sequence (\ref{e-big-BC-diagram}):
\[
0 \to \MO_X(\Delta) \xrightarrow{\Phi^e_{X, \Delta, n}} Q^e_{X, \Delta, n} \to B^e_{X, \Delta, n} \to 0,
\]
it suffices  to prove that $Q^e_{X, \Delta, n}$ is a $W_e\MO_X$-module. 
By (\ref{e-big-BC-diagram}), we have 
\[
W_e\MO_X = \frac{W_n\MO_X}{V^e(F_*^eW_{n-e}\MO_X)}\qquad \text{and}\qquad 
Q^e_{X, \Delta, n} = \frac{F_*^eW_n\MO_X(p^e\Delta)}{VF^e(F_*W_{n-1}\MO_X(p\Delta))}.
\]
Hence 
it is enough  to show that 
\[
(V^eF_*^e\zeta) \cdot ( F^e_* \xi) \in  VF^e(F_*W_{n-1}\MO_X(p\Delta)(U))
\]
for $F_*^e\zeta \in F_*^eW_{n-e}\MO_X(U)$ and $F_*^e\xi \in F_*^eW_n\MO_X(p^e\Delta)(U)$, 
where $U$ is an affine open subset of $X$. 
For $\zeta =(\zeta_0, \zeta_1, ..., \zeta_{n-e-1}) 
\in W_{n-e}\MO_X(U)$, 
we set 
\[
\widetilde{\zeta} := (\zeta_0, \zeta_1, ..., \zeta_{n-e-1}, 0, ..., 0) 
\in W_{n}\MO_X(U).
\]
For $F^e : W_n\MO_X(U) \to F_*^eW_n\MO_X(U)$, we obtain 
 $F^eV^e(F_*^e{\zeta}) = p^eF_*^e\widetilde{\zeta}$. 
 It holds that 
\[
(V^eF_*^e\zeta) \cdot ( F^e_* \xi) = F^e_*( (p^e\widetilde{\zeta}) \cdot \xi) 
= p^e F^e_*( \widetilde{\zeta} \cdot \xi) 
\]
\[
\in p^e F^e_* ( \widetilde{\zeta} \cdot W_n\MO_X(p^e\Delta)(U)) 
\subset p^eF^e_*W_n\MO_X(p^e\Delta)(U) \subset 
VF^e(F_*W_{n-1}\MO_X(p\Delta)). 
\]
Therefore, $Q^e_{X, \Delta, n}$ is a $W_e\MO_X$-module. 
\end{proof}

\begin{rem} 
In the situation of Setting \ref{setting:foundations-of-log-quasi-F-splitting}, assume that $(R,\mathfrak m)$ is a local ring and let $\Delta$ be a $\mathbb Q$-divisor on $X$.
We then get 
\[
H^i_{W_n\m}(X, Q^e_{X, \Delta, n})= H^i_{W_{\ell}\m}(X, Q^e_{X, \Delta, n}) 
= H^i_{W_e\m}(X, Q^e_{X, \Delta, n}) 
\]
for $i \in \Z$, $e \in \Z_{>0}, n \in \Z_{>0}$, and $\ell := \min \{e, n\}$ (Remark \ref{r-loc-coh-W_nm}).
\end{rem}

\begin{lem}\label{l-B-induction}
In the situation of Setting \ref{setting:most-general-foundations-of-log-quasi-F-splitting}, let $\Delta$ be a $\bQ$-divisor on $X$. 
Then we have the following exact sequences
\begin{equation}\label{e-Be-induction}
0 \to B^e_{X, \Delta, n} \to B_{X, \Delta, n}^{e+r} \to 
F_*^eB^r_{X, p^e\Delta, n}  \to 0 
\end{equation}
\begin{equation}\label{e-Bn-induction}
0 \to F_*^nB^e_{X, p^n\Delta, s} \to B_{X, \Delta, n+s}^{e} \to 
 {B^e_{X, \Delta, n}} \to 0 
\end{equation}
for all $e, n, r, s \in \Z_{>0}$. 
\end{lem}

\begin{proof}
The exact sequences 
(\ref{e-Be-induction}) and (\ref{e-Bn-induction}) 
are obtained by applying the snake lemma to the following commutative diagrams in which each horizontal sequence is exact and all the vertical arrows are injective: 
\[
\begin{tikzcd}
 & W_n\MO_X(\Delta) \arrow[r,equal] \arrow[d, "F^e"]& W_n\MO_X(\Delta)\arrow[d, "F^{e+r}"]\\
0 \arrow[r] & F_*^eW_{n}\MO_X(p^e\Delta) \arrow[r, "F^r"] & F_*^{e+r}W_n\MO_X(p^{e+r}\Delta) 
\arrow[r] & F_*^eB^r_{X, p^e\Delta, n} \arrow[r] & 0,
\end{tikzcd}
\]
\[
\begin{tikzcd}
0 \arrow[r] & F_*^nW_s\mathcal{O}_X(p^n\Delta) \arrow[r, "V^{n}"] \arrow[d, "F^e"] & W_{n+s}\mathcal{O}_X(\Delta) \arrow[r, "R^s"] \arrow[d, "F^e"] & W_{n}\mathcal{O}_X(\Delta) \arrow[r] \arrow[d, "F^e"] & 0 \\
0 \arrow[r] & F_*^{n+e}W_s\mathcal{O}_X(p^{n+e}\Delta) \arrow[r, "V^{n}"] & {F^e_*}W_{n+s}\mathcal{O}_X(p^e\Delta) \arrow[r, "R^s"] & F_*^eW_{n}\mathcal{O}_X(p^e\Delta) \arrow[r] & 0
\end{tikzcd}
\]
\end{proof}


\subsection{Alternative definitions via splittings}\label{ss-alternative}

Throughout this subsection, 
we fix an isomorphism $\iota: \MO_X \xrightarrow{\simeq} \omega_X(-K_X)$ 
and identify $\MO_X$ with $\omega_X(-K_X)$ via this fixed isomorphism. 
For a Weil divisor $D$, the equality $\MO_X(D) = \omega_X(D-K_X)$ means 
the fixed isomorphism $\MO_X(D) \xrightarrow{\simeq, \iota} \omega_X(D-K_X)$.

\begin{dfn}\label{d-Q'-Q''}
In the situation of Setting  \ref{setting:most-general-foundations-of-log-quasi-F-splitting}, let $\Delta$ be a $\bQ$-divisor on $X$. 
Fix $e, n \in \Z_{>0}$. 
We define a $W_n\MO_X$-module $Q'^e_{X, \Delta, n}$ 
and a $W_n\MO_X$-module homomorphism 
\[
\Phi'^e_{X, \Delta, n} : W_n\omega_X(\rdown{\Delta}-K_X) \to Q'^e_{X, \Delta, n}
\] 
by the following diagram, where each square is a pushout diagram  of $W_n\MO_X$-modules:  
\begin{equation} \label{e-def-Q'}
\begin{tikzcd}
W_n\cO_X(\Delta) \arrow{r}{F^e} \arrow{d}{R^{n-1}} & F^e_* W_n \cO_X(p^e\Delta)  \arrow{d}\\
\MO_X(\rdown{\Delta}) = \omega_X(\rdown{\Delta}-K_X) \arrow{r}{\Phi^e_{X, \Delta, n}} 
\arrow[d, "\underline{p}^{n-1}"]&
\arrow[lu, phantom, "\usebox\pushoutdr" , very near start, yshift=0em, xshift=0.6em, color=black] Q^e_{X, \Delta, n} \arrow[d]\\
W_n\omega_X(\rdown{\Delta}-K_X) \arrow{r}{\Phi'^e_{X, \Delta, n}}&
\arrow[lu, phantom, "\usebox\pushoutdr" , very near start, yshift=0em, xshift=0.6em, color=black] Q'^e_{X, \Delta, n}.
\end{tikzcd}
\end{equation}
Recall that $\underline{p}^{n-1}$ is the dual of the restriction $R^{{n-1}}$. 
Moreover, $\omega_X(\rdown{\Delta} -K_X) = (\omega_X \otimes \MO_X(\rdown{\Delta}-K_X))^{**}$ 
and 
$W_n\omega_X(\rdown{\Delta} -K_X) = (W_n\omega_X \otimes W_n\MO_X(\rdown{\Delta}-K_X))^{**}$.  
\end{dfn}

{Explicitly, $Q'^e_{X, \Delta, n}$ is defined via the short exact sequence

\begin{equation} \label{eq:explicit-def-Q'}
0 \to \MO_X(\rdown{\Delta}) \to W_n\omega_X(\rdown{\Delta} - K_X) \oplus Q^e_{X, \Delta, n} \to Q'^e_{X, \Delta, n} \to 0.
\end{equation}}

\begin{rem}\label{r-diag-iota}
As explained above, the equality $\MO_X(\rdown{\Delta}) = \omega_X(\rdown{\Delta}-K_X)$ means the isomorphism 
$\iota : \MO_X(\rdown{\Delta}) \xrightarrow{\simeq} \omega_X(\rdown{\Delta}-K_X)$ induced by the fixed isomorphism $\iota: \MO_X \xrightarrow{\simeq} \omega_X(-K_X)$. 
Then the diagram (\ref{e-def-Q'}) is a shortened version of the following diagram: 
\[
\begin{tikzcd}
W_n\cO_X(\Delta) \arrow{r}{F^e} \arrow{d}{R^{n-1}} & F^e_* W_n \cO_X(p^e\Delta)  \arrow{d}\\
\MO_X(\rdown{\Delta})  \arrow{r}{\Phi^e_{X, \Delta, n}} 
\arrow[d, "\iota", "\simeq"']&
\arrow[lu, phantom, "\usebox\pushoutdr" , very near start, yshift=0em, xshift=0.6em, color=black] Q^e_{X, \Delta, n} \arrow[dd]\\
\omega_X(\rdown{\Delta}-K_X)\arrow[d, "\underline{p}^{n-1}"]\\
W_n\omega_X(\rdown{\Delta}-K_X) \arrow{r}{\Phi'^e_{X, \Delta, n}}&
\arrow[lu, phantom, "\usebox\pushoutdr" , very near start, yshift=0em, xshift=0.6em, color=black] Q'^e_{X, \Delta, n}.
\end{tikzcd}
\]
It is easy to see that $Q'^e_{X, \Delta, n}$ and $\Phi'^e_{X, \Delta, n}$ 
are unique up to isomorphisms. 
\end{rem}

\begin{rem}
We use the same notations as in Definition \ref{d-Q'-Q''}. 
Since $F^e$ is injective, $\Phi'^e_{X, \Delta, n}$ 
is injective and 
we have the following exact sequences of coherent $W_n\MO_X$-modules, 
\[
0 \to W_n\omega_X({\rdown{\Delta}}-K_X)\xrightarrow{\Phi'^e_{X, \Delta, n}} Q'^e_{X, \Delta, n} \to B^e_{X, \Delta, n} \to 0, 
\]
because we have similar diagrams to (\ref{e-big-BC-diagram}) for $Q'^e_{X, \Delta, n}$ (cf.\ (\ref{eq:explicit-def-Q'})).
\end{rem}

\begin{dfn}\label{d-Psi-prime}
In the situation of Setting  \ref{setting:most-general-foundations-of-log-quasi-F-splitting}, let $\Delta$ be a $\bQ$-divisor on $X$.  
We define a $W_n\MO_X$-module homomorphism \[(\Phi'^e_{X, \Delta, n})^* : (Q'^e_{X, \Delta, n})^* \to  (W_n\omega_X(\rdown{\Delta}-K_X))^*\]
by applying $(-)^* := \mathcal Hom_{W_n\cO_X}(-, W_n\omega_X(-K_X))$ to the $W_n\cO_X$-module homomorphism \[\Phi'^e_{X, \Delta, n} \colon W_n\omega_X(\rdown{\Delta}-K_X)  \to Q'^e_{X,\Delta,n}.\]
More explicitly, $(\Phi'^e_{X, \Delta, n})^*$ can be written as follows: 
{\small 
\[
\mathcal Hom_{W_n\MO_X}(Q'^e_{X, \Delta, n}, W_n\omega_X(-K_X)) \xrightarrow{(\Phi'^e_{X, \Delta, n})^*} 
{\underbrace{\mathcal Hom_{W_n\MO_X}(W_n\omega_X(\rdown{\Delta}-K_X), W_n\omega_X(-K_X))}_{
{\simeq}\, W_n\MO_X(-\rdown{\Delta})}}.  
\]
}
\end{dfn}


\begin{prop}\label{p-many-defs}
In the situation of Setting  \ref{setting:most-general-foundations-of-log-quasi-F-splitting}, let $\Delta$ be a $\bQ$-divisor on $X$ satisfying $\rdown{\Delta}=0$. 
Fix $e, n \in \Z_{>0}$. 
Then the following are equivalent. 
\begin{enumerate}
\item $(X, \Delta)$ is $n$-quasi-$F^e$-split. 
\item There exists a $W_n\MO_X$-module homomorphism $\alpha : F_*^e W_n\MO_X({p^e\Delta}) \to W_n\omega_X(-K_X)$ that commutes with the following diagram: 
\[\begin{tikzcd}
	W_n\MO_X(\Delta) & {} & F_*^eW_n\MO_X(p^e\Delta) \\
	{\mathcal O_X = \omega_X(-K_X)} \\
	{W_n\omega_X(-K_X)}
	\arrow["R^{n-1}"', from=1-1, to=2-1]
	\arrow["\underline{p}^{n-1}"', from=2-1, to=3-1]
	\arrow["{F^e}", from=1-1, to=1-3]
	\arrow["{\exists \alpha}", dashed, from=1-3, to=3-1]
\end{tikzcd}\]
\item 
There exists a $W_n\MO_X$-module homomorphism 
$\beta : Q^e_{X, \Delta, n} \to W_n\omega_X(-K_X)$ that commutes with the following diagram: 
\[\begin{tikzcd}
	\MO_X = \omega_X(-K_X) && Q^e_{X, \Delta, n} \\
	W_n\omega_X(-K_X)
	\arrow["\underline{p}^{n-1}"', from=1-1, to=2-1]
	\arrow["\Phi^e_{X, \Delta, n}", from=1-1, to=1-3]
	\arrow["{\exists \beta}", dashed, from=1-3, to=2-1]
\end{tikzcd}\]
\item 
$\Phi'^e_{X, \Delta, n} : W_n\omega_X(-K_X) \to Q'^e_{X, \Delta, n}$ 
splits as a $W_n\MO_X$-module homomorphism.
\item The map \[ H^0(X, (\Phi'^e_{X, \Delta, n})^*) : 
H^0(X, (Q'^e_{X, \Delta, n})^*) \to H^0(X, (W_n\omega_X(-K_X))^*)\]
is surjective, where 
$(-)^* := \cHom_{W_n\MO_X}(-, W_n\omega_X(-K_X))$. 
\end{enumerate}
\end{prop}

\noindent Note that we have 
$(W_n\omega_X(-K_X))^* = \cHom_{W_n\MO_X}(W_n\omega_X(-K_X), W_n\omega_X(-K_X)) 
\simeq \cHom_{W_n\MO_X}(W_n\omega_X, W_n\omega_X) \simeq W_n\MO_X$. 

\begin{proof}
First of all, we prove  that $(2) \Leftrightarrow (3) \Leftrightarrow (4) \Leftrightarrow (5)$. 
Since $Q^e_{X,\Delta, n}$ and $Q'^e_{X, \Delta, n}$ are the pushouts (cf. (\ref{e-def-Q'})), we get $(2) \Leftrightarrow (3) \Leftrightarrow (4)$. 
It is clear that (4) $\Rightarrow$ (5). 
It is also easy to see the opposite implication (5) $\Rightarrow$ (4), because the surjectivity in (5) assures the existence of 
$\theta \in \Hom_{W_n\MO_X}(Q'^e_{X, n}, W_n\omega_X(-K_X))$ 
satisfying $(\Phi'^e_{X, \Delta, n})^*(\theta) = {\rm id}$, i.e., 
$\theta \circ \Phi'^e_{X, \Delta, n} = {\rm id}$. 
This completes the proof of $(2) \Leftrightarrow (3) \Leftrightarrow (4) \Leftrightarrow (5)$. 

\medskip


It is obvious that (1) $\Rightarrow$ {(3)}. 
Let us show (3) $\Rightarrow$ (1). 
Assume (3). 
Let $j : X \hookrightarrow W_{n}X$ be the induced closed immersion. 
Recall that the equality $\MO_X = \omega_X(-K_X)$ in (3) means a fixed isomorphism $\iota: \MO_X \xrightarrow{\simeq} \omega_X(-K_X)$. 
We have the following isomorphisms: 
\begin{eqnarray*}
\Hom_{W_n\MO_X}(j_*\MO_X, W_n\omega_X(-K_X))
&=&H^0(X, \mathcal Hom_{W_n\MO_X}(j_*\MO_X, W_n\omega_X(-K_X)) )\\
&\overset{{\rm (i)}}{\simeq} &H^0(X, \mathcal Hom_{\MO_X}(\MO_X, \omega_X(-K_X)) )\\
&=& \Hom_{\MO_X}(\MO_X, \omega_X(-K_X))\\
&\overset{{\rm (ii)}}{\simeq} &H^0(X, \MO_X), 
\end{eqnarray*}
where (i) follows from the Grothendieck duality (as both hand sides are $S_2$) and (ii) is given by  our fixed isomorphism $\iota: \MO_X \xrightarrow{\simeq} \omega_X(-K_X)$. 
It follows from  \[H^0(X, W_n\MO_X) = W_n(H^0(X, \MO_X))\] that  
$\Hom_{W_n\MO_X}(j_*\MO_X, W_n\omega_X(-K_X))$ has the generator as a $H^0(X, W_n\MO_X)$-module 
which is corresponding to $1 \in  H^0(X, \MO_X)$. 
For the trace map $\underline{p}^{n-1}: j_*\omega_X(-K_X) \to W_n\omega_X(-K_X)$ (obtained as the dual of $R^{n-1}:W_n\MO_X \to j_*\MO_X$), 
the first isomorphism (i) is given as follows: 
\begin{eqnarray*}
\Hom_{W_n\MO_X}(j_*\MO_X, W_n\omega_X(-K_X)) &\underset{(i)}{\overset{\simeq}{\longleftarrow}}&\Hom_{\MO_X}(\MO_X, \omega_X(-K_X)) \\
(j_*\MO_X \xrightarrow{j_*\varphi} j_*\omega_X(-K_X) \xrightarrow{\underline{p}^{n-1}} W_n\omega_X(-K_X)) & \longmapsfrom & (\varphi : \MO_X \to \omega_X(-K_X)). 
\end{eqnarray*}
To summarise, we have 
\begin{align*}
\eqmathbox[M1]{\Hom_{W_n\MO_X}(j_*\MO_X, W_{n}\omega_X(-K_X))} &&\overset{{\rm (i)}}{\simeq}& \eqmathbox[M2]{\quad \Hom_{\MO_X}(\MO_X, \omega_X(-K_X))} && 
\overset{{\rm (ii)}}{\simeq} \eqmathbox[M3]{\quad H^0(X, \MO_X)}\\ 
\eqmathbox[M1]{j_*\MO_X \xrightarrow{\underline{p}^{n-1} \circ j_*\iota} W_n\omega_X(-K_X)}
&&\leftrightarrow& \eqmathbox[M2]{\quad \iota : \MO_X \xrightarrow{\simeq} \omega_X(-K_X)}  &&\leftrightarrow \eqmathbox[M3]{\quad 1.}
\end{align*}
\vspace{0.0em}

Pick $\beta : Q^e_{X, \Delta, n} \to W_n\omega_X(-K_X)$ as in (3). 
The commutativity of the diagram in (3) means that
$\underline{p}^{n-1} \circ (j_* \iota) = \beta \circ \Phi^e_{X, \Delta, n}$. 
By 
\begin{eqnarray*}
\Hom_{W_n\MO_X}(Q^e_{X, \Delta, n}, W_n\omega_X(-K_X)) 
&\xrightarrow{H^0(X, (\Phi_{X, \Delta, n}^e)^*)}&
\Hom_{W_n\MO_X}(j_*\MO_X, W_n\omega_X(-K_X))\\
\beta &\longmapsto& \beta \circ \Phi^e_{X, \Delta, n}, 
\end{eqnarray*}
the image of $H^0(X, (\Phi_{X, \Delta, n}^e)^*)$ 
contains $\underline{p}^{n-1} \circ (j_* \iota) =\beta \circ \Phi^e_{X, \Delta, n}$, which is a generator of $\Hom_{W_n\MO_X}(j_*\MO_X, W_n\omega_X(-K_X))$ 
as an $H^0(X, W_n\MO_X)$-module. 
Since $H^0(X,  (\Phi_{X, \Delta, n}^e)^*)$ is a $H^0(X, W_n\MO_X)$-module homomorphism, 
$H^0(X, (\Phi_{X, \Delta, n}^e)^*)$ is surjective. 
\end{proof}

{One of the advantages of the splitting definition above is that it allows for showing that sections of all line bundles are quasi-$F^e$-stable when $X$ is quasi-$F^e$-split (Proposition \ref{prop:qS^0-for-quasi-F-split}). This will be a consequence of the following proposition.}

\begin{prop}\label{p-IQFS-L-surje}
In the situation of Setting \ref{setting:most-general-foundations-of-log-quasi-F-splitting}, 
let $\Delta$ be an effective $\Q$-divisor with $\rdown{\Delta}=0$. 
Take $e, n \in \Z_{>0}$ and suppose that $(X, \Delta)$ is $n$-quasi-$F^e$-split. 
Let $D$ be a Weil divisor on $X$. 
Then  the following induced maps are surjective: 
\[
H^0(X, (\Phi^e_{X, \Delta+D, n})^*) : 
H^0(X, (Q^e_{X, \Delta+D, n})^*) 
\to H^0(X, \MO_X(D)^*) 
\]
\[
H^0(X, (\Phi'^e_{X, \Delta+D, n})^*) : 
H^0(X, (Q'^e_{X, \Delta+D, n})^*) 
\to H^0(X, (W_n\omega_X({-K_X+}D))^*), 
\]
where $(-)^* := \cHom_{W_n\MO_X}(-, W_n\omega_X(-K_X))$ and $W_n\MO_X(D)$ denotes the Witt divisorial sheaf. 
\end{prop}

\noindent Note that 
$\MO_X(D)^* \simeq \MO_X(-D)$  and $W_n\omega_X({-K_X} + D)^* \simeq W_n\MO_X(-D)$. 

\begin{proof}
We have the following pushout diagram (Definition \ref{d-Q'-Q''}): 
\[
\begin{CD}
\omega_X(-K_X) @>\Phi^e_{X, \Delta, n}>> Q^e_{X, \Delta, n}\\
@VV\underline{p}^{n-1} V @VVV\\
W_n\omega_X(-K_X) @>\Phi'^e_{X, \Delta, n} >> Q'^e_{X, \Delta,  n}.
\end{CD}
\]
By applying 
$(-) \otimes_{W_n\MO_X} W_n\MO_X(D)$ and $(-)^* = \mathcal Hom_{W_n\MO_X}(-, W_n\omega_X(-K_X))$, we obtain the following commutative diagram
\[
\begin{CD}
\MO_X(-D) @<(\Phi^e_{X, \Delta+D, n})^*<< (Q^e_{X, \Delta+D, n})^*\\
@AAR^{n-1} A @AAA\\
W_n\MO_X(-D) @<(\Phi'^e_{X, \Delta+D, n})^*<< (Q'^e_{X, \Delta+D, n})^*. 
\end{CD}
\]
Since $\Phi'^e_{X, \Delta, n}$ 
is a split injection (Proposition \ref{p-many-defs}),  
$(\Phi'^e_{X, \Delta+D, n})^*$ is a split surjection. 
Therefore, 
\[
H^0(X, (\Phi'^e_{X, \Delta+D, n})^*) : 
H^0(X, (Q'^e_{X, \Delta+D, n})^*) 
\to H^0(X, W_n\MO_X(-D))
\]
is surjective. Note that 
\[
R^{n-1} : H^0(X, W_n\MO_X(-D)) \to H^0(X, \MO_X(-D))
\]
is surjective, because an arbitrary element $s \in H^0(X, \MO_X(-D))$ has its Teichm\"{u}ller lift $\underline{s} \in  H^0(X, W_n\MO_X(-D))$, which  satisfies $R^{n-1}(\underline{s})=s$. 
By diagram chase, the induced map 
\[
H^0(X, (\Phi^e_{X, \Delta+D, n})^*) : 
H^0(X, (Q^e_{X, \Delta+D, n})^*) 
\to H^0(X, \MO_X(-D))
\]
is surjective. 
\end{proof}

\subsection{Finite covers} \label{ss:finite-covers}



\begin{prop}\label{p-Gal-Q-split}
In the situation of Setting \ref{setting:most-general-foundations-of-log-quasi-F-splitting}, 
let $f:Y \to X$ be a finite surjective morphism 
from an integral normal excellent scheme $Y$ such that {$K(Y)/K(X)$} is a Galois extension 
and its extension degree {$[K(Y):K(X)]$} is not divisible by $p$. 
Let $D$ be a $\Q$-divisor on $X$  and set $D_Y := f^*D$. 
Then the induced $W_n\MO_X$-module homomorphism 
\[
f^* : Q^e_{X, D, n} \to f_*Q^e_{Y, D_Y, n}
\]
splits. 
\end{prop}


\begin{proof}
Set $G := {\rm Gal}(K(Y)/K(X)) = \{ \sigma_1, ..., \sigma_n\}$. 
For every $\Q$-divisor $E$ on $X$, we have a $W_n\MO_X$-module homomorphism 
\[
\pi_{E} : 
{f_*}W_n\MO_Y(f^*E) \to W_n\MO_X(E), \qquad \alpha \mapsto \frac{1}{|G|}\sum_{\sigma \in G} \sigma^*\alpha. 
\]
We then  obtain the following diagram in which each horizontal sequence is exact: 
\[
\begin{CD}
0 @>>> f_*F_*W_{n-1}\MO_Y(pD_Y) @>VF^e>> f_*F^e_*W_n\MO_Y(p^eD_Y) @>>> f_*Q^e_{Y, D_Y, n} @>>> 0\\
@. @VV\pi_{D} V @VVF^e_*\pi_{p^eD} V @VV\overline{\pi} V\\
0 @>>> F_*W_{n-1}\MO_X(pD) @>VF^e>> F^e_*W_n\MO_X(p^eD) @>>> Q^e_{X, D, n} @>>> 0\\
\end{CD}
\]
Indeed, the left square in the above diagram is commutative by the following commutative diagram, and $\overline{\pi}$ is induced from the other vertical arrows 
\[
\begin{CD}
F_*W_{n-1}\MO_Y(pD_Y) @>VF^e>> F^e_*W_n\MO_Y(p^eD_Y)\\
@VV \sigma^* V @VV \sigma^* V\\
F_*W_{n-1}\MO_Y(pD_Y) @>VF^e>> F^e_*W_n\MO_Y(p^eD_Y). 
\end{CD}
\]
Note that we have 
$F^e_*\pi_{p^eD} \circ f^* ={\rm id}$. 
By diagram chase, we have that $\overline{\pi} \circ f^* ={\rm id}$, and hence 
$
f^* : Q^e_{X, D, n} \to f_*Q^e_{Y, D_Y, n}$ splits. 
\qedhere




\end{proof}

\begin{prop}\label{p-Galois-ht}
In the situation of Setting \ref{setting:foundations-of-log-quasi-F-splitting}, 
assume that $(R, \m)$ is a local ring. 
Let $f:Y \to X$ be a finite surjective morphism 
from an integral normal excellent scheme $Y$. 
Let $\Delta_X$ and $\Delta_Y$ be effective $\Q$-divisors on $X$ and $Y$ such that 
$\llcorner \Delta_X \lrcorner =0, \llcorner \Delta_Y \lrcorner =0$, and 
$K_Y +\Delta_Y \sim f^*(K_X+\Delta_X)$.  
\begin{enumerate}
\item Assume that 
$f^* : H^d_{\m}(X, \MO_X(K_X+\Delta_X)) \to H^d_{\m}(Y, \MO_Y(K_Y+\Delta_Y))$ is injective. 
If $(Y, \Delta_Y)$ is $n$-quasi-$F^e$-split, then so is $(X, \Delta_X)$. 
In particular, $h^e(X, \Delta_X) \leq h^e(Y, \Delta_Y)$. 
\item Assume that 
$K(Y)/K(X)$ is a Galois extension, its extension degree $[K(Y):K(X)]$ is not divisible by $p$, 
and $f^* : H^d_{\m}(X, \MO_X(K_X+\Delta_X)) \to H^d_{\m}(Y, \MO_Y(K_Y+\Delta_Y))$ is bijective. 
Then $(X, \Delta_X)$ is $n$-quasi-$F^e$-split if and only if so is $(Y, \Delta_Y)$. 
In particular, $h^e(X, \Delta_X) = h^e(Y, \Delta_Y)$. 
\end{enumerate}
\end{prop}


\begin{proof}
We have the following commutative diagram 
\[
\begin{CD}
f_*W_n\MO_Y(K_Y+\Delta_Y) @>F^e>> f_*F_*^eW_n\MO_Y( p^e(K_Y+\Delta_Y))\\
@AAf^*A @AAf^*A\\
W_n\MO_X(K_X+\Delta_X) @>F^e>> F_*^eW_n\MO_X( p^e(K_X+\Delta_X)).\\
\end{CD}
\]
Taking the pushouts and $H^d_{\m}(Y, -)$, we obtain 
the following commutative diagram 
\[
\begin{tikzcd}[column sep=100]
H^d_{\m}(Y, \mathcal{O}_Y(K_Y+\Delta_Y)) \arrow[r, "{H^d_{\m}(Y, \Phi^e_{Y, K_Y+\Delta_Y, n})}"]  & H^d_{\m}(Y, Q^e_{Y, K_Y+\Delta_Y, n})  \\
H^d_{\m}(X, \mathcal{O}_X(K_X+\Delta_X)) \arrow[r, "{H^d_{\m}(X, \Phi^e_{X, K_X+\Delta_X, n})}"]
\arrow[u, "f^* =:\alpha"']& 
H^d_{\m}(X, {Q^e_{X, K_X+\Delta_X, n}})\arrow[u, "f^* =:\beta"'].
\end{tikzcd}
\]

Let us show (1). 
Assume that $(Y, \Delta_Y)$ is $n$-quasi-$F^e$-split, i.e., 
$H_{\m}^{d}(Y, \Phi^e_{Y, K_Y+\Delta_Y, n})$ is injective (Lemma \ref{lem:cohomological-criterion-for-log-quasi-F-splitting}). 
Since $\alpha$ is injective, 
so is $H^d_{\m}(X, \Phi^e_{X, K_X+\Delta_X, n})$ by diagram chase. 
Then 
$(X, \Delta_X)$ is $n$-quasi-$F^e$-split (Lemma \ref{lem:cohomological-criterion-for-log-quasi-F-splitting}). 
Thus (1) holds. 

The assertion (2) follows from a similar argument by using the fact that $\beta$ is an injection (Proposition \ref{p-Gal-Q-split}).
\end{proof}

\begin{cor}\label{c-Galois-ht}
Let $k$ be an algebraically closed field of characteristic $p>0$ and let 
$f: Y \to X$ be a finite surjective morphism 
of projective normal varieties over $k$. 
Let $\Delta_X$ and $\Delta_Y$ be effective $\Q$-divisors on $X$ and $Y$ such that 
$\llcorner \Delta_X \lrcorner =0, \llcorner \Delta_Y \lrcorner =0$, and 
$K_Y +\Delta_Y \sim f^*(K_X+\Delta_X)$.  
Assume that 
$K(Y)/K(X)$ is a Galois extension and its extension degree $[K(Y):K(X)]$ is not divisible by $p$. 
Then $(X, \Delta_X)$ is $n$-quasi-$F^e$-split if and only if so is $(Y, \Delta_Y)$. 
In particular, $h^e(X, \Delta_X) = h^e(Y, \Delta_Y)$. 
\end{cor}

\begin{proof}
Set $d := \dim X = \dim Y$. 
Since $H^d(X, \MO_X(K_X+\Delta_X)) \to H^d(Y, \MO_Y(K_Y+\Delta_Y))$ is an isomorphism, 
the assertion follows from Proposition \ref{p-Galois-ht}(2). 
\end{proof}

\begin{prop}\label{p-etale-bc}
Let $f:Y \to X$ be an \'etale morphism of normal 
varieties over a perfect field of characteristic $p>0$. 
Take an effective $\Q$-divisor $\Delta_X$ on $X$ {such that $\rdown{\Delta_X}=0$} and set $\Delta_Y :=f^*\Delta_X$. 
If $(X, \Delta_X)$ is $n$-quasi-$F^e$-split, then 
so is $(Y, \Delta_Y)$. In particular, $h^e(Y, \Delta_Y) \leq h^e(X, \Delta_X)$. 
\end{prop}

\begin{proof}
By removing the singular locus of $X$, we may assume that $X$ and $Y$ are smooth (Remark \ref{r-big-open-QFS}).
Since $(X, \Delta_X)$ is $n$-quasi-$F^e$-split, 
there exists a $W_n\MO_X$-module homomorphism $\alpha$ 
which completes the following commutative diagram (Proposition \ref{p-many-defs}(2)): 
\[\begin{tikzcd}
	W_n\MO_X(\Delta_X) & {} & F_*^eW_n\MO_X(p^e\Delta_X) \\
	{\mathcal O_X = \omega_X(-K_X)} \\
	{W_n\omega_X(-K_X).}
	\arrow["R^{n-1}"', from=1-1, to=2-1]
	\arrow["\underline{p}^{n-1}"', from=2-1, to=3-1]
	\arrow["{F^e}", from=1-1, to=1-3]
	\arrow["{\exists \alpha}", dashed, from=1-3, to=3-1]
\end{tikzcd}
\]
Recall that we have the induced \'etale morphism $W_nf: {W_nY \to W_nX}$ \cite[Ch. 0, Proposition 1.5.8]{illusie_de_rham_witt}. 
By applying $(W_nf)^*$ to the above diagram, we see that $(Y, \Delta_Y)$ is $n$-quasi-$F^e$-split. 
Indeed, for every $d \in \Z_{\geq 0}$, we have 
\[
(W_nf)^*(F_*^dW_n\MO_X(p^d\Delta_X)) \overset{{\rm (i)}}{\simeq} 
F^d_*(W_nf)^*(W_n\MO_X(p^d\Delta_X)) 
\overset{{\rm (ii)}}{\simeq} F_*^d W_n\MO_Y(p^df^*\Delta_X) 
\]
and 
\[
(W_nf)^*(W_n\omega_X(-K_X)) 
= (W_nf)^*(W_n\omega_X \otimes W_n\MO_X(-K_X)) \overset{{\rm (iii)}}{\simeq}  W_n\omega_Y \otimes W_n\MO_Y(-K_Y), 
\]
where 
\begin{enumerate}
\item[(i)] holds by \cite[Ch. 0, Proposition 1.5.8]{illusie_de_rham_witt}, 

\item[(ii)] follows from \cite[Lemma 4.6]{tanaka22}, and 

\item[(iii)] is obtained by $(W_nf)^*(W_n\Omega_X^N) \simeq W_n\Omega_Y^N$ 
\cite[Ch. I, Proposition 1.14]{illusie_de_rham_witt} and 
$(W_nf)^*(W_n\MO_X(-K_X)) = W_n\MO_Y(-K_Y)$ \cite[Lemma 4.6]{tanaka22}.  \qedhere
\end{enumerate}
\end{proof}

\begin{prop}\label{p-bc-base-field}
Let $k \subset k'$ be an algebraic separable field extension between 
perfect fields of characteristic $p>0$. 
Take a normal variety $X$ over $k$ such that $X \times_k k'$ is integral. 
Let $\Delta$ be an effective $\Q$-divisor on $X$ with $\rdown{\Delta}=0$. 
Assume that $X$ is projective over an affine variety $V$ over $k$. 
Then $(X, \Delta)$ is $n$-quasi-$F^e$-split if and only if $(X \times_k k', \Delta \times_k k')$ is $n$-quasi-$F^e$-split, 
where $\Delta \times_k k'$ denotes the pullback of $\Delta$ to $X \times_k k'$. 
In particular, $h^e(X, \Delta) = h^e(X \times_k k', \Delta \times_k k')$. 
\end{prop}
\begin{proof}
Set $R := \Gamma(V, \MO_V)$. 
Then all the conditions in Setting \ref{setting:foundations-of-log-quasi-F-splitting} hold (Remark \ref{r-dcpx-in-application}). 
Since $k \subset k'$ can be written as a direct limit of finite separable extension, 
the problem is reduced to the case when $k \subset k'$ is a finite separable extension. 
In particular, $X \times_k k' \to X$ is \'etale. 
If $(X, \Delta)$ is $n$-quasi-$F^e$-split, then 
so is $(X \times_k k', \Delta \times_k k')$  by Proposition \ref{p-etale-bc}.

In order to  prove the  opposite implication, 
it is enough to show that \[f^* : H^d_{\m}(X, \MO_X(K_X)) \to H^d_{\m}(X \times_k k', \MO_{X \times_k k'}(K_{X \times_k k'}))\]
is injective for every maximal ideal $\mathfrak m$ of $R$ (Proposition \ref{p-Galois-ht}(1) applied for $R_\m$), where we used that 
$\MO_X(K_X+\Delta) = \MO_X(K_X)$ and 
$\MO_{X \times_k k'}(K_{X \times_k k'}+\Delta \times_k k') =\MO_{X \times_k k'}(K_{X \times_k k'})$. 
Since $k \hookrightarrow k'$ is a split injection, 
so are $f^* : \MO_X \to f_*\MO_{X \times_k k'}$ and $f^* : \MO_X(K_X) \to f_*\MO_{X \times_k k'}(K_{X \times_k k'})$. 
\end{proof}

\subsection{Quasi-$F^e$-splitting criterion via Cartier operator}
\label{ss:quasi-F^e-Cartier}

{We refer to \cite[Subsection 5.2]{KTTWYY1} for the discussion on the Cartier operator for projective schemes over domains essentially of finite type.}


\begin{lem} \label{l-QFS-criterion}
Let $(R, \m)$ be a local domain essentially of finite type over 
a perfect field of characteristic $p>0$.  
Let $(X,{\Delta)}$ be a $d$-dimensional 
{log pair} which is projective over $R$,  {where} $\rdown{\Delta}=0$. Let $f \colon Y \to X$ be a log resolution of $(X,\Delta)$ and let 
$B_Y$ be a $\bQ$-divisor such that  $\rdown{B_Y}\leq 0$, $-(K_Y+B_Y)$ is ample, and $f_*B_Y = \Delta$. 
{Set} $E := \mathrm{Supp}(B_Y)$. 
Take $e, n \in \Z_{>0}$ and suppose that $(*)$ holds. 
\begin{enumerate}
    \item[$(*)$]  $H^{d-1}_{\m}(Y, B_{Y, p^c(K_Y+B_Y), n}) = 0$ for every $0 \leq c \leq e-1$. 
\end{enumerate}
Then $(X,\Delta)$ is $n$-quasi-$F^e$-split.
\end{lem}

\begin{proof}
Write $K_Y+\Delta_Y = f^*(K_X+\Delta)$. Since $-(K_Y+B_Y)$ is ample and $f_*B_Y=\Delta$, the negativity lemma implies that $K_Y+B_Y \geq K_Y+\Delta_Y$, and so $f_*\cO_Y(K_Y+B_Y) = \cO_X(K_X+\Delta)$. 
Moreover, we obtain a natural map $f^* \colon Q^e_{X,K_X+\Delta,n} \to f_*Q^e_{Y,K_Y+B_Y,n}$ featured inside the following commutative diagram:
\[
\begin{tikzcd}[column sep = huge]
H^d_{\m}(X, \cO_X(K_X+\Delta)) \ar{r}{H^d_{\m}(\Phi^e_{X,K_X+\Delta,n})} \arrow[d,dash,shift left=.1em] \arrow[d,dash,shift right=.1em]  & H^d_{\m}(X, Q^e_{X,K_X+\Delta,n}) \ar{d}
 \ar{d}{f^*}\\
H^d_{\m}(Y, \cO_Y(K_Y+B_Y)) \ar{r}{H^d_{\m}(\Phi^e_{Y,K_Y+B_Y,n})} & 
H^d_{\m}(Y, Q^e_{Y,K_Y+B_Y,n}).
\end{tikzcd}
\]
Note that the injectivity of the upper horizontal arrow is equivalent to $(X,\Delta)$ being $n$-quasi-$F^e$-split  by Lemma \ref{lem:cohomological-criterion-for-log-quasi-F-splitting}. 
Thus, to prove the theorem, it is enough to show that the lower horizontal arrow
\begin{equation} \label{Cartier-criterion-local-cohomology-injectivity}
H^d_{\m}(\Phi^e_{Y,K_Y+B_Y,n}) : 
H^d_{\m}(Y, \mathcal{O}_Y(K_Y+B_Y)) 
\to H^d_{\m}(Q^e_{Y,K_Y+ B_Y,n}) \ \textrm{is  injective}. 
\end{equation}
By the exact sequence $0 \to \MO_Y(K_Y+B_Y) \to Q^e_{Y, K_Y+B_Y, n} \to B^e_{Y, K_Y+B_Y, n} \to 0$ (\ref{e-big-BC-diagram}), 
it suffices to prove that 
$H^{d-1}_{\m}(Y, B^e_{Y, K_Y+B_Y, n})=0$. 
By the following exact sequence (\ref{e-Be-induction}): 
\[
0 \to B^{e-1}_{Y, K_Y+B_Y, n} \to B_{Y, K_Y+B_Y, n}^{e} \to 
F_*^{e-1}B_{Y, p^{e-1}(K_Y+B_Y), n}  \to 0, 
\]
$H^{d-1}_{\m}(B^e_{Y, K_Y+B_Y, n})=0$ holds by  $(*)$. 
\end{proof}

\begin{theorem} \label{thm:higher-Cartier-criterion-for-quasi-F-split}
Let $(R, \m)$ be a local domain essentially of finite type over 
a perfect field of characteristic $p>0$.  
Let $(X,{\Delta)}$ be a $d$-dimensional 
{log pair} which is projective over $R$,  {where} $\rdown{\Delta}=0$. Let $f \colon Y \to X$ be a log resolution of $(X,\Delta)$ and let 
$B_Y$ be a $\bQ$-divisor such that  $\rdown{B_Y}\leq 0$, $-(K_Y+B_Y)$ is ample, and $f_*B_Y = \Delta$. 
{Set} $E := \mathrm{Supp}\, B_Y$. 
Take $e, n \in \Z_{>0}$ and suppose that
\vspace{0.2em}
\begin{enumerate}\setlength\itemsep{0.3em}
    \item  $H^{d-2}_{\m}(Y,\Omega_Y^1(\log E)(p^c(K_Y + B_Y))) = 0$ for every $0 \leq c \leq e-1$, 
    \item $H^{d-2}_{\m}(Y,B_1\Omega^2_Y(\log E)(p^k(K_Y + B_Y))) = 0$ for every $k \geq 1$, and 
    \item $H^1(Y, \Omega^1_Y(\log E)^* \otimes \cO_Y(K_Y-\rdown{p^{n+c}(K_Y + {B_Y})})) =0$
    for every $0 \leq c \leq e-1$, where 
    $\Omega^1_Y(\log E)^* := \mathcal Hom_{\MO_Y}(\Omega^1_Y(\log E), \MO_Y)$. 
\end{enumerate}
Then $(X,\Delta)$ is $n$-quasi-$F^e$-split.
\end{theorem}

\noindent {We refer to Notation \ref{notation:global-local-cohomology} for the definition of $H^i_{\m}(Y, -)$.} 
Also, to avoid confusion, we emphasise that 
\[
B_1\Omega^2_Y(\log E)(p^k(K_Y + B_Y)) \simeq  B_1(\Omega^2_Y(\log E)(p^kB_Y)) \otimes 
\MO_Y(p^{k-1}K_Y). 
\]
Last, we point out that assumption (3) is valid by Serre vanishing as long as $n \gg 0$.
\begin{proof}
By Lemma \ref{l-QFS-criterion}, it is enough to 
show that the condition $(*)$ holds. 
\vspace{0.2em}
\begin{enumerate}
    \item[$(*)$]  $H^{d-1}_{\m}(Y, B_{Y, p^c(K_Y+B_Y), n}) = 0$ for every $0 \leq c \leq e-1$. 
\end{enumerate}
\vspace{0.2em}
Fix an integer $c$ satisfying $ 0\leq c \leq e-1$. 
Recall that (\cite[Lemma 5.9]{KTTWYY1}):
\[
B_{Y, p^c(K_Y+B_Y), n} \simeq 
B_{Y, p^cB_Y, n} \otimes_{\MO_Y} \MO_Y(p^cK_Y) \simeq 
B_n\Omega_Y^1(\log E)(p^{n+c}B_Y)  \otimes_{\MO_Y} \MO_Y(p^cK_Y). 
\]
Now we use the following exact sequence \cite[(5.7.1)]{KTTWYY1}:
\[
0 \to B_n\Omega^1_Y(\log E)(p^{n+c} {B_Y}) \to Z_n\Omega^1_Y(\log E)(p^{n+c} { B_Y}) \xrightarrow{{C^n_{p^c{B_Y}}}}  \Omega^1_Y(\log E)(p^c{B_Y}) \to 0,
\]
By tensoring by {$\MO_Y(p^cK_Y)$}, the problem is reduced to showing  that 
\begin{itemize} \setlength\itemsep{0.2em}
    \item $H^{d-1}_{\m}(Y,Z_n\Omega^1_Y(\log E)(p^{n+c}{B_Y}) \otimes \MO_Y(p^cK_Y)) = 0$, and
    \item $H^{d-2}_{\m}(Y,\Omega^1_Y(\log E)(p^c{B_Y}) \otimes \MO_Y(p^cK_Y))=0$.
\end{itemize}
The latter {assertion} is nothing but our Assumption (1).
 {Thus} it is sufficient to prove the former {assertion}.

Now we use the short exact sequence \cite[Lemma 5.8]{KTTWYY1}:
{\small \[
0 \!\shortrightarrow\! 
F_*^kZ_{n{-}k} \Omega^1_Y(\log E)(p^{n+c}{B_Y}) \!\shortrightarrow\! 
F^{k+1}_*Z_{n{-}k{-}1} \Omega^1_Y(\log E)(p^{n+c}{B_Y}) \!\shortrightarrow\! 
F^k_*B_1\Omega^2_Y(\log E)(p^{k{+}1+c}{B_Y}) \!\shortrightarrow\! 0
\]}
By tensoring by $\MO_Y(p^cK_Y)$ 
and repeatedly applying Assumption (2), 
{we get the following injections
\begin{align*}
H^{d-1}_{\m}(Y,\, Z_n \Omega^1_Y(\log E)(p^{n+c}{B_Y}) \otimes \MO_Y(p^cK_Y)) 
\xhookrightarrow{\hphantom{aa}}\ 
&H^{d-1}_{\m}(Y, F_*Z_{n-1} \Omega^1_Y(\log E)(p^{n+c}{B_Y}) \otimes \MO_Y(p^cK_Y)) \\
\xhookrightarrow{\hphantom{aa}}\  &\cdots \\
\xhookrightarrow{\hphantom{aa}}\
&H^{d-1}_{\m}(Y, F_*^n\Omega^1_Y(\log E)(p^{n+c}{B_Y}) \otimes \MO_Y(p^cK_Y)). 
\end{align*}
}
{Hence it suffices to show that } 
\[
H^{d-1}_{\m}(Y, F_*^n\Omega^1_Y(\log E)(p^{n+c}{B_Y}) \otimes \MO_Y(p^cK_Y)) 
= H^{d-1}_{\m}(Y, \Omega^1_Y(\log E)(p^{n+c}(K_Y + {B_Y})))
\]
is zero. 
By Matlis duality {(Proposition \ref{prop:Matlis-duality})} and exactness of the derived $\mathfrak m$-completion in our setting (\stacksproj{0A06}), it is enough to verify that
\begin{align*}
H^{-d+1}&R\Hom(R\pi_*\Omega^1_Y(\log E)(p^{n+c}(K_Y + {B_Y})), \omega^{\mydot}_R) \\
= H^{1}&R\Hom(\Omega^1_Y(\log E)(p^{n+c}(K_Y + {B_Y})), \omega_Y) \\
= H^1&(Y, \Omega^1_Y(\log E)^* \otimes \cO_Y(K_Y-\rdown{p^{n+c}(K_Y + {B_Y})}))
\end{align*}
is equal to zero, where $\pi \colon Y \to \Spec R$ is the natural projection. 
This is nothing but Assumption (3). 
\end{proof}

\subsection{Definition of pure quasi-$F^e$-splitting}\label{ss:pure-quasi-F-split}

\noindent 
{Throughout this subsection, we work under the assumptions of Setting \ref{setting:most-general-foundations-of-log-pure-quasi-F-splitting}.

\begin{setting}[{general case, with boundary}] \label{setting:most-general-foundations-of-log-pure-quasi-F-splitting}
Let $X$ be a $d$-dimensional 
integral normal {Noetherian} 
$F$-finite 
$\F_p$-scheme such that $W_nX$ is an excellent scheme admitting a dualising complex $W_n\omega^{\mydot}_X$ for every $n \in \Z_{>0}$. 
We fix a 
Weil divisor $S$ on $X$ which {either} is a prime divisor or is equal to zero. 
\end{setting}

\begin{setting}[{projective case, with boundary}] \label{setting:foundations-of-log-pure-quasi-F-splitting}
Let $R$ be {an $F$-finite Noetherian domain of characteristic $p>0$} such that 
$W_nR$ is an excellent ring admitting a dualising complex for every $n \in \Z_{>0}$. 
Let $X$ be a $d$-dimensional integral normal scheme which is projective over $\Spec R$. 
We fix a 
Weil divisor $S$ on $X$ which is either a prime divisor or is equal to zero.
\end{setting}}
We denote the ideal sheaf of $S$ by $\cI_S$. 
We remind the reader that $W_n\cI_S$ and $W_n\cO_X(-S)$ are \emph{not} equal in general. 
Let $\Delta$ be a (non-necessarily effective) $\bQ$-divisor on $X$. 
We define a coherent $W_n\MO_X$-module $Q^{S, e}_{X,\Delta,n}$ and 
a $W_n\MO_X$-module homomorphism $\Phi^{S, e}_{X, \Delta, n}$ by the following pushout diagram of $W_n\MO_X$-modules: 
\begin{equation}\label{diagram:purely-quasi-F-split-definition}
\begin{tikzcd}
W_n\cI_S(\Delta) \arrow{r}{F^e} \arrow{d}{R^{n-1}} & F^e_* W_n \cI_S(p^e\Delta)  \arrow{d}\\
\cI_S(\Delta) \arrow{r}{\Phi^{S, e}_{X, \Delta, n}}& \arrow[lu, phantom, "\usebox\pushoutdr" , very near start, yshift=0em, xshift=0.6em, color=black] Q^{S, e}_{X, \Delta, n}.
\end{tikzcd}
\end{equation}
We remind the reader that
$\cI_S(\Delta) =  \MO_X(\Delta-S)$.

\begin{remark}
In our applications, we are mainly interested in the case {when $\Delta=S+B$ for a prime divisor $S$ and a $\bQ$-divisor $B$ such that $S \not \subseteq \Supp B$ and  $\rdown{B}=0$.}
Under these assumptions, {$B \geq 0$ and}  $\cI_S(\Delta) = \MO_X$.
\end{remark}

We define a coherent $W_n\cO_X$-module 
$B^{S, e}_{X, \Delta, n}$ by 
\begin{align*}
B^{S, e}_{X, \Delta, n} &:= {\rm Coker}(W_n\cI_S(\Delta) \xrightarrow{F^e} F^e_*W_n\cI_S(p^e\Delta)) \\
&\hphantom{:}= F^e_*W_n\cI_S({ p^e}\Delta)/F^e(W_n\cI_S({\Delta})). 
\end{align*} 

\begin{remark}
{The key properties of the construction of $Q^{S, e}_{X,\Delta,n}$,  $B^{S, e}_{X,\Delta,n}$, and $\Phi^{S, e}_{X,\Delta,n}$ may be encapsulated by the following diagram}
\begin{equation}\label{e-big-BC-diagram2}
\begin{tikzcd}
& 0 \arrow{d} & 0 \arrow{d} & & \\
& F_*W_{n-1}\cI_S(p\Delta) \arrow{d}{V} \arrow[r,dash,shift left=.1em] \arrow[r,dash,shift right=.1em] &  F_*W_{n-1}\cI_S(p\Delta) \arrow{d}{ F^eV =VF^e} &&\\
0 \arrow{r} &  W_n\cI_S(\Delta) \arrow{d}{R^{n-1}} \arrow{r}{F^e} & F^e_*W_n\cI_S(p^e\Delta) \arrow{r} \arrow{d} & B^{S, e}_{X, \Delta, n} \arrow[d,dash,shift left=.1em] \arrow[d,dash,shift right=.1em] \arrow{r} & 0 \\
0 \arrow{r} &  \cI_S(\Delta) \arrow{d} \arrow{r}{\Phi^{S, e}_{X, \Delta, n}} & Q^{S, e}_{X, \Delta, n} \arrow[lu, phantom, "\usebox\pushoutdr" , very near start, yshift=0em, xshift=1em, color=black] \arrow{d} \arrow{r} & B^{S, e}_{X, \Delta, n} \arrow{r} & 0.\\
& 0 & 0 &&
\end{tikzcd}
\end{equation}
{All the horizontal and vertical sequences are exact, as} $F^e \colon W_n\cI_S(\Delta) \to F^e_*W_n\cI_S(p^e\Delta)$ is injective, $R^{n-1} \colon W_n\cI_S(\Delta) \to \cI_S(\Delta)$ is surjective, and 
(\ref{diagram:purely-quasi-F-split-definition}) is a pushout diagram. 

\end{remark}

In what follows, we apply exactly the same arguments as in Subsection \ref{ss-log-QFS} and Subsection \ref{ss-alternative}, 
and hence we omit proofs.

\begin{lem}\label{l-B-induction-purely}
In the situation of Setting  \ref{setting:most-general-foundations-of-log-pure-quasi-F-splitting}, let $\Delta$ be a $\bQ$-divisor on $X$. 
Then we have the following exact sequences
\begin{equation}\label{e-Be-induction-purely}
0 \to B^{S, e}_{X, \Delta, n} \to B_{X, \Delta, n}^{S, e+r} \to 
F_*^eB^{S, r}_{X, p^e\Delta, n}  \to 0 
\end{equation}
\begin{equation}\label{e-Bn-induction-purely}
0 \to F_*^nB^{S, e}_{X, p^n\Delta, s} \to B_{X, \Delta, n+s}^{S, e} \to 
{B^{S, e}_{X, \Delta, n}} \to 0 
\end{equation}
for all $e, n, r, s \in \Z_{>0}$. 
\end{lem}

For a Cartier divisor $D$, the following holds: 
\begin{equation}\label{e-C-tensor2}
Q^{S, e}_{X,\Delta+D,n} \simeq Q^{S, e}_{X,\Delta,n} \otimes_{W_n\MO_X} W_n\cO_X(D). 
\end{equation} 

\begin{dfn}\label{d-Psi-S}
In the situation of Setting \ref{setting:most-general-foundations-of-log-pure-quasi-F-splitting}, let $\Delta$ be a $\bQ$-divisor. 
We define $(\Phi^{S, e}_{X, \Delta, n})^* : (Q^{S, e}_{X,\Delta,n})^* \to \cI_S(\Delta)^*$  by applying $(-)^* := \mathcal Hom_{W_n\cO_X}(-, W_n\omega_X(-K_X))$ to the $W_n\cO_X$-module homomorphism $\Phi^{S, e}_{X, \Delta, n} \colon \cI_S(\Delta)  \to Q^{S, e}_{X,\Delta,n}$:
\[
(\Phi^{S, e}_{X, \Delta, n})^* : 
\mathcal Hom_{W_n\MO_X}(Q^{S, e}_{X,\Delta,n}, W_n\omega_X(-K_X)) \to 
\mathcal Hom_{W_n\MO_X}(\cI_S(\Delta), W_n\omega_X(-K_X)). 
\]
Note that $\cI_S(\Delta) = \MO_X(-S+\Delta)$ 
and $\cI_S(\Delta)^* = \mathcal Hom_{W_n\MO_X}(\cI_S(\Delta), W_n\omega_X(-K_X)) 
\simeq \MO_X(\rup{S -\Delta})$. 
\end{dfn}

\begin{dfn}\label{d-pure-IQFS}
In the situation of Setting \ref{setting:most-general-foundations-of-log-pure-quasi-F-splitting}, 
let $B$ be a $\Q$-divisor on $X$ such that $S \not\subseteq \Supp B$ and $\rdown{B}=0$. {Set $\Delta :=S+B$} and 
take $n, e \in \Z_{>0}$. 
We say that $(X, S+B)$ is \emph{purely $n$-quasi-$F^e$-split (along $S$)} if 
the induced map
{
\[
(\Phi^{S, e}_{X, \Delta, n})^* : H^0(X, (Q^{S, e}_{X,\Delta,n})^*) \to H^0(X, \cI_S(\Delta)^*) = H^0(X,\cO_X)
\]}
is surjective. Note that we have 
\[
{\cI_S(\Delta)^* = } \cHom_{W_n\MO_X}(\cI_S(S+B), W_n\omega_X(-K_X)) = 
\cHom_{W_n\MO_X}(\MO_X, W_n\omega_X(-K_X)) \simeq \MO_X. 
\]
We call $(X, \Delta)$ {\em purely quasi-$F^e$-split} if it is purely $n$-quasi-$F^e$-split for some $n \in \Z_{>0}$. 
We say that  $(X, \Delta)$ is 
 {\em purely $n$-quasi-$F^e$-pure} (resp.\ 
 {\em purely quasi-$F^e$-pure}) 
if 
there exists an open cover $X = \bigcup_{i \in I} X_i$ 
such that $(X_i, \Delta|_{X_i})$ is 
 {purely $n$-quasi-$F^e$-split} (resp.\  
 {purely quasi-$F^e$-split}) for every $i \in I$. 
\end{dfn}

\begin{lemma} \label{lem:cohomological-criterion-for-pure-quasi-F-splitting}
In the situation of Setting  \ref{setting:foundations-of-log-pure-quasi-F-splitting}, let  $B$ be a $\bQ$-divisor on $X$ such that $S \not\subseteq \Supp B$ and $\rdown{B}=0$.  Assume that $R$ is a local ring with maximal ideal $\m$. 
Fix $e, n \in \Z_{>0}$. 
Then $(X,S+B)$ is purely $n$-quasi-$F^e$-split if and only if 
the following map 
is injective: 
\[
H^d_{\m}(X, \Phi^{S, e}_{X, K_X+S+B, n}) \colon H^d_{\m}(X, \cI_S(K_X+S+B)) \to H^d_{\m}(X, Q^{S, e}_{X,K_X+S+B, n}).
\]
\end{lemma}

\begin{lemma}\label{l-purely-BQI-ex}
In the situation of Setting \ref{setting:most-general-foundations-of-log-pure-quasi-F-splitting}, 
 let $\Delta$ be a $\bQ$-divisor on $X$. 
Take $e, n \in \Z_{>0}$. 
Then we have the following exact sequences of coherent $W_{n+1}\MO_X$-modules: 
\begin{align}
\label{eq:adjoint-C-restriction-sequence}
{0 \to 
F_*B^{S,e}_{X, p\Delta, n} \to Q^{S,e}_{X,\Delta,n+1} \to F^e_* \cI_S(p^e\Delta) \to 0.} 
\end{align}
\end{lemma}

\begin{proposition} 
In the situation of Setting \ref{setting:most-general-foundations-of-log-pure-quasi-F-splitting},  let $\Delta$ be a $\bQ$-divisor on $X$. 
Then $B^{S, e}_{X, \Delta, n}$ and $Q^{S, e}_{X, \Delta, n}$
are naturally coherent $W_e\MO_X$-modules. 
\end{proposition}

\subsection{Quasi-$F^e$-stable sections $q^eS^0_n(X,\Delta; L)$ and $q^eS^0_{{\rm adj}, n}(X,\Delta; L)$}\hphantom{a}\vspace{0.2em}
\label{ss:quasi-F^e-stable-sections}


\indent In this subsection, we introduce a submodule (subspace) 
\begin{align*} q^eS^0_n(X,\Delta;L) &\subseteq H^0(X, \MO_X(\rup{L-\Delta})), \  \text{and its adjoint variant}\\
q^eS^0_{n, \adj}(X,S+B; L) &\subseteq H^0(X, \MO_X(\rup{L-B})).
\end{align*}

\begin{definition} \label{definition:quasi-F-stable-sections} 
In the situation of Setting \ref{setting:most-general-foundations-of-log-quasi-F-splitting}, let $L$ and $\Delta$ be $\bQ$-divisors on $X$. 
Recall that we have $W_e\MO_X$-module homomorphisms 
\begin{align*}
\Phi^e_{X, \Delta -L, n} \colon\quad &\MO_X(\Delta -L) \longrightarrow Q^e_{X, \Delta -L, n}\\  
(\Phi^e_{X, \Delta -L, n})^* \colon \quad &(Q^e_{X, \Delta -L, n})^* \longrightarrow \MO_X(\Delta -L)^* = \MO_X(\rup{L-\Delta}), 
\end{align*}
where the lower one $(\Phi^e_{X, \Delta -L, n})^*$ is obtained by applying 
$(-)^* =\cHom_{W_n\MO_X}(-, W_n\omega_X(-K_X))$ to the upper one $\Phi^e_{X, \Delta -L, n}$. 
We define $q^eS^0_n(X,\Delta;L)$ by 
\[
q^eS^0_n(X,\Delta;L) := \image\big( 
H^0(X, (Q^e_{X,\Delta-L,n})^*) \xrightarrow{H^0(X, (\Phi^e_{X, \Delta -L, n})^*)} H^0(X, \cO_X(\lceil L-\Delta\rceil))\big).
\]
\end{definition}

When $\rdown{\Delta}=0$ and $L$ is a Weil divisor, we have 
\[
q^eS^0_n(X,\Delta; L) \subseteq H^0(X, \cO_X(\lceil L-\Delta\rceil))  = H^0(X,\cO_X(L)).
\]
However, we do need to consider the most general case of the above definition as non-integral $L$ and non-effective $\Delta$ come up naturally in the context of pullbacks and restrictions of Weil divisors.



\begin{lemma} \label{l:special-cases-of-L}
Under the same notation as in Definition \ref{definition:quasi-F-stable-sections}, 
suppose that $(W,X+B_W)$ is a plt pair, $K_X+\Delta = (K_W+X+B_W)|_X$, and $L=L_W|_X$ for a Weil divisor $L_W$ on $W$. 
Then $\rdown{\Delta}=0$ and $\Delta \geq \{L\}$. 
In particular, $\lceil L - \Delta \rceil = \lfloor L \rfloor$ and 
\[
q^eS^0_{n}(X,\Delta;L) = q^eS^0_{n}(X,\Delta-\{L\};\rdown{L}) \subseteq H^0(X,\cO_X(\rdown{L})).
\]
\end{lemma}

\begin{proof}
See \cite[Lemma 3.36]{KTTWYY1}. 
\end{proof}

\begin{proposition} \label{prop:qS^0-for-quasi-F-split}
In the situation of Setting \ref{setting:most-general-foundations-of-log-quasi-F-splitting}, let $L$ and $\Delta$ be $\bQ$-divisors on $X$. Suppose that $(X,\{\Delta-L\})$ is $n$-quasi-$F^{e}$-split. Then
\[
q^{e}S^0_n(X,\Delta;L) = H^0(X,\cO_X(\rup{L-\Delta})).
\]
\end{proposition}

\begin{proof}
Since $(X,\{\Delta-L\})$ is $n$-quasi-$F^e$-split, 
the following map is surjective by  Proposition \ref{p-IQFS-L-surje}: 
\[
H^0(X, (Q^e_{X, \Delta -L, n})^*) 
\xrightarrow{H^0(X, (\Phi^e_{X, \Delta -L, n})^*)}
H^0(X,\cO_X(\rdown{\Delta -L})^*) =H^0(X, \MO_X( \rup{L-\Delta})). 
\]
Then the required equality follows from Definition \ref{definition:quasi-F-stable-sections}. 
\end{proof}

Analogously, we make an adjoint definition. 

\begin{definition}\label{d-pureQFS-sections}
In the situation of Setting \ref{setting:most-general-foundations-of-log-pure-quasi-F-splitting}, let $L$ and $B$ be $\bQ$-divisors on $X$. 
Recall that we have $W_n\MO_X$-module homomorphisms 
\begin{align*}
\Phi^{S, e}_{X, S+B -L, n}\colon \quad &\cI_S({S+B} -L) \to Q^{S, e}_{X, S+B -L, n}\\  
(\Phi^{S, e}_{X, S+B -L, n})^*\colon \quad &(Q^{S, e}_{X, S+B -L, n})^* \to 
\cI_S(S+B -L)^* = \MO_X(\rup{L-B}), 
\end{align*}
where the lower one $(\Phi^{S, e}_{X, S+B -L, n})^*$ is obtained by applying 
$(-)^* =\cHom_{W_n\MO_X}(-, W_n\omega_X(-K_X))$ to the upper one $\Phi^{S, e}_{X, S+B -L, n}$. 
We define $q^eS^0_{n, \adj}(X, S+B;L)$ by 
\[
q^eS^0_{n, \adj}(X, S+B;L) := \image\big( 
H^0(X, (Q^{S, e}_{X, S+B -L, n})^*) \xrightarrow{H^0(X, (\Phi^{S, e}_{X, S+B -L, n})^*)} H^0(X, \MO_X(\rup{L-B}))\big).
\]
\end{definition}



 \begin{lemma}\label{l-qS^0-adj-nonadj}
In the situation of {Setting \ref{setting:most-general-foundations-of-log-pure-quasi-F-splitting}},  
 let $B$ be a $\Q$-divisor on $X$  such that $S \not\subseteq \Supp B$ and $\rdown{B}=0$. 
Let $L$ be a Weil divisor on $X$. 
Then the following inclusions hold {for every rational number $0 \leq t < 1$}: 
\[
q^{e}S^0_{n, \adj}(X, S+B; L) \subseteq 
q^{e}S^0_n(X, {tS+} B; L) \subseteq H^0(X, \MO_X(L)).
\]
\end{lemma}

\begin{proof}
 We have the following commutative diagram: 
\[
\hphantom{\MO_X({tS+}B-L)} \begin{tikzcd}[row sep=1cm,column sep=2cm]
W_n\MO_X({tS+}B-L) \arrow{r}{ F^e_1}\arrow[hookrightarrow]{d}{i} & F^e_*W_n\MO_X(p^e({tS+}B-L)) \arrow[hookrightarrow]{d}{j}\\
W_n\cI_S(S+B-L)\arrow{r}{F^e_2} \arrow{d}{R^{n-1}}& F^e_*W_n\cI_S({p^e}(S+B-L))\\
\mathllap{\MO_X({tS+}B-L) =}\  \cI_S(S+B-L),
\end{tikzcd}
\]
where $i$ and $j$ are the natural inclusions, 
{and $F^e_1$ and $F^e_2$ denote the $e$-th iterated Frobenii}. 
Since $Q^e_{X, {tS+}B-L, n}$ and $Q^{S, e}_{X, S+B-L, n}$ are the pushouts of 
$(F^e_1, R^{n-1} \circ i)$  and $(F^e_2, R^{n-1})$ respectively, {we get a map $Q^e_{X, {tS+}B-L, n} \to Q^{S, e}_{X, S+B-L, n}$ sitting inside the following diagram}
\begin{center}
\begin{tikzcd}[column sep = huge]
{\MO_X(-L)} \arrow[bend right = 30]{rr}{\Phi^{S, e}_{X, S+B-L, n}} \arrow{r}{\Phi^e_{X, {tS+}B-L, n}} &  Q^e_{X, {tS+}B-L, n} \arrow{r} & Q^{S, e}_{X, S+B-L, n}.
\end{tikzcd}    
\end{center}
Now, by applying $H^0(X, (-)^*) = \Hom_{W_n\cO_X}(-, W_n\omega_X(-K_X))$, we get a factorisation of 
{$H^0(X, (\Phi^{S, e}_{X, S+B-L, n})^*)$} as follows: 
\[
H^0(X, (Q^{S, e}_{X, S+B-L, n})^{*}) \to H^0(X, (Q_{X, {tS+}B-L, n}^{e})^*) 
\xrightarrow{H^0(X, (\Phi^e_{X, tS+B-L, n})^*)} H^0(X, \MO_X(L)), 
\]
which immediately implies the statement of the lemma.
\end{proof}

Finally, we state a reformulation of quasi-$F$-stable sections via local cohomology. 

\begin{proposition} \label{prop:local-cohomology-description-of-S^0} 
Fix  {$n \in \Z_{>0}$ and $e \in \Z_{>0}$}. 
In the situation of Setting \ref{setting:foundations-of-log-quasi-F-splitting},  
{assume that $(R, \m)$ is a local ring. 
Take  $\bQ$-divisors $L$ and $\Delta$ on $X$.} 
Then
\[
q^eS^0_n(X,\Delta;L)^{\wedge} \simeq {\rm Im}\Big(H^d_{\fram}({X}, \mcO_X(\lfloor K_X+\Delta-L \rfloor)) \to H^d_{W_n\fram}({X}, Q^{e}_{X,K_X+\Delta-L,n})\Big)^{\vee},
\]
where $(-)^{\wedge}$ denotes $\fram$-completion, and $(-)^{\vee}$ denotes Matlis duality.

Similarly, in the situation of Setting \ref{setting:foundations-of-log-pure-quasi-F-splitting},
{assume that $(R, \m)$ is a local ring. 
Take  $\bQ$-divisors $L$ and $B$ on $X$.} 
Then 
\[
q^eS^0_{n,\adj}(X,S+B;L)^{\wedge} \simeq {\rm Im}\Big(H^d_{\fram}({X}, \mcO_X(\lfloor K_X+B-L \rfloor)) \to H^d_{W_n\fram}({X}, Q^{S,e}_{X,K_X+S+B-L,n})\Big)^{\vee}.
\]
\end{proposition}
\begin{proof}
This is immediate by Matlis duality (Lemma \ref{lem:Matlis-duality-for-highest-cohomology}) and definitions of quasi-$F$-stable sections (Definition \ref{definition:quasi-F-stable-sections} and Definition \ref{d-pureQFS-sections}). 
\end{proof}

\subsection{Quasi-$F^\infty$-splittings} \label{ss:iterative-and-perfection}

\begin{dfn}\label{d-IQFS2}
In the situation of Setting  \ref{setting:most-general-foundations-of-log-quasi-F-splitting}, let $\Delta$ be a $\bQ$-divisor on $X$ {satisfying $\rdown{\Delta} =0$}. 
We say that $(X, \Delta)$ is 
\begin{align*}
    \text{{\em $n$-quasi-$F^\infty$-split}}\ \ &\text{if\ \  $\forall_{e>0}$ the pair $(X,\Delta)$ is $n$-quasi-$F^e$-split} \\ 
    \text{{\em quasi-$F^\infty$-split}}\ \ &\text{if\ \ $\forall_{e>0}\, \exists_{n>0}$ s.t.\ the pair $(X,\Delta)$ is $n$-quasi-$F^e$-split}\\
    \text{{\em {uniformly} quasi-$F^\infty$-split}}\ \ &\text{if\ \  $\exists_{n>0}\, \forall_{e>0}$ the pair $(X,\Delta)$ is $n$-quasi-$F^e$-split}. 
\end{align*}
We define local variants ({i.e., $n$-quasi-$F^\infty$-pure,  quasi-$F^\infty$-pure, 
uniformly quasi-$F^\infty$-pure}) analogously to Definition \ref{d-IQFS}.
\end{dfn}
As we will see later, 
Calabi-Yau varieties (e.g., elliptic curves)
are never {uniformly}  quasi-$F^{\infty}$-split {unless they are $F$-split} 
(Proposition \ref{p-CY-unif-QFS}). 

\begin{definition} \label{def:all-stable-F^infty-sections}
In the situation of Setting  \ref{setting:most-general-foundations-of-log-quasi-F-splitting}, let $\Delta$ and $L$ be $\bQ$-divisors on $X$. We define the following subspaces of $H^0(X, \cO_X(\lceil L-\Delta\rceil))$:
\begin{align*}
    q^eS^0(X,\Delta;L) &:= \bigcup_{n=1}^{\infty} q^eS^0_n(X,\Delta;L) \\
    q^\infty S^0_n(X,\Delta;L) &:= \bigcap_{e>0} q^e S^0_n(X,\Delta;L) \\
    q^\infty S^0 (X,\Delta;L) &:= \bigcap_{e>0} \bigcup_{n>0} q^e S^0_n(X,\Delta;L) = \bigcap_{e>0} q^e S^0(X,\Delta;L)  \\
    q^\infty_{\rm uni} S^0 (X,\Delta;L) &:= \bigcup_{n>0} \bigcap_{e>0}  q^e S^0_n(X,\Delta;L) = \bigcup_{n>0} q^\infty S^0_n(X,\Delta;L).
\end{align*}
\end{definition}

Note that for $n \leq m$:
\begin{align*}
q^eS^0_n(X,\Delta;L) &\subseteq q^eS^0_m(X,\Delta;L),\quad \text{ and } \\
q^\infty S^0_n(X,\Delta;L) &\subseteq q^\infty S^0_m(X,\Delta;L). 
\end{align*}
Thus, if $X$ is projective over a Noetherian ring $R$, then $H^0(X, \cO_X(\lceil L-\Delta\rceil))$ is {a finitely generated $R$-module}, 
and so
we get the stabilisation:
\begin{align*}
q^eS^0(X,\Delta;L) &= q^eS^0_n(X,\Delta;L),\quad  \text{ and } \\
q^\infty_{\rm uni} S^0(X,\Delta;L) &= q^\infty S^0_n(X,\Delta;L).
\end{align*}
for $n \gg 0$. The stabilisation with respect to $e$ will be discussed in the next subsection.

\begin{remark}\label{r-QFS-summary}
By definition and the above stabilisation, 
the following holds for a $\Q$-divisor $\Delta$ with $\rdown{\Delta}=0$: 
\begin{align*}
 \text{$(X, \Delta)$ is $n$-quasi-$F^e$-split} &\iff q^eS^0_n(X,\Delta;\cO_X) = H^0(X,\cO_X);\\
 \text{$(X, \Delta)$ is quasi-$F^e$-split} &\iff q^eS^0(X,\Delta;\cO_X) = H^0(X,\cO_X);\\
 \text{$(X, \Delta)$ is $n$-quasi-$F^\infty$-split} &\iff q^\infty S^0_n(X,\Delta;\cO_X) = H^0(X,\cO_X);\\
 \text{$(X, \Delta)$ is quasi-$F^\infty$-split} &\iff q^\infty S^0(X,\Delta;\cO_X) = H^0(X,\cO_X);\\
 \text{$(X, \Delta)$ is {uniformly} quasi-$F^\infty$-split} & \iff q^\infty_{\rm uni} S^0(X,\Delta;\cO_X) = H^0(X,\cO_X).
\end{align*}
\end{remark}

\begin{remark}
As in Subsection \ref{ss:pure-quasi-F-split} and Definition \ref{d-pureQFS-sections}, we define pure and adjoint variants of Definition \ref{d-IQFS2} and Definition \ref{def:all-stable-F^infty-sections}. For sake of brevity, we do not spell these definitions out. 
\end{remark}

\begin{remark}\label{r-Q-perf}
In this remark, we sketch an interpretation of quasi-$F^\infty$-splittings using perfection. Define $Q^{\rm perf}_{X,\Delta,n}$ as the pushout sitting in the following diagram:
\begin{equation*}
\begin{tikzcd}
W_n\MO_X(\Delta) \arrow{r}{{\pi^*}} \arrow{d}{R^{n-1}} & \pi_*W_n\MO^{\rm perf}_{X}(\pi^*\Delta) \arrow{d} \\
\MO_X(\Delta) \arrow{r}{\Phi^{{\perf}}_{X, \Delta, n}}  &  \arrow[lu, phantom, "\usebox\pushoutdr" , very near start, yshift=0em, xshift=0.6em, color=black] Q^{\rm perf}_{X, \Delta, n},
\end{tikzcd}
\end{equation*}
where $X^{\rm perf}$ denotes the perfection of $X$ 
and $\pi \colon X^{\rm perf} \to X$ is the induced projection. 
Recall that the perfection $X^{\perf}$ is defined as the scheme $(X, \cO^{\perf}_X)$, 
where 
\[
\cO^{\perf}_X := \varinjlim_F  \MO_X := \varinjlim (\MO_X \xrightarrow{F} F_*\MO_X \xrightarrow{F} F_*^2\MO_X \xrightarrow{F} \cdots).
\]  
{It holds that $Q^{\rm perf}_{X, \Delta, n} \simeq \varinjlim_e Q^e_{X, \Delta, n}$.} 
We define $Q^{\rm perf}_{X, \Delta} := \varprojlim_n Q^{\rm perf}_{X, \Delta,n}$. 
One can check that $Q^{\rm perf}_{X, \Delta}$ is the pushout of  $\MO_X(\Delta) \leftarrow W\MO_X(\Delta) \to \pi_*W\MO_{X^{\rm perf}}(\pi^*\Delta)$, but we shall not need this in our paper.\\



In the situation of Setting \ref{setting:foundations-of-log-quasi-F-splitting}, we assume that $(R, \mathfrak m)$ is a local ring. 
Given $n \in \bZ_{>0}$,  we have that $(X,\Delta)$ is  $n$-quasi-$F^\infty$-split if and only if the natural map 
\[
H^d_{\m}(\cO_X(K_X+\Delta)) \to H^d_{W_n\m}(Q^{\rm perf}_{X,K_X+\Delta,n}) 
\]
is injective. 
Using Artiniaty of local cohomology, one can easily deduce from this that $(X,\Delta)$ is {uniformly} quasi-$F^\infty$-split if and only if the natural map 
\[
H^d_{\m}(\cO_X(K_X+\Delta)) \to H^d_{W\m}(Q^{\rm perf}_{X,K_X+\Delta}) 
\]
is injective.
\end{remark}

{When dealing with non-noetherian objects, one needs to be extra careful about distinguishing theorems stated for local rings and arbitrary rings, as taking infinite intersection may not commute with localisation. This will stop to be an issue after we show stabilisation  for quasi-$F^e$-stable sections in the next subsection {(Corollary \ref{cor-stab-q^e})}.}



\subsection{Stabilisation for quasi-$F^e$-stable sections}

In this subsection, we prove the stabilisation of quasi-$F^e$-stable sections 
{(Corollary \ref{cor-stab-q^e})}. For ease of notation, we will work in the framework of local cohomology, but the proofs can be restated in terms of usual trace maps. As we often twist our divisors by other nef divisors, we shall consider all divisors of the form $rA + kN$ and establish necessary bounds to be independent of the choice of $r>0$ and $k \geq 0$. 

We start with the following lemma generalising Lemma \ref{lemma:stabilisation_S^0} to the {Witt} case. 

\begin{lemma} \label{lem:Witt-stabilisation}
In the situation of Setting \ref{setting:foundations-of-log-quasi-F-splitting}, 
fix $n 
\in \bZ_{>0}$, 
let $A$ be an ample {$\Q$-Cartier} $\bQ$-divisor, and let $N$ be a nef {$\Q$-Cartier} $\bQ$-divisor. Assume that $X$ is divisorially Cohen-Macaulay.
Then there exists $e_0>0$  such that for {every integer} $r \in \Z_{>0}${, every integer} $k \geq 0$, 
and every maximal ideal $\fram$ of $R$,
the kernel of
\[
{F^e \colon} H^d_{W_n\fram}({X}, W_n\cO_X(-(rA + kN))) \to H^d_{W_n\fram}({X}, F^{e}_*W_n\cO_X(-p^{e}(rA + kN)))
\]
is independent of the choice of $e \geq e_0$. 
\end{lemma}

\begin{proof}
{By abuse of notation}, we replace {$rA+kN$ by $A$}  in the following proof. 
All the bounds we {shall} pick are clearly independent of $r$ and $k$, 
as they come from Lemma \ref{lemma:stabilisation_S^0} and {the} Fujita vanishing {theorem}.

Let $e^{\dagger}_0>0$ be the bound coming from Lemma \ref{lemma:stabilisation_S^0}. By Matlis duality (Lemma \ref{lem:Matlis-duality-for-highest-cohomology}), we can assume that the kernel of 
\begin{equation} \label{eq:trace-Witt-Frobenius-stabilise-base-induction}
H^d_\fram(X, \cO_X(-A)) \to H^d_\fram(X, F^{e}_*\cO_X(-p^{e}A))
\end{equation}
{is independent of the choice of} $e \geq e^\dagger_0$. 
By increasing $e^{\dagger}_0$ further, we may also assume that 
\begin{equation} \label{eq:quasi-stabilisation-fujita}
H^{d-1}_\fram({X}, F^{e}_*\cO_X(-p^{e}A)) = 0
\end{equation}
 for all $e \geq e^\dagger_0$. This is possible, because $X$ is divisorially Cohen-Macaulay, and so by Matlis duality we need the vanishing of $H^1(X, \mcO_X(K_X+\rup{p^{e}A}))=0$ which follows by the Fujita vanishing theorem (cf.\ Lemma \ref{l-easy-Serre}).


Fix $n \in \bZ_{>0}$. By induction on $n$, we may assume that the statement of the lemma is valid 
{after replacing $n$ by $n-1$.} 
Specifically, there exists an integer $e^{\dagger\dagger}_0>0$ such that the kernel of
\begin{equation} \label{eq:trace-Witt-Frobenius-stabilise-induction}
F^e \colon H^d_{W_{n-1}\fram}({X}, W_{n-1}\cO_X(-A)) \to H^d_{W_{n-1}\fram}({X}, F^{e}_*W_{n-1}\cO_X(-p^{e}A))
\end{equation}
stabilises for all $e \geq e^{\dagger\dagger}_0$. 
{Note that} the same holds by $A$ replaced by any $rA + kN$ for integers $r>0$ and $k\geq 0$ with the same bound $e^{\dagger\dagger}_0$.




We will prove the statement of the lemma for $e_0 := e^\dagger_0 + e^{\dagger\dagger}_0$. Take $e', e$ such that $e' \geq e \geq e_0$. It is enough to show the following claim.
\begin{claim*} Suppose that 
\[
\zeta \in {\rm Ker}(H^d_{W_n\fram}(X, W_{n}\cO_X(-A)) \to H^d_{W_{n}\fram}(X, F^{e'}_*W_{n}\cO_X(-p^{e'}A))).
\]
Then
\begin{equation} \label{eq:zeta-kernel-bigger}
\zeta \in {\rm Ker}(H^d_{W_n\fram}(X, W_{n}\cO_X(-A)) \to H^d_{W_{n}\fram}({X}, F^{e}_*W_{n}\cO_X(-p^{e}A))).
\end{equation}
\end{claim*}

In what follows, we prove {the above claim}. 
{
We set $H^d_{\m}(-) := H^d_{W_n\m}(X, -)$ by abuse of notation.} 
{The exact sequence 
\[
0 \to F_*W_{n-1}\MO_X(pD) \xrightarrow{V} W_n\MO_X(D) \xrightarrow{R^{n-1}} \MO_X(D) \to 0 
\]
induces} the following {commutative} diagram:
{\small\[
\begin{tikzcd}[column sep = tiny]
0 \arrow{r} & H^d_{\fram}(F^{e'+1}_*W_{n-1}\mcO_X(-p^{e'+1}A)) \arrow{r} &  H^d_{\fram}(F^{e'}_*W_n\mcO_X(-p^{e'}A)) \arrow{r} & H^d_{\fram}(F^{e'}_*\mcO_X(-p^{e'}A)) \arrow{r} & 0\\
0 \arrow{r} & H^d_{\fram}(F^{e+1}_*W_{n-1}\mcO_X(-p^{e+1}A)) \arrow{r} \arrow{u} &  H^d_{\fram}(F^{e}_*W_n\mcO_X(-p^{e}A)) \arrow{r}  \ar{u} & H^d_{\fram}(F^{e}_*\mcO_X(-p^{e}A)) \ar{u} \arrow{r} & 0\\
0 \arrow{r} & H^d_{\fram}(F^{e^\dagger_0+1}_*W_{n-1}\mcO_X(-p^{e^\dagger_0+1}A))  \arrow{u}{(\dagger\dagger)} \arrow{r}{\psi} &  H^d_{\fram}(F^{e^\dagger_0}_*W_n\mcO_X(-p^{e^\dagger_0}A)) \arrow{u} \arrow{r} & H^d_{\fram}(F^{e^\dagger_0}_*\mcO_X(-p^{e^\dagger_0}A)) \arrow{u} \arrow{r} & 0\\
 & H^d_{\fram}(F_*W_{n-1}\mcO_X(-pA)) \arrow{r} \arrow{u} &  H^d_{\fram}(W_n\mcO_X(-A)) \arrow{u}{\theta} \arrow{r}{\varphi} & H^d_{\fram}(\mcO_X(-A)) \arrow{u}{(\dagger)} \arrow{r} & 0.
\end{tikzcd}
\]}
\!{Note that each horizontal sequence is exact by (\ref{eq:quasi-stabilisation-fujita})}
{the above claim} 
now follows immediately by diagram chase. Specifically, by (\ref{eq:trace-Witt-Frobenius-stabilise-base-induction}), $\varphi(\zeta)$ lies in the kernel of $(\dagger)$. Therefore, $\theta(\zeta) = \psi(\zeta')$ for some
\[
\zeta' \in H^d_{\fram}(F^{e^\dagger_0+1}_*W_{n-1}\mcO_X(-p^{e^\dagger_0+1}A)).
\]
Then, by (\ref{eq:trace-Witt-Frobenius-stabilise-induction}) {and 
$(e+1)-(e^\dagger_0 +1)  = e -e^\dagger_0 \geq e^{\dagger\dagger}_0$}, 
$\zeta'$ lies in the kernel of $(\dagger\dagger)$. 
Therefore, (\ref{eq:zeta-kernel-bigger}) holds, which concludes 
{the proofs of {the above claim} 
and Lemma \ref{lem:Witt-stabilisation}}. 
\end{proof}

\begin{proposition} \label{prop:qeS0-stabilisation} \label{l-stab-q^e}
In the situation of Setting \ref{setting:foundations-of-log-quasi-F-splitting}, fix $n\in \bZ_{>0}$, let $A$ be an ample {$\Q$-Cartier} $\bQ$-divisor, and let $N$ be a nef {$\Q$-Cartier} $\bQ$-divisor. 
Assume that $X$ is divisorially Cohen-Macaulay. 
{Then there exists $e_0 \in \Z_{>0}$ such that for
\begin{itemize}
    \item {every} $\bQ$-divisor $\Delta$, and 
    \item all integers $r,k,e$
satisfying $r \geq 1$, $k \geq 0$, and $e \geq e_0$,
\end{itemize}
the following {equality} 
holds:
\[
q^{e}S^0_n(X,\Delta; L) =q^{e_0}S^0_n(X,\Delta; L),
\]
where $L := K_X + \Delta + rA + kN$.}
\end{proposition}

\noindent {In other words, $q^{\infty}S^0_n(X,\Delta; L) =q^{e}S^0_n(X,\Delta; L)$ for $e \geq e_0$.}

\begin{proof}
{By abuse of notation,} we replace $rA + kN$ by $A$. 
In particular, $L=K_X+\Delta+A$.
All the bounds in this proof are clearly independent of the choice of integers $r \geq 1$ and $k \geq 0$.

Let $\fram$ be a maximal ideal of $R$. We shall prove the following statement: 
there exists $e_0>0$ such that the kernel of
\begin{equation} \label{eq:stabilisation-auxiliary-qinfty}
H^d_{\fram}(X, \cO_X(-A))  \to H^d_{W_n\fram}(X, Q^{e}_{X,-A,n}).
\end{equation}
is independent of the choice of $e\geq e_0$. By Matlis duality (Proposition \ref{prop:local-cohomology-description-of-S^0}) and Definition \ref{definition:quasi-F-stable-sections}, this implies the statement of the proposition: {$q^eS^0_n(X, \Delta; L) = q^{e_0}S^0_n(X,\Delta; L)$} {for the case when $(R, \m)$ is a local ring}. 
Since the bound $e_0$ in Lemma \ref{lem:Witt-stabilisation} 
is independent of $\mathfrak m$, the general case follows
by localisation at all maximal ideals of $R$.


We pick $e_0>0$ as in Lemma \ref{lem:Witt-stabilisation} and {fix} 
$e \geq e_0$. 
Consider the following {commutative} diagram {in which each horizontal sequence is exact (Remark \ref{r-Q-B-big-diag})}:
\[
\begin{tikzcd}
0 \arrow{r} & F_*W_{n-1}\mcO_X(-pA) \arrow[r, "{V}"] \arrow{d}{=} &  W_n\mcO_X(-A) \arrow[d, rightarrow, "{F^{e}}"] 
\arrow[r, "{R^{n-1}}"] & \cO_X(-A) \arrow[d, rightarrow] \arrow{r} & 0\\
0 \arrow{r} & F_*W_{n-1}\mcO_X(-pA) \arrow[r, "{VF^{e}}"] &  F^{e}_*W_n\mcO_X(-p^{e}A) \arrow{r} & Q^{e}_{X,-A,n} \arrow{r} & 0. 
\end{tikzcd}
\]
By applying 
cohomology {$H^i_{W_n\m}(-) := H^i_{W_n\m}(X, -)$}, we get
{\small
\[
\begin{tikzcd}[column sep=small]
0 \arrow{r} & H^d_{W_{n-1}\fram}(F_*W_{n-1}\mcO_X(-pA))/M \arrow{r} 
\arrow[d, rightarrow, two heads] &  H^d_{W_n\fram}(W_n\mcO_X(-A)) \arrow[d] \arrow{r} & H^d_{\fram}(\cO_X(-A)) \arrow[d] \arrow{r} & 0\\
0 \arrow{r} & H^d_{W_{n-1}\fram}(F_*W_{n-1}\mcO_X(-pA))/N \arrow{r} &  H^d_{W_n\fram}(F^{e}_*W_n\mcO_X(-p^{e}A)) \arrow{r} & H^d_{W_n\fram}(Q^{e}_{X,-A,n}) \arrow{r} & 0.
\end{tikzcd}
\]
}
\!\!Here $M$ and $N$ are the images of $H^{d-1}_{\fram}(\cO_X(-A))$ and $H^{d-1}_{W_n\fram}(Q^{e}_{X,-A,n})$, respectively. 
By the snake lemma, 
\[
K^e := \Ker(H^d_{W_n\fram}(W_n\mcO_X(-A)) \to H^d_{W_n\fram}(F^{e}_*W_n\mcO_X(-p^{e}A)))
\]
maps surjectively onto the kernel of (\ref{eq:stabilisation-auxiliary-qinfty}).  
Since 
${K^e}$
stabilises for $e \geq e_0$ by Lemma \ref{lem:Witt-stabilisation}, 
the kernel of (\ref{eq:stabilisation-auxiliary-qinfty}) stabilises as well for 
${e \geq e_0}$. 
\end{proof}

In particular, we obtain the following stabilisation result. We use quantifiers to emphasise the subtlety of how various indices depend on one another. We do not know a similar statement for $q^\infty S^0(X, \Delta; L)$.


\begin{corollary}
\label{cor-stab-q^e}
In the situation of {Setting  \ref{setting:foundations-of-log-quasi-F-splitting}}, let $\Delta$ and $L$ be $\bQ$-divisors on $X$ such that $L-(K_X+\Delta)$ is an ample {$\Q$-Cartier} $\bQ$-divisor. 
{Assume that $X$ is divisorially Cohen-Macaulay.} 
Then {the following hold.}
\begin{alignat}{3}
\tag{1} \forall_{n>0}\ \exists_{e_0>0}\ \forall_{e \geq e_0}\quad  && q^\infty S^0_n(X, \Delta; L)  &\,=\, q^eS^0_n(X, \Delta; L). \\
\tag{2} \exists_{n_0>0}\ \forall_{n\geq n_0}\ \exists_{e_0>0}\ \forall_{e\geq e_0} \quad  && q^\infty_{\rm uni} S^0(X, \Delta; L)  &\,=\, q^eS^0_n(X, \Delta; L).
\end{alignat}
\end{corollary}

\begin{proof}
Statement (1) is clear by 
{Definition \ref{def:all-stable-F^infty-sections} and} 
Proposition \ref{prop:qeS0-stabilisation}.


{Let us show (2). By Definition \ref{def:all-stable-F^infty-sections}, we have an ascending chain 
\[
q^{\infty}S^0_1(X, \Delta; L) \subseteq 
q^{\infty}S^0_2(X, \Delta; L) \subseteq \cdots \subseteq H^0(X, \MO_X(\rup{L-\Delta}))
\]
of $R$-submodules of $H^0(X, \MO_X(\rup{L-\Delta}))$. 
Since $X$ is projective over an Noetherian ring $R$,} 
there exists $n_0>0$ such that for all $n \geq n_0$,  
\begin{equation*} 
q^\infty_{\rm uni} S^0(X, \Delta; L) = \bigcup_{n=1}^{\infty} 
q^\infty S^0_n(X, \Delta; L)
= q^\infty S^0_n(X, \Delta; L). 
\end{equation*}
{Fix $n \geq n_0$.} 
Then there exists $e_0$ (dependent on $n$) such that for all $e \geq e_0$, 
\begin{equation*} 
q^\infty S^0_n(X, \Delta; L) = q^e S^0_n(X, \Delta; L)
\end{equation*}
by (1).
\end{proof}

\begin{proposition} \label{prop:Serre-type-result-qS0} 
{In the situation of Setting \ref{setting:foundations-of-log-quasi-F-splitting}, fix $n\in \bZ_{>0}$, let $A$ be an ample {$\Q$-Cartier} $\bQ$-divisor, and let $N$ be a nef {$\Q$-Cartier} $\bQ$-divisor. 
Assume that $X$ is divisorially Cohen-Macaulay. 
Then there exists $r_0 \in \Z_{>0}$ such that for
\begin{itemize}
    \item all integers $r,k,e$
satisfying $r \geq r_0$, $k \geq 0$, and $e >0 $, and
    \item {every} $\bQ$-divisor $\Delta$ such that 
    $(X,\Delta)$ is $n$-quasi-$F^\infty$-pure and $L := K_X + \Delta + rA + kN$ is {a \textcolor{black}{Weil} divisor}, 
\end{itemize}
the following 
{equality} holds:
\[
q^{e}S^0_n(X,\Delta; L) = H^0(X, \cO_X(L)).
\]}
\end{proposition}
\noindent {
{Since} we are assuming that $(X,\Delta)$ is $n$-quasi-$F^\infty$-pure, this implicitly requires that $\rdown{\Delta}=0$, and so $\Delta = \{-(rA+kN)\}$. Moreover, $L = K_X + \lceil rA + kN \rceil$. }
\begin{proof}
By Proposition \ref{prop:qeS0-stabilisation}, we may pick $e_0$ such that 
\[
{q^{e}S^0_n(X,\Delta; L) = q^{e_0}S^0_n(X, \Delta; L)}
\]
for all integers $e \geq e_0$, $r \geq 1$, and $k \geq 0$. 
{On the other hand, we have 
$q^{e}S^0_n(X,\Delta; L) \supseteq q^{e_0}S^0_n(X, \Delta; L)$ when $e \leq e_0$. Hence} it is enough to prove the statement of the proposition for $e = e_0$.

By Definition \ref{definition:quasi-F-stable-sections},
\begin{multline*}
q^{e_0}S^0_n(X,\Delta;{L}) \\ := \image\Big( 
\Hom_{{W_n\MO_X}}
(Q^{e_0}_{X,-({K_X+} rA + kN),n}, W_n\omega_X{{(-K_X)}}) 
\to H^0(X, \cO_X(K_X+\lceil rA + kN \rceil))\Big).  
\end{multline*}
Define
\[
\cG_{r,k} := {\rm Ker}\Big( 
\cHom_{{W_n\MO_X}}(Q^{e_0}_{X,-({K_X+} rA + kN),n}, W_n\omega_X{{(-K_X)}}) 
\xrightarrow{(\dagger)} \cO_X(K_X+\lceil rA + kN \rceil)\Big).  
\]

\begin{claim} There exists an integer $r_0 > 0$ such that $H^1(X, \cG_{r,k})=0$ for all $r \geq r_0$ and $k \geq 0$.
\end{claim}

{ Assuming this claim, 
we now finish the proof of Proposition \ref{prop:Serre-type-result-qS0}.} 
Note that $(\dagger)$ is surjective as {$(X,\Delta)$ is $n$-quasi-$F^{e_0}$-pure (see Proposition \ref{p-IQFS-L-surje} for $D= -(K_X + \lceil rA + kN \rceil)$ and $\Delta + D = -(K_X+rA + kN)$)}. Therefore, the claim immediately implies that $q^{e_0}S^0_n(X,\Delta;K_X+\lceil rA + kN \rceil) = H^0(X, \cO_X(K_X+\lceil rA + kN \rceil))$ for $r \geq r_0$, so it is enough to prove 
{the above} claim.

To this end, take 
{$m_0 \in \Z_{>0}$ such that $m_0A$ and $m_0 N$ are Cartier}.
Take the integers $r'$ and $k'$ defined by 
$r' := r \ {\rm mod}\ {m_0}$ and $k' := k \ {\rm mod}\ {m_0}$ 
{(cf.\ Subsection \ref{ss:notation}({\ref{ss-nota-mod}}))}. 
{Since $(r-r')A$ and $(k-k')N$ are Cartier, we get} 
\[
\cG_{r,k} = \cG_{r',k'} \otimes_{W_n\MO_X} W_n\cO_X((r-r')A + (k-k')N).
\]
Thus there exists $r_0$ (independent of {$(r, k)$}) 
such that
\[
H^1(X, \cG_{r',k'} \otimes_{W_n\MO_X} W_n\cO_X((r-r')A + (k-k')N))=0
\]
for {every} $r\geq r_0$ 
{and every $k\geq 0$}  
by {the} Fujita vanishing {theorem} (\cite[Theorem 1.5]{keeler03}),  
{because} there are only finitely many {possibilities} 
for $\cG_{r',k'}$.
\qedhere


\end{proof}

Finally, we state a reformulation of quasi-$F$-stable sections via local cohomology. 

\begin{proposition} \label{prop:local-cohomology-description-of-S^0-perf} 
Fix   {$n \in \Z_{>0}$}. 
In the situation of Setting \ref{setting:foundations-of-log-quasi-F-splitting},  
{assume that $(R, \m)$ is a local ring. 
Take  $\bQ$-divisors $L$ and $\Delta$ on $X$.} 
Then 
\[
H^d_{\m}(X, Q^{\perf}_{X, K_X+\Delta-L, n})^{\vee} \simeq \varprojlim_e   \left(H^0(X, (Q^{e}_{X, \Delta-L, n})^*)^{\wedge}\right) 
\]
Furthermore, if $X$ is divisorially Cohen-Macaulay and $L-(K_X+\Delta)$ is ample, then we have
\[
{\q} 
S^0_n(X,\Delta;L)^{\wedge} \simeq {\rm Im}\Big(H^d_{\fram}({X}, \mcO_X(\lfloor K_X+\Delta-L \rfloor)) \to H^d_{W_n\fram}({X}, Q^{\perf}_{X,K_X+\Delta-L,n})\Big)^{\vee},
\]
where $(-)^{\wedge}$ denotes $\fram$-completion and $(-)^{\vee} 
:= \Hom_{W_nR}(-, E)$ denotes Matlis duality.
\end{proposition}



{\begin{proof}
 
The first assertion holds by the following argument: 
\begin{align}\label{e2-qinfty-Matlis}
H^d_{\m}(X, Q^{\perf}_{X, K_X+\Delta-L, n})^{\vee}
&= \Hom_{W_nR}(H^d_{\m}(X, \varinjlim_e Q^{e}_{X, K_X+\Delta-L, n}), E)\notag\\
&\simeq \Hom_{W_nR}( \varinjlim_e H^d_{\m}(X, Q^{e}_{X, K_X+\Delta-L, n}), E)\notag\\
&\simeq \varprojlim_e  \Hom_{W_nR}(H^d_{\m}(X, Q^{e}_{X, K_X+\Delta-L, n}), E)\\
&\simeq \varprojlim_e  \left( \Hom_{W_n\MO_X}(Q^{e}_{X, K_X+\Delta-L, n}, W_n\omega_X)^{\wedge}\right)\notag\\
&\simeq \varprojlim_e   \left(H^0(X, (Q^{e}_{X, \Delta-L, n})^*)^{\wedge}\right)\notag. 
\end{align}

Let us show the second assertion.
Take the image: 
\begin{equation} \label{e1-desc-S^0-perf}
H^d_{\m}(X, \MO_X(K_X+\Delta-L))
\twoheadrightarrow 
I_n 
\hookrightarrow 
H^d_{W_n\m}(X, Q^{\perf}_{X, K_X+\Delta-L, n}). 
\end{equation}
It suffices 
to show that $I_n^{\vee} \simeq   \q S^0_n(X,\Delta;L)^{\wedge}$. Since Matlis duality turns colimits into limits, we get that
\[
I_n^{\vee} = (\varinjlim_e I^e_n)^{\vee}  = \varprojlim_e (I^e_n)^{\vee},
\]
where $I^e_n$ is defined as the image in 
\begin{equation} \label{e1-desc-S^0-perf2}
H^d_{\m}(X, \MO_X(K_X+\Delta-L))
\twoheadrightarrow 
I^e_n 
\hookrightarrow 
H^d_{W_n\m}(X, Q^{e}_{X, K_X+\Delta-L, n}).
\end{equation}
By Lemma \ref{lem:Matlis-duality-for-highest-cohomology}, the Matlis duality functor {$(-)^{\vee} = \Hom_{W_nR}(-, E)$} yields   
\begin{align*}
H^d_{\m}(X, \MO_X(K_X+\Delta-L))^{\vee} &\simeq H^0(X, \MO_X(\rup{L-\Delta}))^{\wedge} \text{ and } \\
H^d_{\m}(X, Q^{e}_{X, K_X+\Delta-L, n})^{\vee} 
&= \Hom_{W_nR}(H^d_{\m}(X, Q^{e}_{X, K_X+\Delta-L, n}), E) \\ 
&\simeq  \Hom_{{W_n\MO_X}}(Q^{e}_{X, K_X+\Delta-L, n}, W_n\omega_X)^{\wedge} \\ 
&\simeq H^0(X, (Q^{e}_{X, \Delta-L, n})^*)^{\wedge}, 
\end{align*}
{where $(-)^* :=\mathcal Hom_{W_n\MO_X}(-, W_n\omega_X(-K_X))$.} 
{By applying $(-)^{\vee} = \Hom_{W_nR}(-, E)$ to (\ref{e1-desc-S^0-perf2}), we get} 
\[
H^0(X, \cO_X(\rup{L-\Delta}))^{\wedge} \hookleftarrow (I^e_n)^{\vee} 
\twoheadleftarrow 
{\displaystyle H^0(X, (Q^{e}_{X, \Delta-L, n})^*)^{\wedge}}.
\]

Now by 
Corollary \ref{cor-stab-q^e}(1), there is an integer $e_0$ such that $q^{e}S_n^0(X, \Delta; L) = q^{e_0}S_n^0(X, \Delta; L)$ for any $e \geq e_0$.
By the above argument, we get $(I_n^{e})^{\vee}=(I_n^{e_0})^{\vee}$ for any $e \geq e_0$.
Hence we have
\[
I_n^{\vee}=\bigcap_{e>0}(I_n^e)^{\vee}=(I_n^{e_0})^{\vee}=q^{e_0}S_n^0(X, \Delta; L)^{\wedge}=q^{\infty}S_n^0(X, \Delta; L)^{\wedge}.
\]
This completes the proof.
\qedhere
\end{proof}

\section{Quasi-$F$-regularity and quasi-+-regularity}

\subsection{Quasi-$F$-regularity}



\begin{dfn}\label{d-QFR}
In the situation of Setting  \ref{setting:most-general-foundations-of-log-quasi-F-splitting}, let $\Delta$ be a $\bQ$-divisor on $X$. 
Take $n \in \Z_{>0}$. 
\begin{enumerate}
\item 
We say that $(X, \Delta)$ is {\em globally $n$-quasi-$F$-regular} 
if 
\begin{enumerate}
\item $\rdown{\Delta}=0$, and 
\item 
given an effective $\Q$-divisor $E$ on $X$, 
there exist $\epsilon \in \Q_{>0}$  
such that $(X, \Delta + \epsilon E)$ is $n$-quasi-$F^e$-split 
for every $e \in \Z_{>0}$. 
\end{enumerate}
\item 
We say that $(X, \Delta)$ is {\em globally quasi-$F$-regular} 
if $(X, \Delta)$ is globally $m$-quasi-$F$-regular for some $m \in \Z_{>0}$. 
\item 
We say that $(X, \Delta)$ is {\em {locally} $n$-quasi-$F$-regular} 
if there exists an  open cover $X = \bigcup_{i \in I} X_i$ such that 
$(X_i, \Delta|_{X_i})$ is globally $n$-quasi-$F$-regular for every $i \in I$. 
\item 
We say that $(X, \Delta)$ is {\em {locally} quasi-$F$-regular} 
if $(X, \Delta)$ is {locally} $m$-quasi-$F$-regular for some $m \in \Z_{>0}$. 
This condition is equivalent to the following: 
there exists an  open cover $X = \bigcup_{i \in I} X_i$ such that $(X_i, \Delta|_{X_i})$ is globally quasi-$F$-regular for every $i \in I$. 
\end{enumerate}
\end{dfn}

\begin{dfn}\label{d-weak-QFR}
In the situation of Setting  \ref{setting:most-general-foundations-of-log-quasi-F-splitting}, let $\Delta$ be a $\bQ$-divisor on $X$. 
We say that $(X, \Delta)$ is {\em {feebly} globally quasi-$F$-regular} 
if 
\begin{enumerate}
\item $\rdown{\Delta}=0$, and 
\item 
given an effective $\Q$-divisor $E$ on $X$, 
there exist 
$n \in \Z_{>0}$ and $\epsilon \in \Q_{>0}$  
such that $(X, \Delta + \epsilon E)$ is $n$-quasi-$F^e$-split 
for every $e \in \Z_{>0}$. 
\end{enumerate}
We say that $(X, \Delta)$ is {\em {feebly} {locally} quasi-$F$-regular} 
if there exists an  open cover $X = \bigcup_{i \in I} X_i$ such that 
$(X_i, \Delta|_{X_i})$ is 
feebly
globally quasi-$F$-regular for every $i \in I$. 
\end{dfn}



\begin{rem}\label{r-infinite-vs-every}
If $(X, \Delta + \epsilon E)$ is $n$-quasi-$F^{e+1}$-split, 
then $(X, \Delta + \epsilon E)$ is $n$-quasi-$F^{e}$-split 
(Remark \ref{r-e-to-e+1}). 
Therefore, $(X, \Delta + \epsilon E)$ is $n$-quasi-$F^e$-split 
for every $e \in \Z_{>0}$ if and only if 
$(X, \Delta + \epsilon E)$ is $n$-quasi-$F^e$-split 
for infinitely many $e \in \Z_{>0}$. 
\end{rem}

In order to compare our definition of global quasi-$F$-regularity (Definition \ref{d-QFR}) with the global $F$-regularity, 
let us recall its definition. 

\begin{dfn}
Let $X$ be an integral normal $F$-finite Noetherian scheme 
and let $\Delta$ be an effective $\Q$-divisor. 
We say that $(X, \Delta)$ is {\em globally $F$-regular} 
if given an effective $\Q$-divisor $E$, 
there exists $e \in \Z_{>0}$ such that 
\[
F^e : \MO_X \to F^e_*\MO_X( (p^e-1)\Delta + E)
\]
splits.
\end{dfn}

Note that some authors use the splitting of 
\[
F^e : \MO_X \to F^e\MO_X( \ulcorner  (p^e-1)\Delta + E\urcorner)
\]
for the definition of the global $F$-regularity (cf.\ \cite[Definition 3.1]{schwedesmith10}). 
It is easy to see that two definitions are equivalent. 
The following result 
shows that our definition of global quasi-$F$-regularity (Definition \ref{d-QFR}) can be considered 
as an analogue of the global $F$-regularity. 

\begin{prop}\label{p-GFR-compare}
In the situation of Setting  \ref{setting:most-general-foundations-of-log-quasi-F-splitting}, let $\Delta$ be a $\bQ$-divisor on $X$ such that $\rdown{\Delta}=0$. 
Then the following are equivalent. 
\begin{enumerate}
\item $(X, \Delta)$ is globally $F$-regular, i.e., 
given an effective $\Q$-divisor $E$, 
there exists $e \in \Z_{>0}$ such that 
\[
F^e : \MO_X \to F_*^e\MO_X( (p^e-1)\Delta +E)
\]
splits. 
\item {$(X, \Delta)$ is globally $1$-quasi-$F$-regular. In other words,} given an effective $\Q$-divisor $E$, there exists $\epsilon \in \Q_{>0}$ such that 
$(X, \Delta + \epsilon E)$ is $1$-quasi-$F^e$-split 
for every $e \in \Z_{>0}$, i.e., 
\begin{equation}\label{e1-GFR-compare}
\Phi^e_{X, \Delta, 1}: 
\MO_X \to F_*^e\MO_X( p^e(\Delta +\epsilon E))
\end{equation}
splits for every $e \in \Z_{>0}$. 
\end{enumerate}
\end{prop}

\begin{proof}
Let us show $(2) \Rightarrow (1)$. 
Assume (2). 
Fix an effective $\Q$-divisor $E$. 
Take $\epsilon \in \Q_{>0}$ as in (2). 
Then we can find $e \in \Z_{>0}$ satisfying $1 \leq p^e \epsilon$, which implies 
\[
(p^e-1)\Delta +E \leq p^e(\Delta +\epsilon E). 
\]
Hence the composite $\MO_X$-module homomorphism 
\[
\Phi^e_{X, \Delta, 1}:  \MO_X  \xrightarrow{F^e =:\alpha} F_*^e\MO_X(
(p^e-1)\Delta +E) 
\hookrightarrow F_*^e\MO_X(p^e(\Delta +\epsilon E))
\]
splits, and hence also $\alpha$ splits. 
Thus (1) holds. 

\medskip

Let us show $(1) \Rightarrow (2)$. 
Assume (1). 
Fix $E$ an effective $\Q$-divisor. 
Then 
\[
(X, (1+\delta)(\Delta + \epsilon E))
\]
is globally $F$-regular for some $\epsilon, \delta \in \Q_{>0}$ \cite[Corollary 6.1]{schwedesmith10}. 
In particular, $(X, (1+\delta)(\Delta + \epsilon E))$ is globally sharply $F$-split in the sense of \cite[Definition 3.1]{schwedesmith10}. 
Then there exists  $d_1 \in \Z_{>0}$ such that 
\[
F^{dd_1} : \MO_X \to F_*^{dd_1}\MO_X( 
\rup{(p^{dd_1}-1)(1+\delta)(\Delta +\epsilon E)})
\]
splits for every $d \in \Z_{>0}$ \cite[Proposition 3.8(b)]{schwedesmith10}. 
Since we have 
$(p^{dd_1} -1)(1+\delta) \geq p^{dd_1}$ 
and 
$(p^{dd_1} -1)(1+\delta)\epsilon \geq 1$ for $d \gg 0$, 
we can find infinitely many $e \in \Z_{>0}$ such that 
(\ref{e1-GFR-compare}) splits. 
Therefore, 
(\ref{e1-GFR-compare}) splits for every $e \in \Z_{>0}$ 
(Remark \ref{r-infinite-vs-every}). 
\end{proof}

\subsection{Quasi-+-regularity}

\begin{dfn}\label{d-Q-for-plus}
In the situation of Setting \ref{setting:most-general-foundations-of-log-quasi-F-splitting}, 
let  $\Delta$ be a $\Q$-divisor. 
Let $f:Y \to X$ be a finite surjective morphism from an integral normal scheme $Y$. 
Recall that $f^*\Delta$ can be naturally defined as in \cite[the proof of Proposition 5.20]{KM98}. 
We define a $W_n\MO_X$-module homomorphism $\Phi^f_{X, \Delta, n}$ and 
a coherent $W_n\MO_X$-module $Q^f_{X, \Delta, n}$  on $X$ by the following pushout diagram: 
\begin{equation}\label{e1-Q-for-plus}
\begin{tikzcd}
W_n\MO_X(\Delta) \arrow[r, "f^*"] \arrow[d, "R^{n-1}"'] & f_*W_n\MO_Y(f^*\Delta) \arrow[d] \\
\MO_X(\Delta) \arrow[r, "\Phi^f_{X, \Delta, n}"] & Q^f_{X, \Delta, n}.
\end{tikzcd}
\end{equation}
We define $(\Phi^f_{X, \Delta, n})^*$ by applying 
$(-)^* := \cHom_{W_n\MO_X}(-, W_n\omega_X(-K_X))$ to $\Phi^f_{X, \Delta, n}$: 
\[
(\Phi^f_{X, \Delta, n})^*: 
(Q^f_{X, \Delta, n})^* \to 
\MO_X(\Delta)^*
\]
Note that $\MO_X(\Delta)^* = \cHom_{W_n\MO_X}(\MO_X(\Delta), W_n\omega_X(-K_X)) \simeq 
\cHom_{\MO_X}(\MO_X(\Delta), \omega_X(-K_X)) \simeq \MO_X(\rup{-\Delta})$. 
Set 
$B^f_{X,\Delta,n} := \Coker(W_n\MO_X(\Delta)\xrightarrow{f^*} f_*W_n\MO_Y(f^*\Delta))$. 
\end{dfn}

\begin{dfn}\label{d-Q+R}
In the situation of Setting \ref{setting:most-general-foundations-of-log-quasi-F-splitting}, 
let  $\Delta$ be a $\Q$-divisor on $X$. 
\begin{enumerate}
\item 
Take $n \in \Z_{>0}$ and a finite surjective morphism $f: Y \to X$ from an integral normal excellent scheme $Y$. 
We say that $(X, \Delta)$ is {\em globally $n$-quasi-$f$-regular} if 
$\rdown{\Delta} =0$ and the induced map 
\[
\hphantom{aaaaaa}\Hom_{W_n\MO_X}(Q^f_{X, \Delta, n}, W_n\omega_X(-K_X)) \xrightarrow{H^0(X, (\Phi^f_{X, \Delta, n})^*)}
\Hom_{W_n\MO_X}(\MO_X, W_n\omega_X(-K_X))
\]
is surjective. 
\item 
Given $n \in \Z_{>0}$, 
we say that $(X, \Delta)$ is {\em globally $n$-quasi-+-regular} if 
$\rdown{\Delta} =0$ and it is globally $n$-quasi-$f$-regular for every finite surjective morphism $f: Y \to X$ from an integral normal excellent scheme $Y$. 
We say that $(X, \Delta)$ is {\em globally quasi-+-regular} if 
$(X, \Delta)$ is globally $n$-quasi-+-regular for some $n \in \Z_{>0}$. 
\item 
We say that $(X, \Delta)$ is {\em feebly globally quasi-+-regular} if 
$\rdown{\Delta} =0$ and, for every finite surjective morphism $f: Y \to X$ from an integral normal excellent scheme $Y$, 
there exists $n \in \Z_{>0}$ such that $(X, \Delta)$ is $n$-quasi-
$f$
-regular. 
\end{enumerate}
We say that $(X, \Delta)$ is {\em $($feebly$)$ locally quasi-$+$-regular} 
(resp.\ {\em $($feebly$)$ locally $n$-quasi-$+$-regular}) 
if there exists an open cover $X = \bigcup_{i \in I} X_i$ such that 
$(X_i, \Delta|_{X_i})$ is {(feebly)  quasi-$+$-regular} 
(resp.\ {(feebly) $n$-quasi-$+$-regular}) for every $i \in I$. 
\end{dfn}



\begin{lem}\label{l-QFR-to-Q+R}
Let $g:Y \to Z$ be a finite surjective morphism of integral normal Noetherian $F$-finite $\F_p$-schemes. 
Assume that the induced field extension $j : K(Z) \hookrightarrow K(Y)$ is separable. 
Fix $n \in \Z_{>0}$. 
Then the following hold. 
\begin{enumerate}
\item $W_n(K(Z))$ is a local Artinian ring. 
\item 
The induced map $W_nj : W_n(K(Z)) \to W_n(K(Y))$ splits as a $W_n(K(Z))$-module homomorphism, 
i.e, there exists a $W_n(K(Z))$-module homomorphism $\varphi: W_n(K(Y)) \to W_n(K(Z))$ such that the composition 
\[
W_n(K(Z)) \xrightarrow{W_nj} W_n(K(Y)) \xrightarrow{\varphi} W_n(K(Z))
\]
is the identity map {\rm id}. 
\item 
Let $\varphi :  W_n(K(Y)) \to  W_n(K({Z}))$ 
be as in (2), i.e., 
a $W_n(K(Z))$-module homomorphism such that $\varphi \circ W_nj ={\rm id}$. 
Take a $\Q$-divisor $\Delta_Z$.  
Then there exists an effective Weil divisor $D_Z$ on $Z$  
which induces the following commutative diagram consisting of $W_n\MO_Z$-module homomorphisms for every $m \in \Z$
\[
\begin{tikzcd}
W_n(K(Z)) \arrow[r, "W_nj"] \arrow[rr, bend left, "{\rm id}"] &W_n(K(Y))\arrow[r, "\varphi"] & W_n(K(Z))\\
W_n\MO_Z(m\Delta_Z)  \arrow[r, "g^*"] \arrow[u, hook] & g_*W_n\MO_Y(mg^*\Delta_Z)  \arrow[u, hook] \arrow[r, "\psi_m"] & W_n\MO_Z(m\Delta_Z+D_Z), \arrow[u, hook] 
\end{tikzcd}
\]
where $W_n(K(Y))$ and $W_n(K(Z))$ denote the corresponding constant sheaves on $Z$, 
and $\psi_m : g_*W_n\MO_Y(mg^*\Delta_Z) \to W_n\MO_Z(m\Delta_Z + D_Z)$ is a $W_n\MO_Z$-module homomorphism. 
\end{enumerate}
\end{lem}

\begin{proof}
Let us show (1). 
Since $K(Z)$ is an $F$-finite Noetherian $\F_p$-algebra, 
$W_n(K(Z))$ is a Noetherian ring. 
As $\Spec W_n(K(Z))   \simeq \Spec K(Z)$  consists of one point, $W_n(K(Z))$ is a local artinian ring. 
Thus (1) holds.


Let us show (3) assuming (2). 
By removing the non-regular locus of $Z$, we may assume that $Z$ is regular. 
In particular, there exists $m_0 \in \Z_{>0}$ such that $m_0\Delta_Z$ is Cartier. 
Since $g_*W_n\MO_Y(mg^*\Delta_Z)$ is a coherent $W_n\MO_Z$-submodule of $W_n(K(Y))$, 
it is easy to find an effective Weil divisor $D_{Z, m}$, depending on $m$, satisfying 
\[\varphi(g_*W_n\MO_Y(mg^*\Delta_Z)) \subseteq W_n\MO_Z(m\Delta_Z + D_{Z, m}).
\]
Then the assertion holds for an effective Weil divisor $D_Z$ satisfying 
$D_{Z, m} \leq D_Z$ for every $0 \leq m \leq m_0-1$, 
because $\psi_m : g_*W_n\MO_Y(mg^*\Delta_Z) \to W_n\MO_Z(m\Delta_Z+D_Z)$ induces 
$\psi_{m \pm m_0} :  g_*W_n\MO_Y( (m\pm m_0) g^*\Delta_Z) \to W_n\MO_Z( (m\pm m_0) \Delta_Z+D_Z)$ 
by applying $(-) \otimes_{W_n\MO_Z} W_n\MO_Z(\pm m_0\Delta_Z)$.


Therefore, it is now enough to show (2). 
Set 
\[
K := K(Z)\qquad \text{and} \qquad L := K(Y). 
\]
Fix a $K$-linear basis $\alpha_1, \alpha_2, ..., \alpha_d$ of $L$ with $\alpha_1 = 1$, 
so that $L = \bigoplus_{i=1}^d K \alpha_i$. 
It suffices to show that 
$W_n(L) = \bigoplus_{i=1}^d W_n(K) \underline{\alpha_i}$, 
where each $\underline{\alpha_i} \in W_n(L)$ denotes the Teichm\"{u}ller lift $(\alpha_i, 0, ..., 0)$ of $\alpha_i$. 

First, we prove that 
\begin{equation}\label{e1-l-QFR-to-Q+R}
W_n(L) = \sum_{i=1}^d W_n(K) \underline{\alpha_i}. 
\end{equation}
To this end, it follows from (1) that $(W_n(K), V(W_{n-1}K))$ is a Noetherian complete local ring. 
By \cite[Theorem 8.4]{Matsumura},  
it is enough to show that 
\[
W_n(L)/ (V(W_{n-1}K) \cdot W_n(L))
\]
is generated by the images 
of $\underline{\alpha_1}, ..., \underline{\alpha_d}$. 
This follows from 
\[
W_n(L)/ (V(W_{n-1}K)\cdot W_n(L))
\simeq W_n(L) \otimes_{W_n(K)} (W_n(K)/V(W_{n-1}(K))) 
\]
\[
\simeq 
W_n(L) \otimes_{W_n(K)} W_1(K) \simeq W_1(L) =L
\]
where the last isomorphism follows from the fact that $K \to L$ is \'etale \cite[Ch. 0, Proposition 1.5.8]{illusie_de_rham_witt}. 
This completes the proof of (\ref{e1-l-QFR-to-Q+R}). 

Fix $e \in \Z_{>0}$ and set $q :=p^e$. 
It suffices to show that 
$\alpha_1^q, ..., \alpha_d^q$ are linearly independent over $K$. 
Consider the following commutative diagram consisting of field extensions: 
\[
\begin{tikzcd}
K \arrow[r, hook] \arrow[d, hook] & L \arrow[d, hook] \\
K^{1/q} \arrow[r, hook] & K^{1/q} \otimes_K L\arrow[r, hook, "\theta"] & L^{1/q}
\end{tikzcd}
\]
where $K^{1/q} \otimes_K L$ is a field, 
because $\Spec$ of the horizontal arrows are \'etale and $\Spec$ of the vertical arrows are universal homeomorphisms. 
Since $1 \otimes \alpha_1, ..., 1 \otimes \alpha_d \in K^{1/q} \otimes_K L$ is a $K^{1/q}$-linear basis of $K^{1/q} \otimes_K L$,  
it follows from $\theta(1 \otimes \alpha_i) =\alpha_i$ that 
$\alpha_1, ..., \alpha_d \in L^{1/q}$ are linearly independent over $K^{1/q}$. 
Take an equation 
\[
c_1 \alpha_1^q  + \cdots + c_d \alpha_d^q =0
\]
with $c_i \in K$. 
Apply $(-)^{1/q}$ inside $L^{1/q}$: 
\[
c_1^{1/q} \alpha_1  + \cdots + c_d^{1/q} \alpha_d =0. 
\]
Note that $c_i^{1/q} \in K^{1/q}$. 
Since $\alpha_1, ..., \alpha_d \in L^{1/q}$ is linearly independent over $K^{1/q}$, 
we obtain $c_1^{1/q} = \cdots =c_d^{1/q}=0$, which implies $c_1 = \cdots =c_d =0$. 
Thus (2) holds. 
\end{proof}

\begin{prop}\label{p-QFR-to-Q+R}
In the situation of Setting \ref{setting:most-general-foundations-of-log-quasi-F-splitting}, 
let $\Delta$ be a  $\Q$-divisor on $X$. 
{If $(X, \Delta)$ is globally quasi-$F$-regular, then $(X, \Delta)$ is globally quasi-$+$-regular.}
\end{prop}


\begin{proof}
Assume that $(X, \Delta)$ is 
globally quasi-$F$-regular. 
In particular, $\rdown{\Delta}=0$. 
{Fix $n \in \Z_{>0}$ such that $(X, \Delta)$ is globally $n$-quasi-$F$-regular.} 
It suffices to show the following: 
\begin{enumerate}
\item[$(\star)$]
Let $f: Y \to X$ be a finite surjective morphism from an integral normal excellent scheme $Y$. 
Then 
\[
H^0(X, (\Phi^f_{X, \Delta, n})^*) : 
H^0(X, \mathcal Hom(Q^f_{X, \Delta, n}, W_n\omega_X(-K_X)))
\to H^0(X, \MO_X)
\]
is surjective. 
\end{enumerate}

{In order to} prove {$(\star)$}, 
we may take the normal closure $L$ of $K(Y)/K(X)$ and replace $Y$ by its normalisation in $L$. Then 
there exists an intermediate field $K(X) \subseteq M \subseteq K(Y)$ 
such that $M/K(X)$ is purely inseparable and $K(Y)/M$ is separable {\cite[Proposition 6.11]{Lan02}}.






In particular, there is a factorisation of 
finite surjective morphisms of integral normal excellent schemes  
\[
f: Y \xrightarrow{g} Z \xrightarrow{h} X, 
\]
where $h : Z \to X$ is purely inseparable and $g: Y \to Z$ is separable. 
We have the induced homomorphisms: 
\[
W_n \MO_X(\Delta) \xrightarrow{h^*} h_*W_n\MO_Z(h^*\Delta) \xrightarrow{g^*} h_*g_*W_n\MO_Y(g^*h^*\Delta) = f_*W_n\MO_Y(f^*\Delta). 
\]
By Lemma \ref{l-QFR-to-Q+R}, there exist an effective Weil divisor $D_Z$ on $Z$ 
and a $W_n(K(Z))$-module homomorphism $\varphi : W_n(K(Y)) \to W_n(K(Z))$ 
{which induces, for every $e \in \Z_{>0}$,} 
 a $W_n\MO_Z$-module homomorphism $\varphi : g_*W_n\MO_Y(p^ef^*\Delta) \to W_n\MO_Z(p^eh^*\Delta+D_Z)$ 
such that $\varphi \circ g^*(1) =1$:  
\begin{align*}
W_n\MO_Z(p^eh^*\Delta) &&\xrightarrow{g^*}&& g_*W_n\MO_Y(p^ef^*\Delta) &&\xrightarrow{\varphi}&&  W_n\MO_Z(p^eh^*\Delta+D_Z), \\
1 &&\mapsto&& g^*(1) &&\mapsto&& \varphi \circ g^* (1) = 1. 
\end{align*}
Pick an effective Weil divisor $D_X$ on $X$ such that 
$h^*D_X \geq D_Z$. 
By enlarging $D_Z$, we may assume that $D_Z= h^*D_X$. 
Since $h : Z \to X$ is a finite purely inseparable surjective morphism, 
we can find $e_1 \in \Z_{>0}$ such that the 
$e_1$-th iterated absolute Frobenius morphism $F^{e_1} : X \to X$ 
factors through $h$: 
\[
F^{e_1}: X \xrightarrow{\alpha} Z \xrightarrow{h} X.  
\]
In particular, we get the induced 
$W_n\MO_X$-module homomorphism: 
\[
\alpha^* : h_*{W_n}\MO_Z(h^*(\Delta+D_X)) \to F^{e_1}_*{W_n}\MO_X( (F^{e_1})^*(\Delta+D_X)) = F^{e_1}_*{W_n}\MO_X( p^{e_1}(\Delta+D_X)).
\]
Since $(X, \Delta)$ is globally $n$-quasi-$F$-regular, 
there exists $\epsilon \in \Q_{>0}$ 
such that $(X, \Delta + \epsilon D_X)$ is $n$-quasi-$F^e$-split for every $e \in \Z_{>0}$ 
(Definition \ref{d-QFR}). 
Fix $e_2 \in \Z_{>0}$ such that $p^{e_1+e_2} \epsilon  \geq p^{e_1}$, i.e., 
$p^{e_2} \epsilon \geq 1$. 
By $h_*g_*W_n\MO_Y   = f_*W_n\MO_Y$, 
we can factor the following composition through
$f^*:W_n\MO_X(\Delta)\to f_*W_n\MO_Y(f^*\Delta)$: 
\begin{align*}
W_n\MO_X(\Delta) &\xrightarrow{F^{e_2}} F_*^{e_2}W_n\MO_X(p^{e_2}\Delta)
\xrightarrow{h^*}  F_*^{e_2}h_*W_n\MO_Z(p^{e_2}h^*\Delta)  \xrightarrow{g^*} F_*^{e_2}h_*g_*W_n\MO_Y(p^{e_2}f^*\Delta) \\ 
&\xrightarrow{\hphantom{a}\varphi\hphantom{a}}  F_*^{e_2}h_*W_n\MO_Z(p^{e_2}h^*\Delta + h^*D_X) 
\xrightarrow{\alpha^*} F_*^{e_1+e_2}W_n\MO_X(p^{e_1+e_2}\Delta+ p^{e_1}D_X) \\
&\xhookrightarrow{\hphantom{a\varphi a}} F_*^{e_1+e_2}W_n\MO_X( p^{e_1+e_2}(\Delta+ \epsilon D_X)). 
\end{align*}
By $\varphi \circ g^*(1) = 1$ and $F^{e_1} = h \circ \alpha$,
this composition coincides with the following canonical one: 
\[
W_n\MO_X(\Delta) \xrightarrow{F^{e_1+e_2}} F_*^{e_1+e_2}W_n\MO_X(p^{e_1+e_2}\Delta) \hookrightarrow F_*^{e_1+e_2}W_n\MO_X(p^{e_1+e_2}(\Delta+ \epsilon D_X)). 
\]
Taking the pushouts, we obtain the following factorisation: 
\[
\Phi^{e_1+e_2}_{X, \Delta +\epsilon D_X, n} : \MO_X \xrightarrow{\Phi^f_{X, \Delta, n}} Q^f_{X, \Delta, n} 
\to Q^{e_1+e_2}_{X, \Delta+\epsilon D_X, n}. 
\]
Since $(X, \Delta +\epsilon D_X)$ is $n$-quasi-$F^{e_1+e_2}$-split, 
\[
H^0(X, (\Phi^{e_1+e_2}_{X, \Delta+\epsilon D_X, n})^*): 
\Hom(Q^{e_1+e_2}_{X, \Delta+\epsilon D_X, n}, W_n\omega_X(-K_X)) 
\to 
\Hom(\MO_X, W_n\omega_X(-K_X)) 
\]
is surjective. 
Therefore, also 
\[
H^0(X, (\Phi^f_{X, \Delta, n})^*): 
\Hom(Q^f_{X,\Delta, n}, W_n\omega_X(-K_X)) 
\to 
\Hom(\MO_X(\Delta), W_n\omega_X(-K_X)) 
\]
is surjective. 
\end{proof}

\begin{rem}
In \cite{KTTWYY3}, 
it will be shown that in the $\bQ$-Gorenstein case $X$ is locally quasi-$F$-regular if and only if 
$X$ is locally quasi-+-regular. 
\end{rem}

\subsection{Quasi-$+$-stable sections and stabilisation}

\begin{definition}\label{d-def-+-section}
In the situation of Setting \ref{setting:most-general-foundations-of-log-quasi-F-splitting},  
let $\Delta$ and $L$ be $\Q$-divisors on $X$.
Given $n \in \Z_{>0}$ and a finite cover $f: Y \to X$ (i.e., a finite surjective morphism from an excellent integral normal scheme $Y$), we set 
\[
q^fB^0_n(X, \Delta; L) := \Im 
\big( (\Phi^f_{X, \Delta-L, n})^*: 
H^0(X,(Q^f_{X, \Delta -L,n})^*) \to 
H^0(X, \cO_X(\Delta-L)^*)\big), 
\]
where $(-)^* := \cHom_{W_n\MO_X}(-, W_n\omega_X(-K_X))$ and 
$\MO_X(\Delta-L)^* \simeq \MO_X( \rup{L-\Delta})$. 
We define 
\begin{eqnarray*}
qB^0_n(X, \Delta;L) &:=& \bigcap_{f: Y \to X}q^fB^0_n(X, \Delta;L),\\ 
qB^0(X, \Delta;L) &:=& \bigcup_{n\geq 1} qB^0_n(X, \Delta;L), 
\end{eqnarray*}
\[
Q^+_{X, \Delta, n} := \varinjlim_{f: Y \to X} Q^f_{X, \Delta, n} \qquad {\text{and} \qquad 
B^+_{X, \Delta, n} := \varinjlim_{f: Y \to X} B^f_{X, \Delta, n},}
\]
where $f: Y \to X$ runs over all finite covers. 
Note that if $X$ is projective over a Noetherian ring $R$, 
we have $\bigcup_{n\geq 1} qB^0_n(X, \Delta;L) =qB^0_N(X, \Delta;L)$ for a sufficiently large integer $N$ (note that we have $qB^0_n(-) \subseteq qB^0_{n+1}(-))$). 
\end{definition}

\begin{proposition} \label{prop:qeB0-stabilisation} 
In the situation of Setting \ref{setting:foundations-of-log-quasi-F-splitting}, fix $n\in \bZ_{>0}$, let $A$ be an ample {$\Q$-Cartier} $\bQ$-divisor, and let $N$ be a nef {$\Q$-Cartier} $\bQ$-divisor. 
Then there exists a finite surjective {morphism} 
$f \colon Y \to X$ from a normal integral scheme $Y$ such that for
\begin{itemize}
    \item 
     {every} $\bQ$-divisor $\Delta$, and 
    \item all integers $r,k$
satisfying $r \geq 1$ and $k \geq 0$, 
\end{itemize}
the following  
{equality} 
holds:
\[
qB^0_n(X,\Delta; L) = \Im 
\big( (\Phi^f_{X, \Delta-L, n})^*: 
H^0(X,(Q^f_{X, \Delta -L,n})^*) \to 
H^0(X, \cO_X(\Delta-L)^*)\big),
\]
where $L := K_X + \Delta + rA + kN$.
\end{proposition}

\begin{proof}
In view of Proposition \ref{prop:stabilisation_B^0}, 
the proof is analogous to that of Lemma \ref{lem:Witt-stabilisation} and Proposition \ref{prop:qeS0-stabilisation}. Therefore, we only sketch the key ideas and leave the verification of the details to the reader. Let $\pi \colon X^+ \to X$ denote the natural projection. 
By abuse of notation, we replace {$rA+kN$ by $A$},
but all the 
{finite surjective morphisms} 
we {shall} construct in the proof are clearly independent of $r$ and $k$.

In what follows, we may assume that $R$ is local, since the construction of the desired $f$ relies only on Proposition \ref{prop:stabilisation_B^0}, which is independent of the choice of maximal ideal.

First, take any finite surjective {morphism} $f \colon Y \to X$ from a normal integral scheme $Y$. By the same argument as in the proof of Proposition \ref{prop:qeS0-stabilisation}, looking at the {commutative} diagram 
{in which each horizontal sequence is exact} 
\[
\begin{tikzcd}[column sep=small]
H^d_{W_{n-1}\fram}(F_*W_{n-1}\mcO_X(-pA)) \arrow{r} 
\arrow[d, equal]
&  H^d_{W_n\fram}(W_n\mcO_X(-A)) \arrow[d] \arrow{r} & H^d_{\fram}(\cO_X(-A)) \arrow[d] \arrow{r} & 0\\
H^d_{W_{n-1}\fram}(F_*W_{n-1}\mcO_X(-pA)) \arrow{r} &  H^d_{W_n\fram}(f_*W_n\mcO_Y(-f^*A)) \arrow{r} & H^d_{W_n\fram}(Q^{f}_{X,-A,n}) \arrow{r} & 0,
\end{tikzcd}
\]
we get that
\begin{align*}
\Ker \Big(H^d_{W_n\m}(X,W_n\cO_X(-A)) \to H^d_{W_n\m}(X,f_*W_n\cO_Y(-f^*A))\Big)
\end{align*}
surjects onto 
\[
\Ker\Big(H^d_\m(X,\cO_X(-A)) \to H^d_{W_n\m}(X,Q^f_{X,-A,n})\Big) 
\]
for every finite surjective {morphism} $f \colon Y \to X$ from a normal integral scheme $Y$. In particular, it is enough to show that 
\begin{multline*}
\Ker \Big(H^d_{W_n\m}(X,W_n\cO_X(-A)) \to H^d_{W_n\m}(X^+,W_n\cO_{X^+}(-\pi^*A))\Big)\\
= \Ker \Big(H^d_{W_n\m}(X,W_n\cO_X(-A)) \to H^d_{W_n\m}(Y,W_n\cO_Y(-f^*A))\Big)
\end{multline*}
for some fixed finite surjective {morphism} $f$.

To this end, by Proposition \ref{prop:stabilisation_B^0}, we pick $f \colon Y \to X$ such that
\begin{multline*}
\Ker \Big(H^d_{\m}(X,\cO_X(-A)) \to H^d_{\m}(X^+,\cO_{X^+}(-\pi^*A))\Big)\\
= \Ker \Big(H^d_{\m}(X,\cO_X(-A)) \to H^d_{\m}(Y,\cO_Y(-f^*A))\Big).
\end{multline*}
Now by induction on $n$ applied to $Y$, we can find $f' \colon Z \to Y$ such that
\begin{multline*}
\Ker \Big(H^d_{W_{n-1}\m}(Y,W_{n-1}\cO_Y(-f^*A)) \to H^d_{W_{n-1}\m}(Y^+,W_{n-1}\cO^+_Y(-\pi^*A))\Big)\\
= \Ker \Big(H^d_{W_{n-1}\m}(Y,W_{n-1}\cO_Y(-f^*A)) \to H^d_{W_{n-1}\m}(Z,W_{n-1}\cO_Z(-g^*A))\Big),
\end{multline*}
where $g := f \circ f'$.


Note that $Z^+ = Y^+ = X^+$. We claim that the statement of the proposition holds true for $f$ replaced by $g$. 
To this end, consider the following diagram
{\small\[
\begin{tikzcd}[column sep = small]
H^d_{\fram}(F_*\pi_*W_{n-1}\mcO_{X^+}(-p\pi^*A)) \arrow{r}{(\dagger)} &  H^d_{\fram}(\pi_*W_n\mcO_{X^+}(-\pi^*A)) \arrow{r} & H^d_{\fram}(\pi_*\mcO_{X^+}(-\pi^*A)) \arrow{r} & 0\\
 H^d_{\fram}(F_*g_*W_{n-1}\mcO_Z(-pg^*A)) \arrow{r} \arrow{u} &  H^d_{\fram}(g_*W_n\mcO_Z(-g^*A)) \arrow{r}  \ar{u} & H^d_{\fram}(g_*\mcO_Z(-g^*A)) \ar{u} \arrow{r} & 0\\
 H^d_{\fram}(F_*f_*W_{n-1}\mcO_Y(-pf^*A))  \arrow{u}{} \arrow{r}{} &  H^d_{\fram}(f_*W_n\mcO_Y(-f^*A)) \arrow{u} \arrow{r} & H^d_{\fram}(f_*\mcO_Y(-f^*A)) \arrow{u} \arrow{r} & 0\\
 H^d_{\fram}(F_*W_{n-1}\mcO_X(-pA)) \arrow{r} \arrow{u} &  H^d_{\fram}(W_n\mcO_X(-A)) \arrow{u}{} \arrow{r}{} & H^d_{\fram}(\mcO_X(-A)) \arrow{u}{} \arrow{r} & 0.
\end{tikzcd}
\]}
Note that $H^{d-1}_{\fram}(X,\pi_*\mcO_{X^+}(-\pi^*A))=0$ by \cite[Corollary 3.7]{BMPSTWW20}, and so the map ($\dagger$) is injective. Then
\begin{multline*}
\Ker\Big(H^d_{\fram}(W_n\mcO_X(-A)) \to  H^d_{\fram}(\pi_*W_n\mcO_{X^+}(-\pi^*A)) \Big) \\
= \Ker\Big(H^d_{\fram}(W_n\mcO_X(-A)) \to  H^d_{\fram}(g_*W_n\mcO_{Z}(-g^*A))\Big) 
\end{multline*}
by a diagram chase similar to that in the proof of Proposition \ref{prop:qeS0-stabilisation}.
\end{proof}

\begin{proposition} \label{prop:Serre-type-result-qB0} 
In the situation of Setting \ref{setting:foundations-of-log-quasi-F-splitting}, fix $n\in \bZ_{>0}$, let $A$ be an ample {$\Q$-Cartier} $\bQ$-divisor, and let $N$ be a nef {$\Q$-Cartier} $\bQ$-divisor. 
Assume that $X$ is divisorially Cohen-Macaulay. 
Then there exists $r_0 \in \Z_{>0}$ such that for
\begin{itemize}
    \item all integers $r$, $k$
satisfying $r \geq r_0$ and $k \geq 0$, and
    \item 
    {every} 
    $\bQ$-divisor $\Delta$ such that 
    $(X,\Delta)$ is locally $n$-quasi-$+$-regular and $L := K_X + \Delta + rA + kN$ is 
    {a \textcolor{black}{Weil} divisor},
\end{itemize}
the following 
{equality} 
holds:
\[
qB^0_n(X,\Delta; L) = H^0(X, \cO_X(L)).
\]
\end{proposition}
\begin{proof}
In view of Proposition \ref{prop:qeB0-stabilisation}, the proof is analogous to that of Proposition \ref{prop:Serre-type-result-qS0}. We leave the verification of the details to the reader.
\end{proof}

\begin{proposition} \label{prop:local-cohomology-description-of-B^0-new} 
Fix   {$n \in \Z_{>0}$}. 
In the situation of Setting \ref{setting:foundations-of-log-quasi-F-splitting},  
{assume that $(R, \m)$ is a local ring. 
Take  $\bQ$-divisors $L$ and $\Delta$ on $X$.} 
Assume that $X$ is divisorially Cohen-Macaulay and $L-(K_X+\Delta)$ is ample.
Then
\[
qB^0_n(X,\Delta;L)^{\wedge} \simeq {\rm Im}\Big(H^d_{\fram}({X}, \mcO_X(\lfloor K_X+\Delta-L \rfloor)) \to H^d_{W_n\fram}({X}, Q^{+}_{X,K_X+\Delta-L,n})\Big)^{\vee},
\]
where $(-)^{\wedge}$ denotes $\fram$-completion, and $(-)^{\vee}$ denotes Matlis duality.
\end{proposition}
\begin{proof}
The proof is the same as that of Proposition \ref{prop:local-cohomology-description-of-S^0-perf} by using Proposition \ref{prop:qeB0-stabilisation} 
instead of Corollary \ref{cor-stab-q^e}(1).
\end{proof}


In order to be as explicit as possible in the proof of the proposition below, we introduce two more {non-standard} definitions. Let ${\varphi} \colon \{A_n\}_{n>0} \to \{B_n\}_{n>0}$ be a {homomorphism} 
of 
{projective systems} 
induced by ${\varphi}_n \colon A_n \to B_n$. We denote the 
{projections} 
in the first {projective} system by $\pi_{m,n} \colon A_m \to A_n$ for $m \geq n$. 
For an integer {$n_0 >0$},  
we say that ${\varphi}$ is \emph{$2$-injective at $n_0$} 
if $\pi_{2n_0,n_0}({\rm Ker}({\varphi}_{2n_0} \colon 
{A_{2n_0} \to B_{2n_0}}))=0$.  
We say that $\{A_n\}_{n>0}$ is \emph{$2$-zero at $n_0$} if $\pi_{2n_0,n_0} = 0$.


In what follows, we show that quasi-$+$-stable sections {$qB^0(-)$} agree with {uniformly} quasi-$F^\infty$-stable sections ${\qS^0(-)}$ 
when 
each $(X, \{p^s\Delta\})$ is locally quasi-$+$-regular.

\begin{theorem}\label{t-qB^0-vs-qS^0}
In the situation of Setting \ref{setting:foundations-of-log-quasi-F-splitting}, let $\Delta$ be an effective $\bQ$-divisor on $X$ such that $\rdown{\Delta}=0$ and 
$(X, \{p^s\Delta\})$ is locally quasi-$+$-regular for every integer $s \geq 0$. Let $L$ be a Weil divisor such that $L-(K_X+\Delta)$ is ample. Assume that $X$ is divisorially Cohen-Macaulay. Then
\[
qB^0(X,\Delta; L) = q^\infty_{\rm uni} S^0(X,\Delta; L).
\]
\end{theorem}

\begin{proof}
Set $A := L-(K_X+\Delta)$. 
{Since $X$ is projective over a Noetherian ring $R$}, we can find $n_0>0$ such that
\begin{align*}
qB^0(X,\Delta; L) &=  qB^0_n(X,\Delta; L) \qquad \text{ and }\\
q^\infty_{\rm uni} S^0(X,\Delta; L) &=  q^\infty S^0_n(X,\Delta; L)
\end{align*}
for all $n  \geq n_0$. 
By increasing $n_0$, we may assume that $(X, \{p^s\Delta\})$ is locally $n$-quasi-$+$-regular for all $s \geq 0$ and $n \geq n_0$. 
In view of Proposition \ref{prop:qeS0-stabilisation} and Proposition \ref{prop:qeB0-stabilisation}, 
the question is stable under localisation, and so we may assume that $R$ 
is local. 
{Let $\m$ be the maximal ideal of $R$.}\\



\noindent \textbf{Step 1.} Reduction to the following (\ref{eq:injection-of-pro-systems-Sperf=B}):\\

\begin{equation} \label{eq:injection-of-pro-systems-Sperf=B}
\begin{aligned}
&\text{The homomorphism of projective system}\\
&\{\hspace{-0.05em}H^d_{W_n\fram}({X},\hspace{-0.1em} Q^{e}_{X,{-}A,n}) \}_{n >0} 
\to \{\hspace{-0.05em}H^d_{W_n\fram}({X},\hspace{-0.1em} Q^{+}_{X,{-}A,n})\}_{n > 0}\! \\
&\text{ is $2$-injective at $n_0$ {for $e \gg 0$}.}
\end{aligned}
\end{equation}

\noindent By Proposition \ref{prop:local-cohomology-description-of-S^0-perf} and Proposition \ref{prop:local-cohomology-description-of-B^0-new}, it is enough to prove that the natural map
\[
H^d_{W_n\fram}({X}, Q^{\perf}_{X,-A,n}) \to H^d_{W_n\fram}({X}, Q^{+}_{X,-A,n})
\]
is an injection. Unfortunately, this is too much to hope for, but we will show the following weaker (
{but} sufficient) statement: 
the {homomorphism} 
of 
{projective} systems 
\begin{equation*} 
\{H^d_{W_n\fram}({X}, Q^{\perf}_{X,-A,n})\}_{n>0} \to \{H^d_{W_n\fram}({X}, Q^{+}_{X,-A,n})\}_{n>0} \ \text{ is injective}.
\end{equation*}
Specifically, {we will show} that {this homomorphism} is $2$-injective at $n_0$, which immediately implies the statement of the proposition by 
{chasing} the following diagram:
\[
\begin{tikzcd}
H^d_{\fram}({X}, \cO_X(-A)) \ar{r} \ar{d}{=} &   H^d_{W_{2n}\fram}({X}, Q^{\perf}_{X,-A,2n}) \ar{r} \ar{d} & H^d_{W_{2n}\fram}({X}, Q^{{+}}_{X,-A,2n}) \ar{d}\\
H^d_{\fram}({X}, \cO_X(-A)) \ar{r} &   H^d_{W_n\fram}({X}, Q^{\perf}_{X,-A,n}) \ar{r} & H^d_{W_n\fram}({X}, Q^{{+}}_{X,-A,n})
\end{tikzcd}
\]
for $n = n_0$. 
Since $Q^{\perf}_{X,-A,n} = \varinjlim_e Q^e_{X,-A,n}$, the reduction step follows.


\vspace{1em}

\noindent \textbf{Step 2.} The proof of (\ref{eq:injection-of-pro-systems-Sperf=B}).\\

\noindent To prove (\ref{eq:injection-of-pro-systems-Sperf=B}), consider the natural projection 
$\pi \colon X^+ \to X$ and {the factorisation $\pi \colon X^+ \xrightarrow{\pi'} X \xrightarrow{F^e} X$.
Although $\pi$ and $\pi'$ are isomorphic to each other 
via $F^e : X^+ \xrightarrow{\simeq} X^+$, 
we shall distinguish them in what follows.} 
{We have} the following {commutative} 
diagram with horizontal and vertical sequences exact: 
\[
\begin{tikzcd}
 & F_*W_{n-1}\cO_X(-pA) \ar[hook]{d}{F^e \circ V} \ar{r}{=} 
 & F_*W_{n-1}\cO_X(-pA) \ar[hook]{d}{\pi^* \circ V} &  & \\
0 \ar{r} & F^e_*W_n\cO_X(-p^eA) \ar[two heads]{d} \ar{r}{{\pi'^*}} & \pi_*W_n\cO_{{X^+}}(-\pi^*A) \ar[two heads]{d} \ar{r} & \cF_n \ar{d}{=} \ar{r} & 0. \\
0 \ar{r} & Q^{e}_{X,-A,n} \ar{r} & Q^{+}_{X,-A,n} \ar{r} & \cF_n \ar{r} & 0.
\end{tikzcd}
\] 
It is enough to show that $\{H^{d-1}_{W_n\fram}(X, \cF_n)\}_{n > 0}$ is $2$-zero at $n_0$. We also see that $\cF_n$ features in the following {commutative} diagram with 
{horizontal and vertical sequences} exact:
\[
\begin{tikzcd}
 & F^{e+1}_*W_{n-1}\cO_X(-p^{e+1}A) \ar[hook]{d}{V} \ar{r}{=} & F^{e+1}_*W_{n-1}\cO_X(-p^{e+1}A) \ar[hook]{d}{{\pi'^* \circ V}} &  & \\
0 \ar{r} & F^e_*W_n\cO_X(-p^eA) \arrow[d, two heads, "{R^{n-1}}"] 
\ar{r}{{\pi'^*}} & \pi_*W_n\cO_{{X^+}}(-\pi^*A) \ar[two heads]{d} \ar{r} & \cF_n \ar{d}{=} \ar{r} & 0. \\
0 \ar{r} & F^e_*\cO_X(-p^eA) \ar{r} & (\star_n) \ar{r} & \cF_n \ar{r} & 0.
\end{tikzcd}
\]
Since $
{F^e \circ {\pi'}} = \pi \colon X^+ \to X$, we get that 
\[
\pi_*W_n\cO_{{X^+}}(-\pi^*A) \simeq F^e_*{\pi'}_*W_n\cO_{{X^+}}(-{\pi'^*p^eA}). 
\]
Then $(\star_n)$ {is isomorphic to $F^e_*Q^+_{X,-p^eA,n}$, as $\pi$ and $\pi'$ are isomorphic.} 
Thus to show that $\{H^{d-1}_{W_n\fram}(X, \cF_n)\}_{n>0}$ is $2$-zero at $n_0$, it is sufficient to prove that for $e \gg 0$,
\begin{enumerate} \setlength\itemsep{0.5em}
\item $\{H^{d-1}_{W_n\fram}({X}, Q^+_{X,-p^eA,n})\}_{n>0}$ is $2$-zero at $n_0$, and
\item $H^{d}_{W_n\fram}({X}, \cO_X(-p^eA)) \to H^{d}_{W_n\fram}({X}, Q^+_{X,-p^eA,n})$ is injective for $n=2n_0$.
\end{enumerate}
Now (1) follows from Lemma \ref{lem:towards-CM-of-Q+} {below} in view of Proposition \ref{prop:Serre-type-result-qB0}, whilst (2) {holds by} 
Proposition \ref{prop:Serre-type-result-qB0} thanks to Matlis duality (Proposition \ref{prop:local-cohomology-description-of-B^0-new}). This concludes the proof of the theorem.
\qedhere




\end{proof}

{\begin{lemma} \label{l-BQO+}
In the situation of Setting \ref{setting:most-general-foundations-of-log-quasi-F-splitting}, 
fix $n \in \Z_{>0}$ and let $\Delta$ be a $\bQ$-divisor on $X$. 
\begin{enumerate}
\item 
For a finite surjective morphism $f \colon Y \to X$, there is a canonical exact sequence 
\[
0 \to F_*B^f_{X, p\Delta, n} \to Q^f_{X, \Delta, n+1} \to f_*\MO_Y(f^*\Delta) \to 0.  
\]
\item 
For the induced morphism $\pi \colon X^+ \to X$, there is a canonical exact sequence 
\[
0 \to F_*B^+_{X, p\Delta, n} \to Q^+_{X, \Delta, n+1} \to \pi_*\MO_{X^+}(\pi^*\Delta) \to 0.  
\]
\item There is a canonical exact sequence 
\[
0 \to \MO_X(\Delta) \to Q^+_{X, \Delta, n} \to B^+_{X, \Delta, n} \to 0. 
\]
\end{enumerate}


\end{lemma}

\begin{proof}
The assertion (1) holds by applying the snake lemma to  the following commutative diagram in which the horizontal sequence is exact: 
\vspace{0.3em}
\[
\begin{tikzcd}
& F_*W_{n}\MO_X(p\Delta) \arrow{d}{f^*} \arrow[r,dash,shift left=.1em] \arrow[r,dash,shift right=.1em] & F_*W_{n}\MO_X(p\Delta) \arrow{d}{V \circ f^*} & &\\
0 \arrow{r}&  F_*f_*W_{n}\MO_Y(pf^*\Delta) \arrow{r}{V} & f_*W_{n+1}\MO_Y(f^*\Delta) \arrow{r}{R^n} & f_*\MO_Y(f^*\Delta) \arrow{r} & 0.
\end{tikzcd}
\]
Then (1) implies (2) by applying $\varinjlim_f$, where 
$f$ runs over the finite surjective morphisms $f \colon Y \to X$ from  normal integral schemes $Y$. 
Similarly, (3) is obtained by applying $\varinjlim_f$ to $0 \to \MO_X(\Delta) \to Q^f_{X, \Delta, n} \to B^f_{X, \Delta, n} \to 0$. 
\end{proof}
}

\begin{lemma} \label{lem:towards-CM-of-Q+}
In the situation of Setting \ref{setting:foundations-of-log-quasi-F-splitting}, 
fix $n \in \Z_{>0}$ and suppose that $(R,\mathfrak m)$ is a local ring and $X$ is divisorially Cohen-Macaulay. 
Let $A$ be an ample $\bQ$-Cartier $\bQ$-divisor such that 
\begin{equation} \label{eq:ass-towards-CM-of-Q+}
qB^0_{2n-i}(X,
\{-p^iA\}
, K_X + \lceil p^i A \rceil) = H^0(X, \cO_X(K_X + \lceil p^i A \rceil))
\end{equation}
for all $0 \leq i \leq n$. Then
\[
H^{d-1}_{W_{2n}\fram}(X, Q^+_{X,-A,2n}) \to H^{d-1}_{W_n\fram}(X, Q^+_{X,-A,n})
\]
is a zero map.
\end{lemma}

\begin{proof}
By abuse of notation, we set $H^i_{\m}(-) := H^i_{W_{2n}\m}(-)$. 
More generally, we will show by descending induction on $0 \leq i \leq {n-1}$ that
\begin{equation} \label{eq:induction_Q_2n-A-bigg}
H^{d-1}_{\fram}({X}, Q^+_{X,-p^iA,2n-i}) \to H^{d-1}_{\fram}({X}, Q^+_{X,-p^iA,n-i})
\text{\quad is\quad zero.}
\end{equation}
{The base case $i=n-1$ of this induction follows from 
\[
H^{d-1}_{\m}(X, Q^+_{X, -p^{n-1}A, 1}) = H^{d-1}_{\m}(X^+, \MO_{X^+}(-p^{n-1}\pi^*A)) =0, 
\]
where the latter equality is a consequence of \cite[Corollary 3.7]{BMPSTWW20}.} 
Thus we may assume that (\ref{eq:induction_Q_2n-A-bigg}) holds for some $0 < i \leq {n-1}$ and aim for showing that:
\begin{equation} \label{eq:induction_-assumptionQ_2n-A-bigg}
H^{d-1}_{\fram}({X}, Q^+_{X,-p^{i-1}A,2n-i+1}) \to H^{d-1}_{\fram}({X}, Q^+_{X,-p^{i-1}A,n-i+1}) {\text{\quad is zero}}.
\end{equation}

First, consider the following map of short exact sequences ({Lemma \ref{l-BQO+}(2)}):
\[
\begin{tikzcd}
0 \arrow{r} & F_*({B^+_{X,-p^iA,2n-i}}) \arrow{d} \arrow{r} & Q^+_{X,-p^{i-1}A,2n-i + 1} \arrow{r} \arrow{d} & \pi_*\mcO_{X^+}(-p^{i-1}\pi^*A) \arrow{d}{=} \arrow{r} & 0\\
0 \arrow{r} & F_*({B^+_{X,-p^iA,n-i}}
) \arrow{r} & Q^+_{X,-p^{i-1}A,n-i+1} \arrow{r} & \pi_*\mcO_{X^+}(-p^{i-1}\pi^*A) \arrow{r} & 0.
\end{tikzcd}
\]
Note that $H^j_{\fram}(X, \pi_*\mcO_{X^+}(-p^{i-1}\pi^*A)) = 0$ for $j<d$ by \cite[Corollary 3.7]{BMPSTWW20}.
Thus, it is enough to show that
{\small \begin{equation} \label{eq:basic_lifting_auxiliary}
H^{d-1}_{\fram}(X, {B^+_{X,-p^iA,2n-i}}) \to 
H^{d-1}_{\fram}(X,  {B^+_{X,-p^iA,n-i}})
\end{equation}}
is zero.

Consider the following diagram ({Lemma \ref{l-BQO+}(3)})  
\begin{center}
\begin{tikzcd}
0 \arrow{r} & \mcO_X(-p^iA) \arrow{d}{=} \arrow{r} & 
Q^+_{X,-p^iA,2n-i} \arrow{r} \arrow{d} & 
{B^+_{X,-p^iA,2n-i}}
\arrow{d} \arrow{r} & 0\\
0 \arrow{r} & \mcO_X(-p^iA) \arrow{r} & Q^+_{X,-p^iA,n-i} \arrow{r} & 
{B^+_{X,-p^iA,n-i}}
\arrow{r} & 0.
\end{tikzcd}
\end{center}
By our assumption (\ref{eq:ass-towards-CM-of-Q+}), Proposition \ref{prop:local-cohomology-description-of-B^0-new} gives us that
\[
H^d_{\fram}(X, \mcO_X(-p^iA)) \to H^d_{\fram}(X, Q^+_{X,-p^iA,2n-i})
\]
is injective. Moreover, \[
H^{d-1}_{\fram}(X,Q^+_{X,-p^iA,2n-i}) \to H^{d-1}_{\fram}(X,Q^+_{X,-p^iA,n-i})
\]
is a zero map by {the} induction {hypothesis}. Therefore, a diagram chase shows that (\ref{eq:basic_lifting_auxiliary}) is a zero map as well.
\end{proof}

\begin{remark}
Note that the key assumption (\ref{eq:ass-towards-CM-of-Q+}) of Lemma \ref{lem:towards-CM-of-Q+} is satisfied when $n \gg 0$ if 
\begin{itemize}
    \item $(X,\{-p^iA\})$ is locally quasi-$+$-regular and $A$ is sufficiently ample, or
    \item $(X,\{-p^iA\})$ is globally quasi-$+$-regular,
\end{itemize} 
for $0 \leq i \leq n$.
    \end{remark}

{\begin{cor}\label{c-Fano-QFS-vs-+}
In the situation of Setting \ref{setting:foundations-of-log-quasi-F-splitting}, let $\Delta$ be an effective $\bQ$-divisor on $X$ with $\rdown{\Delta}=0$. 
Assume that $X$ is divisorially Cohen-Macaulay, 
$-(K_X+\Delta)$ is ample, and 
$(X, \{p^s\Delta\})$ is locally quasi-$+$-regular for every integer $s \geq 0$. 
Then $(X, \Delta)$ is globally quasi-$+$-regular if and only if 
$(X, \Delta)$ is uniformly quasi-$F^{\infty}$-split. 
\end{cor}

\begin{proof}
The assertion holds by applying Theorem \ref{t-qB^0-vs-qS^0} for $L:=0$. 
\end{proof}
}

Finally, we give definitions in the adjoint setting.
\begin{definition}\label{d purely globally Q+R}
We work in the situation of  Setting \ref{setting:most-general-foundations-of-log-pure-quasi-F-splitting}. 
We say that $f: (Y, S_Y) \to (X, S)$ is a {\em finite cover} (of $(X, S)$) if 
$f : Y \to X$ is a finite surjective morphism from an integral normal scheme $Y$ 
and $S_Y$ is a prime divisor on $Y$ satisfying $f(S_Y) = S$. 
By abuse of notation, the induced morphism $S_Y \to S$ is also denoted by $f$.  

For a finite cover $f:(Y,S_Y)\to (X,S)$ and a $\bQ$-divisor $\Delta$ on $X$, 
we define a $W_n\MO_X$-module homomorphism
$\Phi^{S,f}_{X,\Delta,n}$ and a coherent $W_n\MO_X$-module
$Q^{S,f}_{X,\Delta,n}$ on $X$ by the following pushout diagram:
\[
\begin{tikzcd}
W_n\cI_S(\Delta)
  \arrow[r, "f^*"]
  \arrow[d, "R^{n-1}"']
&
f_*W_n\cI_{S_Y}(f^*\Delta)
  \arrow[d]
\\
\cI_S(\Delta)
  \arrow[r, "\Phi^{S,f}_{X,\Delta,n}"]
&
Q^{S,f}_{X,\Delta,n}.
\end{tikzcd}
\]
We define $(\Phi^{S,f}_{X,\Delta,n})^*$ by applying
$(-)^*:=\cHom_{W_n\MO_X}(-,W_n\omega_X(-K_X))$ to $\Phi^{S,f}_{X,\Delta,n}$:
\[
(\Phi^{S,f}_{X,\Delta,n})^*:
(Q^{S,f}_{X,\Delta,n})^*
\to
\cI_S(\Delta)^*.
\]
Note that
{\small 
\[
\cI_S(\Delta)^*
=
\cHom_{W_n\MO_X}
\bigl(\cI_S(\Delta),W_n\omega_X(-K_X)\bigr)
\simeq
\cHom_{\MO_X}
\bigl(\cI_S(\Delta),\omega_X(-K_X)\bigr)
\simeq
\MO_X(\rup{S-\Delta}).
\]}
Let $B$ be a $\bQ$-divisor on $X$ such that
$S\not\subseteq\Supp B$ and $\rdown{B}=0$, and let
$L$ be a Weil divisor on $X$ such that
$S\not\subseteq\Supp L$.
For a finite cover $f:(Y, S_Y) \to (X, S)$, 
we set 
\begin{align*}
q^fB^0_{n,\adj}(X,S+B;L) &:= \image\big(
H^0(X,(Q^{S, f}_{X,S+B-L,n})^*) \to  
H^0(X, \cO_X(L))\big),\ 
\end{align*}

We define 
\begin{align*}
qB^0_{n, \adj}(X,S+B;L) \ &:= \quad\  \bigcap_{\mathclap{f: (Y, T) \to (X, S)}} q^fB^0_{n,\adj}(X,S+B;L),\\
qB^0_{\adj}(X,S+B;L) \ &:= \quad\ \bigcup_{n\geq 1} qB^0_{n, \adj}(X,S+B;L) =   qB^0_{N, \adj}(X,S+B;L),
\end{align*}
where $f: (Y, T) \to (X, S)$ 
runs over all finite covers. 
Note that if $X$ is projective over a Noetherian ring $R$, 
we have $\bigcup_{n\geq 1} qB^0_{n, \adj}(X,S+B;L) =   qB^0_{N, \adj}(X,S+B;L)$ for a sufficiently large integer $N$ (note that we have $qB^0_{n, \adj}(-) \subseteq qB^0_{n+1, \adj}(-))$).

We say that $(X,S+B)$ is \emph{purely globally
$n$-quasi-$+$-regular} if
\[
qB^0_{n,\adj}(X,S+B;0)=H^0(X,\cO_X).
\]
We say that $(X,S+B)$ is \emph{purely globally
quasi-$+$-regular} if it is purely globally  $n$-quasi-$+$-regular for some $n\in\Z_{>0}$.
\end{definition}

\section{Inversion of adjunction}

\subsection{The restriction map for $q^eS^0$}

We start by showing that the restriction of $Q^e_{X,\Delta,n}$ to $S$ is equal to $Q^e_{S,\Delta|_S,n}$.

\begin{proposition} \label{prop:restrictiong-ses-for-C}
In the situation of Setting \ref{setting:most-general-foundations-of-log-quasi-F-splitting}, suppose that $X$ is divisorially Cohen-Macaulay. Let $S$ be a normal prime divisor on $X$ and 
let $\Delta$ be a {$\bQ$-Cartier} $\bQ$-divisor on $X$ such that 
{$S \not\subseteq \Supp\, \Delta$} and  
$(X,S+\{p^i\Delta\})$ is plt for every $i \geq 0$.
Then restricting $Q^e_{X,\Delta, n}$ to $S$ yields the following short exact sequence:
\[
0 \to Q^{S, e}_{X,\Delta,n} \to Q^e_{X,\Delta,n} \to Q^e_{S,\Delta_S,n} \to 0,
\]
where $\Delta_S := \Delta|_S$ (cf. Definition \ref{definition:restricting-divisor}).
\end{proposition}


\begin{proof}
By \cite[the first paragraph of the proof of Proposition 4.1]{KTTWYY1}
we have the following exact sequence
\begin{equation} \label{eq:witt-restriction}
0 \to W_n \cI_S(\Delta) \to W_n\cO_X(\Delta) \xrightarrow{\res} W_n\cO_S(\Delta_S) \to 0. 
\end{equation}
Note that the $W_n\cO_X$-module homomorphism $\res \colon W_n\cO_X(\Delta) \to W_n\cO_S(\Delta_S)$ is defined by the formula 
\[
(a_0, \ldots, a_{n-1}) \mapsto (a_0|_S, \ldots, a_{n-1}|_S),
\]
where $a_i \in F^i_*\cO_X(p^i\Delta)(U)$ for some open subset $U \subseteq X$, and $a_i|_S := \res(a_i)$ for {the} surjective homomorphism $\res \colon F^i_*\cO_X(p^i\Delta) \to F^i_*\cO_S(p^i\Delta_S)$ induced from Proposition \ref{prop:basic_restriction}. 
By (\ref{eq:witt-restriction}), 
we get the following commutative diagram in which each horizontal sequence is exact 
\begin{center}
\begin{tikzcd}
0 \arrow{r} & F_*W_{n-1} \cI_S({p}\Delta) \arrow{d}{F^eV} \arrow{r} & 
F_*W_{n-1}\cO_X({ p}\Delta) \arrow{d}{F^eV} \arrow{r}{\res} & 
F_*W_{n-1}\cO_S({ p}\Delta_S) \arrow{d}{F^eV} \arrow{r} & 0 \\
0 \arrow{r} & F_*^eW_n\cI_S({ p^e}\Delta) \arrow{r} & 
F_*^eW_n\cO_X({ p^e}\Delta) \arrow{r}{\res} & 
F_*^eW_n\MO_S({p^e}\Delta_S) \arrow{r} & 0.
\end{tikzcd}
\end{center} 
Since all the vertical arrows are injective, the snake lemma yields the required exact sequence (cf.\ (\ref{e-big-BC-diagram2})).
\end{proof}

We now construct a restriction map for $q^eS^0$, 
which will be a key tool for inversion of adjunction. 
With exactly the same assumptions as in Proposition \ref{prop:restrictiong-ses-for-C},  
consider the following commutative diagram 
in which each horizontal sequence is exact:
\begin{equation}\label{e-adj-diag}
\begin{tikzcd}
0 \arrow{r} & Q^{S, e}_{X,\Delta,n}  \arrow{r} & Q^e_{X,\Delta,n} \arrow{r} & Q^e_{S,\Delta_S,n} \arrow{r} & 0 \\
0 \arrow{r} & \cO_X(\Delta-S) \arrow{r} 
\arrow{u}{{\Phi^{S, e}_{X, \Delta, n}}} & \cO_X(\Delta) \arrow{r} \arrow{u}{{\Phi^e_{X, \Delta, n}}} & \cO_S(\Delta_S) \arrow{r} \arrow{u}{{\Phi^e_{S, \Delta_S, n}}} & 0.
\end{tikzcd}
\end{equation}
Applying 
$\cHom_{{W_n\MO_X}}(- , W_n\omega_X) = 
(-\otimes_{W_n\MO_X} W_n\MO_X(-K_X))^*$ (recall that $(-)^* := \cHom_{W_n\MO_X}(-, W_n\omega_X(-K_X))$), 
we get the following commutative diagram 
in which each horizontal sequence is exact:
\begin{equation}\label{e-C-adj-diag}
\begin{tikzcd}[column sep = small]
0 \arrow{r} & (Q^e_{X,\Delta-K_X,n})^* \arrow{d}{} \arrow{r} & (Q^{S, e}_{X, \Delta-K_X, n})^* \arrow{d}{}  \arrow{r} & (Q^e_{S,\Delta_S-K_S, n})^* \arrow{d}{}  &   \\
0 \arrow{r} & \cO_X(K_X\!-\!\rdown{\Delta})  \arrow{r} & \cO_X(K_X\!+\!S\!-\!\rdown{\Delta})  \arrow{r} & \cO_S(K_S\!-\!\rdown{\Delta_S})  \arrow{r} & 0.
\end{tikzcd}
\end{equation}
{The existence of this diagram follows by the same argument as in \cite[p.47]{KTTWYY1}. Note that the lower sequence is exact, because the cokernel of   $\cO_X(K_X\!-\!\rdown{\Delta}) \to \cO_X(K_X\!+\!S\!-\!\rdown{\Delta})$ is reflexive, due to $X$ being divisorially Cohen-Macaulay.} 
We are ready to define the restriction map for $q^eS^0$.

\begin{definition} \label{definitin:restriction-for-qS^0} 
In the situation of {Setting \ref{setting:most-general-foundations-of-log-quasi-F-splitting}},
let $(X,S+B)$ be a  divisorially Cohen-Macaulay plt pair, 
where
$S$ is a normal prime divisor and 
$B$ is a $\Q$-divisor such that $\rdown{B}=0$. 
Suppose that $(X,S+\{p^iB\})$ is plt for every $i\geq 0$ and write $K_S + B_S = (K_X+S+B)|_S$. 

Let $L$ be a $\bQ$-Cartier  Weil divisor {such that $S \not\subseteq \Supp L$}. Recall that 
\begin{align*}
q^eS^0_{n,\adj}(X,S+B;L) &= \image\big(
H^0(X,(Q^{S, e}_{X,S+B-L,n})^*) \to  
H^0(X, \cO_X(L))\big),\ {\rm and } \\
q^eS^0_n(S, B_S;L|_S) &= \image\big(
H^0(S,(Q^e_{S,B_S -L|_S,n})^*) \to H^0(S, \cO_{S}(\lceil L|_S-B_S\rceil))\big).
\end{align*}
Then we {obtain} a map
\[
q^eS^0_{n, \adj}(X, S+B; L) \to q^eS^0_n(S, B_S; L|_S)
\]
induced by plugging $\Delta=K_X+S+B-L$ into the above diagram (\ref{e-C-adj-diag}) to get
\begin{equation} \label{eq:diagram-res-C}
\begin{tikzcd}[column sep = small]
0 \arrow{r} & (Q^e_{X,S+B-L,n})^* \arrow{d}{} \arrow{r} & (Q^{S, e}_{X,S+B-L, n})^* \arrow{d}{}  \arrow{r} & (Q^e_{S,B_S-L|_S, n})^* \arrow{d}{}\\
0 \arrow{r} & \cO_X(L-S)  \arrow{r} & \cO_X(L)  \arrow{r} & \cO_S(\rup{L|_S - B_S})  \arrow{r} & 0.
\end{tikzcd}
\end{equation}
By taking the union $\bigcup_{n \geq 1}$, we obtain a map 
\[
q^eS^0_{\adj}(X, S+B; L) \to q^eS^0(S, B_S; L|_S)
\]
\end{definition}


Under our assumptions, $\rup{L|_S- B_S} = \rdown{L|_S}$ (see Lemma  \ref{l:special-cases-of-L}), and so the lower horizontal row in Diagram (\ref{eq:diagram-res-C}) does not contradict Proposition \ref{prop:basic_restriction}.


{
\begin{remark} \label{remark:implicit-replacement-of-L}
In the above definition, in order to make $(K_X+S+B)|_S$ well-defined, we always, implicitly, pick a canonical divisor $K_X$ so that $S \not \subseteq \Supp\,(K_X+S+B)$. Specifically, this is achieved by using prime avoidance to pick $K_X=-S+D$ for a Weil divisor $D$ such that $S \not \subseteq \Supp D$.

More importantly, we will apply Definition \ref{definitin:restriction-for-qS^0} even when $S \subseteq \Supp L$. In that case, we implicitly find $L' \sim L$ such that $S \not \subseteq \Supp L'$ and then replace $L$ by $L'$. 
As making such replacements in the body of the text would make the proofs and  notation impenetrable, we elected \emph{not} to do it. 
\end{remark}}

\begin{lem}\label{l-why-Matlis}
In the situation of {Setting \ref{setting:foundations-of-log-quasi-F-splitting}, 
assume that  $(R, \m)$ is a local ring.}
Let $(X,S+B)$ be a  divisorially Cohen-Macaulay plt pair, 
where
$S$ is a normal prime divisor and 
$B$ is a $\Q$-divisor such that $\rdown{B}=0$. 
Suppose that $(X,S+\{p^iB\})$ is plt for every $i\geq 0$ and write $K_S + B_S = (K_X+S+B)|_S$. 
Let $L$ be a $\bQ$-Cartier  Weil divisor on $X$. {Fix an integer $e>0$}. 
Then the following hold. 
\begin{enumerate}
\item If $H^{d-1}_{\m}(X, Q^e_{X, K_X+S+B-L, n})=0$, 
then 
\[
q^eS^0_{n, \adj}(X, S+B; L) \to q^eS^0_n(S, B_S; L|_S)
\]
is surjective. 
\item Assume that for every $n_0$, there exist integers 
$m$ and $n$ such that $m > n \geq n_0$  and  
\[
H^{d-1}_{\m}(X, Q^e_{X, K_X+S+B-L, m }) \to H^{d-1}_{\m}(X, Q^e_{X, K_X+S+B-L, n})
\]
is zero. Then 
\[
q^eS^0_{\adj}(X, S+B; L) \to q^eS^0(S, B_S; L|_S)
\]
is surjective. 
\end{enumerate}
\end{lem}

\begin{proof}
By definition, we obtain the following commutative diagram: 
\begin{equation} \label{e1-why-Matlis}
\begin{tikzcd}
H^0(X, (Q^{S, e}_{X,S+B-L, n})^*) \arrow[d, twoheadrightarrow]  
\arrow[r, "\alpha_n"] & 
H^0(S, (Q^e_{S,B_S-L|_S, n})^*)  \arrow[d, twoheadrightarrow]\\
q^eS^0_{n,\adj}(X,S+B;L)  \arrow[d, hook]  \arrow[r, "\beta_n"] & 
q^eS^0_n(S, B_S;L|_S) \arrow[d, hook] \\
H^0(X, \cO_X(L))  \arrow[r, "\gamma_n"] & H^0(S, \cO_S(\rup{L|_S - B_S})).
\end{tikzcd}
\end{equation}
By applying Matlis duality $(-)^{\vee} := \Hom_{W_nR}(-, E)$, we obtain 
\begin{equation} \label{e2-why-Matlis}
\begin{tikzcd}
H^d_{\m}(X, Q^{S, e}_{X, K_X+S+B-L, n}) \arrow[d, hookleftarrow]  \arrow[r, leftarrow, "{\alpha_n^{\vee}}"] & 
H^{d-1}_{\m}(S, Q^e_{S, K_S+B_S-L|_S, n})  \arrow[d, hookleftarrow, "j_n"]\\
q^eS^0_{n,\adj}(X,S+B;L)^{\vee}  \arrow[d,  twoheadleftarrow]  \arrow[r, leftarrow, "{\beta_n^{\vee}}"] & 
q^eS^0_n(S, B_S;L|_S)^{\vee} \arrow[d, twoheadleftarrow] \\
H^d_{\m}(X, \cO_X(K_X-L)) \arrow[r, leftarrow, "{\gamma_n^{\vee}}"] & 
H^{d-1}_{\m}(S, \cO_S(K_S+B_S-L|_S)), 
\end{tikzcd}
\end{equation}
where we used $H^d_{\m}(X, \mathcal F) \simeq H^0(X, \cHom_{W_n\MO_X}(\mathcal F, W_n\omega_X))^{\vee}$ 
{for a coherent sheaf $\cF$ on $W_n X$ (see Lemma \ref{lem:Matlis-duality-for-highest-cohomology} applied to $W_n X$)}. 
Note that the diagram (\ref{e1-why-Matlis}) coincides with 
the one obtained from (\ref{e-adj-diag}).


Let us show (1). 
By $H^{d-1}_{\m}(X, Q^e_{X, K_X+S+B-L, n})=0$, 
 $\alpha_n^{\vee}$ is injective. 
By chasing the diagram (\ref{e2-why-Matlis}),  $\beta_n^{\vee}$ is injective, and hence $\beta_n$ is surjective. 
 Thus (1) holds. 

Let us show (2). 
Pick $n_0 \gg 0$ and take integers $m$ and $n$ as in the statement of (2). 
In particular, $m > n \geq n_0 \gg 0$. 
Recall that we have 
\begin{align*}
q^eS^0_1(S, B_S; L|_S) \subseteq 
\, \cdots\,  &\subseteq 
q^eS^0_n(S, B_S; L|_S) \\
&= q^eS^0_{n+1}(S, B_S; L|_S) = \, \cdots\, =
q^eS^0(S, B_S; L|_S). 
\end{align*}
Applying $(-)^{\vee} \,(=\Hom_{W_eR}(-, E))$, we obtain 
\begin{align*}
q^eS^0_1(S, B_S; L|_S)^{\vee} 
\twoheadleftarrow \,\cdots\,  &\twoheadleftarrow 
q^eS^0_n(S, B_S; L|_S)^{\vee} \\
&\xleftarrow{\simeq} q^eS^0_{n+1}(S, B_S; L|_S)^{\vee} \xleftarrow{\simeq} \, \cdots \, \xleftarrow{\simeq}
q^eS^0(S, B_S; L|_S)^{\vee}. 
\end{align*}
Similarly, we may assume that 
$q^eS^0_{n,\adj}(X,S+B;L)^{\vee}
 \xleftarrow{\simeq}
q^eS^0_{m,\adj}(X,S+B;L)^{\vee}$. 
It is enough to show that $\beta^{\vee}_n$ is injective.\\

To this end, pick $\zeta_n \in q^eS^0_n(S, B_S; L|_S)^{\vee}$ 
such that $\beta_n^{\vee}(\zeta_n)=0$. 
Via the induced isomorphism \[
q^eS_n^0(S, B_S; L|_S) \xleftarrow{\simeq}
q^eS_m^0(S, B_S; L|_S),
\]
there is an element $\zeta_m \in q^eS_m^0(S, B_S; L|_S)^{\vee}$ corresponding to $\zeta_n \in q^eS_n^0(S, B_S; L|_S)^{\vee}$. 
Therefore, $\beta^{\vee}_m(\zeta_m)=0$. 
Then $j_m(\zeta_m)$ is the image of some element 
\[
\xi_m \in H^{d-1}_{\m}(X, Q^e_{X, K_X+S+B-L, m }).
\]
By our assumption, the image of $\xi_m$ in $H^{d-1}_{\m}(X, Q^e_{X, K_X+S+B-L, n })$ is zero. 
Then $j_n(\zeta_n)=0$, and hence $\zeta_n=0$. 
Thus (2) holds. 
\end{proof}

\begin{remark}
Assertion (2) can be replaced by a more precise statement. 
Specifically, pick an integer $n>0$ such that
\[
q^eS^0_{n}(S, B_S; L|_S) =  q^eS^0_{n+1}(S, B_S; L|_S) = \ldots = q^eS^0(S, B_S; L|_S),
\]
and assume that $H^{d-1}_{\m}(X, Q^e_{X, K_X+S+B-L, m }) \to H^{d-1}_{\m}(X, Q^e_{X, K_X+S+B-L, n})$ is zero for some $m>n$. Then
\[
q^eS^0_{\adj}(X, S+B; L) \to q^eS^0(S, B_S; L|_S)
\]
is surjective. We emphasise that the choice of $n$ depends on $e$!
\end{remark}

\begin{remark}
With some work involving the Mittag-Leffler condition, one can check that the assumption in (2) is equivalent to the vanishing  $H^{d-1}_{\m}(Q^{e}_{X,K_X+S+B-L})=0$, where $Q^{e}_{X,K_X+S+B-L} := \varprojlim_n Q^{e}_{X,K_X+S+B-L,n}$, which one can use to give a unified proof of both (1) and (2).
\end{remark}

\subsection{Inversion of adjunction in the log Calabi-Yau case}

  
\begin{theorem} \label{thm:log-trivial-inversion-of-adjunction}
Fix $n \in \Z_{>0}$ and $e \in \Z_{>0}$. 
Let $k$ be an $F$-finite field of characteristic $p>0$ and 
let $(X,S+B)$ be a  divisorially Cohen-Macaulay plt pair, 
where $X$ is projective over $k$, $S$ is a normal prime divisor, and 
$B$ is a $\Q$-divisor such that $\rdown{B}=0$. 
Set $d := \dim X$ and define the effective $\Q$-divisor $B_S$ on $S$ by adjunction $K_S + B_S = (K_X+S+B)|_S$. 
{Suppose that the following hold.}
\begin{enumerate}
\item $(X,S+\{p^iB\})$ is plt for every $i\geq 0$. 
    \item $(S,B_S)$ is $n$-quasi-$F^e$-split. 
    \item $H^{d-1}(X, \MO_X(p^m(K_X+S+B)))=H^{d}(X, \MO_X(p^m(K_X+S+B)))=0$ for every integer $m\geq 1$. 
\end{enumerate}
Then $(X,S+B)$ is   {purely} $n$-quasi-$F^e$-split.
\end{theorem}

\begin{proof}
After replacing $k$ by $H^0(X, \MO_X)$, we may assume that $H^0(X, \MO_X)=k$. 
By (2), we have $q^eS_{{n}}^0(S, B_S; 0) =H^0(S, \MO_S)$, and hence $q^eS^0(S, B_S; 0) \neq 0$. 
It suffices to show 
\begin{equation}\label{e1-log-CY}
    H^{d-1}(X, Q^e_{X, K_X+S+B, m})=0\qquad 
\text{for \quad every }\quad m \in \{1, 2,..., n\}, 
\end{equation}
because 
the equality $H^{d-1}(X, Q^e_{X, K_X+S+B, n})=0$, together with Lemma \ref{l-why-Matlis}, implies $q^eS^0_{{n}, \adj}(X, S+B; 0) \neq 0$, 
which in turn yields $q^eS^0_{{n}, \adj}(X, S+B; 0) = H^0(X, \MO_X)$, i.e., $(X,S+B)$ is   {purely} $n$-quasi-$F^e$-split.  

Let us prove (\ref{e1-log-CY}) by induction on $m$. 
The base case $m=1$ is settled by 
\[
H^{d-1}(X, Q^e_{X, K_X+S+B, 1}) 
=
H^{d-1}(X, F_*^e\MO_X(p^e(K_X+S+B))) =0, 
\]
where the latter equality follows from Assumption (3). 
We have the following exact sequence (\ref{e-Bn-induction}): 
\[
0 \to F_*^{m}B^e_{X, p^{m}(K_X+S+B), 1}
\to 
Q^e_{X, K_X+S+B, m+1} \to 
Q^e_{X, K_X+S+B, m} \to 0, 
\]
and hence the problem is reduced to 
$H^{d-1}(X, B^e_{X, p^m(K_X+S+B), 1})=0$ for every $m \geq 1$. 
 
As a more generalised statement, it suffices to show 
\begin{equation}\label{e2-log-CY}
H^{d-1}(X, B^f_{X, p^m(K_X+S+B), 1}) =0\qquad 
\text{if \quad }\quad f \geq 0 \quad \text{and}\quad m \geq 1. 
\end{equation}
We prove (\ref{e2-log-CY}) by induction on $f$. 
The base case  $f=0$ follows from 
$B^0_{X, p^m(K_X+S+B), 1} = \Coker( F^0: \MO_X \to F_*^0\MO_X) = \Coker( {\rm id} :  \MO_X \to \MO_X) =0$. 
By the following exact sequence  (\ref{e-Be-induction}): 
\[
0 \to B^{f-1}_{X, p^m(K_X+S+B), 1} \to B^f_{X, p^m(K_X+S+B), 1} \to F_*^{f-1}B^1_{X, p^{m+f-1}(K_X+S+B), 1} \to 0, 
\]
it is enough to prove 
\begin{equation}\label{e3-log-CY}
H^{d-1}(X, B^1_{X, p^m(K_X+S+B), 1})=0 \qquad \text{ for \quad every \quad $m  \geq 1$}. 
\end{equation}
By the definition of 
$B^1_{X, p^m(K_X+S+B), 1}$ (\ref{eq:definition-of-log-B}), we have an exact sequence 
\[
0 \to \MO_X(p^m(K_X+S+B)) \to F_*\MO_X(p^{m+1}(K_X+S+B)) \to 
B^1_{X, p^m(K_X+S+B), 1} \to 0. 
\]
Then (\ref{e3-log-CY}) follows from Assumption (3). 
\end{proof}

\begin{rem}\label{r-log-CY}
We use the same notation as in the statement of Theorem \ref{thm:log-trivial-inversion-of-adjunction}. 
\begin{itemize}
\item Assumption (1) holds 
if $X$ is $\Q$-factorial and $B$ has standard coefficients. 
Indeed, as we are assuming that $(X, S+B)$ is plt, 
also $(X, S+\{p^i B\})$ is plt for every $i \geq 0$, because 
$K_X+S+\{p^iB\}$ is $\Q$-Cartier and 
$\{ p^iB\} \leq B$.
\item 
The equality $H^{d}(X, \MO_X(p^m(K_X+S+B)))=0$
in Assumption (3) holds when $K_X+S+B \equiv 0$ and $B$ has standard coefficients, because we have 
\[
h^{d}(X, \MO_X(p^m(K_X+S+B))) = h^0(X,\MO_X( K_X - \rdown{p^m(K_X+S+B)})) 
\]
and 
\[
K_X - \rdown{p^m(K_X+S+B)}
\equiv K_X + \{p^{m}B\} \leq K_X+B \equiv -S.
\]
\end{itemize}  
\end{rem}


\begin{cor}\label{c-log-CY}
Fix $n \in \Z_{>0}$ and $e \in \Z_{>0}$. 
Let $k$ be an $F$-finite field of characteristic $p>0$ and 
let $X$ be a $d$-dimensional regular projective variety over $k$. 
Take a regular prime divisor $S$ on $X$. 
Assume that 
\begin{enumerate}
\renewcommand{\labelenumi}{(\roman{enumi})}
\item $H^{d-1}(X, \MO_X)=0$, 
\item $K_X+S \sim 0$, and 
\item $S$ is $n$-quasi-$F^e$-split. 
\end{enumerate}
Then  $(X, S)$ is  {purely} $n$-quasi-$F^e$-split. 
In particular, $X$ is $n$-quasi-$F^e$-split. 
\end{cor}

\begin{proof}
We apply Theorem \ref{thm:log-trivial-inversion-of-adjunction} by setting $B := 0$. 
It is enough to verify its assumptions. 
By Remark \ref{r-log-CY}, 
Assumption (1) of  Theorem \ref{thm:log-trivial-inversion-of-adjunction} holds. 
Assumption (2) of  Theorem \ref{thm:log-trivial-inversion-of-adjunction} follows from (iii). 
Finally, Assumption (3) of  Theorem \ref{thm:log-trivial-inversion-of-adjunction} holds by (i) and (ii). 
\end{proof}



\subsection{Inversion of adjunction for anti-semi-ample divisors}

\begin{setting}\label{s-IOA-anti-ample}
Let $R$ be {an $F$-finite Noetherian local domain of characteristic $p>0$} such that 
$W_nR$ is an excellent ring admitting a dualising complex for every $n \in \Z_{>0}$. 
{Set $\m$ to be the maximal ideal of $R$.} 
Let $X$ be a $d$-dimensional integral normal scheme 
such that {$X$ admits a projective morphism $\pi \colon X \to \Spec R$}  and $H^0(X, \MO_X)=R$. 
Let $(X,S+B)$ be a divisorially Cohen-Macaulay 
$\Q$-factorial plt pair, where 
$S$ is a normal prime divisor and  
$B$ is a $\Q$-divisor such that $\rdown{B}=0$ {and $(X, S + \{ p^iB\})$ is plt for every integer $i\geq 0$}. 
We define an effective $\Q$-divisor $B_S$ on $S$ by $K_S + B_S = (K_X+S+B)|_S$. 
Let $L$ be a Weil divisor on $X$ such that  $A := L-(K_X+S+B)$ is ample. 
For $n \in \Z_{>0}$ and a coherent $W_n\MO_X$-module $\cF$, we set $H^i_{\m}(X, \cF) := H^i_{W_n\m}(X, \cF)$ by abuse of notation. 

We also introduce two conditions on an integer $e_0>0$:
\begin{align} \label{eq:cond-star} \tag{$\star$}
H^j_{\m}(X, \MO_X(-sA))&=0 \quad \text{ for all } j < d \text{ and } s \geq p^{e_0},
\end{align}
and in the case when $-S$ is semiample over $\Spec R$:
\begin{align}
 \label{eq:cond-star2} \tag{$\star\star$}
H^j_{\m}(X, \MO_X(-sA+kS))&=0 \quad \text{ for all } j < d, k \geq 0, \text{and } s \geq p^{e_0}.
\end{align}
Such $e_0$ may be chosen by Lemma \ref{l-easy-Serre} (more precisely, its variant using the Fujita vanishing theorem).
\end{setting}
{The goal of this section is to establish that quasi-$F^e$-stable sections of $L$ lift from $S$, when $-S$ is semiample.}

\begin{remark}
{The main difficulty with working with iterative quasi-$F$-splittings is that one needs to be very careful about the order of quantifiers between the power of Frobenius $e$, and the height $n$. In this subsection we deal with non-uniform version of quasi-$F$-splittings, and so the choice of $e$ cannot depend on $n$.

Assume that $-S$ is semiample and fix an integer $e_0$ satisfying Condition (\ref{eq:cond-star2}). 
We will establish the following surjectivity as long as $e \geq e_0$ (see Theorem \ref{t-IOA-anti-ample}):
\begin{equation*}
q^eS^0_{n,\adj}(X,S+B; L) \to q^eS^0_n(S, B_S; L|_S)
\end{equation*}
assuming that $n$ is large enough. We emphasise that the choice of large enough $n$ depends on $e$ itself, and not just $e_0$. For example, we need $n$ large enough so that:
\[
q^eS^0_{n,\adj}(X,S+B; L) = q^eS^0_{n+1,\adj}(X,S+B; L) = \ldots = q^eS^0_{\adj}(X,S+B; L), 
\]
and similarly for $q^eS^0_n(S, B_S; L|_S)$. In fact, what we need is that  $n$ is large enough so that the above equalities hold for $L$ replaced by $K_X + \lceil p^i A \rceil - kS$ and $B$ replaced by $\{p^iB\}$ for all $i, k \geq 0$\footnote{or to be more precise, for all $0 \leq k \leq k_0$ where $k_0$ is bounded in terms of $i$ (see Remark \ref{remark:jakub-choice-of-k})}. The existence of such $n$ follows by Noetherianity and the stabilisation of quasi-$F^e$-stable sections for every $i \gg 0$ and $k\geq 0$ (see the lemma below).

Last we point out that there are only finitely many possibilities for $\{p^iB\}$ when $i$ is an integer, and so the choice of $i$ will not affect the order of quantifiers.}
\end{remark}

\begin{lem}\label{l-why-q-F-pure}
In the situation of Setting \ref{s-IOA-anti-ample}, fix an integer $e>0$ and a $\pi$-nef Weil divisor $N$. 
Let $n$ be a positive integer such that $(X,\{p^iB\})$ is $n$-quasi-$F^e$-pure for every integer $i \geq 0$.
Then there exists $i_0 \in \Z_{>0}$ such that 
\begin{equation}\label{e1-why-q-F-pure}
q^{e}S^0_n(X,\{p^iB\}; K_X + \lceil p^iA \rceil +kN) =H^0(X,\cO_X(K_X + \lceil p^iA \rceil +kN))
\end{equation}
for all integers $i \geq i_0$ and $k \geq 0$.
\end{lem}
\noindent {We immediately get that \[q^{e}S^0_{n'}(X,\{p^iB\}; K_X + \lceil p^iA \rceil +kN) =H^0(X,\cO_X(K_X + \lceil p^iA \rceil +kN))\]
for all $n' \geq n$, $i \geq i_0$, and $k \geq 0$. The following proof is the same as that of Proposition \ref{prop:Serre-type-result-qS0}. For the 
reader's convenience, 
we repeat it below.}

\begin{proof}
By Definition \ref{definition:quasi-F-stable-sections}, $q^{e}S^0_n(X,\{p^iB\};K_X+\lceil p^iA \rceil + kN)$ is equal to 
\[
\image\Big( 
H^0(X, 
(Q^{e}_{X,-({K_X+} p^iA + kN),n})^*)  
\to H^0(X, \cO_X(K_X+\lceil p^iA + kN \rceil))\Big), 
\]
where $(-)^* := \cHom_{W_n\cO_X}(-, W_n\omega(-K_X))$. 
Define
\[
\cG_{r,k} := {\rm Ker}\Big( (Q^{e}_{X,-({K_X+} rA + kN),n})^* 
\xrightarrow{(\dagger)} \cO_X(K_X+\lceil rA + kN \rceil)\Big).  
\]

\begin{claim} There exists an integer $r_0 > 0$ such that $H^1(X, \cG_{r,k})=0$ for all $r \geq r_0$ and $k \geq 0$.
\end{claim}

Assuming this claim, 
we now finish the proof of Lemma \ref{l-why-q-F-pure}. Pick an integer $i_0$ such that $p^{i_0}\geq r_0$ and set $r = p^i$.
Note that $(\dagger)$ is surjective as {$(X,\{p^iB\})$ is $n$-quasi-$F^{e}$-pure (see Proposition \ref{p-IQFS-L-surje} for $D= -(K_X + \lceil p^iA + kN \rceil)$ and $\{p^iB\} + D = -(K_X+p^iA + kN)$)}. Therefore, the claim immediately implies that $q^{e}S^0_n(X,\{p^iB\};K_X+\lceil p^iA + kN \rceil) = H^0(X, \cO_X(K_X+\lceil p^iA + kN \rceil))$ for all $i \geq i_0$, so it is enough to prove 
{the above} claim.

\medskip

To this end, take 
{$m_0 \in \Z_{>0}$ such that $m_0A$ and $m_0 N$ are Cartier}.
Take integers $r'$ and $k'$ defined by 
$r' := r \ {\rm mod}\ {m_0}$ and $k' := k \ {\rm mod}\ {m_0}$ 
{(cf.\ Subsection \ref{ss:notation}({\ref{ss-nota-mod}}))}. 
{Since $(r-r')A$ and $(k-k')N$ are Cartier, we get} 
\[
\cG_{r,k} = \cG_{r',k'} \otimes_{W_n\MO_X} W_n\cO_X((r-r')A + (k-k')N).
\]
Thus there exists $r_0$ (independent of {$(r, k)$}) 
such that
\[
H^1(X, \cG_{r',k'} \otimes_{W_n\MO_X} W_n\cO_X((r-r')A + (k-k')N))=0
\]
for {every} $r\geq r_0$ 
{and every $k\geq 0$}  
by {the} Fujita vanishing {theorem} (\cite[Theorem 1.5]{keeler03}),  
{because} there are only finitely many {possibilities} 
for $\cG_{r',k'}$.
\end{proof}

The following proposition allows one to run an inductive proof that quasi-$F^e$-stable sections lift.


\begin{prop} \label{p-B-to-A}
In the situation of Setting \ref{s-IOA-anti-ample}, 
take integers $e$ and $e_0$ satisfying $e \geq e_0>0$ 
and Condition $($\ref{eq:cond-star}$)$. 
Suppose that  the following hold. 
\begin{enumerate}
\item 
$(X, {\{p^i B\}})$ is 
quasi-$F^e$-pure for 
every integer $i\geq 0$.  
\item 
$q^eS^0(X,\{p^iB\}; K_X + \lceil p^iA \rceil) =H^0(X,\cO_X(K_X + \lceil p^iA \rceil))$ for {every integer $i \geq 1$}. 
\end{enumerate}
Then 
the map 
\[
q^eS^0_{\adj}(X,S+B;L) \to q^eS^0(S,B_S; L|_S)
\]
from Definition \ref{definitin:restriction-for-qS^0} is surjective 
(cf.\ Remark \ref{remark:implicit-replacement-of-L}).
\end{prop}



\noindent Note that $K_X + \lceil p^iA \rceil = K_X + \{p^iB\} + p^iA$.
\begin{proof}
{Pick an integer $\ell >0$ such that $(X,\{p^iB\})$ is $\ell$-quasi-$F^e$-pure for every $i \geq 0$ and}
pick $i_1 \in \Z_{>0}$ as in Lemma \ref{l-why-q-F-pure}. 
Then
we have 
\[
 q^eS^0_{\ell}(X,\{p^iB\}; K_X + \lceil p^iA \rceil) =H^0(X,\cO_X(K_X + \lceil p^iA \rceil)) 
\]
for every $i>i_1$. 
For every $1 \leq i \leq i_1$, 
(2) enables us to  find $n_i \in \Z_{>0}$ such that 
\[ 
q^eS^0_{n_i}(X,\{p^iB\}; K_X + \lceil p^iA \rceil) =H^0(X,\cO_X(K_X + \lceil p^iA \rceil)). 
\]
Fix integers $n$ and $m$ satisfying 
$n \geq \max \{ \ell, n_1, ..., n_{i_1}\}$ and $m \geq 2n$. 
In particular, 
\begin{equation}\label{e2-B-to-A}
\begin{split}
 q^eS^0_{n'}(X,\{p^iB\}; K_X + \lceil p^iA \rceil) &=H^0(X,\cO_X(K_X + \lceil p^iA \rceil))\\
&\text{if }\qquad n' \geq n, \qquad {i\geq 1}. 
\end{split}
\end{equation}
By Lemma \ref{l-why-Matlis}(2), it is enough to show  that the map 
\[
H_{\m}^{d-1}(X, Q^e_{X,-A,m}) \to 
H_{\m}^{d-1}(X, Q^e_{X,-A,n})
\]
is zero. 
We shall prove a more general statement that 
\[
H_{\m}^{d-1}(X, Q^e_{X,-p^{i-1}A,m-i+1}) \to 
H_{\m}^{d-1}(X,  Q^e_{X,-p^{i-1}A,n-i+1})
\]
is zero by descending induction on $1 \leq i \leq n$. 
This holds when $i=n$, because {Condition (\ref{eq:cond-star})} 
implies 
\[
H_{\m}^{d-1}(X,  Q^e_{X,-p^{n-1}A, 1}) 
\simeq 
H_{\m}^{d-1}(X,  F_*^e\MO_X(-p^{e+n-1}A)) =0. 
\]

Now consider the following  commutative diagram in which each horizontal sequence is exact {
(Lemma \ref{l-BQO})}: 
\[
\begin{tikzcd}
0 \arrow{r} & F_*(B^e_{X,-p^iA,m-i}) \arrow[d, "\alpha^e_{i, m, n}"] \arrow{r} & Q^e_{X,-p^{i-1}A,m-i + 1} \arrow{r} \arrow[d, "\beta^e_{i-1, m, n}"] & F^e_*\mcO_{X}(-p^{e+i-1}A) \arrow[d, equal] \arrow{r} & 0\\
0 \arrow{r} & F_*(B^e_{X,-p^iA,n-i}) \arrow{r} & Q^e_{X,-p^{i-1}A,n-i+1} \arrow{r} & 
F^e_*\mcO_{X}(-p^{e+i-1}A) \arrow{r} & 0.
\end{tikzcd}
\]
It suffices to show the implication 
$H^{d-1}_{\m}(X, \beta^e_{i, m, n}) =0 \Rightarrow 
H^{d-1}_{\m}(X, \beta^e_{i-1, m, n})=0$, 
where 
\[
H^{d-1}_{\m}(X, \beta^e_{i, m, n}): 
H^{d-1}_{\m}(X, Q^e_{X,-p^{i}A,m-i}) \to 
H^{d-1}_{\m}(X, Q^e_{X,-p^{i}A,n-i}). 
\]
By $H^{d-2}_{\m}(X, \MO_X(-p^{e+i-1}A))= H^{d-1}_{\m}(X, \MO_X(-p^{e+i-1}A))=0$ (\ref{eq:cond-star}), 
the above diagram implies that $H^{d-1}_{\m}(X, \beta^e_{i-1, m, n})=0 \Leftrightarrow 
H^{d-1}_{\m}(X, \alpha^e_{i, m, n})=0$. 
Consider the following commutative diagram in which each horizontal sequence is exact {(Remark \ref{r-Q-B-big-diag})}: 
\begin{center}
\begin{tikzcd}[row sep=large, column sep=50]
0 \arrow{r} & \mcO_X(-p^iA) \arrow[d, equal] \arrow[r, "\Phi^e_{X, -p^iA, m-i}"] & Q^e_{X,-p^iA,m-i} \arrow[r, "{\rho}"]  \arrow[d, "\beta^e_{i, m, n}"] & 
B^e_{X,-p^iA,m-i} \arrow[d, "\alpha^e_{i, m, n}"] \arrow{r} & 0\\
0 \arrow{r} & \mcO_X(-p^iA) \arrow[r, "\Phi^e_{X, -p^iA, n-i}"] & Q^e_{X,-p^iA,n-i} \arrow{r} & B^e_{X,-p^iA,n-i} \arrow{r} & 0.
\end{tikzcd}
\end{center}
By $m-i \geq 2n -n = n$, 
(\ref{e2-B-to-A}) implies that $q^eS_{m-i}^0(X,\{p^iB\}; K_X + \lceil p^iA \rceil) =H^0(X, \MO_X(K_X + \lceil p^iA \rceil))$, i.e., 
\[
H^0(X, (Q^e_{X, -K_X-p^iA, m-i})^*) \xrightarrow{H^0(X, (\Phi^e_{X, -K_X-p^iA, m-i})^*)} H^0(X, 
\MO_X(K_X + \lceil p^iA \rceil))
\]
is surjective. 
We get that
\[
H^d_{\m}(X, \Phi^e_{X, -p^iA, m-i}): 
H^{d}_{\m}(X, \mcO_X(-p^iA)) \to H^{d}_{\m}(X, Q^e_{X,-p^iA,m-i})
\]
is injective by Lemma \ref{lem:cohomological-criterion-for-log-quasi-F-splitting}. 
Since $H^{d-1}_{\m}({\rho})$ is surjective, 
we obtain 
$H^{d-1}_{\m}(X, \beta^e_{i, m, n}) =0 \Rightarrow H^{d-1}_{\m}(X, \alpha^e_{i, m, n})=0$. 
To summarise, 
\[
H^{d-1}_{\m}(X, \beta^e_{i, m, n}) =0 \Rightarrow H^{d-1}_{\m}(X, \alpha^e_{i, m, n})=0 
\Leftrightarrow 
H^{d-1}_{\m}(X, \beta^e_{i-1, m, n})=0, 
\]
as required. 
\end{proof}

\begin{prop} \label{p-A-to-B}
Fix $e \in \Z_{>0}$. 
In the situation of Setting \ref{s-IOA-anti-ample}, 
suppose that  the following hold. 
\begin{enumerate}
\item $-S$ is {semiample}.

\item \vspace{0.05em}$(S,B_S)$ is quasi-$F^e$-split. 
\item \vspace{0.05em}
$q^eS^0_{\adj}(X,S+{B}; L-kS) \!\to\! q^eS^0(S,{B_S}; (L-kS)|_S)$ 
is surjective for all $k \geq 0$. 
\end{enumerate}
Then
\begin{equation}\label{e1-A-to-B}
q^eS^0_{{\adj}}(X,{{S+B}}; {L}) = H^0(X,\cO_X(L)).
\end{equation}
{In particular, $q^eS^0(X,B; {L}) = H^0(X,\cO_X(L))$ (cf.\ Lemma \ref{l-qS^0-adj-nonadj}).} 
\end{prop}


\noindent 
We omit some details of the following proof, 
as it is identical to that of \cite[Lemma 4.9]{KTTWYY1}. 


\begin{proof}
Let $\m$ be the maximal ideal of $R$. 
We have $\m \in \pi(S)$. 
First of all, we reduce the problem to showing  that 
 \begin{equation} \label{eq:second-induction}
 H^0(X,\cO_X(L-kS)) + q^eS^0_{{\adj}}(X, {S+B}; {L}) = H^0(X,\cO_X({L}))
 \end{equation}
for every integer {$k \geq 0$}. 
To this end, we first observe that in order to prove (\ref{e1-A-to-B}), 
it suffices to show 
\[
\m H^0(X,\cO_X({L})) + q^eS^0_{{\adj}}(X, {S+B}; {L}) = H^0(X,\cO_X({L}))
\]
by Nakayama's lemma. 
One can show that 
$\m H^0(X,\cO_X({L})) \supseteq H^0(X,\cO_X(L-kS))$ for a large integer $k \gg 0$, and hence (\ref{e1-A-to-B}) holds if 
(\ref{eq:second-induction}) holds for every $k \geq 0$. 
For more details, see \cite[Reduction to (4.9.1) and Reduction to (4.9.3) in the proof of Lemma 4.9]{KTTWYY1}. 

\medskip

We shall prove (\ref{eq:second-induction}) for every integer {$k \geq 0$}  by increasing induction on $k$. Since the base case of the induction {$k = 0$} is trivial, we assume that  (\ref{eq:second-induction})
is true for some $k\geq {0}$ and aim for showing that
\begin{equation} \label{eq:l-IOA-antiample-k+1}
H^0(X,\cO_X(L-(k+1)S)) + q^eS^0_{{\adj}}(X, {S+B}; {L}) = H^0(X,\cO_X({L})).
\end{equation}
{We start by making two observations. First, }note that
\begin{equation} \label{eq:inclusion_of_im_psi-new}
q^eS^0_{\adj}(X,S+B; L) \supseteq q^eS^0_{\adj}(X,S+B; L-kS).    
\end{equation}
Second, observe that 
Proposition \ref{prop:qS^0-for-quasi-F-split} {yields}
\begin{equation}\label{e-IOAlem-QFS-S}
q^eS^0(S, {B_{S}}; (L-kS)|_S)
 = H^0(S, \cO_S(\rup{(L-kS)|_S - {B_{S}}})). 
\end{equation}
Here we used that $(S, \{{B_{S}} - (L-kS)|_S\})$ is quasi-$F^e$-split, which is true by Assumption (2) and $\{{B_{S}} - (L-kS)|_S\} \leq B_S$ 
\cite[Claim 4.10 in the proof of Lemma 4.9]{KTTWYY1}.

Now pick a section $\gamma \in H^0(X,\cO_X(L))$. 
{By the induction hypothesis (\ref{eq:second-induction}), we have $\gamma= \gamma_1 + \gamma_2$ 
for some 
}
$\gamma_1 \in H^0(X,\cO_X(L-kS))$ and 
$\gamma_2 \in q^eS^0_{\adj}(X, S+B;{L})$. 
Then 
\[
\gamma_1|_S \in H^0(S, \cO_S(\rup{(L-kS)|_S-B_{S}}))
\]
(cf.\ the lower {horizontal sequence} {in the diagram (\ref{eq:diagram-res-C})}). 
By (\ref{e-IOAlem-QFS-S}) and Assumption (3),  
there exists a section 
\[
\sigma \in q^eS^0_{\adj}(X,S+{B}; L-kS) \subseteq H^0(X,\cO_X(L-kS)) 
\]
such that $\sigma|_S = \gamma_1|_S$. Hence  
\[
\gamma_1 = (\gamma_1-\sigma) + \sigma,
\]
where $\gamma_1-\sigma \in H^0(X,\cO_X(L-(k+1)S))$ and 
$\sigma \in q^eS^0_{\adj}(X,S+{B}; L)$ 
by (\ref{eq:inclusion_of_im_psi-new}). 
This concludes the proof of the claim 
{(\ref{eq:l-IOA-antiample-k+1})}.
\end{proof}

\begin{remark} \label{remark:jakub-choice-of-k} 
{The above proposition looks very subtle in that for different choices of $k$ we may be forced to pick a different height in assumption (3). Thus, at first glance, it may seem that one cannot work with a fixed height depending on $e$. 
However, this is not really the case as we only need to check this condition for $k\leq k_0$, where $k_0$ is chosen so that
\[
\m H^0(X,\cO_X({L})) \supseteq H^0(X,\cO_X(L-k_0S)).
\]}
\end{remark}

\begin{theorem} \label{t-IOA-anti-ample}
In the situation of Setting \ref{s-IOA-anti-ample}, 
take integers $e$ and $e_0$ satisfying $e \geq e_0>0$ 
and Condition $($\ref{eq:cond-star2}$)$. 
Suppose that  the following hold. 
\begin{enumerate}
\item 
$(X,B)$ is quasi-$F^e$-pure. 
\item $B \geq \{ p^iB\}$ for every $i \in \Z_{\geq 0}$.  
\item $-S$ is semiample. 
\item $(S, B_S)$ is quasi-$F^e$-split. 
\end{enumerate}
Then $q^eS^0_{\adj}(X,S+B; L) \to q^eS^0(S, B_S; L|_S)$ is surjective.  
\end{theorem}


\begin{proof}
For $i \in \Z_{\geq 0}$ and $k \in \Z_{\geq 0}$, we set
\begin{itemize}
\item $B_i := \{ p^iB\}$, 
\item $K_S + B_{S, i} := (K_X+S+B_i)|_S$, 
\item $A_{i, k} := p^iA-kS$, and 
\item $L_{i, k} := \rup{K_X+S+p^iA-kS} = 
\rup{K_X+S+A_{i, k}} = 
K_X+S+ B_i +A_{i, k}$. 
\end{itemize}
Note that 
$B = B_0$, 
$B_{S, 0} = B_S$, $A = A_{0, 0}$, and 
$L = K_X+S+B+A = \rup{K_X+S+A} =L_{0, 0}$. 
All the conditions assumed in Setting \ref{s-IOA-anti-ample}  
hold after replacing $(B, B_S, A, L)$ 
by $(B_i, B_{S, i}, A_{i, k}, L_{i, k})$. 
{After such replacement, 
Condition (\ref{eq:cond-star}) is guaranteed by Condition $($\ref{eq:cond-star2}$)$, too.} {However, note that it is not necessarily  true that $B_i \geq \{p^jB_i\}$.}

By (1), (2), and (4), the following hold. 
\begin{enumerate}
\item[(1')] $(X, B_i)$ is quasi-$F^e$-pure. 
    \item[(4')] $(S, B_{S, i})$ is quasi-$F^e$-split. 
\end{enumerate}
For $i \geq 0$, consider the following two assertions. 
\begin{enumerate}
\item[$({\rm a}_{i})$] 
The following map is surjective for every $k \geq 0$:  
\[
q^eS^0_{\adj}(X,S+B_i; L_{i, k}) \!\to\! q^eS^0(S, B_{S, i}; L_{i, k}|_S).
\]
\item[$({\rm b}_{i})$] 
 $q^eS^0(X, B_i; L_{i, k}) = H^0(X, \MO_X(L_{i, k}))$ 
 for all $k \geq 0$. 
\end{enumerate}
We prove $({\rm a}_i)$ and $({\rm b}_i)$ by descending induction on $i$. 
We can find $i_1\in \Z_{>0}$ such that $({\rm b}_i)$ holds for every $i \geq i_1$ 
(Lemma \ref{l-why-q-F-pure}).
Fix 
an integer $i \geq 0$. 
It suffices to show $(\alpha)$ and $(\beta)$ below. 
\begin{enumerate}
\item[$(\alpha)$] If $({\rm a}_i)$ holds, then $({\rm b}_i)$ holds. 
\item[$(\beta)$] If $({\rm b}_{i'})$ holds for every $i' \geq i+1$, then $({\rm a}_{i})$ holds. 
\end{enumerate}


Let us show $(\alpha)$. 
Assume $({\rm a}_i)$. Fix $k \geq 0$. 
We now apply Proposition \ref{p-A-to-B} 
after replacing $(B, B_S, A, L)$ 
by $(B_i, B_{S, i}, A_{i, k}, L_{i, k})$. 
The assumptions (1), (2), and (3) of Proposition \ref{p-A-to-B} 
follow from (3), (4'), and $({\rm a}_i)$, respectively. 
We then obtain $q^{e}S^0(X,{B_i}; L_{i, k}) = H^0(X,\cO_X(L_{i, k}))$, i.e., $({\rm b}_i)$ holds. 
Thus $(\alpha)$ holds. 

\medskip

Let us show $(\beta)$. 
Assume that $({\rm b}_{i'})$ holds for every $i' \geq i+1$. 
Fix $k \geq 0$. 
In order to apply Proposition \ref{p-B-to-A} after replacing $(B, B_S, A, L)$ 
by $(B_{i}, B_{S, {i}}, A_{i, k}, L_{i, k})$, 
we now verify  its assumptions. 
First, Condition (\ref{eq:cond-star}) is satisfied as mentioned earlier. Further, 
for every $j \in \Z_{\geq 0}$, it holds that 
$\{p^j B_{i}\} = \{ p^j \{ p^{i} B\}\} = \{ p^{i+j} B\} = B_{i+j}$ and 
\[
p^jA_{i, k} = p^j (p^{i}A-kS) = p^{i+j}A -kp^j S = A_{i+j, kp^j}.  
\]
Thus {Assumption (1) in} Proposition \ref{p-B-to-A} holds by (1'). 
As for {Assumption (2) in} Proposition \ref{p-B-to-A}, it is enough to check the following. 
\begin{enumerate}
\item[($\dagger$)] 
$q^eS^0(X, \{ p^jB_{i}\}; K_X+\rup{p^jA_{i, k}})= H^0(X, \MO_X(K_X+\rup{p^jA_{i, k}}))$ 
for every $j\geq 1$. 
\end{enumerate}
We have $K_X + \rup{p^jA_{i, k}} = K_X + \rup{A_{i+j, kp^j}} = L_{i+j,  kp^j} -S {= L_{i+j,  kp^j+1}}$. Therefore, ($\dagger$) is equivalent to $(\ddagger)$ below. 
\begin{enumerate}
\item[$(\ddagger)$] 
$q^eS^0(X, B_{i+j}; {L_{i+j,  kp^j+1}})= H^0(X, \MO_X({L_{i+j,  kp^j+1}}))$ 
for every $j>0$. 
\end{enumerate}
Then $(\ddagger)$ holds, because we are assuming that $({\rm b}_{i+j})$ holds for every $j \geq 1$. 
This completes the proof of  $(\beta)$. 
\qedhere

\end{proof}

\begin{corollary} \label{cor:qSadjoint-is-everything-semiample}
With the same assumptions as in Theorem \ref{t-IOA-anti-ample}, we have that $q^eS^0_{\adj}(X,S+B; L) = H^0(X,\cO_X(L))$.
\end{corollary}
\begin{proof}
By applying Theorem \ref{t-IOA-anti-ample} with $L$ replaced by $L - kS$, we get that
\[
q^eS^0_{\adj}(X,S+B; L-kS) \to q^eS^0(S, B_S; (L-kS)|_S)
\]
is surjective for every $k \geq 0$. Then we can conclude by Proposition \ref{p-A-to-B}.
\end{proof}

\begin{corollary} \label{c-IOA-anti-ample}
In the situation of Setting \ref{setting:foundations-of-log-quasi-F-splitting},  suppose that $R$ is a local ring and $H^0(X, \MO_X)=R$. 
Let $(X,S+B)$ be a divisorially Cohen-Macaulay $\Q$-factorial 
plt pair, where 
$S$ is a normal prime divisor and 
$B$ is a $\Q$-divisor such that $\rdown{B}=0$. 
Suppose that $-(K_X+S+B)$ is ample and the following hold. 
\begin{enumerate}
\item 
$(X, B)$ is quasi-$F^e$-pure for every $e \in \Z_{>0}$. 
\item $B \geq \{ p^i B\}$ for every $i \in \Z_{>0}$. 
\item $-S$ is semiample.  
\item $(S, B_S)$ is quasi-$F^e$-split for every $e \in \Z_{>0}$,  where $K_S + B_S = (K_X+S+B)|_S$. 
\end{enumerate}
Then $(X, S+B)$ is purely quasi-$F^e$-split for every $e \in \Z_{>0}$. 
\end{corollary}


\begin{proof}Set $A := -(K_X+S+B)$. Take  $e_0 \in \Z_{>0}$ such that  
\[ 
H^j_{\m}(X, \MO_X(-sA+kS))=0 \quad  \text{for all integers } j<d, k \geq 0, \text{and } s \geq p^{e_0}, 
\]
whose existence is guaranteed by {a variant} of Lemma \ref{l-easy-Serre} {using the Fujita vanishing theorem}. 
Fix an integer $e  \geq e_0$. Then 
it follows from Corollary \ref{cor:qSadjoint-is-everything-semiample} 
that
\[
q^eS^0_{\rm adj}(X,S+B;0) = H^0(X,\cO_X),
\]
and so $(X,S+B)$ is purely quasi-$F^e$-split (see Definition \ref{d-pure-IQFS}). 
\qedhere

\end{proof}

\section{Inversion of adjunction for uniform quasi-$F^{\infty}$-splitting}





\begin{lem}\label{l-why-Matlis-uni}
In the situation of {Setting \ref{s-IOA-anti-ample}}, 
assume that $S$ is divisorially Cohen-Macaulay.
Then the following hold. 
\begin{enumerate}
\item 
If $H^{d-1}_{\m}(X, Q^{\perf}_{X, K_X+S+B-L, n })=0$ for some integer $n>0$, then 
\[
\q S^0_{n, \adj}(X, S+B; L) \to \q S^0_n(S, B_S; L|_S)
\]
is surjective. 
\item 
Suppose that for every $n_0$, there exist integers 
$m$ and $n$ such that $m > n \geq n_0$  and  
\[
H^{d-1}_{\m}(X, Q^{\perf}_{X, K_X+S+B-L, m }) \to H^{d-1}_{\m}(X, Q^{\perf}_{X, K_X+S+B-L, n})
\]
is zero. Then 
\[
\qS^0_{\adj}(X, S+B; L) \to \qS^0(S, B_S; L|_S)
\]
is surjective. 
\end{enumerate}
\end{lem}


\begin{proof}

As in the proof of Lemma \ref{l-why-Matlis}, we
have the following commutative diagram: 
\begin{equation} \label{e1-why-Matlis-uni}
\begin{tikzcd}
H^d_{\m}(X, Q^{S, \perf}_{X, K_X+S+B-L, n}) \arrow[d, hookleftarrow]  \arrow[r, leftarrow, "{\alpha_n}"] & 
H^{d-1}_{\m}(S, Q^{\perf}_{S, K_S+B_S-L|_S, n})  \arrow[d, hookleftarrow, "j_n"]\\
I_n 
\arrow[d,  twoheadleftarrow]  \arrow[r, leftarrow, "{\beta_n}"] & 
J_n  \arrow[d, twoheadleftarrow] \\
H^d_{\m}(X, {\cO_X(K_X-L)}) \arrow[r, leftarrow, "{\gamma_n}"] & 
H^{d-1}_{\m}(S, \cO_S(K_S+B_S-L|_S)), 
\end{tikzcd}
\end{equation}
where $I_n$ and $J_n$ are defined as the images of the vertical arrows. Note that by Proposition \ref{prop:local-cohomology-description-of-S^0-perf} and an analogous statement in the adjoint case, we have that
\begin{align*}I_n^{\vee}  &\simeq {\q} 
S^0_{n,\adj}(X,S+B;L)^{\wedge},\quad \text{ and } \\
J_n^{\vee}  &\simeq {\q} 
S^0_n(S,B_S;L|_S)^{\wedge}.
\end{align*}
Let us show (1). 
By $H^{d-1}_{\m}(X, Q^{\perf}_{X, K_X+S+B-L, n})=0$,  we have that
 $\alpha_n$ is injective. 
By chasing the diagram (\ref{e1-why-Matlis-uni}),  $\beta_n$ is injective. 
Thus $\beta_n^{\vee} \colon \q S^0_{n, \adj}(X, S+B; L) \to \q S^0_n(S, B_S; L|_S)$ is surjective 
and hence (1) holds. 


Let us show (2). 
Pick $n_0 \gg 0$ and take integers $m$ and $n$ as in the statement of (2). 
In particular, $m > n \geq n_0 \gg 0$. 
Recall that we have 
\begin{align*}
{\q}S^0_1(S, B_S; L|_S) \subseteq 
\, \cdots\,  &\subseteq 
\q S^0_n(S, B_S; L|_S) \\
&= {\q}S^0_{n+1}(S, B_S; L|_S) = \, \cdots\, =
\qS^0(S, B_S; L|_S). 
\end{align*}
Applying $(-)^{\vee} \,(=\Hom_{W_eR}(-, E))$, we obtain 
\begin{align*}
{\q}S^0_1(S, B_S; L|_S)^{\vee} 
\twoheadleftarrow \,\cdots\,  &\twoheadleftarrow 
{\q}S^0_n(S, B_S; L|_S)^{\vee} \\
&\xleftarrow{\simeq} {\q}S^0_{n+1}(S, B_S; L|_S)^{\vee} \xleftarrow{\simeq} \, \cdots \, \xleftarrow{\simeq}
\qS^0(S, B_S; L|_S)^{\vee}. 
\end{align*}
Similarly, we may assume that 
${\q}S^0_{n,\adj}(X,S+B;L)^{\vee}
 \xleftarrow{\simeq}
{\q}S^0_{m,\adj}(X,S+B;L)^{\vee}$. 
It is enough to show that $\beta_n$ is injective.\\

To this end, pick $\zeta_n \in {\q}S^0_n(S, B_S; L|_S)^{\vee}$ 
such that $\beta_n(\zeta_n)=0$. 
Via the induced isomorphism \[
{\q}S_n^0(S, B_S; L|_S) \xleftarrow{\simeq}
{\q}S_m^0(S, B_S; L|_S),
\]
there is an element $\zeta_m \in {\q}S_m^0(S, B_S; L|_S)^{\vee}$ corresponding to $\zeta_n \in {\q}S_n^0(S, B_S; L|_S)^{\vee}$. 
Therefore, $\beta_m(\zeta_m)=0$. 
Then $j_m(\zeta_m)$ is the image of some element 
\[
\xi_m \in H^{d-1}_{\m}(X, Q^\perf_{X, K_X+S+B-L, m }).
\]
By our assumption, the image of $\xi_m$ in $H^{d-1}_{\m}(X, Q^\perf_{X, K_X+S+B-L, n })$ is zero. 
Then $j_n(\zeta_n)=0$, and hence $\zeta_n=0$. 
Thus (2) holds. 
\qedhere

\end{proof}

\subsection{Uniformly quasi-$F^{\infty}$-split case}



\begin{lem}\label{l-why-q-F-pure-uni}
In the situation of Setting \ref{s-IOA-anti-ample}, 
fix a $\pi$-nef Weil divisor $N$. 
Let $n$ be a positive integer such that $(X,\{p^iB\})$ is $n$-quasi-$F^\infty$-pure for every integer $i \geq 0$.
Then there exists $i_1 \in \Z_{>0}$ such that 
\begin{equation}\label{e1-why-q-F-pure-uni}
q^{\infty}S^0_n(X,\{p^iB\}; K_X + \lceil p^iA \rceil +kN) =H^0(X,\cO_X(K_X + \lceil p^iA \rceil +kN))
\end{equation}
for all integers $i \geq i_1$ and $k \geq 0$.
\end{lem}
\begin{proof}
This follows immediately from Proposition \ref{prop:Serre-type-result-qS0}.
\end{proof}

\begin{prop} \label{p-B-to-A-uni}
In the situation of Setting \ref{s-IOA-anti-ample}, 
assume that $S$ is divisorially Cohen-Macaulay.
Suppose that  the following hold. 
\begin{enumerate}
\item 
$(X, {\{p^i B\}})$ is 
uniformly quasi-$F^\infty$-pure for 
every integer $i\geq 0$.  
\item 
$\qS^0(X,\{p^iB\}; K_X + \lceil p^iA \rceil) =H^0(X,\cO_X(K_X + \lceil p^iA \rceil))$ for {all integers $i \geq 1$}. 
\end{enumerate}
Then 
the map 
\[
\qS^0_{\adj}(X,S+B;L) \to \qS^0(S,B_S; L|_S)
\]
induced from Definition \ref{definitin:restriction-for-qS^0} is surjective 
(cf.\ Remark \ref{remark:implicit-replacement-of-L}).
\end{prop}



\noindent Note that $K_X + \lceil p^iA \rceil = K_X + \{p^iB\} + p^iA$.
\begin{proof}
{Pick an integer $\ell >0$ such that $(X,\{p^iB\})$ is $\ell$-quasi-$F^\infty$-pure for every $i \geq 0$ and}
pick $i_1 \in \Z_{>0}$ as in Lemma \ref{l-why-q-F-pure-uni}. 
Then
we have 
\[
 q^\infty S^0_{\ell}(X,\{p^iB\}; K_X + \lceil p^iA \rceil) =H^0(X,\cO_X(K_X + \lceil p^iA \rceil)) 
\]
for every $i>i_1$. 
For every $1 \leq i \leq i_1$, 
(2) enables us to  find $n_i \in \Z_{>0}$ such that 
\[ 
q^\infty S^0_{n_i}(X,\{p^iB\}; K_X + \lceil p^iA \rceil) =H^0(X,\cO_X(K_X + \lceil p^iA \rceil)). 
\]
Fix integers $n$ and $m$ satisfying 
$n \geq \max \{ \ell, n_1, ..., n_{i_1}\}$ and $m \geq 2n$. 
In particular, 
\begin{equation}\label{e2-B-to-A-uni}
\begin{split}
 q^\infty S^0_{n'}(X,\{p^iB\}; K_X + \lceil p^iA \rceil) &=H^0(X,\cO_X(K_X + \lceil p^iA \rceil))\\
&\text{if }\qquad n' \geq n, \qquad {i\geq 1}. 
\end{split}
\end{equation}
By Lemma \ref{l-why-Matlis-uni}(2), it is enough to show  that the map 
\[
H_{\m}^{d-1}(X, Q^\perf_{X,-A,m}) \to 
H_{\m}^{d-1}(X, Q^\perf_{X,-A,n})
\]
is zero. 
We shall prove a more general statement that 
\[
H_{\m}^{d-1}(X, Q^\perf_{X,-p^{i-1}A,m-i+1}) \to 
H_{\m}^{d-1}(X,  Q^\perf_{X,-p^{i-1}A,n-i+1})
\]
is zero by descending induction on $1 \leq i \leq n$. 
This holds when $i=n$, because
\[
H_{\m}^{d-1}(X,  Q^\perf_{X,-p^{n-1}A, 1}) 
\simeq 
H_{\m}^{d-1}(X,  \theta_*\MO_{X^{\rm perf}}(-p^{n-1}A)) =0 
\]
by Lemma \ref{l-easy-Serre} and the fact that colimits commute with local cohomology, {where $\theta \colon X^{\perf} \to X$ denotes the canonical projection.}

Now consider the following  commutative diagram in which each horizontal sequence is exact (built from Lemma \ref{l-BQO} by taking 
{the colimit $\varinjlim_e$}):  
\[
\begin{tikzcd}
0 \arrow{r} & F_*(B^\perf_{X,-p^iA,m-i}) \arrow[d, "\alpha_{i, m, n}"] \arrow{r} & Q^\perf_{X,-p^{i-1}A,m-i + 1} \arrow{r} \arrow[d, "\beta_{i-1, m, n}"] & \theta_*\mcO_{X^\perf}(-p^{i-1}A) \arrow[d, equal] \arrow{r} & 0\\
0 \arrow{r} & F_*(B^\perf_{X,-p^iA,n-i}) \arrow{r} & Q^\perf_{X,-p^{i-1}A,n-i+1} \arrow{r} & 
\theta_*\mcO_{X^\perf}(-p^{i-1}A) \arrow{r} & 0.
\end{tikzcd}
\]
It suffices to show the implication 
$H^{d-1}_{\m}(X, \beta_{i, m, n}) =0 \Rightarrow 
H^{d-1}_{\m}(X, \beta_{i-1, m, n})=0$, 
where 
\[
H^{d-1}_{\m}(X, \beta_{i, m, n}): 
H^{d-1}_{\m}(X, Q^\perf_{X,-p^{i}A,m-i}) \to 
H^{d-1}_{\m}(X, Q^\perf_{X,-p^{i}A,n-i}). 
\]
By $H^{d-2}_{\m}(X, \theta_*\MO_{X^\perf}(-p^{i-1}A))= H^{d-1}_{\m}(X, \theta_*\MO_{X^\perf}(-p^{i-1}A))=0$ (Lemma \ref{l-easy-Serre}), 
the above diagram implies that $H^{d-1}_{\m}(X, \beta_{i-1, m, n})=0 \Leftrightarrow 
H^{d-1}_{\m}(X, \alpha_{i, m, n})=0$. 
Consider the following commutative diagram in which each horizontal sequence is exact {(built from Remark \ref{r-Q-B-big-diag} by taking {the colimit $\varinjlim_e$})}: 
\begin{center}
\begin{tikzcd}[row sep=large, column sep=50]
0 \arrow{r} & \mcO_X(-p^iA) \arrow[d, equal] \arrow[r, "\Phi^\perf_{X, -p^iA, m-i}"] & Q^\perf_{X,-p^iA,m-i} \arrow[r, "{\rho}"]  \arrow[d, "\beta_{i, m, n}"] & 
B^\perf_{X,-p^iA,m-i} \arrow[d, "\alpha_{i, m, n}"] \arrow{r} & 0\\
0 \arrow{r} & \mcO_X(-p^iA) \arrow[r, "\Phi^\perf_{X, -p^iA, n-i}"] & Q^\perf_{X,-p^iA,n-i} \arrow{r} & B^\perf_{X,-p^iA,n-i} \arrow{r} & 0.
\end{tikzcd}
\end{center}
By $m-i \geq 2n -n = n$, 
(\ref{e2-B-to-A-uni}) implies that $
{q^{\infty}S}_{m-i}^0(X,\{p^iB\}; K_X + \lceil p^iA \rceil) =H^0(X, \MO_X(K_X + \lceil p^iA \rceil))$.
By Proposition \ref{prop:local-cohomology-description-of-S^0-perf}, we get that
\[
H^d_{\m}(X, \Phi^\perf_{X, -p^iA, m-i}): 
H^{d}_{\m}(X, \mcO_X(-p^iA)) \to H^{d}_{\m}(X, Q^\perf_{X,-p^iA,m-i})
\]
is injective. 
Since $H^{d-1}_{\m}({\rho})$ is surjective, 
we obtain 
$H^{d-1}_{\m}(X, \beta_{i, m, n}) =0 \Rightarrow H^{d-1}_{\m}(X, \alpha_{i, m, n})=0$. 
To summarise, 
\[
H^{d-1}_{\m}(X, \beta_{i, m, n}) =0 \Rightarrow H^{d-1}_{\m}(X, \alpha_{i, m, n})=0 
\Leftrightarrow 
H^{d-1}_{\m}(X, \beta_{i-1, m, n})=0, 
\]
as required. 
\end{proof}

\begin{prop} \label{p-A-to-B-uni}
In the situation of Setting \ref{s-IOA-anti-ample}, 
suppose that  the following hold. 
\begin{enumerate}
\item $-S$ is {semiample}.

\item \vspace{0.05em}$(S,B_S)$ is uniformly quasi-$F^\infty$-split. 
\item \vspace{0.05em}
$\qS^0_{\adj}(X,S+{B}; L-kS) \!\to\! \qS^0(S,{B_S}; (L-kS)|_S)$ 
is surjective for all $k \geq 0$. 
\end{enumerate}
Then
\begin{equation}\label{e1-A-to-B-uni}
\qS^0_{{\adj}}(X,{{S+B}}; {L}) = H^0(X,\cO_X(L)).
\end{equation}
{In particular, $\qS^0(X,B; {L}) = H^0(X,\cO_X(L))$  (cf.\ Lemma \ref{l-qS^0-adj-nonadj}, Definition \ref{def:all-stable-F^infty-sections}).} 
\end{prop}
\begin{proof}
{The same argument as in Proposition \ref{p-A-to-B} works 
after replacing $q^eS^0(-)$ by $\qS^0(-)$.} 
\end{proof}

\begin{theorem} \label{t-IOA-anti-ample-uni}
In the situation of Setting \ref{s-IOA-anti-ample}, 
suppose that  the following hold. 
\begin{enumerate}
\item 
$(X,B)$ is uniformly quasi-$F^\infty$-pure. 
\item $B \geq \{ p^iB\}$ for every $i \in \Z_{\geq 0}$.  
\item $-S$ is semiample. 
\item $(S, B_S)$ is uniformly quasi-$F^\infty$-split 
and $S$ is divisorially Cohen-Macaulay.
\end{enumerate}
Then $\qS^0_{\adj}(X,S+B; L) \to \qS^0(S, B_S; L|_S)$ is surjective.  
\end{theorem}


\begin{proof}
For $i \in \Z_{\geq 0}$ and $k \in \Z_{\geq 0}$, we set
\begin{itemize}
\item $B_i := \{ p^iB\}$, 
\item $K_S + B_{S, i} := (K_X+S+B_i)|_S$, 
\item $A_{i, k} := p^iA-kS$, and 
\item $L_{i, k} := \rup{K_X+S+p^iA-kS} = 
\rup{K_X+S+A_{i, k}} = 
K_X+S+ B_i +A_{i, k}$. 
\end{itemize}
Note that 
$B = B_0$, 
$B_{S, 0} = B_S$, $A = A_{0, 0}$, and 
$L = K_X+S+B+A = \rup{K_X+S+A} =L_{0, 0}$. 
All the conditions assumed in Setting \ref{s-IOA-anti-ample}  
hold after replacing $(B, B_S, A, L)$ 
by $(B_i, B_{S, i}, A_{i, k}, L_{i, k})$. {However, note that it is not necessarily  true that $B_i \geq \{p^jB_i\}$.}

By (1), (2), and (4), the following hold. 
\begin{enumerate}
\item[(1')] $(X, B_i)$ is uniformly quasi-$F^\infty$-pure. 
    \item[(4')] $(S, B_{S, i})$ is uniformly quasi-$F^\infty$-split. 
\end{enumerate}
For $i \geq 0$, consider the following two assertions. 
\begin{enumerate}
\item[$({\rm a}_{i})$] 
The following map is surjective for every $k \geq 0$:  
\[
\qS^0_{\adj}(X,S+B_i; L_{i, k}) \!\to\! \qS^0(S, B_{S, i}; L_{i, k}|_S).
\]
\item[$({\rm b}_{i})$] 
 $\qS^0(X, B_i; L_{i, k}) = H^0(X, \MO_X(L_{i, k}))$ 
 for all $k \geq 0$. 
\end{enumerate}
We prove $({\rm a}_i)$ and $({\rm b}_i)$ by descending induction on $i$. 
We can find $i_1\in \Z_{>0}$ such that $({\rm b}_i)$ holds for every $i \geq i_1$ 
(Lemma \ref{l-why-q-F-pure-uni}).
Fix 
an integer $i \geq 0$. 
It suffices to show $(\alpha)$ and $(\beta)$ below. 
\begin{enumerate}
\item[$(\alpha)$] If $({\rm a}_i)$ holds, then $({\rm b}_i)$ holds. 
\item[$(\beta)$] If $({\rm b}_{i'})$ holds for every $i' \geq i+1$, then $({\rm a}_{i})$ holds. 
\end{enumerate}


Let us show $(\alpha)$. 
Assume $({\rm a}_i)$. Fix $k \geq 0$. 
We now apply Proposition \ref{p-A-to-B-uni} 
after replacing $(B, B_S, A, L)$ 
by $(B_i, B_{S, i}, A_{i, k}, L_{i, k})$. 
The assumptions (1), (2), and (3) of Proposition \ref{p-A-to-B-uni} 
follow from (3), (4'), and $({\rm a}_i)$, respectively. 
We then obtain $\qS^0(X,{B_i}; L_{i, k}) = H^0(X,\cO_X(L_{i, k}))$, i.e., $({\rm b}_i)$ holds. 
Thus $(\alpha)$ holds. 

\medskip

Let us show $(\beta)$. 
Assume that $({\rm b}_{i'})$ holds for every $i' \geq i+1$. 
Fix $k \geq 0$. 
In order to apply Proposition \ref{p-B-to-A-uni} after replacing $(B, B_S, A, L)$ 
by $(B_{i}, B_{S, {i}}, A_{i, k}, L_{i, k})$, 
we now verify  its assumptions. 
Further, 
for every $j \in \Z_{\geq 0}$, it holds that 
$\{p^j B_{i}\} = \{ p^j \{ p^{i} B\}\} = \{ p^{i+j} B\} = B_{i+j}$ and 
\[
p^jA_{i, k} = p^j (p^{i}A-kS) = p^{i+j}A -kp^j S = A_{i+j, kp^j}.  
\]
Thus {Assumption (1) in} Proposition \ref{p-B-to-A-uni} holds by (1'). 
As for {Assumption (2) in} Proposition \ref{p-B-to-A-uni}, it is enough to check the following. 
\begin{enumerate}
\item[($\dagger$)] 
$\qS^0(X, \{ p^jB_{i}\}; K_X+\rup{p^jA_{i, k}})= H^0(X, \MO_X(K_X+\rup{p^jA_{i, k}}))$ 
for every $j\geq 1$. 
\end{enumerate}
We have $K_X + \rup{p^jA_{i, k}} = K_X + \rup{A_{i+j, kp^j}} = L_{i+j,  kp^j} -S {= L_{i+j,  kp^j+1}}$. Therefore, ($\dagger$) is equivalent to $(\ddagger)$ below. 
\begin{enumerate}
\item[$(\ddagger)$] 
$\qS^0(X, B_{i+j}; {L_{i+j,  kp^j+1}})= H^0(X, \MO_X({L_{i+j,  kp^j+1}}))$ 
for every $j>0$. 
\end{enumerate}
Then $(\ddagger)$ holds, because we are assuming that $({\rm b}_{i+j})$ holds for every $j \geq 1$. 
This completes the proof of  $(\beta)$. 
\qedhere

\end{proof}

\begin{corollary} \label{cor:qSadjoint-is-everything-semiample-uni}
With the same assumptions as in Theorem \ref{t-IOA-anti-ample-uni}, we have that \[
\qS^0_{\adj}(X,S+B; L) = H^0(X,\cO_X(L)).
\]
\end{corollary}
\begin{proof}
By applying Theorem \ref{t-IOA-anti-ample-uni} with $L$ replaced by $L - kS$, we get that
\[
\qS^0_{\adj}(X,S+B; L-kS) \to \qS^0(S, B_S; (L-kS)|_S)
\]
is surjective for every $k \geq 0$. Then we can conclude by Proposition \ref{p-A-to-B-uni}.
\end{proof}

\begin{corollary} \label{c-IOA-anti-ample-uni}
In the situation of Setting \ref{setting:foundations-of-log-quasi-F-splitting},  suppose that $R$ is a local ring and $H^0(X, \MO_X)=R$. 
Let $(X,S+B)$ be a divisorially Cohen-Macaulay $\Q$-factorial 
plt pair, where 
$S$ is a normal prime divisor and 
$B$ is a $\Q$-divisor such that $\rdown{B}=0$. 
 {Suppose that $-(K_X+S+B)$ is ample} and the following hold. 
\begin{enumerate}
\item $(X, B)$ is uniformly quasi-$F^{\infty}$-pure. 
\item $\{ p^i B\} \leq B$ for every $i \in \Z_{>0}$. 
\item $-S$ is semiample.  
\item $(S, B_S)$ is uniformly quasi-$F^{\infty}$-split,  where $K_S + B_S = (K_X+S+B)|_S$ 
and $S$ is divisorially Cohen-Macaulay.
\end{enumerate}
Then $(X, S+B)$ is purely uniformly quasi-$F^{\infty}$-split. 
\end{corollary}

\begin{proof}
{By the same argument} as in Corollary \ref{c-IOA-anti-ample}, this follows immediately from Corollary \ref{cor:qSadjoint-is-everything-semiample-uni}. 
\end{proof}

\subsection{Quasi-$+$-regular case}


In what follows, we shall work in the situation of  Setting \ref{s-IOA-anti-ample} until the end of this {subsection}. 
Recall that 
\begin{eqnarray*}
qB^0_{n, \adj}(X,S+B;L) &:=& \bigcap_{f: (Y, T) \to (X, S)} q^fB^0_{n,\adj}(X,S+B;L),\\
qB^0_{\adj}(X,S+B;L) &:=& \bigcup_{n\geq 1} qB^0_{n, \adj}(X,S+B;L) =   qB^0_{N, \adj}(X,S+B;L),
\end{eqnarray*}
where $f: (Y, T) \to (X, S)$ 
runs over all finite covers and 
$N$ is a sufficiently large integer (note that $qB^0_n(-) \subseteq qB^0_{n+1}(-))$). 
As before, we shall apply the above definitions even when $S \subseteq \Supp L$ (cf.\ Remark \ref{remark:implicit-replacement-of-L}). 

Recall that $(X, S+B)$ is purely globally $n$-quasi-$+$-regular if and only if 
$qB^0_{n, \adj}(X,S+B;0) = H^0(X, \MO_X)$. 
Similarly, $(S, B_S)$ is globally $n$-quasi-$+$-regular if and only if 
$qB^0_n(S, B_S;0) = H^0(S, \MO_S)$. 



\begin{theorem} \label{t-IOA-anti-ample+}
In the situation of Setting \ref{s-IOA-anti-ample}, 
suppose that  the following hold. 
\begin{enumerate}
\item $(X, B)$ is locally quasi-$+$-regular. 
\item $B \geq \{ p^iB\}$ for every $i \in \Z_{\geq 0}$.  
\item 
 $S$ is divisorially Cohen-Macaulay and
 $-S$ is semiample. 
\item $(S, B_S)$ is globally quasi-+-regular. 
\end{enumerate}
Then
\[
qB^0_{\adj}(X, S+B; L) \to qB^0(S, B_S; L|_S)
\]
is surjective.
Moreover,
$\qS^0_{\adj}(X, S+B; L)=H^0(X, \cO_X(L)).$
\end{theorem}

\begin{proof}
 If 
 \[
 \qS^0(S, B_S; L|_S) \overset{(\star)}{=} qB^0(S, B_S; L|_S) \,\,\,
 {\rm and} \,\,\,
 \qS^0_{\adj}(X,S+B; L) \overset{(\star\star)}{=} 
 qB^0_{\adj}(X,S+B; L),  
 \]
then the assertion holds by Theorem \ref{t-IOA-anti-ample-uni} {and Corollary \ref{cor:qSadjoint-is-everything-semiample-uni}}. 
The adjoint version of Theorem \ref{t-qB^0-vs-qS^0} implies 
$(\star\star)$.

It is enough to show ($\star$).
The following hold: 
\begin{itemize}
\item $\qS^0(S, B_S; L|_S) = \qS^0(S, B_S-\{L|_S\};\rdown{L|_S})$ (Lemma \ref{l:special-cases-of-L}, Definition \ref{def:all-stable-F^infty-sections}).
\item 
$qB^0(S, B_S; L|_S)= qB^0(S, B_S-\{L|_S\};\rdown{L|_S})$
(Definition \ref{d-def-+-section}). 
\item 
$B_S - \{L|_S\} \geq 0$ (Lemma \ref{l:special-cases-of-L}). 
\item The $\Q$-divisor 
$\rdown{L|_S}-(K_S+B_S-\{L|_S\})=A|_S$ is ample. 
\end{itemize}
Fix $s \geq 0$. 
By Theorem \ref{t-qB^0-vs-qS^0}, it suffices to show that 
$(S, \{ p^s(B_S-\{L|_S\})\})$ is locally quasi-$+$-regular. 
In view of Assumption (4), 
the problem is reduced to proving $\{ p^s(B_S-\{L|_S\})\} \leq B_S$. 
For $K_S+B_{S,s}:=(K_X+S+\{p^sB\})|_S$ and $L_s:=K_X+S+ \rup{p^sA}$, we have 
$0\leq B_{S,s}\leq B_S$ by Assumption (2). 
It holds that 
\[
L_s|_S = (K_X+S+ \rup{p^sA})|_S= 
(K_X+S+ p^sA+\{p^sB\})|_S= K_S+B_{S, s} +p^s(A|_S), 
\]
where the second equality follows from $-p^sA = \rdown{-p^sA} + \{ -p^sA\} 
= -\rup{p^sA} + \{ p^sB\}$. 
Lemma \ref{l:special-cases-of-L} implies $\{L_s|_S\}\leq B_{S,s}$. 
Therefore, we obtain 
\[
\{p^s(B_S-\{L|_S\})\}
=\{-p^s(A|_S)\} = \{ B_{S,s}-L_s|_S\}
\]
\[
= \{ B_{S,s}- \{L_s|_S\} \}=B_{S,s}-\{L_s|_S\}
\leq B_{S,s}\leq B_S, 
\]
as required.
\end{proof}

\begin{corollary} \label{c-IOA-anti-ample+}
In the situation of Setting \ref{setting:foundations-of-log-quasi-F-splitting},  suppose that $R$ is a local ring and $H^0(X, \MO_X)=R$. 
Let $(X,S+B)$ be a divisorially Cohen-Macaulay $\Q$-factorial 
plt pair, where 
$S$ is a normal prime divisor and 
$B$ is a $\Q$-divisor such that $\rdown{B}=0$. 
 {Suppose that $-(K_X+S+B)$ is ample} and the following hold. 
\begin{enumerate}
\item $(X, B)$ is locally quasi-$+$-regular. 
\item $\{ p^i B\} \leq B$ for every $i \in \Z_{>0}$. 
\item 
 $S$ is divisorially Cohen-Macaulay and
 $-S$ is semiample.  
\item $(S, B_S)$ is globally quasi-+-regular,  where $K_S + B_S = (K_X+S+B)|_S$. 
\end{enumerate}
Then $(X, S+B)$ is purely globally quasi-+-regular. 
\end{corollary}

\begin{proof}
By applying Theorem \ref{t-IOA-anti-ample+} with $L=0$, 
we get 
$qB^0_{\adj}(X,S+B;0)=H^0(X,\cO_X)$, and hence $(X,S+B)$ is purely globally quasi-$+$-regular by Definition \ref{d purely globally Q+R}.
\end{proof}

\section{Examples}

\subsection{Calabi-Yau varieties}



In this paragraph, we recall our constructions and results 
for the case when $k$ is an algebraically closed field of characteristic $p>0$ and $X$ is a  smooth proper variety over $k$. 
Set $d := \dim X$. 
We have the induced morphism $\pi_n \colon W_nX \to \Spec W_n(k)$ and $W_n\omega_X^{\mydot} :=\pi_n^!(W_n(k)) \simeq W_n\omega_X[d]$.
Note that $W_n(k)$ itself is the injective hull of the residue field $k$ and hence the dualising module of the Artin local ring $W_n(k)$. 
Let $\cF$ be a coherent $W_n\cO_X$-module.
By the Grothendieck duality for  $\pi_n \colon W_nX \to \Spec W_n(k)$, we have a canonical isomorphism 
\[
H^i(X, R\cHom_{W_n\cO_X}(\cF, W_n\omega_X)) \simeq \Hom_{W_n(k)} (H^{d-i}(X, \cF), W_n(k)).
\]
In particular, we have
\[
\Hom_{W_n\cO_X}(\cF, W_n\omega_X) \simeq \Hom_{W_n(k)} (H^d(X, \cF), W_n(k)).
\]
We now apply this functorial  isomorphism for $\Phi^e_{X, n} : \MO_X \to Q^e_{X, n} := Q^e_{X, 0, n}$. 
Then 
\[
\Hom_{W_n\cO_X}(Q_{X, n}^e, W_n\omega_X(-K_X)) \xrightarrow{H^0(X, (\Phi^e_{X, n})^*)} 
\Hom_{W_n\cO_X}(\cO_X, W_n\omega_X(-K_X))
\]
is surjective if and only if its $W_n(k)$-dual
\[
H^d(X, \Phi^e_{X, K_X, n}) : H^d(X, \cO_X(K_X)) \to H^d(X, Q_{X, K_X, n}^e)
\]
is injective. 

Now assume $K_X=0$.
Then the above morphism fits into the following diagram
\[
\begin{tikzcd}[column sep=large, row sep=large]
    H^d(X, W_n\cO_X) \arrow[r, "F^e"]\arrow[d, "R^{n-1}", twoheadrightarrow] & H^d(X, W_n\cO_X)\arrow[d, twoheadrightarrow] \\
    H^d(X, \cO_X) \arrow[r, "{H^d(X, \Phi^e_{X, n})}"] & H^d(X, Q_{X, n}^e).
\end{tikzcd}
\]
By definition, we have 
\[
Q_{X, n}^e \simeq W\cO_X/(V^nW\cO_X+VF^eW\cO_X).
\]
Set $H:=H^d(X, W\cO_X)$.
Note that $H$ is not necessarily finitely generated as a $W(k)$-module but equipped with induced operators $F$ and $V$.
Since $H^d(X, -)$ is right exact, we have
\[
H^d(X, Q_{X, n}^e) \simeq H/(V^nH+VF^eH).
\]

Further assume that $H=H^d(X, W\MO_X)$ is a finitely generated free $W(k)$-module of rank $h \geq 1$. 
By the structure theorem of one-dimensional $p$-divisible group over an algebraically closed field, there is a $W(k)$-basis $v_1, \cdots, v_h$ of $H$ such that
\begin{align*}
    &Vv_1=v_2,\,\,\, Vv_2=v_3,\,\,\, \cdots,\,\,\, Vv_{h-1}=v_h,\,\,\, Vv_h=pv_1, \\
    &Fv_1=v_h,\,\,\, Fv_2=pv_1,\,\,\, \cdots,\,\,\, Fv_{h-1}=pv_{h-2},\,\,\, Fv_h=pv_{h-1}.
\end{align*}
When $h=1$, these equations mean $Vv_1=pv_1, {Fv_1=v_1}$. 

Observe that $Fv_i=V^{h-1}v_i$ for each $i=1, \cdots, h$.
Then  it holds that 
\begin{align*}
    V^nH+VF^eH &= V^nH+V^{eh-e+1}H\\
    &=
    \begin{cases}
      V^nH \hspace{19mm} \text{if} \qquad  eh-e+1> n, \\
      V^{eh-e+1}H \hspace{10mm} \text{if}\qquad eh-e+1 \leq n.
    \end{cases}
\end{align*}

\begin{thm}\label{t-CY}
    Let $k$ be an algebraically closed field of characteristic $p>0$.
    Let $X$ be a $d$-dimensional  smooth proper variety  over $k$ with $\omega_X \simeq \cO_X$. 
    Assume that $H^d(X, W\cO_X)$ is a finitely generated free $W(k)$-module of rank $h \geq 1$.
    Take  $e \in \Z_{>0}$. 
    Then $X$ is quasi-$F^e$-split of height $eh-e+1$.
\end{thm}

\begin{proof}
We use the same notation as above.
The image of $v_1$ under the natural surjection $H \twoheadrightarrow H^d(X, \cO_X)$ is a basis of the $k$-vector space $H^d(X, \cO_X)$. 
Therefore, $X$ is $n$-quasi-$F^e$-split (i.e., $h^e(X) \leq n$) if and only if 
$F^ev_1 \not\in V^nH+VF^eH$. 
Since $Fv_1=V^{h-1}v_1$, we have
\[
F^ev_1=V^{eh-e}v_1.
\]

First assume that $n \leq eh-e$.
Then 
\[
F^ev_1 = V^{eh-e}v_1 \in V^{eh-e}H \subseteq V^nH \subseteq V^nH+F^eVH=V^nH. 
\]
Hence the image of $F^ev_1$ in $H/(V^nH+VF^eH)$ is zero. 
Then it holds that $h^e(X) > eh-e$. 

Next assume $n = eh-e+1$.
In this case, we have $V^nH+F^eVH=V^{eh-e+1}H$.
Since $V : H \to H$ is injective and $v_1 \notin VH$, 
the element $F^ev_1=V^{eh-e}v_1$ is not contained in $V^{eh-e+1}H$. 
Thus $h^e(X) \leq  eh-e+1$. 
\end{proof}


\begin{cor}\label{c-str-CY}
    Let $k$ be an algebraically closed field of characteristic $p>0$.
    Let $X$ be a $d$-dimensional smooth proper variety  over $k$ with $\omega_X \simeq \cO_X$. 
    {Assume that $H^{d-1}(X, \cO_X)=0$ when $d \geq 2$.} 
    Take  $e \in \Z_{>0}$. 
    Then $h^e(X) = eh(X)-e+1$. 
\end{cor}

\begin{proof}
If $h(X)=\infty$, then we have $h^e(X)=\infty$ by definition, and hence 
$h^e(X) = eh(X)-e+1$ holds. 
In what follows, we assume $h(X) < \infty$. 
By our assumption, $X$ has the associated Artin-Mazur formal group, which is a one-dimensional formal group \cite[II. Proposition 1.8]{artin-mazur77}. 
Note that we have $h(X) = \dim_K H^d(X, W\MO_X) \otimes_{W(k)} K$ for 
$K := \text{Frac}\, W(k)$ \cite[Theorem 4.5]{yobuko19}. 
Moreover, $h(X)< \infty$ implies that $H^d(X, W\MO_X)$ is a finitely generated free $W(k)$-module. 
Then the assertion follows from Theorem \ref{t-CY}. 
\end{proof}

\begin{cor}\label{c-abel-var}
Let $k$ be a perfect field of characteristic $p>0$.
Let $X$ be a $g$-dimensional abelian variety over $k$. 
Set $f(X)$ to be its $p$-rank. 
\begin{enumerate}
\item If $f(X)=g$, then $h^e(X)=1$. 
\item If $f(X)=g-1$, then $h^e(X)=e+1$
\item If $f(X) \leq g-2$, then $h^e(X) = \infty$. 
\end{enumerate}
\end{cor}
\noindent Recall that $f(X)$ is equal to the $\F_p$-dimension of the $p$-torsion subgroup $X(\overline{k})[p]$.
\begin{proof}
When $f(X)=g$ (resp. $f(X)=g-1$), $H^g(X, W\cO_X)$ is a free $W(k)$-module of rank one (resp. two) \cite[Theorem 3.2 and its proof]{yobuko23}. 
Hence (1) and (2) follow from Theorem \ref{t-CY}. 
If $f(X) \leq g-2$, then $h(X)=\infty$ \cite[Theorem 3.2]{yobuko23}, and hence $h^e(X) =\infty$. 
Thus (3) holds. 
\end{proof}



\begin{proposition}\label{p-CY-unif-QFS}
    Let $X$ be a smooth projective variety over a perfect field $k$ of characteristic $p>0$. 
Assume that $K_X$ is pseudo-effective (e.g., $\kappa(X) \geq 0$). 
Then $X$ is $F$-split if and only if $X$ is uniformly quasi-$F^{\infty}$-split. 
\end{proposition}

\begin{proof}
By Proposition \ref{p-bc-base-field}, we may assume that $k$ is an algebraically closed field. 
As the \lq\lq only-if" part is clear, let us prove the \lq\lq if" part. 
In what follows, we assume that $X$ is uniformly quasi-$F^{\infty}$-split. 
The proof consists of the following two steps. 
\begin{enumerate}
\item[(A)] Reduction to the case when $K_X \sim 0$. 
\item[(B)] The proof for the case when $K_X \sim 0$. 
\end{enumerate}

(A) Let us reduce the proof to the case when $K_X \sim 0$. 
Since $X$ is quasi-$F$-split, 
there exists an integer $n>0$ such that $H^0(X, -(p^n-1)K_X) \neq 0$ \cite[the proof of Proposition 3.14]{KTTWYY1}. 
Since $K_X$ is pseudo-effective, we get $(p^n-1)K_X \sim 0$. 
Let $\pi : Y \to X$ be the associated cyclic \'etale cover from a smooth projective variety $Y$ with $K_Y \sim 0$. 
It follows from Corollary \ref{c-Galois-ht} that $h^e(X) = h^e(Y)$, and hence the problem is reduced to the case when $K_X \sim 0$. 

\medskip

(B) 
Assume $K_X \sim 0$. 
Since $X$ is uniformly quasi-$F^{\infty}$-split, 
there exists $e_0 \geq 1$ such that 
\[
h^1(X) \leq h^2(X) \leq \cdots \leq h^{e_0}(X) = h^{e_0+1}(X) = \cdots =: n <\infty. 
\]
    Fix an integer $e \geq e_0+1$.
    By assumption, there exists a nonzero $W_n\MO_X$-module homomorphism
$\alpha \colon F_*^eW_n\mcO_X \to W_n\omega_X(-K_X)$ 
    such that $\alpha \circ F^e = \underline{p}^{n-1} \circ R^{n-1}$: 
\[\begin{tikzcd}
	W_n\MO_X & {} & F_*^eW_n\MO_X \\
	{\mathcal O_X = \omega_X(-K_X)} \\
	{W_n\omega_X(-K_X)}
	\arrow["R^{n-1}"', from=1-1, to=2-1]
	\arrow["\underline{p}^{n-1}"', from=2-1, to=3-1]
	\arrow["{F^e}", from=1-1, to=1-3]
	\arrow["{\exists \alpha}", dashed, from=1-3, to=3-1]
\end{tikzcd}\]
{Consider the composition 
    \[
    \beta := \alpha \circ V^{n-1} \colon F_*^{e+n-1}\mcO_X \to F_*^eW_n\mcO_X \to W_n\omega_X(-K_X). 
    \]
    Then $\beta$ is nonzero by $n= h^e(X)$ and the exact sequence 
\[
0 \to F_*^{e+n-1}\mcO_X \xrightarrow{V^{n-1}} F_*^eW_n\mcO_X \xrightarrow{R} F_*^eW_{n-1}\MO_X \to 0. 
\]
    Also consider the composition
    \[
    \alpha' :=\alpha \circ F \colon F_*^{e-1}W_n\mcO_X 
    \xrightarrow{F}  F_*^{e}W_n\mcO_X \xrightarrow{\alpha} W_n\omega_X(-K_X),
    \]
    which is nonzero, because 
    $\alpha' \circ F^{e-1} = \alpha \circ F^e = \underline{p}^{n-1} \circ R^{n-1} \neq 0$.
    By $n = h^{e-1}(X)$,  $\beta':=\alpha' \circ V^{n-1}$ is nonzero.
    By definition, under the induced $W_n\MO_X$-module homomorphism
    \[
    F^* \colon \Hom(F_*^{e+n-1}\mcO_X, W_n\omega_X(-K_X)) \to \Hom(F_*^{e+n-2}\mcO_X, W_n\omega_X(-K_X)), 
    \]
    we have 
    \[
    F^*(\beta)=\beta \circ F = \alpha \circ V^{n-1} \circ F  = \alpha \circ F \circ V^{n-1}= \alpha' \circ V^{n-1} = \beta'.
    \]
        Via Grothendieck duality, the above morphism is isomorphic to $H^0$ of 
    \[
    \Psi \colon F_*^{e+n-1}\mcO_X((1-p^{e+n-1})K_X) \to F_*^{e+n-2}\mcO_X((1-p^{e+n-2})K_X),
    \]
    which we know to be nonzero on global sections.
    By tensoring $\Psi$ with $\mcO_X(K_X)$, we get
    \[
    F_*^{e+n-1}\mcO_X(K_X) \to F_*^{e+n-2}\mcO_X(K_X)
    \]
    which is the pushforward $F_*^{e+n-2}$ of the Cartier operator $C \colon F_*\mcO_X(K_X) \to \mcO_X(K_X)$.
    Now assume that $K_X$ is trivial.
    Then $\Psi$ and $C$ are isomorphic to each other.
    So we see that the Cartier operator $C$ on $H^0(X, \mcO_X(K_X))$ is nonzero, which implies that $X$ is $F$-split.
   }
\end{proof}


\subsection{Log Fano curves}

\begin{thm}\label{t-fano-curve-weak}
Let $\kappa$ be a field {which is} finitely generated over a perfect field $k$ of characteristic $p>0$. 
Let $X$ be a regular projective curve over $\kappa$ and 
let $\Delta$ be an effective $\bQ$-divisor on $X$ such that 
{$\rdown{\Delta} =0$} 
and  $-(K_X+\Delta)$ is ample. Then {the following hold.}
\begin{enumerate}
    \item $(X,\Delta)$ is {feebly} globally quasi-$F$-regular. 
    \item 
    {$(X, \Delta)$ is globally quasi-$+$-regular}. 
\end{enumerate}
\end{thm}

\begin{proof}
{Let us show (1).} 
Fix an effective $\Q$-divisor $E$ on $X$. 
We can find a sufficiently small $\epsilon \in \Q_{>0}$ such that $-(K_X+\Delta +\epsilon E)$ is ample. 
It suffices to find $n \in \Z_{>0}$ such that 
$(X, \Delta +\epsilon E)$ is $n$-quasi-$F^e$-split for every $e>0$. 
This follows from Theorem \ref{thm:higher-Cartier-criterion-for-quasi-F-split}, which can be applied for ${\rm id} \colon Y =X \to X$, because $\kappa$ is essentially of finite type over $k$. 
{Thus (1) holds.} 
{By Definition \ref{d-IQFS2} and Definition \ref{d-weak-QFR}, (1) implies that $(X, \Delta)$ is uniformly quasi-$F^{\infty}$-split. 
Then  (2)  holds by Corollary \ref{c-Fano-QFS-vs-+}, 
which is applicable because $(X, \{p^s\Delta\})$ 
is strongly $F$-regular for every $s \geq 0$.} 
\end{proof}

{At first, one could think that the above proof shows global quasi-$F$-regularity of log Fano curves in the strong sense thanks to 
{the Fujita vanishing theorem}. 
Unfortunately, as one imposes no bounds on the support of $E$, it is not easy to control the fractional parts of divisors that come up in the above proof. Thus, to prove global quasi-$F$-regularity of curves in the strong sense, one needs a more careful argument, which is explained below. We start with the following lemma.} 


\begin{lem}\label{l-e-to-e+1}
Let $X$ be a normal projective variety over an $F$-finite field of characteristic $p>0$. 
Take an effective $\Q$-divisor $\Delta$ on $X$ satisfying $\rdown{\Delta} =0$.  
Fix $e \in \Z_{>0}$ and $n \in \Z_{>0}$. 
Assume that 
\begin{enumerate}
\item $(X, \Delta)$ is $n$-quasi-$F^e$-split. 
\item $H^{\dim X-1}(X, \MO_X(p^f(K_X+\Delta))) =0$ for every integer $f \geq e$. 
\item 
$
F: \MO_X \to F_*\MO_X(p \{ p^f\Delta\})
$ 
splits as an $\MO_X$-module homomorphism  for every integer $f \geq e$. 
\end{enumerate}
Then $(X, \Delta)$ is $n$-quasi-$F^{e+1}$-split. 
\end{lem}

\begin{proof}
Set $d:=\dim X$ and $D := K_X+\Delta$. 
By construction, we have the induced $W_n\MO_X$-module homomorphisms 
\begin{equation}\label{e1-e-to-e+1}
\Phi^{e+1}_{X, D, n} : \MO_X(D) \xrightarrow{\Phi^e_{X, D, n}} Q^e_{X, D, n} \xrightarrow{\psi} Q^{e+1}_{X, D, n}. 
\end{equation}
Recall that the exact sequence 
\begin{equation}\label{e2-e-to-e+1}
0\to F_*^eW_n\MO_X(p^eD)  \xrightarrow{F} 
 F_*^{e+1}W_n\MO_X(p^{e+1}D) 
 \to F_*^e B^1_{X, p^eD, n} \to 0 
\end{equation}
induces another exact sequence by taking pushouts: 
\begin{equation}\label{e3-e-to-e+1}
0\to Q^e_{X, D, n} \xrightarrow{\psi} Q^{e+1}_{X, D, n} \to 
F_*^e B^1_{X, p^eD, n}\to 0. 
\end{equation}
By (1), (\ref{e1-e-to-e+1}), and Lemma \ref{lem:cohomological-criterion-for-log-quasi-F-splitting}, it suffices to show that  
\[
H^d(X, \psi) : H^d(X,  Q^e_{X, D, n}) 
\to 
H^d(X, Q^{e+1}_{X, D, n})
\]
is injective. 
It follows from (\ref{e3-e-to-e+1}) that the problem is reduced to 
$H^{d-1}(X, B^1_{X, p^eD, n})=0$. 
Then, by (\ref{e2-e-to-e+1}), it is enough to prove that 
\begin{enumerate}
\item[(i)] $H^{d-1}(X, F_*^{e+1}W_n\MO_X(p^{e+1}D))=0$, and 
\item[(ii)] $F: H^d(X, W_n\MO_X(p^eD))
\to H^d(X, F_*W_n\MO_X(p^{e+1}D))$ is injective. 
\end{enumerate}
The assertion (i) follows from (2) and an exact sequence 
\[
0 \to F^{n-1}_*\MO_X(p^{e+n}D) \xrightarrow{V^{n-1}} W_n\MO_X(p^{e+1}D)  \xrightarrow{R} W_{n-1}\MO_X(p^{e+1}D) \to 0.  
\]
It suffices to show (ii). Consider the following commutative diagram in which each horizontal sequence is exact: 
\[
\begin{tikzcd}[column sep = small]
 0 \arrow[r] & 
H^d(F^{n-1}_*\MO_X(p^{e+{n-1}}D)) \arrow[r, "V^{n-1}"] \arrow[d, "F"]& 
H^d(W_n\MO_X(p^eD)) \arrow[r, "R"] \arrow[d, "F"]& 
H^d(W_{n-1}\MO_X(p^eD)) \arrow[r]\arrow[d, "F"] & 0\\
 0 \arrow[r] & 
H^d(F^n_*\MO_X(p^{e+{n}}D)) \arrow[r, "V^{n-1}"] & 
H^d(F_*W_n\MO_X(p^{e+1}D)) \arrow[r, "R"] & 
H^d(F_*W_{n-1}\MO_X(p^{e+1}D)) \arrow[r] & 0,
\end{tikzcd}
\]
where $H^d(-) :=H^d(X, -)$. 
By induction on $n$, it is enough to show that 
\[
F : H^d(X, \MO_X(p^fD)) \to H^d(X, F_*\MO_X(p^{f+1}D))
\]
is injective for every $f \geq e$. 
This holds, because (3) implies that 
$\MO_X(p^fD) \to F_*\MO_X(p^{f+1}D)$ splits for every $f \geq e$. This completes the proof of (ii). 
\end{proof}

\begin{thm}\label{t-fano-curve-strong}
Let  $k$ be an algebraically closed field of characteristic $p>0$. 
Let $X$ be a smooth projective curve over $k$ and 
let $\Delta$ be an effective $\bQ$-divisor on $X$ such that $\rdown{\Delta}=0$ and  $-(K_X+\Delta)$ is ample. 
Then $(X,\Delta)$ is globally quasi-$F$-regular.
\end{thm}

\begin{proof}
By enlarging the coefficients of $\Delta$, 
the problem is reduced to the case when 
\[
\Delta = \sum_{i=1}^r \frac{\alpha_i}{p^{\nu}} P_i, \qquad \text{where} \qquad \nu, \alpha_1, ..., \alpha_r \in \Z_{>0}. 
\]
Replacing $\nu$ by a larger integer, we may assume that 
$\rdown{\Delta'}=0$ and 
$-(K_X+\Delta')$ is ample for 
\[
\Delta' := \sum_{i=1}^r \frac{\alpha_i +1}{p^{\nu}} P_i. 
\]
Then there is $n \in \Z_{>0}$ such that $(X, \Delta')$ is $n$-quasi-$F^e$-split for every $e >0$ (Theorem \ref{t-fano-curve-weak}(2)).

Fix an effective $\Q$-divisor $E$ on $X$. 
Take the decomposition $E = E_1 + E_2$ into the effective $\Q$-divisors $E_1$ and $E_2$ 
such that $\Supp E_1 \subseteq \Supp \Delta$ and $\Supp E_2 \cap \Supp \Delta = \emptyset$. 
Fix $\epsilon_1 \in \Q_{>0}$ satisfying 
\[
\Delta + \epsilon_1 E_1 \leq \Delta'. 
\]
Since $n$ and $\nu$ are fixed, 
we have $Q^{\nu}_{X, \Delta'  +\epsilon_2 E_2, n} = Q^{\nu}_{X, \Delta', n}$ for some $0 < \epsilon_2 \ll 1$, 
which implies that  
\[
H^1(X, \Phi^{\nu}_{X, K_X+\Delta'  +\epsilon_2 E_2, n}) : H^1(X, \MO_X(K_X+\Delta'  +\epsilon_2 E_2))  \to H^1(X, Q^{\nu}_{X, K_X+\Delta'  +\epsilon_2 E_2, n})
\]
is injective, i.e., 
$(X, \Delta' + \epsilon_2 E_2)$ is $n$-quasi-$F^{\nu}$-split.  
Let  
$E_2 = b_1 Q_1 + \cdots + b_s Q_s$ be the irreducible decomposition of $E_2$. 
Pick positive integers $\mu_1, ..., \mu_s$ satisfying 
\begin{itemize}
\item $\nu < \mu_1 < \mu_2 < \cdots < \mu_s$, and 
\item $\overline{E}_2 \leq \epsilon_2 E_2$ for 
\[
\overline{E}_2 := \frac{1}{p^{\mu_1}} Q_1 + \cdots + \frac{1}{p^{\mu_s}}Q_s. 
\]
\end{itemize}
By $\Delta' +\overline{E}_2 \leq \Delta'+\epsilon_2 E_2$, $(X, \Delta' + \overline{E}_2)$ is $n$-quasi-$F^{\nu}$-split. 
We may assume that 
{$-(K_X+\Delta' + \overline{E}_2)$ is ample and $\rdown{  \Delta' + \overline{E}_2}=0$.}

We can find $\epsilon \in \Q_{>0}$ such that 
$\epsilon \leq \epsilon_1$ and $\Delta + \epsilon E \leq \Delta' + \overline E_2$. 
Then it suffices to show that $(X, \Delta'+\overline{E}_2)$ is $n$-quasi-$F^e$-split for every integer $e \geq \nu$. 
Let us prove this by induction on $e$. 
The base case (i.e., the case when $e=\nu$) has been settled already. 
Since $H^0(X, \MO_X(p^f(K_X+\Delta'+\overline{E}_2)))=0$ for every $f \in \Z_{\geq 0}$, 
it is enough, by Lemma \ref{l-e-to-e+1}, to show that 
\[
F: \MO_X \to F_*\MO_X(p \{ p^f(\Delta'+\overline{E}_2)\})
\]
splits for every $f \geq \nu$. 
If $f \geq \nu$, then we have that 
\[
\{ p^f (\Delta' + \overline{E}_2)\} \overset{{\rm (a)}}{=} \{ p^f \overline{E}_2\} 
=
\{ \frac{p^f}{p^{\mu_1}} Q_1 + \cdots + \frac{p^f}{p^{\mu_s}}Q_s\} 
\overset{{\rm (b)}}{\leq }
\frac{1}{p}Q_j + \frac{1}{p^2} (Q_{j+1} + \cdots + Q_s)
\]
for some $j$, 
where 
(a) holds by $p^f \cdot \frac{\alpha_i+1}{p^{\nu}} = p^{f-\nu} (\alpha_i+1) \in \Z$ and (b) follows from $\mu_1 < \cdots < \mu_s$. 
Hence the problem is reduced to the splitting of 
\[
F : \MO_X \to 
F_*\MO_X(p \cdot (\frac{1}{p}Q_j + \frac{1}{p^2}(Q_{j+1} + \cdots +Q_s))) =F_*\MO_X(Q_j), 
\]
which is well known. 
\qedhere



\end{proof}

\subsection{Log Calabi-Yau curves}

\begin{dfn}\label{d-logCY}
We work over an algebraically closed field $k$  of characteristic $p>0$. 
We say that $(X=\mathbb P^1, \Delta)$ is a log Calabi-Yau pair (over $k$) with standard coefficients 
if $\Delta$ is an effective $\Q$-divisor whose coefficients are contained in $\{1 \} \cup 
\{ 1- 1/n \,|\, n \in \Z_{>0}\}$, and $\deg \Delta =2$. 
\end{dfn}

It is easy to see that 
the classification of log Calabi-Yau curves $(X =\P^1, \Delta)$ with standard coefficients (except for elliptic curves) is given as follows, where $P_1, P_2, P_3, P_4$ are distinct points: 
\begin{enumerate}
\renewcommand{\labelenumi}{(\roman{enumi})}
\item $\Delta= \frac{2}{3}P_1 + \frac{2}{3}P_2 +\frac{2}{3}P_3$. 
\item $\Delta= \frac{1}{2}P_1 + \frac{3}{4}P_2 +\frac{3}{4}P_3$. 
\item $\Delta= \frac{1}{2}P_1 + \frac{2}{3}P_2 +\frac{5}{6}P_3$. 
\item $\Delta= \frac{1}{2}P_1 + \frac{1}{2}P_2 +\frac{1}{2}P_3 + \frac{1}{2}P_4$. 
\item $\Delta= P_1 + \frac{1}{2}P_2 +\frac{1}{2}P_3$. 
\item {$\Delta= P_1 + P_2$.}
\end{enumerate}
Since we are interested in the case when $\rdown{\Delta} =0$, we shall treat (i)-(iv). 

\begin{prop}\label{p-logCY-Yes}
Let $(X = \mathbb P^1, \Delta)$ be a log Calabi-Yau pair with standard coefficients. 
Assume that {$\rdown{\Delta}=0$ and} 
$(p^s-1)\Delta$ is Cartier for some $s \in \Z_{>0}$. {Then $(X,\Delta)$ is quasi-$F^e$-split for every $e>0$. Moreover:} 
\begin{enumerate}
\item 
if $(X, \Delta)$ is quasi-$F$-split of height $1$, 
then $(X, \Delta)$ is  quasi-$F^e$-split of height $1$ for every $e>0$;
\item 
if $(X, \Delta)$ is not quasi-$F$-split of height $1$, 
then $(X, \Delta)$ is quasi-$F^e$-split of height $e+1$ for every $e>0$. 
\end{enumerate}
\end{prop}
\noindent {Recall that $(X, \Delta)$ is quasi-$F$-split of height $1$ if and only if $\MO_X \to F_*\MO_X(p\Delta)$ splits as an $\MO_X$-module homomorphism. This condition is called \emph{naively keenly $F$-split} in \cite[{Definition 2.19}]{KTTWYY1}.}


{The assumption that $(p^s-1)\Delta$ is Cartier is equivalent to: {$p\neq 3$ in Case (i), $p\neq 2$ in Cases (ii) and (iv), and $p \neq 2, 3$ in Case (iii)}}.
\begin{proof}
By the proof of {\cite[Proposition 3.21]{KTTWYY1}} 
there is a finite Galois cover $f: Y \to X$ 
from an elliptic curve $Y$ such that $\deg f$ is not divisible by $p$ and $K_Y \sim f^*(K_X+\Delta)$ (note that the resulting field extension $K(Y)/K(X)$ is a Kummer extension, because $K(X)$ contains an algebraically closed subfield). {We point out that the proof therein works even in the case of $p=2$ and $3$, as we are assuming that $(p^s-1)\Delta$ is Cartier.}
By Corollary \ref{c-Galois-ht}, we have that $h^e(X, \Delta) = h^e(Y)$. 

(1) Assume that $(X, \Delta)$ is quasi-$F$-split of height $1$. 
Then $h^1(Y) = h^1(X, \Delta)=1$, i.e., $Y$ is an ordinary elliptic curve. 
Therefore, $h^e(X, \Delta) = h^e(Y)=1$ (Corollary \ref{c-abel-var}). 

(2) Assume that $(X, \Delta)$ is not quasi-$F$-split of height $1$. 
Then $h^1(Y) = h^1(X, \Delta)\neq 1$, i.e.,  $Y$ is a supersingular elliptic curve. 
Therefore, $h^e(X, \Delta) = h^e(Y)=e+1$ (Corollary \ref{c-abel-var}). 
\end{proof}

\begin{prop}\label{p-logCY-No}
Let $(X = \mathbb P^1, \Delta)$ be a log Calabi-Yau pair with standard coefficients such that {$\lfloor \Delta \rfloor = 0$}. 
\begin{enumerate}
\item Assume that $H^1(X, \MO_X(p^r(K_X+\Delta)))=0$ for every $r \in \Z_{>0}$. 
Then $(X, \Delta)$ is not quasi-$F$-split. 
\item 
Assume that $p$ and $(X, \Delta)$ satisfy one of the following. 
\begin{enumerate}
\item $p=2, \Delta = \frac{1}{2}P_1+\frac{1}{2}P_2+\frac{1}{2}P_3+\frac{1}{2}P_4$.
\item $p=2, \Delta =  \frac{1}{2}P_1+\frac{3}{4}P_2+\frac{3}{4}P_3$.  
\item $p=2, \Delta = \frac{1}{2}P_1+\frac{2}{3}P_2+\frac{5}{6}P_3$.
\item $p=3, \Delta = \frac{2}{3}P_1+\frac{2}{3}P_2+\frac{2}{3}P_3$.
\item $p=3, \Delta = \frac{1}{2}P_1+\frac{2}{3}P_2+\frac{5}{6}P_3$.
\end{enumerate}
Then $H^1(X, \MO_X(p^r(K_X+\Delta)))=0$ holds for every $r \in \Z_{>0}$ 
and $(X, \Delta)$ is not quasi-$F$-split. 
\end{enumerate}
\end{prop}

\begin{proof}
Let us show (1).  
Fix $n \in \Z_{>0}$. 
It is enough to prove that 
the map $H^1(X, \Phi_{X, K_X+\Delta, n})$, appearing below, is zero for every $n \in \Z_{>0}$ (Lemma \ref{lem:cohomological-criterion-for-log-quasi-F-splitting}):  
\[
\begin{tikzcd}
H^1(X, W_n\mathcal{O}_X(K_X+\Delta)) \arrow[r, "F"] \arrow[d, "R^{n-1}"] & H^1(X,F_*W_n\mathcal{O}_X(p(K_X+\Delta))) \arrow[d] \\
H^1(X, \mathcal{O}_X(K_X+\Delta)) \arrow[r] & H^1(X, Q_{X,K_X+\Delta, n})
\end{tikzcd}
\]
By an exact sequence 
\[
0 \to F_*W_{n-1}\MO_X ( p^{r+1}(K_X+\Delta))  \xrightarrow{V} W_n\MO_X ( p^r(K_X+\Delta)) \xrightarrow{R^{n-1}} 
\MO_X(p^r(K_X+\Delta)) \to 0, 
\]
our assumption  $H^1(X, \MO_X(p^r(K_X+\Delta)))=0$ $(r>0)$ implies that 
 \[
 H^1(X, F_* W_n\MO_X(p^r(K_X+\Delta))) =0
 \]
for every $r \in \Z_{>0}$ and every $n \in \Z_{>0}$. 
Hence the composite map $H^1(X, \Phi_{X, K_X+\Delta, n}) \circ R^{n-1}$ is zero. 
Since the left vertical arrow $R^{n-1}$   is surjective, 
$H^1(X, \Phi_{X, K_X+\Delta, n})$ is also zero. 
Thus (1) holds.\\

Let us show (2).  
If $p^r(K_X+\Delta) \sim 0$, then $H^1(X, \MO_X ( p^r(K_X+\Delta))) \simeq 
 H^1(X, \MO_X)=0$. 
Thus the cases (a) and (d) are settled by (1). 
\vspace{1em}

\noindent \textbf{Case (b).} Assume that $p=2$ and $\Delta = \frac{1}{2}P_1+\frac{3}{4}P_2+\frac{3}{4}P_3$. 
Then it suffices to show 
\[
H^1(X, \MO_X(p^r(K_X+\Delta))) =0
\]
for every integer $r \geq 1$. 
If $r \geq 2$, then we have $p^r(K_X+\Delta) \sim 0$. 
We may assume $r=1$. 
In this case, we have 
\begin{equation}\label{e1-logCY-No}
h^1(X, 2(K_X+\Delta)) = h^0(X, K_X - \llcorner 2(K_X+\Delta) \lrcorner)
=h^0(X, -K_X - \llcorner 2 \Delta \lrcorner)=0. 
\end{equation}
\vspace{0.0em}

\noindent \textbf{Case (c).}
Assume that $p=2$ and $\Delta =\frac{1}{2}P_1+\frac{2}{3}P_2+\frac{5}{6}P_3$. 
Since we have \[(2^{s+2}-2^{s})(K_X+\Delta)  =2^{s-1} \cdot 6(K_X+\Delta) \sim 0\] for every $s \in \Z_{> 0}$, 
$2^r(K_X+\Delta)$ is linearly equivalent to 
either $2(K_X+\Delta)$ or $4(K_X+\Delta)$. 
It suffices to show that $H^1(X, \MO_X(2(K_X+\Delta)))= H^1(X, \MO_X(4(K_X+\Delta)))=0$. 
The former one follows from the same computation as in (\ref{e1-logCY-No}). 
The latter one holds by  
\[
h^1(X, 4(K_X+\Delta))=h^0(X, -3K_X -\llcorner 4\Delta \lrcorner)=0
\]
where the last equality follows from $\deg\,  \llcorner 4\Delta \lrcorner=
\rdown{4\cdot \frac{1}{2}} + \rdown{4\cdot \frac{2}{3}} + \rdown{4\cdot \frac{5}{6}}=2+2+3=
7$. 
\vspace{0.0em}

\noindent \textbf{Case (e).}
Assume that $p=3$ and $\Delta =\frac{1}{2}P_1+\frac{2}{3}P_2+\frac{5}{6}P_3$. 
Since we have 
\[(3^{s+1}-3^s)(K_X+\Delta) \sim 3^{s-1} \cdot 6(K_X+\Delta) \sim 0\] for $s \in \Z_{>0}$, 
it holds that $3^r(K_X+\Delta) \sim 3(K_X+\Delta)$. 
Thus it suffices to show that $H^1(X, \MO_X(3(K_X+\Delta)))=0$. 
This follows from  
\[
h^1(X, 3(K_X+\Delta))=h^0(X, -2K_X -\llcorner 3\Delta \lrcorner)=0
\]
where the last equality holds by 
$\deg\, \llcorner 3\Delta \lrcorner = \rdown{3\cdot \frac{1}{2}} +  \rdown{3\cdot \frac{2}{3}} +  \rdown{3\cdot \frac{5}{6}}= 1+2+2=5$.  
\end{proof}

\begin{thm}\label{t-logCY1}
Let $(X=\mathbb P^1, \Delta)$ be a log Calabi-Yau pair with standard coefficients. 
Assume $\llcorner \Delta \lrcorner =0$. 
Fix $e \in \Z_{>0}$. 
Then the following hold. 
\begin{enumerate}
\item $(X, \Delta)$ is quasi-$F$-split if and only if 
there exists $s \in \Z_{>0}$ such that $(p^s-1)\Delta$ is Cartier. 
\item If $h(X, \Delta)=1$, then $h^e(X, \Delta)=1$. 
\item If $1 < h(X, \Delta) < \infty$, then $h^e(X, \Delta)=e+1$. 
\end{enumerate}
\end{thm}

\begin{proof}
The assertion follows from Proposition \ref{p-logCY-Yes} and Proposition \ref{p-logCY-No}. 
\end{proof}

The following theorem is an explicit version of the above theorem.

\begin{thm}\label{t-logCY2}
Let $(X=\mathbb P^1, \Delta)$ be a log Calabi-Yau pair with standard coefficients. 
Fix $e \in \Z_{>0}$ and take distinct closed points $P_1, P_2, P_3, P_4 \in \P^1$.  
Then the following hold. 
\begin{enumerate}
\item Assume $\Delta=\frac{2}{3}P_1+\frac{2}{3}P_2+\frac{2}{3}P_3$. 
\begin{enumerate}
\item If $p=3$, then $(X, \Delta)$ is not quasi-$F^e$-split. 
\item If $p \equiv 1 \mod 3$, then $(X, \Delta)$ is quasi-$F^e$-split of height $1$. 
\item If $p \equiv 2 \mod 3$, then $(X, \Delta)$ is quasi-$F^e$-split of height $e+1$. 
\end{enumerate}
\item Assume $\Delta= \frac{1}{2}P_1+\frac{3}{4}P_2+\frac{3}{4}P_3$. 
\begin{enumerate}
\item If $p=2$, then $(X, \Delta)$ is not quasi-$F^e$-split. 
\item If $p \equiv 1 \mod 4$, then $(X, \Delta)$ is quasi-$F^e$-split of height $1$. 
\item If $p \equiv 3 \mod 4$, then $(X, \Delta)$ is quasi-$F^e$-split of height $e+1$. 
\end{enumerate}
\item Assume $\Delta= \frac{1}{2}P_1+\frac{2}{3}P_2+\frac{5}{6}P_3$. 
\begin{enumerate}
\item If $p \in \{2, 3\}$, then $(X, \Delta)$ is not quasi-$F^e$-split. 
\item If $p \equiv 1 \mod {3}$, then $(X, \Delta)$ is quasi-$F^e$-split of height $1$. 
\item If {$p \neq 2$ and $p \equiv 2 \mod 3$}, then $(X, \Delta)$ is quasi-$F^e$-split of height $e+1$. 
\end{enumerate}
\item Assume $\Delta=\frac{1}{2}P_1+\frac{1}{2}P_2+\frac{1}{2}P_3+\frac{1}{2}P_4$. 
\begin{enumerate}
\item If $p=2$, then $(X, \Delta)$ is not quasi-$F$-split. 
\item Assume $p\neq 2$. Then $(X, \Delta)$ is quasi-$F^e$-split of height $1$ or $e+1$ 
(in this case, $h^e(X, \Delta)$ depends on $\Supp \Delta$). 
\end{enumerate}
\end{enumerate}
\end{thm}

\begin{proof}
The assertion follows from Theorem \ref{t-logCY1} and \cite[Theorem 4.2(2)]{watanabe91}. 
Recall that the following are equivalent. 
\begin{enumerate}
\renewcommand{\labelenumi}{(\roman{enumi})}
    \item $(X, \Delta)$ is quasi-$F$-split of height $1$. 
    \item $ F : \MO_X \to F_*\MO_X(p\Delta)$ splits as an $\MO_X$-module homomorphism. 
    \item $F: H^{\dim X}(X, \MO_X(K_X+\Delta))  \to H^{\dim X}(X, \MO_X(p(K_X+\Delta)))$ is injective. 
\end{enumerate}
Then, for $R$ as in the statement of 
\cite[Theorem 4.2]{watanabe91}, 
$(X, \Delta)$ is quasi-$F$-split of height $1$ if and only if $R$ is $F$-pure 
\cite[Theorem 3.3(i)]{watanabe91}.
\end{proof}


\subsection{Two-dimensional klt singularities}

\begin{thm}\label{t-2dim-klt-QFR}
Let $R$ be a ring essentially of finite type over a perfect field 
of characteristic $p>0$. 
Let $X$ be a two-dimensional integral normal scheme and 
let $\pi \colon X\to \Spec R$ be a projective morphism with $\dim \pi(X) \geq 1$. 
Let $\Delta$ be an effective $\bQ$-divisor on $X$ such that $(X,\Delta)$ is klt and $-(K_X+\Delta)$ is ample. 
Then  $(X,\Delta)$ is {feebly} globally quasi-$F$-regular. 
\end{thm}

\noindent In particular, 
if $X$ is a Noetherian scheme essentially 
of finite type over a perfect field $k$ of characteristic $p>0$ (i.e., 
there is a finite open affine cover $X = \bigcup_{i=1}^r \Spec\,R_i$ such that 
each $R_i$ is a ring essentially of finite type over $k$) 
and 
$(X, \Delta)$  is a two-dimensional  klt pair, 
then $(X,\Delta)$ is {feebly} locally  quasi-$F$-{regular}. 

\begin{proof} 
{By} taking the Stein factorisation of $\pi \colon X \to \Spec R$ and a localisation of $R$, 
we may assume that $H^0(X, \MO_X)=R$ and $R$ is a local domain. 
Take an effective $\Q$-divisor $E$ on $X$. 
Fix $\epsilon \in \Q_{>0}$ such that 
$(X, \Delta' := \Delta + \epsilon E)$ is klt and $-(K_X+\Delta')$ is ample. 
Let $f \colon Y \to X$ be a log resolution of $(X,\Delta')$.
 We have the induced morphisms: 
\[
g \colon Y \xrightarrow{f} X \xrightarrow{\pi} \Spec R. 
\]
Set $K_Y+\Delta'_Y = f^*(K_X+\Delta')$. 
As $X$ is $\Q$-factorial, 
there exists an $f$-exceptional effective $\Q$-divisor $F$ such that $-F$ is $f$-ample. 
Set $B_Y := \Delta'_Y + \delta F$ for some $0 < \delta \ll 1$. 
We may assume that $\rdown{B_Y} \leq 0$, $-(K_Y+B_Y)$ is ample, and $f_*B_Y = \Delta'$.

It is enough to find $n \in \Z_{>0}$ such that $(X, \Delta')$ is 
$n$-quasi-$F^e$-split for every $e \in \Z_{>0}$ (Definition \ref{d-QFR}). 
To this end, it suffices to find $n \in \Z_{>0}$ such that 
(1)-(3) of Theorem  \ref{thm:higher-Cartier-criterion-for-quasi-F-split} 
hold for every $e \in \Z_{>0}$. 
{Recall that}
\[
H^0 R\Gamma_{\m} Rg_* \cF  = H^0_{\m}(g_*\mathcal{F}) =0
\]
for a torsion-free coherent $\MO_Y$-module $\cF$, where the first equality holds by $R\Gamma_{\m} Rg_* = R(\Gamma_{\m} \circ g_*)$ and the second one follows from the fact that  $g_*\mathcal{F}$ is also torsion-free. 
Hence (1) and (2) 
{of Theorem  \ref{thm:higher-Cartier-criterion-for-quasi-F-split}} 
hold automatically for all $n, e \in \Z_{>0}$. 
By the Serre vanishing theorem, 
we can find $n \in \Z_{>0}$ such that 
\[
H^1(Y, \Omega^1_Y(\log (B_Y)_{\red})^* \otimes \cO_Y(K_Y- \rdown{p^{n+c}(K_Y + {B_Y})})) =0
\]
holds for every $c \geq 0$, i.e.,  (3) {of Theorem  \ref{thm:higher-Cartier-criterion-for-quasi-F-split}} holds for every $e >0$. 
\qedhere
\end{proof}

{When the boundary has standard coefficients, we can prove a stronger result.}
\begin{thm}\label{t-2dim-klt-QFR2}
Let $R$ be a ring essentially of finite type over a perfect field 
of characteristic $p>0$. 
Let $X$ be a two-dimensional integral normal scheme and 
let $\pi \colon X\to \Spec R$ be a projective morphism with $\dim \pi(X) \geq 1$. 
Let $B$ be an effective $\bQ$-divisor on $X$ such that 
$(X, B)$ is klt, $-(K_X+B)$ is ample, and $B$ has standard coefficients. 
Then  $(X, B)$ is globally quasi-$+$-regular. 
\end{thm}

\begin{proof}
By Theorem \ref{t-2dim-klt-QFR}, $(X, B)$ is feebly globally quasi-$F$-regular, and hence 
uniformly quasi-$F^{\infty}$-split (Definition \ref{d-IQFS2}, Definition \ref{d-weak-QFR}). 
Then it is enough to show that $(X, B)$ is locally quasi-$+$-regular (Corollary \ref{c-Fano-QFS-vs-+}). 
In other words, the problem is reduced to the case when $X=\Spec R$. 

For a log resolution $\mu \colon W \to X$ of $(X, B)$, we run a $(K_W+ \mu_*^{-1}B + \Ex(\mu))$-MMP over $X$ \cite[Theorem 1.1]{tanaka16_excellent}. 
Let 
\[
f: Y \to X
\]
be the last step of this MMP. 
By construction, 
\begin{itemize}
\item  $E:=\Ex(f)$ is a prime divisor, 
\item $(Y, E+B_Y)$ is plt for $B_Y := f_*^{-1}B$, and 
\item $-E$ and $-(K_Y+E+B_Y)$ are ample. 
\end{itemize}
We define the effective $\Q$-divisor $B_E$ on $E$ by adjunction: $(K_Y+E+B_Y)|_E = K_E+B_E$. 
We now finish the proof assuming that 
\begin{enumerate}
\item $(E, B_E)$ is globally quasi-$+$-regular, and 
\item  $(Y, B_Y)$ is locally quasi-$+$-regular.  
\end{enumerate}
By (1) and (2), we may apply inversion of adjunction (Corollary \ref{c-IOA-anti-ample+}), and hence 
$(Y, E+B_Y)$ is purely globally quasi-$+$-regular. 
Then its pushforward $(X, f_*(E+B_Y)) = (X, B)$ is globally quasi-$+$-regular. 
Indeed, let $g \colon Z \to X$ be a normal finite cover and $W$ be the normalisation of $Z \times _XY$.
Denote by $h \colon W \to Y$ the induced finite cover.
Since $(Y, B_Y)$ is globally quasi-$+$-regular, for some $n$, there exists a morphism $\psi \colon Q_{Y, B_Y, n}^h \to W_n\omega_Y(-K_Y)$ that “reduces to $1$”, in the sense that the composite with $\Phi \colon \cO_Y \to Q_{Y, B_Y, n}$ coincides with the canonical map $\cO_Y \to W_n\omega_Y(-K_Y)$.
Now let $U \subset X$ be an open subset with complement of codimension two and such that $f^{-1}(U) \to U$ is an isomorphism.
 Via this identification,the restriction of $\psi$ to $f^{-1}(U)$ gives a morphism over $U$.
 Since the relevant sheaves are $S_2$, this morphism extends to a morphism
$Q_{X, B, n}^g \to W_n\omega_X(-K_X)$.
Moreover, the extended morphism also “reduces to $1$”, since this holds over $U$. 

It is enough to show (1) and (2).
The assertion (1) follows from Theorem \ref{t-fano-curve-weak}. 
Since $(E, B_E)$ is strongly $F$-regular, 
so is $(Y, (1-\epsilon)E + B_Y)$ for every $0<\epsilon \leq 1$ by \cite[Theorem 4.1]{Das15} 
(although this reference works over an algebraically closed field of characteristic $p>0$, 
the same argument works under our setting). Hence (2) holds. 
\end{proof}

\section{Appendix: stabilisation for $S^0$ and $B^0$ in characteristic $p>0$}
\label{s:appendix-stabilisation}

In this appendix, we study the stabilisation for $F$-splittings. 
The results in this section 
are well known {to} 
the experts (cf.\ \cite{BST15}), 
but, as far as we know, they have not been written with the exact assumptions we need.


\subsection{Stabilisation for $S^0$}



\begin{lemma} \label{lemma:stabilisation_S^0_local}
In the situation of Setting \ref{setting:most-general-foundations-of-log-quasi-F-splitting} let $D$ be a {$\Q$-Cartier} $\bQ$-divisor on $X$. Then the image of the trace map
\[
{\tr^e \colon} F^e_* \mcO_X(K_X + \lceil p^e D \rceil) \to \mcO_X(K_X + \lceil D \rceil)
\]
stabilises for $e \gg 0$, 
{i.e., 
there exists $e_0>0$ such that $\Im\,\tr^e = \Im\,\tr^{e_0}$ for all $e \geq e_0$}. 
{We denote this image by $\sigma_{\rm keen}(K_X+\lceil D \rceil)$.}
\end{lemma}

\begin{proof}
Fix $n \in \Z_{>0}$ such that $nD$ is Cartier.


We now reduce the problem to the case when $n \not\in p\Z$. 
Otherwise, we can write $n = p^dn'$ for some $d \in \Z_{>0}$ and $n' \in \Z_{>0} \setminus p\Z$. 
For $e \geq d$, we get the following factorisation: 
\[
 \tr^e \colon
F^e_* \mcO_X(K_X + \lceil p^e D \rceil) \xrightarrow{F_*^d\tr^{e-d}} 
F^d_* \mcO_X(K_X + \lceil p^d D \rceil) \xrightarrow{\tr^d}
\mcO_X(K_X + \lceil D \rceil). 
\]
As we are assuming that the statement holds for $p^dD$, 
the image of 
\[
F_*^d\tr^{e-d} \colon 
F^e_* \mcO_X(K_X + \lceil p^e D \rceil) \to 
F^d_* \mcO_X(K_X + \lceil p^{d} D \rceil)
\]
stabilises for $e \gg 0$, and hence so does the image of $\tr^e$.

In what follows, we assume $n \not\in p\Z$.
By taking an open cover $X = \bigcup_{i \in I} X_i$ 
which trivialises $\MO_X(nD)$, 
the problem is reduced to the case when 
there exists an $\MO_X$-module isomorphism $\theta \colon  \cO_X(nD) \xrightarrow{\simeq} \cO_X$.  
For a positive integer $m$ satisfying $p^m \equiv 1 \mod n\Z$, we have that
$\theta$ induces an $\MO_X$-module isomorphism
\[
{\Theta : \MO_X(K_X +\rup{D}) \xrightarrow{\simeq} \MO_X(K_X+ \rup{p^mD}).}
\]
Then the statement follows from \cite[Lemma 13.1]{Gab04} (stating that the image of a high enough power of a $p^{-1}$-{endomorphism} 
$\lambda \colon M \to M$ of a coherent $\cO_X$-module $M$ stabilises) applied to 
\[
M = \bigoplus_{i=0}^{m-1}\underbrace{\mcO_X(K_X + \lceil p^iD \rceil)}_{M_i}. 
\]
Here the 
{$\MO_X$-module homomorphism}
$\lambda \colon F_*M \to M$ is defined {componentwisely} 
for direct summands $F_*M_i$ as the compositions
\begin{itemize}
    \item (for ${0 < i \leq m-1}$): $F_*M_i \to M_{i-1} \hookrightarrow M$, where the first {homomorphism} 
    is the trace map {$\tr \colon F_*M_i=F_* \mcO_X(K_X + \lceil p^i D \rceil) \to \mcO_X(K_X + \lceil p^{i-1}D \rceil) = M_{i-1}$}, 
    and the second {homomorphism} 
    is {the} natural inclusion;
    \item (for $i = 0$): $F_*M_0 {\xrightarrow{F_*\Theta, \simeq}} F_*\cO_X(K_X + \lceil p^m D \rceil) \xrightarrow{{\tr}} \cO_X(K_X + \lceil p^{m-1}D \rceil) = M_{m-1} \hookrightarrow M$, 
    where {the first isomorphism is given by   
    $F_*M_0 = F_*\MO_X(K_X+\rup{D}) \xrightarrow{F_*\Theta, \simeq} F_*\MO_X(K_X+\rup{p^mD})$}, 
    the 
    {second homomorphism $\tr$} 
    is the trace map, 
    and the last {homomorphism} 
    is the natural inclusion. 
\end{itemize}
\end{proof}

With the notation of the above lemma, 
we {get the following} sequence of surjective {$\MO_X$-module homomorphisms induced by the trace map $\tr \colon F^{e+1}_*\MO_X(K_X+\rup{p^{e+1}D}) 
\to F^{e}_*\MO_X(K_X+\rup{p^{e}D})$:}
\[
{\cdots}  
\overset{{\tr}}{\twoheadrightarrow} F^e_*\sigma_{\rm keen}(K_X + \lceil p^eD \rceil) 
\overset{{\tr}}{\twoheadrightarrow} {\cdots} 
\overset{{\tr}}{\twoheadrightarrow}  \sigma_{\rm keen}(K_X + \lceil D \rceil).
\]
The key property of $\sigma_{\rm keen}(K_X + \lceil D \rceil)$ is that for a Cartier divisor $E$, we have that
\begin{equation} \label{eq:key-property-of-sigma-keen}
\sigma_{\rm keen}(K_X + \lceil D + E \rceil) = \sigma_{\rm keen}(K_X + \lceil D \rceil) \otimes \cO_X(E).
\end{equation}

\begin{lemma} \label{lemma:stabilisation_S^0}
In the situation of Setting \ref{setting:foundations-of-log-quasi-F-splitting}, let $A$ be an ample {$\Q$-Cartier} $\bQ$-divisor and let $N$ be a nef {$\Q$-Cartier} $\bQ$-divisor.  Then there exists $e_0>0$ such that 
for {every integer} $r \geq 1$ and  {every integer} $k \geq 0$, the image of
\begin{equation} \label{eq:trace-Frobenius-stabilise}
{\Tr^e \colon}   
H^0(X, F^{e}_* \mcO_X( K_X + \lceil p^{e}(rA + kN) \rceil)) \to H^0(X,\mcO_X( K_X +
\lceil rA + kN \rceil))
\end{equation}
is independent of the choice of $e \geq e_0$, {i.e., $\Im\, \Tr^e = \Im \Tr^{e_0}$ for every $e \geq e_0$}. 
\end{lemma}

\begin{proof}
Set $A_{r,k} := rA + kN$ and 
\[
{M^e_{r, k}} := H^0(X, F^{e}_*\sigma_{\rm keen}(K_X + \lceil p^{e}A_{r,k} \rceil)).
\]
As mentioned earlier, the trace map 
${F^{e+1}_*\MO_X(K_X+\rup{p^{e+1}A_{r, k}}) 
\xrightarrow{\tr} F^{e}_*\MO_X(K_X+\rup{p^{e}A_{r, k}})}$ induces 
{the following commutative diagram for every integer $e \geq 0$: 
\[
\begin{tikzcd}
H^0(X, F^{e+1}_* \mcO_X( K_X + \lceil p^{e+1}A_{r, k} \rceil))
\arrow[r, "\Tr"] &
H^0(X, F^{e}_* \mcO_X( K_X + \lceil p^{e}A_{r, k} \rceil))\\
M^{e+1}_{r, k} \arrow[u, hook]
\arrow[r, "\widetilde{\Tr}"] &
M^{e}_{r, k}, \arrow[u, hook]
\end{tikzcd}
\]
where the vertical arrows are the induced inclusions. 
}

\begin{claim} \label{claim:trace-Frobenius-stabilise} 
There exists $e_1>0$ such that 
{for every $r\geq 1$ and $k \geq 0$,} 
the image ${\rm Im}({\widetilde{\Tr}^e} \colon {M^e_{r, k}} \to {M^0_{r, k}})$ is independent of the choice of $e \geq e_1$. 
\end{claim}

\begin{proof}[Proof of Claim \ref{claim:trace-Frobenius-stabilise}]
{We set} 
\[
{K^e_{r,k} := 
\Ker}( F^{e+1}_*\sigma_{\rm keen}(K_X + \lceil p^{e+1}A_{r,k} \rceil) \to F^{e}_*\sigma_{\rm keen}(K_X +\lceil  p^{e}A_{r,k} \rceil)).
\]
By the long exact sequence {in cohomology}, it is enough to find 
${e_1 \in \Z_{> 0}}$ such that 
\begin{equation} \label{eq:Kerk-vanishing}
H^1(X, {K^e_{r,k}}) = 0
\end{equation}
for all $e \geq e_1$, $r \geq 1$, and $k \geq 0$. 
{Fix $m_0 \in \Z_{>0}$ such that $m_0A$ and $m_0N$ are Cartier. 
Pick $e \geq 0$. 
{ 
Take the integers $r'$ and $k'$ defined by 
$r' := p^e r \mod m_0$ and $k' := p^ek \mod m_0$ (cf.\ Notation \ref{ss:notation}(\ref{ss-nota-mod})).} 
As $(p^er - r')A$ and $(p^ek - k')N$ are Cartier, it holds that} 
\[
{K^e_{r,k}} = F^e_*{K^0_{p^er, p^ek}} = 
F^e_*(K^0_{r', k'} {\otimes \MO_X}((p^er-r')A +(p^ek-k')N)). 
\]
Since there are only finitely many 
{possibilities} for $r'$ and $k'$, 
{the Fujita vanishing theorem enables us to} 
find $e_1 > 0$, independent of $r'$ and $k'$, 
such that the vanishing (\ref{eq:Kerk-vanishing}) holds for all $e \geq e_1$, $r \geq 1$, and $k \geq 0$. This concludes the proof of {Claim \ref{claim:trace-Frobenius-stabilise}}.
\qedhere

\end{proof}

We {now} observe that by Lemma \ref{lemma:stabilisation_S^0_local} applied to $D=p^eA_{r,k}$, 
{it holds that} for every $e>0$, $r>0$ and $k \geq 0$, there exists ${e_2 > 0}$  
{(possibly dependent on $(e, r, k)$)} such that the image of 
\[
{\tr^{e'-e}} \colon F^{e'}_* \mcO_X( K_X + \lceil p^{e'}A_{r,k} \rceil) \to F^e_*\mcO_X( K_X +\lceil p^eA_{r,k} \rceil)
\]
is equal to $F^e_*\sigma_{\rm keen}(K_X+ \lceil p^eA_{r,k} \rceil)$ for all ${e' -e \geq e_2}$. 
Moreover, we may pick ${e_2}$ to be independent of {$(e, r, k)$}. Indeed, by replacing $r$ and $k$ by $p^er$ and $p^ek$, respectively,
we may assume that $e=0$. 
Then, by (\ref{eq:key-property-of-sigma-keen}), 
{it is enough to consider finitely many possibilities for $(r, k)$; 
more explicitly, if $m_0 A$ and $m_0 N$ are Cartier for some $m_0 \in \Z_{>0}$, 
then we may assume that $0 \leq r < m_0$ and $0 \leq k < m_0$.} 
Hence a uniform ${e_2}$ exists.

In particular, we get the inclusions of images:
\begin{multline*}
{\rm Im}({M^e_{r, k}} \xrightarrow{{\widetilde{\Tr}}^e} {M^0_{r, k}})\\
\subseteq {\rm Im}\Big(H^0(X, F^{e}_*\mcO_X( K_X +\lceil p^{e}A_{r,k} \rceil)) \xrightarrow{{\Tr^e}} H^0(X,\mcO_X(K_X + \lceil A_{r,k} \rceil))\Big) \\
\subseteq {\rm Im}(M^{e-e_2}_{r, k} \xrightarrow{{\widetilde{\Tr}^{e-e_2}}} M^0_{r, k})
\end{multline*}
for every $e \geq {e_2}$. By Claim \ref{claim:trace-Frobenius-stabilise}, 
if we take $e-{e_2} \geq e_1$, then the left hand side is equal to the right hand side, and hence all these inclusions are equalities. 
In particular, the images of (\ref{eq:trace-Frobenius-stabilise}) stabilise for all $e \geq {e_0 := e_1 + e_2}$. 
This completes the proof of Lemma \ref{lemma:stabilisation_S^0}.
\qedhere
\end{proof}


\subsection{Stabilisation for  $B^0$}


We fix the assumptions and notation as in Setting \ref{setting:foundations-of-log-quasi-F-splitting}. For a finite surjective {morphism} 
$f \colon Y \to X$ from a normal integral scheme $Y$ and a 
$\bQ$-divisor $D$ {on $X$}, we consider the trace map:
\[
\Tr^f_D \colon H^0(Y, \mcO_{Y}(K_Y + \lceil f^*D\rceil)) \to H^0(X,\mcO_X( K_X +
\lceil D \rceil)).
\]
For $\bQ$-divisors $\Delta$ and $L$ {on $X$}, we define 
\[
B^0(X, \Delta; L) := \bigcap_{f \colon Y \to X} \Im(\Tr^{f}_{L-(K_X+\Delta)}),
\]
where the intersection is taken over all finite surjective {morphism}s $f \colon Y \to X$ from normal integral schemes $Y$. 
{Recall that the pullback $f^*D$ is naturally defined even if $D$ is not $\Q$-Cartier.} 
The following result immediately implies that the above intersection stabilises.

\begin{proposition} \label{prop:stabilisation_B^0}
In the situation of Setting \ref{setting:foundations-of-log-quasi-F-splitting}, let $A$ be an ample {$\Q$-Cartier} $\bQ$-divisor and let $N$ be a nef {$\Q$-Cartier} $\bQ$-divisor.

Then there exists a finite surjective {morphism} $f \colon Y \to X$ from a normal integral scheme $Y$ such that 
\begin{itemize}
    \item for every {integer} $r \geq 1$, {every integer} $k \geq 0$, and
    \item every $\bQ$-divisor $\Delta$ such that $K_X+\Delta$ is $\bQ$-Cartier, 
\end{itemize}
the following equality holds:
\[
B^0(X,\Delta; L) = \Im\big(\Tr^{f}_{L-(K_X+\Delta)}\big),
\]
where $L := K_X+\Delta+rA+kN$.
\end{proposition}
\begin{proof}
By definition, 
\[
B^0(X,\Delta; L) = \bigcap_{g \colon Z \to X}\Im\big(\Tr^{{g}}_{rA+kN}\big),
\]
{where} 
$g \colon Z \to X$ {runs over all} 
finite surjective {morphism}s from normal integral schemes {$Z$}. 
Now, pick {a finite surjective morphism} $f \colon Y \to X$ as in Lemma \ref{lemma:stabilisation_B^02} below. Then
\[
B^0(X,\Delta; L) = \bigcap_{g \colon Z \to X}\Im\big(\Tr^{{g}}_{rA+kN}\big) \supseteq \bigcap_{e>0} \Im\big(\Tr^{F^e \circ f}_{rA+kN}\big) = \Im\big(\Tr^{F^{e_0} \circ f}_{rA+kN}\big)
\]
for some $e_0>0$, where the last equality follows immediately from Lemma \ref{lemma:stabilisation_S^0} applied to $Y$. Clearly,
\[
B^0(X,\Delta; L) \subseteq \Im\big(\Tr^{F^{e_0} \circ f}_{rA+kN}\big),
\]
and so the conclusion of the proposition holds 
{after replacing $F^{e_0} \circ f$ by $f$}. 
\qedhere

\end{proof}

\begin{lemma} \label{lemma:stabilisation_B^02}
In the situation of Setting \ref{setting:foundations-of-log-quasi-F-splitting}, let $A$ be an ample {$\Q$-Cartier} $\bQ$-divisor and let $N$ be a nef {$\Q$-Cartier} $\bQ$-divisor. 

Then there exists a finite surjective {morphism} $f \colon Y \to X$ from a normal integral scheme $Y$ such that the following holds: for every 
finite surjective {morphism} $g \colon Z \to X$ from a normal integral scheme $Z$, there exists $e_0 >0$ such that 
\begin{equation} \label{eq:trace-finite-stabilise}
\Im(\Tr^g_{rA + kN}) \supseteq \Im(\Tr^{F^e \circ f}_{rA+kN}).
\end{equation}
for all $e \geq e_0$, $r>0$, and $k \geq 0$.
\end{lemma}

\begin{proof}
Set $A_{r,k} := rA + kN$. 
For a $\bQ$-divisor $D$ {on $X$}, define
\[
\tau_+(D) := \bigcap_{g \colon Z \to X} \Im(\tr^g_D),
\]
where $g \colon Z \to X$ are finite surjective {morphism}s from normal integral schemes $Z$, and $\tr^g_D$ is the trace map:
\[
\tr^g_D \colon g_*\mcO_Z(K_Z + \lceil g^*D\rceil) \to \mcO_X(K_X +
\lceil D \rceil).
\]
\begin{claim} \label{claim:stabilisation_B^02}  $\tau_+(A_{r,k}) = \Im(\tr^f_{A_{r,k}})$ 
for some fixed finite surjective {morphism} $f \colon Y \to X$ from a normal integral scheme $Y$, which is independent of $r >0$ and $k \geq 0$.
\end{claim}

\begin{proof}[Proof of {Claim \ref{claim:stabilisation_B^02}}]
Fix ${m_0 \in \Z_{>0}}$ such that ${m_0}A$ and ${m_0}N$ are Cartier. 
{Take the integers $r'$ and $k'$ defined by} $r' := r \mod {m_0}$ and $k' := k \mod {m_0}$ 
{(cf.\ Notation \ref{ss:notation}(\ref{ss-nota-mod}))}. In particular, $(r-r')A$ and $(k-k')N$ are Cartier. Then 
\[
\tau_+(A_{r,k}) = \tau_+(A_{r',k'}) \otimes \cO_X((r-r')A + (k-k')N).
\]
Moreover, {if $r'$ and $k'$ are fixed, then we get }
\[
\tau_+(A_{r',k'}) = \Im(\tr^f_{A_{r',k'}})
\]
for some 
finite surjective {morphism} $f \colon Y \to X$ {(dependent on $r'$ and $k'$)} 
from a normal integral scheme $Y$ by \cite[Proposition 4.2, Theorem 4.6]{BST15}. Since there are only finitely many {possibilities} 
for $r'$ and $k'$, we can pick such $f$ uniformly for all the choices. 
{This completes the proof of Claim \ref{claim:stabilisation_B^02}.}
\end{proof}

We fix $f$ as in {Claim \ref{claim:stabilisation_B^02}}. 
Define $\overline{\Tr}^f_{A_{r,k}}$ as the induced map {as follows}: 
\[
\begin{tikzcd}[column sep = 4em]
H^0(Y, \tau_+(f^*A_{r,k})) \ar[hook]{d} \ar{r}{\overline{\Tr}^f_{A_{r,k}}} & H^0(X, \tau_+(A_{r,k})) \ar[hook]{d} \\ 
H^0(Y, \cO_Y(K_Y + \lceil f^*A_{r,k} \rceil)) \ar{r}{{\Tr}^f_{A_{r,k}}} & H^0(X, \cO_X(K_X + \lceil A_{r,k} \rceil)).
\end{tikzcd}
\]
\begin{claim} \label{claim:trace-finite-stabilise_trbar}
For every other finite surjective {morphism} $g \colon Z \to X$ from a normal integral scheme $Z$, there exists $e_0 >0$ such that
\[
\Im(\overline{\Tr}^g_{rA + kN}) \supseteq \Im(\overline{\Tr}^{F^e \circ f}_{rA+kN})
\]
for all $e \geq e_0$, $r >0$, and $k \geq 0$.
\end{claim}

\begin{proof}[Proof of {Claim \ref{claim:trace-finite-stabilise_trbar}}] 
Pick {a} 
finite surjective {morphism} $g \colon Z \to X$ as in the statement of 
{Claim \ref{claim:trace-finite-stabilise_trbar}}. 
{After replacing $g \colon Z \to X$ by a higher one $g' \colon Z' \to X$,} 
{we may assume} 
that ${g}$ 
factors through $f$, that is, we have a diagram:
\[
g \colon Z \xrightarrow{h} Y \xrightarrow{f} X.
\]
Consider the following {commutative} diagram:
\[
\begin{tikzcd}[column sep = 4.5em]
H^0(Z, \tau_+(p^eg^*A_{r,k}))  \ar{r}{\overline{\Tr}^h_{p^ef^*A_{r,k}}} \ar{d}{\overline{\Tr}^{{F^e}}} & H^0(Y, \tau_+(p^ef^*A_{r,k})) \ar{d}{\overline{\Tr}^{{F^e}}} \ar[xshift=2em, bend left = 60]{dd}{\Tr^{f \circ F^e}_{A_{r,k}}} \\
H^0(Z, \tau_+(g^*A_{r,k})) \ar{dr}[swap]{\overline{\Tr}^g_{A_{r,k}}} \ar{r}{\overline{\Tr}^h_{f^*A_{r,k}}} & H^0(Y, \tau_+(f^*A_{r,k})) \ar{d}{\overline{\Tr}^f_{A_{r,k}}}\\
& H^0(X, \tau_+(A_{r,k})).
\end{tikzcd}
\]
For $e \gg 0$, we have that 
\[
\Im\big(\overline{\Tr}^g_{A_{r,k}}\big) \supseteq \Im\big(\overline{\Tr}^{g\circ F^e}_{A_{r,k}}\big) = \Im\big(\overline{\Tr}^{f\circ F^e}_{A_{r,k}}\big),
\]
where the last equality follows from the fact that $\overline{\Tr}^h_{p^ef^*A_{r,k}}$ is surjective. 
{The surjectivity of $\overline{\Tr}^h_{p^ef^*A_{r,k}}$ 
holds, because  
the} Fujita vanishing {theorem assures} 
\[
H^1(Y, {K^e_{r,k}}) = 0,
\]
where
\[
{K^e_{r,k}} := \Ker \big(h_*\tau_+(p^eg^*A_{r,k}) \dhxrightarrow{\tr^h_{f^*A_{p^er,p^ek}}} \tau_+(p^ef^*A_{r,k})\big).
\] 
This is analogous to the proof of Claim \ref{claim:trace-Frobenius-stabilise}, and so we leave the details to the reader. 
{This completes the proof of Claim \ref{claim:trace-finite-stabilise_trbar}.}
\end{proof}

By Claim \ref{claim:stabilisation_B^02} applied to $Y$, 
{there is} a finite surjective {morphism} $h \colon {\widetilde{Y}} 
\to Y$ from a normal integral scheme ${\widetilde Y}$ 
such that
{
\[
\tau_+(f^*A_{r,k}) = \Im(\tr^h_{f^*A_{r,k}}), 
\]
which implies 
\[
H^0(Y, \tau_+(f^*A_{r,k})) = H^0(Y, \Im(\tr^h_{f^*A_{r,k}})) \supseteq 
\Im(H^0(\tr^h_{f^*A_{r,k}})) = \Im (\Tr^h_{f^*A_{r,k}}). 
\]
Taking the images of $H^0(Y, \tau_+(p^ef^*A_{r,k}))  \supseteq \Im (\Tr^h_{p^ef^*A_{r,k}})$  under $\Tr^{F^e \circ f}_{A_{r, k}}$, we get }
\begin{equation}\label{eq-fFe-fhFe}
\Im\big(\overline{\Tr}^{F^e \circ f}_{rA+kN}\big) \supseteq  \Im\big(\Tr^{F^e \circ f \circ h}_{rA+kN}\big)
\end{equation}
for every $e \geq 0$.

Now pick a finite surjective {morphism} $g \colon Z \to X$ from a normal integral scheme $Z$. Then
\[
\Im\big(\Tr^g_{rA+kN}\big) \supseteq \Im(\overline{\Tr}^g_{rA+kN})\overset{{\rm (i)}}{\supseteq}  \Im\big(\overline{\Tr}^{F^e \circ f}_{rA+kN}\big) \overset{{\rm (ii)}}{\supseteq} \Im\big(\Tr^{F^e \circ f \circ h}_{rA+kN}\big).
\]
for $e \geq e_0$, {where (i) follows from Claim \ref{claim:trace-finite-stabilise_trbar} and (ii) holds by (\ref{eq-fFe-fhFe})}. 
Thus the statement of the lemma (that is, (\ref{eq:trace-finite-stabilise})) follows 
{after replacing $f \circ h$ by $f$.} 
\end{proof}

\bibliographystyle{skalpha}
\bibliography{bibliography.bib}

\end{document}